\documentclass[11pt]{article}
\usepackage{paper,asb}

\usepackage{amsmath}
\usepackage{amsfonts}
\usepackage{amssymb}
\usepackage{amsthm}
\usepackage{mathtools}
\usepackage{hyperref}
\usepackage{makecell}

\usepackage[numbers, sort&compress ]{natbib}

\numberwithin{equation}{section}

\newtheorem{innercustomoracle}{Oracle}
\newenvironment{customoracle}[1]
  {\renewcommand\theinnercustomoracle{#1}\innercustomoracle}
  {\endinnercustomoracle}

\definecolor{darkgreen}{rgb}{0,0.7,0}

\usepackage{subcaption}

\usepackage{graphicx} 
\usepackage{adjustbox} 

\usepackage{longtable}
\usepackage{booktabs}
\usepackage{arydshln}
\usepackage{makecell}

\usepackage{chngcntr}

\newtheorem{theorem}{Theorem}[section]
\newtheorem{lemma}[theorem]{Lemma}
\newtheorem{assumption}[theorem]{Assumption}

\newtheorem{corollary}[theorem]{Corollary}
\newtheorem{proposition}[theorem]{Proposition}

\newcommand{\papertitle}{Stochastic and Noisy Oracles in Smooth Continuous Optimization}
\newcommand{\paperauthora}{Albert S. Berahas}
\newcommand{\paperauthoraaffiliation}{Department of Industrial and Operations Engineering, University of Michigan}
\newcommand{\paperauthoraemail}{albertberahas@gmail.com}
\newcommand{\paperauthorb}{Liyuan Cao}
\newcommand{\paperauthorbaffiliation}{School of Mathematics, Nanjing University}
\newcommand{\paperauthorbemail}{caoliyuan@nju.edu.cn}
\newcommand{\paperauthorc}{Katya Scheinberg}
\newcommand{\paperauthorcaffiliation}{School of Industrial and Systems Engineering, Georgia Institute of Technology}
\newcommand{\paperauthorcemail}{katyascheinberg@gmail.com}

\DeclarePairedDelimiter{\ceil}{\lceil}{\rceil}

\usepackage{tikz}
\usepackage{pgfplots}
\usetikzlibrary{fillbetween}
\usepackage{amsmath}
\usepackage{bbm}
			
\begin{document}
\title{\papertitle}

\author{\paperauthora\footnotemark[1] \ \and
    \paperauthorb\footnotemark[2]\ \and   \paperauthorc\footnotemark[3]\  \footnotemark[4]}

\maketitle

\renewcommand{\thefootnote}{\fnsymbol{footnote}}
\footnotetext[1]{\paperauthoraaffiliation. (\url{\paperauthoraemail})}
\footnotetext[2]{\paperauthorbaffiliation. (\url{\paperauthorbemail})}
\footnotetext[3]{\paperauthorcaffiliation. (\url{\paperauthorcemail})}
\footnotetext[4]{Corresponding author.}
\renewcommand{\thefootnote}{\arabic{footnote}}

\begin{abstract}{
We review models of inexact computational oracles in smooth continuous optimization and their implications for the convergence, computational complexity, and attainable accuracy of optimization algorithms. We adopt a framework in which zeroth-, first-, and second-order oracles are characterized by quantitative error bounds and may accept additional control parameters that influence both their accuracy and computational cost. We introduce deterministic and stochastic oracle models and illustrate their properties through a variety of examples, including finite-difference approximations, stochastic gradients, and randomized derivative estimators. Particular attention is given to the distinction between absolute and relative errors, bias and variance, and stochastic function evaluations with independent or coupled randomness.

A central theme of the review is the use of common analytical tools to recover and compare convergence and complexity guarantees under different oracle models. For deterministic and stochastic first-order methods, we show how a general descent inequality leads to complexity bounds under various assumptions on the objective function and all proposed oracle models.  In the stochastic setting, we develop bounds on expected stopping times and the accuracy that can be guaranteed under different oracle conditions. We also examine the total computational complexity of algorithms relying on specific oracle constructions and briefly review models of inexact second-order information and their use in optimization algorithms.

}
\end{abstract}

\tableofcontents

\newpage

\section{Introduction}\label{sec.intro}
Continuous optimization has a long history, with some of its key elements dating back to Newton, Lagrange, and Cauchy. The second half of the 20th century saw major advances in modern optimization theory and practice. In the deterministic optimization community, these advances included the development of computational complexity as an important concept that became widely used in the analysis of optimization algorithms. Much of the classical complexity analysis for deterministic optimization was developed under the assumption that the objective function and its derivatives could be evaluated exactly. While some research addressed the presence of computational errors, such work was less widespread, often relied on assumptions that the errors diminished as an optimal solution was approached, and did not always provide explicit computational complexity bounds.

At the same time, the stochastic and simulation optimization communities developed a significant body of work on optimizing functions expressed as expectations over probability distributions and approximated via sample averages. In other words, they focused on problems where computational errors are stochastic and, in many classical settings, unbiased. Historically, theoretical analyses focused primarily on the asymptotic behavior of function and gradient estimators and the optimization algorithms that use them, with less emphasis on finite-time complexity bounds.

During the first quarter of the 21st century, we have witnessed a substantial expansion of research on continuous optimization with inexact, often stochastic, function values and gradients. The range of such problems and applications has expanded enormously, driven by developments in machine learning, reinforcement learning, and artificial intelligence (AI),  among other engineering innovations of this century such as digital twins.  While many optimization approaches have been developed in these domains, much of this research shares two common objectives: first, to model an optimization problem through an idealized function with desirable properties that cannot be computed exactly but can be approximated; and second, to develop an algorithm with {\em finite-time} complexity guarantees for optimizing this function. The types of functions considered, the forms of inexact information available about them, and the ways this information is used in algorithms give rise to a wide variety of settings to consider and analyze.

\paragraph{What is this review for?}
The goal of this article is to identify common patterns that have emerged in past research and to provide a simple, unified framework through which key results can be recovered, their relationships clarified, and, in some cases, new results obtained using similar analytical tools.
Specifically, we propose to characterize inexact information about the (possibly idealized) objective function through the properties of computable zeroth-, first-, and second-order oracles.
 These oracles have certain properties, which we will define and illustrate through many examples. We analyze basic algorithms that use these oracles, such as gradient descent, and establish a clear connection between the properties of the oracles and those of the algorithms.

Overall, this review seeks to address three fundamental questions:
\begin{enumerate}
    \item How does inexact information arise in continuous optimization, and how can different sources and types of inexactness be characterized through computational oracle models?

    \item How do the properties of inexact oracles influence the design of optimization algorithms, their convergence guarantees, computational complexity, and attainable accuracy?

    \item What common analytical tools can be used to  unify theoretical guarantees across different oracle models and optimization algorithms?
\end{enumerate}

\paragraph{What is an oracle?}
An oracle, in a broad sense, refers to a source of wisdom or prophecies. Historically, oracles played significant roles in various cultures, with one of the most famous being the Oracle of Delphi in ancient Greece, where the priestess Pythia was said to have channeled the god Apollo and delivered enigmatic predictions. In modern usage, the term oracle encompasses any person or system that provides authoritative answers or predictions, and it is often used in computing to refer to a ``black-box'' operation that takes an input and produces an output.

In optimization, an oracle has traditionally been understood as an operation that computes a function value (zeroth-order oracle), its gradient (first-order oracle), or its Hessian (second-order oracle). These oracles can be exact or approximate, but they are usually considered black-box operations whose output is not directly controlled by the algorithm designer. However, a closer look at the optimization literature reveals examples in which gradient approximations with controllable properties are constructed using available first- or zeroth-order oracles. As a basic example, choosing a minibatch of a certain size in empirical risk minimization controls the variance of the stochastic gradient used by a stochastic gradient descent algorithm. This, in turn, may affect the choice of step size in stochastic gradient descent. What is an oracle here? Is it the gradient of the loss function computed from a single data point, or the average of such gradients over a minibatch? The traditional optimization view of oracles suggests the former, because the output of this operation is not under the control of the algorithm designer. Here we propose to apply the term to the latter, because this allows us, in some cases, to separate the algorithm more clearly from the computational procedures it relies on, while establishing a more direct connection between the properties of these objects (oracles) and the behavior of the algorithm.
In other words, we introduce additional input parameters in our oracle definitions that allow the algorithm developer to control the oracle properties to some extent. This perspective allows us to study not only how oracle accuracy affects the convergence and complexity of an optimization algorithm, but also how the computational cost of obtaining that accuracy influences the overall complexity.

We focus on unconstrained nonlinear optimization problems of the form
\begin{align} \label{eq.prob}
    \min_{x \in \mathbb{R}^n} \phi(x),
\end{align}
where $\phi: \mathbb{R}^n \rightarrow \mathbb{R}$ is a smooth objective function that may or may not be convex.

Optimization algorithms for this problem typically generate iterates approaching a stationary point or a minimizer, using information about function values, gradients, and Hessians at or near the current iterate. In our framework, computational oracles provide approximations to this information and accept additional control parameters, such as minibatch sizes, discretization levels, or accuracy tolerances. These parameters may affect both the accuracy and the computational cost of the oracle.

In practice, oracles do not always return exact values of the desired information.
Apart from unavoidable rounding errors in floating-point arithmetic, direct evaluation of the function or its derivatives can be expensive or even impossible, depending on the nature of the problem.
Such inexact oracles arise in many fields, including Derivative-Free Optimization (DFO) \cite{ConnScheVice08c,LarsMeniWild2019,brent2013algorithms,ghadimi2013stochastic,kiefer1952stochastic,spall2005introduction,berahas2018derivative}, Simulation Optimization \cite{pasupathy2018sampling,shashaani2018astro}, and Machine Learning \cite{choro,fazel2018global,liu2018zeroth,TRPO,shamir2017optimal,duchi2015optimal,jamieson2012query,bollapragada2019adaptive,bogolubsky2016learning}, to name just a few areas and examples. The list of optimization applications involving inexact function information, i.e., inexact oracles, continues to grow and is too extensive to discuss exhaustively in this article. We will provide many examples and hope that they inspire others to formulate further applications within the  proposed framework.

\paragraph{Model of error}
Different mathematical models of errors in inexact oracles have been proposed in the literature, reflecting the nature of the underlying computational procedures and the information available about their accuracy.
For zeroth-order oracles (function values), for example, the noise may be deterministic or stochastic. When it is stochastic, it may be biased or unbiased. Moreover, a random sample in a stochastic zeroth-order oracle may produce a single random value (which we refer to as a {\em single-point} oracle) or an entire random function (which can subsequently be evaluated at multiple points and is therefore referred to as a {\em multi-point} oracle). These distinctions lead to different oracle classes, whose properties, in turn, affect the design of optimization algorithms and their complexity guarantees. As we will show, the ability to evaluate a random function at multiple points using the same realization of randomness can substantially improve the accuracy of finite-difference approximations and the resulting complexity bounds.

We distinguish these key properties and, in each case, bound oracle errors either deterministically or through their moments, by a quantity that characterizes the accuracy of the corresponding oracle and ultimately contributes to the convergence properties of optimization algorithms. Specifically, we use $\epsilon_f$ to bound the absolute value of the function-evaluation bias (or the error itself in the deterministic case), and $\sigma_f^2$ to bound the variance of a single-point oracle. For a multi-point oracle, the notation requires further setup and will be introduced in due course.

The types of inexactness in first-order oracles can be even more varied and we propose four error parameters:  $\epsilon_g$ and $\kappa_g$ characterize   absolute and relative errors, respectively, while $\sigma_c$ and $\sigma_v$ characterize  absolute and relative variance components. Depending on the oracle definition, the first two parameters may bound deterministic errors, stochastic bias, or expected squared errors.
Because the relative noise terms scale with $\|\nabla \phi(x)\|^2$, they diminish as the algorithm approaches a stationary point. As we will see, these terms often have a relatively benign effect on optimization algorithms, provided their magnitude is sufficiently controlled. In contrast, absolute errors generally impose a limit on the accuracy that can be guaranteed. We illustrate the sources of all these types of error for both zeroth- and first-order oracles and show that they frequently arise in the literature.

\paragraph{Scope limitation}
Our detailed algorithmic analysis focuses primarily on methods based on zeroth- and first-order oracles. Although we also introduce models of inexact second-order information, illustrate them through examples, and review their use in optimization, we do not attempt a comprehensive analysis of second-order methods.

Furthermore, in each analysis we restrict the assumptions on the global properties of the objective function $\phi$ to a subset of commonly used conditions. In particular, we do not consider nonsmooth objective functions $\phi$. That said, it is important to note that, for a nonsmooth (and possibly discontinuous) function $f(x)$, this review provides useful analysis whenever such a poorly behaved function can be approximated by a smooth function $\phi(x)$ to some prescribed accuracy. In this case, our analysis establishes bounds on how well and how quickly $\phi(x)$ can be optimized, which, under suitable approximation conditions, can in turn imply guarantees for $f(x)$.

\paragraph{Organization and key results}
We begin with Section \ref{sec:assump}, where we list various assumptions that we use and refer to throughout the paper. These include smoothness conditions, convexity, strong convexity, and other structural properties that allow us to establish different convergence and complexity guarantees.

Section \ref{sec:deterministic} is dedicated to deterministic oracles and includes methods such as fixed-step gradient descent, gradient descent with line search, and trust-region methods. We present the oracle models and several examples in which these models apply readily and conveniently. We also establish connections between our first-order oracle model and several other commonly used models in the literature.
We then show how algorithms can be analyzed under the proposed oracle models. For gradient descent, we base our analysis on a fundamental descent inequality that provides a convenient and intuitive way to account for the influence of noise. In the case of convex $\phi(x)$, we also apply the alternative, well-known analysis based on telescoping sums. 

In Section \ref{sec:stochastic}, we extend the definitions and analysis to stochastic oracles with bounded second moments. The section includes alternative types of zeroth- and first-order stochastic oracles and multiple examples illustrating each type. Our main theoretical analysis focuses on fixed-step stochastic gradient descent under various oracle conditions and assumptions on $\phi(x)$. As in the deterministic case, a single general descent inequality allows us to derive many of the key convergence and complexity results for stochastic gradient descent.
 While we deliberately keep the algorithm simple, by considering only a fixed step-size stochastic gradient descent, we demonstrate that the choice of the step-size parameter is dictated by the oracle errors and show how final complexity and these errors are connected. The complexity results we derive are comparable to those in the literature; however, as we will see, we express them in terms of expected complexity, namely, bounds on the expected stopping time at which the algorithm reaches the desired accuracy. Importantly, we also state clear and explicit bounds on the best accuracy that can be guaranteed under the inexact oracle models. We further consider telescoping-sum analysis under convexity. 
We summarize the key results of each (sub)section at the beginning of the appropriate section.

We finish the paper with an overview of inexact second-order oracles and their use in the literature in Section \ref{sec:second_order}.

 \section{Assumptions}\label{sec:assump}

In this section, we state the main assumptions and introduce several preliminary concepts. We first present regularity assumptions on the objective function and its derivatives, followed by additional geometric conditions used to establish particular convergence and complexity results. The first two assumptions hold throughout the manuscript, whereas the remaining assumptions are invoked only when needed.

\begin{assumption}[\textbf{Lower boundedness of $\pmb{\phi}$}] \label{assum:low_bound}
The function $\phi$ is bounded below on $\mathbb{R}^n$, i.e., there exists a scalar $\hat{\phi}\in\mathbb{R}$ such that $\phi(x)\geq\hat{\phi}$ for all $x\in\mathbb{R}^n$.
\end{assumption}

\begin{assumption}[\textbf{Lipschitz continuity of the gradient of $\phi$}] \label{assum:lip_cont}
The function $\phi$ is continuously differentiable on $\mathbb{R}^n$, and its gradient is $L_1$-Lipschitz continuous for some constant $L_1>0$, i.e., $\|\nabla\phi(y)-\nabla\phi(x)\|\leq L_1\|y-x\|$ for all $x,y\in \mathbb{R}^n\times\mathbb{R}^n$.
\end{assumption}

\begin{assumption}[\textbf{Lipschitz continuity of the Hessian of $\phi$}] \label{assum:lip_cont2}
The function $\phi$ is twice continuously differentiable on $\mathbb{R}^n$, and its Hessian is $L_2$-Lipschitz continuous for some constant $L_2>0$, i.e., $\|\nabla^2\phi(y) - \nabla^2\phi(x)\| \le L_2 \|y-x\|$ for all $x,y\in \mathbb{R}^n\times\mathbb{R}^n$.
\end{assumption}

The next four assumptions impose additional geometric conditions on the objective function that are used to establish particular convergence and complexity results.

\begin{assumption}[\textbf{Convexity and boundedness of the initial level set of $\phi$}] \label{assum:convex} The function $\phi$ is convex, i.e., 
\begin{align*}
	\phi(y)\geq \phi(x)+\nabla\phi(x)^T(y-x) 
\end{align*}
for all $x,y\in\mathbb{R}^n$, and attains its minimum value $\phi^* = \phi(x^*)$ on a set $\mathcal{X}^*$, where $\mathcal{X}^*:=\{x\in\mathbb{R}^n:\phi(x)=\phi^*\}$ denotes the set of minimizers. Moreover, for a given point $x_0 \in \mathbb{R}^n$ there exists a constant $D>0$ such that for some $x^*\in \mathcal{X}^*$,
\begin{align*} 
	\sup\left\{\|x-x^*\|:\phi(x)\leq\phi(x_0)\right\}\leq D. 
\end{align*}
\end{assumption}

\begin{assumption}[\textbf{Strong convexity of $\phi$}] \label{assum:str.convex}
The function $\phi$ is $\mu$-strongly convex for some constant $\mu>0$, i.e.,
\begin{align*} 
	\phi(y)\geq \phi(x)+\nabla\phi(x)^T(y-x)+\frac{\mu}{2}\|y-x\|^2
\end{align*} 
for all $x,y\in\mathbb{R}^n\times\mathbb{R}^n$. Consequently, $\phi$ has a unique global minimizer $x^*$, and we define $\phi^*:=\phi(x^*)$.
\end{assumption}

The next two assumptions introduce the Kurdyka-{\L}ojasiewicz (KL) and Polyak-{\L}ojasiewicz (PL) conditions in the particular forms used throughout this manuscript. These conditions serve as weaker alternatives to the corresponding convexity assumptions in parts of our analysis, as discussed below.

\begin{assumption}[\textbf{Linear case of the KL condition of $\phi$}] \label{assum:KLcond}
Let $\mathcal{S}\subseteq\mathbb{R}^n$, and let $\phi^*=\inf \phi(x)$. The function $\phi$ satisfies the linear case of the Kurdyka-{\L}ojasiewicz (KL) condition on $\mathcal{S}$ for some constant $c>0$, i.e.,
\begin{align*} 
	\|\nabla\phi(x)\|\geq c\left(\phi(x)-\phi^*\right)
\end{align*} 
for all $x\in\mathcal{S}$.
\end{assumption}

Assumption \ref{assum:KLcond} is a less common variation of a power-type Kurdyka--{\L}ojasiewicz condition. A commonly used exponent-based formulation is
\begin{align*} 
	\|\nabla\phi(x)\|\geq c\left(\phi(x)-\phi^*\right)^\theta,
\end{align*} 
where $\theta\in[0,1)$; see, e.g., \cite{attouch2013convergence,li2018calculus}.  The case $\theta=\frac{1}{2}$ is equivalent, up to a rescaling of the constant, to the Polyak-{\L}ojasiewicz (PL) condition introduced below. Assumption \ref{assum:KLcond} corresponds to the case $\theta=1$, which we refer to as the linear case because the gradient norm is bounded below by a linear function of the objective-function gap and extended case of Kurdyka--{\L}ojasiewicz condition, because $\theta=1$ is not included in the original Kurdyka--{\L}ojasiewicz condition. While less common than $\theta\in[0,1)$, the linear case, nevertheless, arises naturally in the analysis of first-order methods as a consequence of convexity and boundedness of the initial level set. Specifically, Assumption \ref{assum:convex} implies Assumption \ref{assum:KLcond} on the set $\mathcal{S}:=\{x\in\mathbb{R}^n:\phi(x)\leq\phi(x_0)\}$ with $c=D^{-1}$. Indeed, for any $x\in\mathcal{S}$, convexity of $\phi$ and the definition of $D$ yield
\begin{align}\label{eq:convex_to_KL}
\phi(x)-\phi^*
\leq \nabla\phi(x)^T(x-x^*)
\leq \|\nabla\phi(x)\|\|x-x^*\|
\leq D\|\nabla\phi(x)\|.
\end{align}
Thus, Assumption \ref{assum:KLcond} provides a convenient way to analyze first-order methods under a condition that is weaker than convexity. The reverse implication does not hold, since a function satisfying Assumption \ref{assum:KLcond} need not be convex.

The case $\theta=1$ is studied in detail in \cite{fatkhullin2022sharp}, which provides several examples of functions satisfying Assumption \ref{assum:KLcond} globally, i.e., with $\mathcal{S}=\mathbb{R}^n$. One such example is the convex function $\phi(x)=\frac{e^x+e^{-x}}{2}-1$, which satisfies Assumption \ref{assum:KLcond} with $c=1$. The authors also construct nonconvex variants of this example that satisfy the condition with smaller values of $c$. Although some of these examples are primarily illustrative and somewhat artificial, \cite{fatkhullin2022sharp} also discusses applications in which the linear case arises, including policy optimization in reinforcement learning and problems in operations management.

More generally, the strength of Assumption~\ref{assum:KLcond}
depends on the set $\mathcal{S}$ and the corresponding constant
$c$. To illustrate this dependence, consider the convex function
$\phi(x)=x^4$, whose minimum is attained at $x^*=0$ with
$\phi^*=0$. For any $x_0\neq 0$, the corresponding level set is
$\mathcal{S}=[-|x_0|,|x_0|]$, and, for all $x\in\mathcal{S}$,
\begin{align*}
	\frac{|x_0|}{4}|\phi'(x)|=|x_0||x|^3 \geq |x|^4=\phi(x)-\phi^*=x^4,
\end{align*}
Thus, Assumption~\ref{assum:KLcond} holds on $\mathcal{S}$
with $c=\frac{4}{|x_0|}$, which can be arbitrarily small
as the initial level set grows larger.


Assumption~\ref{assum:KLcond} is particularly useful for analyzing gradient-descent-type methods in the presence of inexact information. For general smooth objectives, convergence to stationarity can be analyzed using a descent inequality that holds at each iteration. Under strong convexity, a similar approach yields convergence rates in terms of the function value gap. In contrast, standard convergence analyses of gradient descent under convexity typically rely on telescoping-sum arguments involving distances to a minimizer. Such arguments can be sensitive to the presence of noise and may require additional boundedness assumptions, even in relatively simple settings. In particular, convexity together with boundedness of the relevant level set implies a condition of the form given in Assumption~\ref{assum:KLcond}. Under Assumption~\ref{assum:KLcond}, the descent inequality can instead be used directly to obtain useful complexity bounds through relatively simple arguments that remain applicable in the presence of gradient errors. However, these bounds depend on the constant $c$ and must therefore be interpreted relative to the set $\mathcal{S}$ on which the assumption holds. Both approaches have been used in the deterministic setting and, as we demonstrate in this manuscript, yield comparable complexity bounds. Interestingly, in the stochastic setting, we observe a gap in the worst-case complexity bounds: the expected complexity bound obtained under Assumption~\ref{assum:KLcond} is worse than the corresponding bound obtained under convexity. This observation is consistent with existing results in the literature and reflects, at least in part, the differences between the underlying assumptions and analytical approaches. On the other hand, telescoping-sum analyses under convexity can be sensitive to the presence of relative bias, whereas analyses based on Assumption~\ref{assum:KLcond} can accommodate such errors. These complementary properties motivate our consideration of both approaches throughout this manuscript.

\begin{assumption}[\textbf{PL condition of $\phi$}]
\label{assum:PLcond}
Let $\mathcal{S}\subseteq\mathbb{R}^n$, and let $\phi^*=\inf_{x\in\mathbb{R}^n}\phi(x)$. The function $\phi$ satisfies the Polyak-{\L}ojasiewicz (PL) condition on $\mathcal{S}$ if there exists a constant $\mu>0$ such that
\begin{align*}
    \frac{1}{2}\|\nabla\phi(x)\|^2\geq \mu\left(\phi(x)-\phi^*\right)
\end{align*}
for all $x\in\mathcal{S}$.
\end{assumption}

The PL condition, introduced in \cite[Condition C]{polyak1963gradient}, is a well-known alternative to strong convexity that is sufficient to establish linear convergence of first-order methods under appropriate smoothness conditions. In particular, Assumption~\ref{assum:str.convex} implies Assumption~\ref{assum:PLcond} on $\mathcal{S}=\mathbb{R}^n$ with the same constant $\mu$, whereas the converse does not hold. Consequently, the complexity results established in this manuscript under Assumption~\ref{assum:PLcond} also apply to strongly convex functions.

%


\paragraph{Notation and conventions}
A list of notation and definitions used throughout the manuscript is provided in Appendix~\ref{app:notation}. We also adopt the following convention: for any expression of the form $\frac{a}{b}$, where $a>0$ and $b$ may equal zero, we interpret $\frac{a}{0}$ as infinity. This convention allows us to avoid treating degenerate cases separately while preserving the validity of our theoretical results.

\section{Deterministic Inexact Oracles}\label{sec:deterministic}

In this section, we focus on zeroth- and first-order oracles with deterministic errors and on optimization algorithms that employ these oracles. For each oracle type, we first introduce a general definition and then present several representative examples. We subsequently study how these oracles are employed within standard optimization frameworks, including fixed-step gradient descent, line-search methods, and trust-region methods, and establish convergence and iteration complexity guarantees.

\subsection{Deterministic Inexact Zeroth-Order Oracles}

We begin with deterministic inexact zeroth-order oracles, i.e., oracles that, given $x \in \mathbb{R}^n$, return an estimate of the objective function value. We present the definition of a deterministic inexact zeroth-order oracle and then provide examples.

\begin{customoracle}{D-0}[Deterministic Inexact Zeroth-Order Oracle] \label{def:oracle_det_0}
Given a point $x\in \mathbb{R}^n$ and  a control parameter $\lambda$  the oracle generates output $f(x)$ that satisfies
\begin{equation} \label{eq:0_oracle_det}
	| f(x) - \phi(x)| \leq \epsilon_f(\lambda) 
\end{equation} 
where $\epsilon_{f}(\lambda)$ is some nonnegative deterministic function of $\lambda$. The computational cost of this oracle is some function ${\cal C}_{f}(\lambda)$. 
\end{customoracle}

\begin{figure}[h!]
\centering
\includegraphics[width=0.475\linewidth]{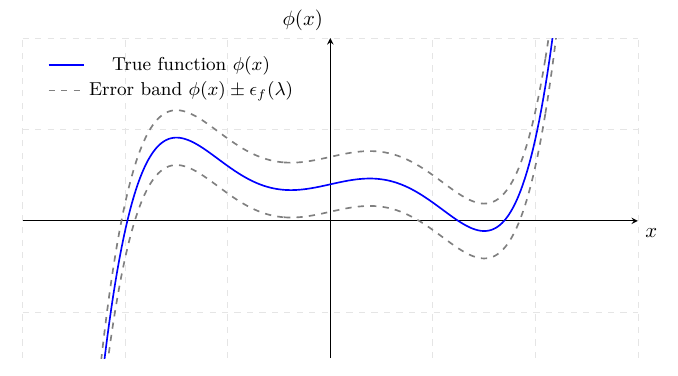}
\includegraphics[width=0.475\linewidth,trim={100 400 100 100},clip]{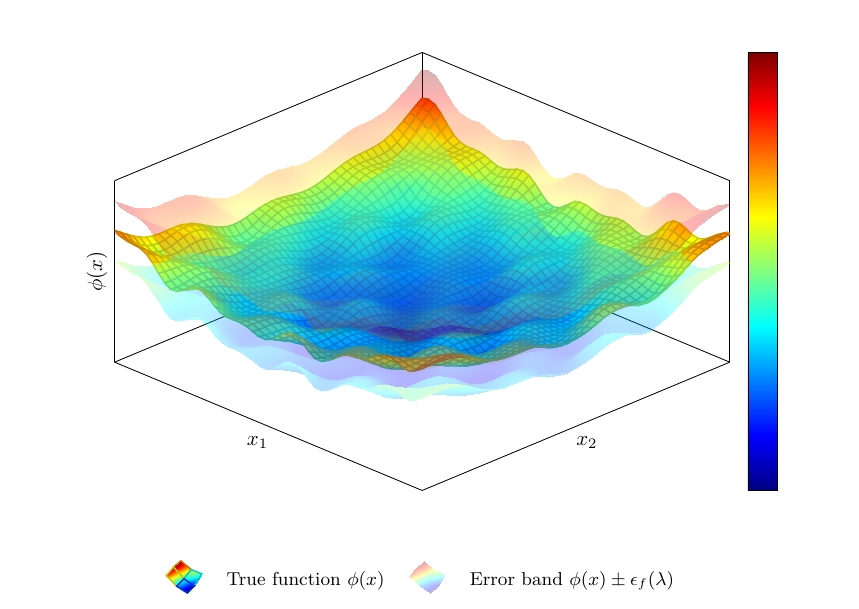}
\caption{Illustration of the error bands for deterministic inexact zeroth-order oracles in one and two dimensions}
\label{fig:error-band}
\end{figure}

Figure~\ref{fig:error-band} shows a one-dimensional and a two-dimensional illustration of a function $\phi$ and the associated error bands that satisfy the condition of Oracle~\ref{def:oracle_det_0} for some particular fixed value $\epsilon_f(\lambda)  = \epsilon_f$.

In what follows, we will often omit the dependence of $\epsilon_f$ on $\lambda$, except when this dependence must be made explicit, as this dependence will be easily understood from the context of  specific examples. 

\subsection{Examples of Deterministic Inexact Zeroth-Order Oracles}

Objective functions in many scientific and engineering applications cannot be evaluated exactly because each evaluation may require solving a system of differential equations \cite{ziems2011adaptive,kouri2014inexact}, solving an embedded system of equations to high accuracy \cite{cartis2022evaluation}, or performing calculations that are available only through numerical discretization \cite{kouri2014inexact}. In such settings, the control parameter $\lambda$ may represent, for example, a mesh size, discretization level, or iterative-solver tolerance. Increasing the accuracy of the underlying computation generally decreases the oracle error $\epsilon_f(\lambda)$ but increases the computational cost ${\cal C}_f(\lambda)$. Depending on the application, $\epsilon_f(\lambda)$ may approach zero as computational fidelity increases or remain bounded away from zero due to irreducible approximation errors that cannot be eliminated by increasing computational effort.

\subsubsection{Inexact function evaluations due to discretization of continuous objects} 
In many applications, the objective function is defined in terms of continuous objects that cannot be evaluated or represented exactly. Consequently, function evaluations rely on discretized approximations of these objects, introducing errors whose magnitude depends on the discretization accuracy. We illustrate this setting with three examples.

\begin{figure}[h!]
\centering
\includegraphics[width=0.375\linewidth]{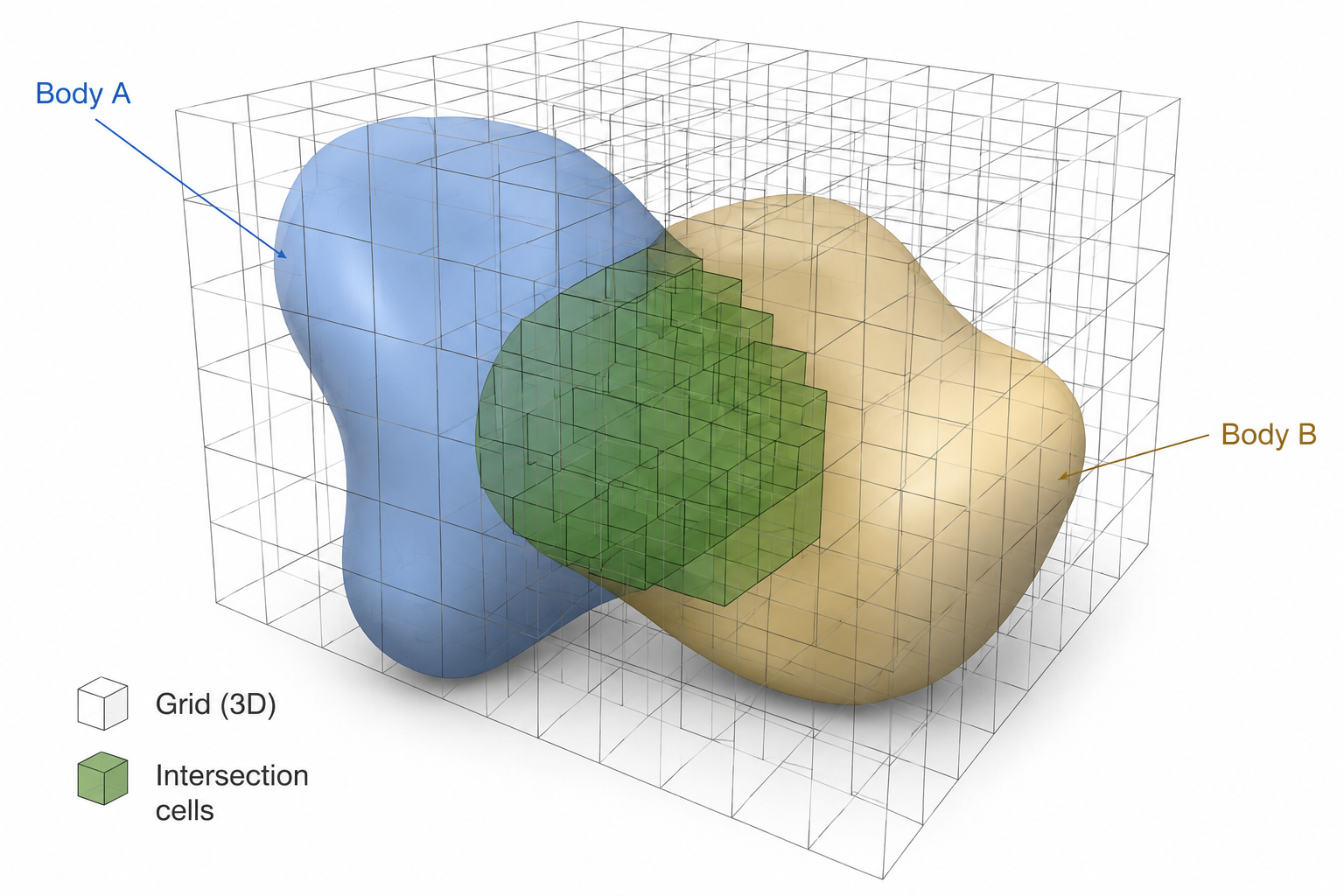}
\includegraphics[width=0.375\linewidth]{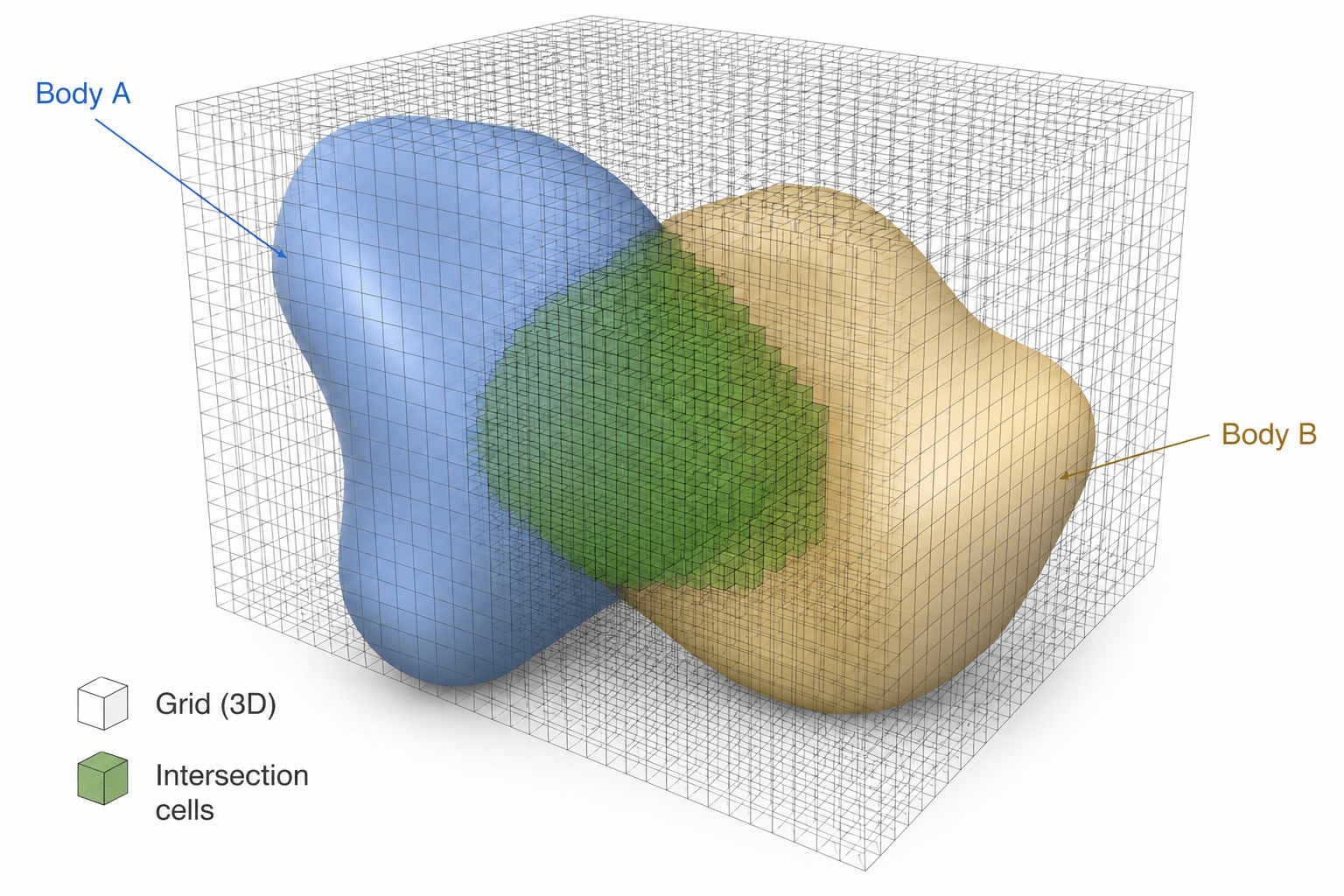}
\caption{Estimating volume of intersection via discretization.}
\label{fig:volume_intersect}
\end{figure}

\paragraph{Volume computation by discretization} 
Consider the problem of maximizing the volume of the intersection of two three-dimensional shapes over possible rotations and translations. The intersection volume can be approximated by discretizing a bounded region of space into cubes, classifying each voxel as inside or outside the intersection, and summing the corresponding voxel volumes \cite{klette2004digital,lachaud2017robust}. Here, the control parameter $\lambda$ may denote the side length of the voxels, with smaller values of $\lambda$ corresponding to finer discretization (see Figure \ref{fig:volume_intersect}). The parameter may be scalar for a uniform grid or vector-valued when different grid resolutions are used along the coordinate directions. For sufficiently regular shapes, only voxels intersecting a boundary layer of thickness proportional to $\lambda$ contribute to the approximation error, and therefore one may expect an error bound of the form $\epsilon_f(\lambda)=O(\lambda)$. On the other hand, examining all voxels in a fixed three-dimensional bounding region requires $O(\lambda^{-3})$ voxel evaluations, and hence ${\cal C}_f(\lambda)=O(\lambda^{-3})$ for a simple full-grid implementation. The precise dependencies may differ for more sophisticated volume-estimation procedures.

\paragraph{PDE-constrained optimization} A common example of inexact function evaluations arising from discretization occurs in PDE-constrained optimization, where computing the objective function requires solving a partial differential equation \cite{ziems2011adaptive,kouri2014inexact}. Suppose that $\phi(x)=J(x,u(x))$, where $u(x)$ denotes the exact solution of a PDE depending on $x$. In practice, the PDE is solved numerically using a discretization with resolution parameter $\lambda$, producing an approximate solution $u_\lambda(x)$ and the oracle output $f(x)=J(x,u_\lambda(x))$. If the discretization satisfies the uniform error estimate $|f(x)-\phi(x)|=|J(x,u_\lambda(x))-J(x,u(x))|\leq C\lambda^p$, then $\epsilon_f(\lambda)=C\lambda^p$, while ${\cal C}_f(\lambda)$ is the cost of assembling and solving the discretized PDE and typically increases as the discretization is refined.

\paragraph{Piecewise constant approximation}
A simple example of inexact function evaluations due to discretization
arises when a continuous objective function is approximated by a
piecewise constant function constructed from its values on a discrete
grid \cite{deBoor2001}.
Let $\phi:\mathbb R\to\mathbb R$ be continuously differentiable
and suppose that its gradient is $L_1$-Lipschitz continuous
(Assumption~\ref{assum:lip_cont}) and that
$|\phi'(x)|\leq G$ for all $x\in\mathbb R$.
For a given $\rho>0$, partition $\mathbb R$ into intervals
$I_k=[k\rho,(k+1)\rho)$, $k\in\mathbb Z$, and define
the piecewise constant function
\[
f(x)=\frac{\phi(k\rho)+\phi((k+1)\rho)}{2},
\qquad x\in I_k.
\]
By the mean value theorem,
\[
|f(x)-\phi(x)|
\le \frac{G}{2}
\bigl(|x-k\rho|+|x-(k+1)\rho|\bigr)
=\frac{G\rho}{2}.
\]
Thus, it follows that,
\[
\sup_{x\in\mathbb R}|f(x)-\phi(x)|
\le \frac{G\rho}{2}.
\]
This example illustrates settings in which a continuous function
can only be evaluated at discrete grid points. If the grid spacing
is controllable, it can be treated as an input to the oracle,
with $f(x)$ representing the oracle output. The computational
cost of each oracle call is that of evaluating $\phi$ at the two
grid points bounding $x$.

\subsubsection{Inexact evaluations arising from approximate subproblem solvers}

In some applications, evaluating the objective function requires
solving an auxiliary numerical problem, whose approximate solution
introduces errors in the function evaluations.

\paragraph{Oracles using iterative solvers}
For example, evaluating the objective function may require solving
a linear system $A(x)u(x)=b(x)$, as is common when the function
evaluation contains an embedded iterative numerical process
\cite{dembo1982inexact,cartis2022evaluation}. Instead of solving the system exactly, an iterative method may be terminated once $\|A(x)\widetilde{u}(x;\lambda)-b(x)\|\leq\lambda$, and the oracle returns $f(x)=J(x,\widetilde{u}(x;\lambda))$. If $J$ is uniformly Lipschitz continuous in its second argument with constant $L_J$ and $\|A(x)^{-1}\|\leq\kappa_A$ over the region of interest, then $|f(x)-\phi(x)|\leq L_J\kappa_A \lambda$. Thus, $\epsilon_f(\lambda)=L_J\kappa_A\lambda$, while ${\cal C}_f(\lambda)$ is the cost of the iterative linear solver and generally increases as the residual tolerance $\lambda$ decreases.

\subsubsection{Smoothing of nonsmooth and discontinuous functions}

Another large class of examples arises when $f(x)$ is a nonsmooth function that can be approximated by a smooth $\phi(x)$  \cite{Nesterov05,beck2012smoothing}. 
These examples are of a different nature from the preceding ones in that $\phi(x)$ is not the objective function, but itself is a convenient surrogate to which the theoretical
results developed in this paper can be applied. Under suitable approximation conditions, the resulting convergence and complexity guarantees can be used to draw conclusions about the original objective function $f(x)$.

The parameter $\lambda$ controls how closely $\phi_\lambda(x)$ approximates $f(x)$ and typically also affects the smoothness properties of $\phi_\lambda$. In this setting, $\lambda$ can still be viewed as an oracle input, but its role differs from that in the preceding examples: rather than controlling the accuracy of an approximation to a fixed objective function, it determines the smooth surrogate itself. If $f(x)$ is viewed as the zeroth-order oracle output for $\phi_\lambda(x)$, the corresponding oracle error is $|f(x)-\phi_\lambda(x)|$, and the cost of each oracle call is simply the cost of evaluating $f(x)$.

\paragraph{Smoothing a maximum of smooth functions}
The first simple illustrative example is where $f(x)$ is a maximum of several linear functions, which includes 
many functions in modern machine learning, such as ReLU, but also other more general objectives when the worst output among multiple objectives is to be minimized. 
Let
\[
    f(x)=\max_{1\leq i\leq m} f_i(x),
\]
where $f_i(x)=a_i^Tx+b_i$. Then $f$ is generally nonsmooth at points where the maximizing index changes.

One can define a  smoothing of $f$ as 
\[
    \phi_\lambda(x)
    =\lambda\log\left(
        \sum_{i=1}^m e^{f_i(x)/\lambda}
    \right),
    \qquad \lambda>0,
\]
which satisfies 
\[
    f(x)
    \leq \phi_\lambda(x)     \leq f(x)+\lambda{\log m},
    \qquad \forall x\in\mathbb{R}^n.
\]
This $\phi_\lambda$ is $L_1$-smooth, with $L_1= \frac{\max_i \|a_i\|^2}{\lambda}$. 
Hence, $f(x)$ can be seen as the inexact zeroth-order oracle for $\phi_\lambda(x)$ with
$\epsilon_f(\lambda)=\lambda\log m$. This example is purely illustrative, since 
$\phi_\lambda(x)$ itself is computable and can be optimized without the use of inexact oracles. Our next example is 
more compelling. 

\paragraph{Moreau envelope} 
The Moreau envelope \cite{moreau1962fonctions} of a convex,
possibly nonsmooth, $L_0$-Lipschitz continuous function
$f:\mathbb{R}^n\to\mathbb{R}$ is defined as
\[
    \phi_\lambda(x)
    =\min_{y\in\mathbb{R}^n}
    \left\{f(y)+\frac{1}{2\lambda}\|y-x\|^2\right\},
\]
where $\lambda>0$. Then $\phi_\lambda(x)$ is $L_1$-smooth with
$L_1=\frac{1}{\lambda}$ and satisfies
\[
0\leq f(x)-\phi_\lambda(x)\leq \frac{\lambda L_0^2}{2}.
\]
Thus we can view $f(x)$ as the inexact first order oracle of a smooth function $   \phi_\lambda(x)$ with 
$\epsilon_f=\frac{\lambda L_0^2}{2}$. 

In \cite{DavisDrusvyatskiy2019}  Moreau envelopes for nonconvex, nonsmooth functions are considered and 
$\| \nabla \phi_\lambda(x)\|$ is explicitly used as the measure of stationarity, thus essentially they consider optimizing 
$\phi_\lambda(x)$ instead of $f(x)$. The nonconvex case is more complicated than the convex case  presented here, because 
the related smoothness of $\phi_\lambda(x)$ and expression for $\epsilon_f$ depend on the weak-convexity constant as well as $L_0$. Nevertheless that setting too can fit into our inexact oracle framework.

\paragraph{Gaussian smoothing} 
Let $f:\mathbb{R}^n\to\mathbb{R}$ be an $L_0$-Lipschitz continuous function that is not necessarily smooth. For $\lambda>0$, the Gaussian smoothing of $f$ \cite{duchi2012randomized,nesterov2017random} is defined as
\[
\phi_\lambda(x)
=\mathbb{E}[f(x+\lambda Z)],
\qquad Z\sim\mathcal{N}(0,I_n).
\]
The function $\phi_\lambda$ is infinitely differentiable, with gradient given by
\[
\nabla\phi_\lambda(x) =\frac{1}{\lambda}\mathbb{E}[f(x+\lambda Z)Z].
\]
Moreover, its gradient is $L_1$-Lipschitz continuous with $L_1=\frac{L_0}{\lambda}$, and the approximation error satisfies
\[
|f(x)-\phi_\lambda(x)| \leq \sqrt{n}\lambda L_0.
\]
Thus $\phi(x)$ is a smooth function and $f(x)$ can be seen as its  approximation\footnote{In fact $ \nabla\phi(x)$ is $\sqrt{\frac{2}{\pi}}\frac{L}{\lambda}$-Lipschitz and $\sqrt{n}$ in the error bound improves to $\sqrt{\frac{2}{\pi}}$ when $n=1$ but such minor constant changes are omitted for brevity.}

Now consider the piecewise constant approximation example discussed above,  but this time consider 
$\phi(x)$ that is not smooth, but $L_0$-Lipschitz. 
Let us define Gaussian smoothings of $f$ and $\phi$
\[
f_{\lambda}(x)
=\mathbb E[f(x+{\lambda}Z)],
\qquad \phi_{\lambda}(x)
=\mathbb E[\phi(x+{\lambda}Z)], \qquad Z\sim {\cal N}(0,1),\quad {\lambda}>0.
\]
From our previous derivations we have the following bound. 
\[
\begin{aligned}
|f_{\lambda}(x)-f(x)|& \leq |f_{\lambda}(x)-\phi_h(x)|+ |\phi_{\lambda}(x)-\phi(x)|+|\phi(x)-f(x)|
\le {L_0\rho}+L_0{\lambda}
\end{aligned}
\]

$f_{\lambda}(x)$ is infinitely smooth and we next bound $|f_h^{\prime \prime}(x)|$. We already established that
$|\phi_{\lambda}^{\prime \prime}(x)|\leq \frac{L}{\lambda}$. To bound $|f_h^{\prime \prime}(x)-\phi_{\lambda}^{\prime \prime}(x)|$
we observe that 
\[
f_{\lambda}^{\prime \prime}(x)-\phi_{\lambda}^{\prime \prime}(x)
=\frac{1}{{\lambda}^2} \mathbb E[(f(x+{\lambda}Z)-\phi(x+_{\lambda}Z))(Z^2-1)].
\]
Consequently,
\[|f_{\lambda}^{\prime \prime}(x)-\phi_{\lambda}^{\prime \prime}(x)|
\le \frac{L_0\rho}{2{\lambda}^2}
\mathbb E|Z^2-1|\leq \frac{L_0\rho}{{\lambda}^2}.
\]
Combining the bounds gives
\[
|f_{\lambda}^{\prime \prime}(x)|
\le
\frac{L_0}{{\lambda}}
+\frac{L_0\rho}{{\lambda}^2}.
\]
In particular, if $\rho={\lambda}$, then $f(x)$ can be seen as a noisy approximation with noise bounded by $2L_0\rho$ of a 
$\frac{2L_0}{\rho}$-smooth function $f_{\lambda}(x)$.

\subsection{Deterministic Inexact First-Order Oracles}

Next, we introduce deterministic inexact first-order oracles, i.e., oracles that given $x \in \mathbb{R}^n$ return an estimate of the gradient of the objective function. We present the definition of a deterministic inexact first-order oracle and list some of the key properties. 
We then describe some common examples where these oracles arise, and investigate the connections to other first-order oracles from the literature.

\begin{customoracle}{D-1}[Deterministic Inexact First-Order Oracle] \label{def:oracle_det_1}
Given a point $x\in \mathbb{R}^n$ and a control parameter $\lambda$  the oracle generates output  $g(x)$ that satisfies
\begin{equation} \label{eq:1_oracle_det}
	\|g(x) - \nabla \phi(x)\|^2 \leq \epsilon_g(\lambda)^2+\kappa_g(\lambda)^2\|\nabla \phi(x)\|^2,   
\end{equation} 
where $\epsilon_g(\lambda)$ and $\kappa_g(\lambda)$ are some nonnegative deterministic functions of $\lambda$. The computational cost of this oracle is some function ${\cal C}_{g}(\lambda)$. 
\end{customoracle}

Unlike the zeroth-order oracle, the  first-order oracle is characterized by two error terms  $\epsilon_g(\lambda)$ and $\kappa_g(\lambda)^2\|\nabla \phi(x)\|^2$ denoting the \emph{absolute error} and \emph{relative error}, respectively. These terms are also often referred to as the \emph{additive} and \emph{multiplicative error components}. 
The error constants $\epsilon_g(\lambda)$ and $\kappa_g(\lambda)$ are often, but not always, functions of $\lambda$, and, thus, can be controlled to some degree by the appropriate choice of the oracle input. We will motivate this dependence when presenting examples. 

 The key reason for identifying the two error terms separately  is that each terms has distinct influence  on the behavior of optimization algorithms which utilize such an oracle. Specifically,  since we are interested in unconstrained optimization,  we aim at driving $\|\nabla \phi(x)\|$  to zero as the algorithm progresses, thereby affecting absolute and relative error terms differently. We illustrate this distinction through analyses of basic optimization algorithms later in this section and throughout the survey.
 
\begin{figure*}[h!]
    \centering
    \begin{subfigure}[t]{0.24\textwidth}
        \centering
        \includegraphics[width=1\linewidth,trim={10 0 0 0},clip]{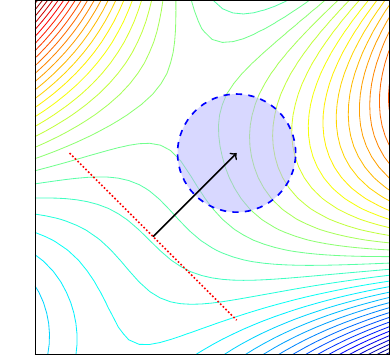}
        \caption{$\epsilon_g=0$, $\kappa_g=0.5$}
    \end{subfigure}%
    ~ 
     \begin{subfigure}[t]{0.24\textwidth}
        \centering
        \includegraphics[width=1\linewidth,trim={10 0 0 0},clip]{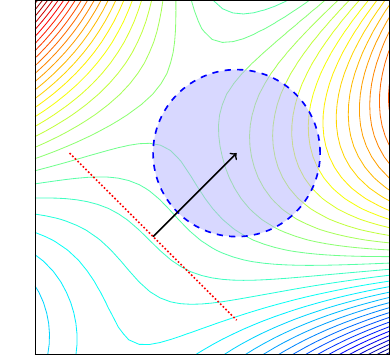}
        \caption{$\epsilon_g=0.5$, $\kappa_g=0.5$}
    \end{subfigure}%
    ~
     \begin{subfigure}[t]{0.24\textwidth}
        \centering
        \includegraphics[width=1\linewidth,trim={10 0 0 0},clip]{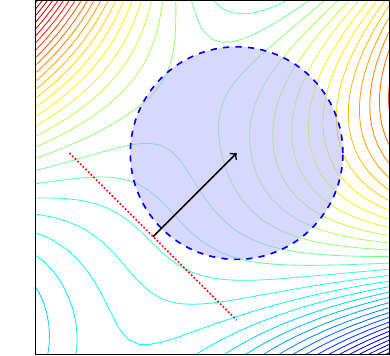}
        \caption{$\epsilon_g=0$, $\kappa_g=0.9$}
    \end{subfigure}%
    ~
     \begin{subfigure}[t]{0.24\textwidth}
        \centering
        \includegraphics[width=1\linewidth,trim={10 0 0 0},clip]{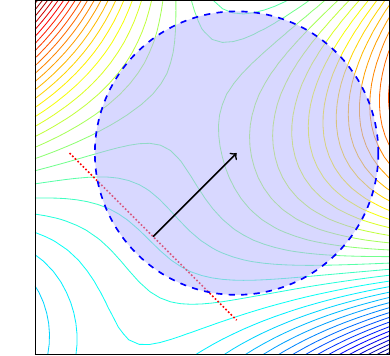}
        \caption{$\epsilon_g=0$, $\kappa_g=1.2$}
    \end{subfigure}%
\caption{Illustration of Oracle~\ref{def:oracle_det_1}. Oracle~\ref{def:oracle_det_1} requires that 
$\|g(x) - \nabla \phi(x)\|^2 \leq \epsilon_g^2 + \kappa_g^2 \|\nabla \phi(x)\|^2$. 
In the figure, the black arrow represents the true gradient $\nabla \phi(x)$. The red dashed line is orthogonal to the true gradient and indicates separation of all $g(x)$ for which $-g(x)$ is a descent direction. 
The shaded blue region in each subplot depicts the set of vectors $g(x)$ that satisfy the oracle condition.}
\label{fig:first_order_oracle}
\end{figure*}

 Figure~\ref{fig:first_order_oracle} is an illustration of Oracle~\ref{def:oracle_det_1} for different values of the oracle parameters ($\epsilon_g$ and $\kappa_g$). What is interesting to note is the relation between the oracle parameters and whether the set of vectors $g(x)$ that satisfy the oracle condition are descent directions. 
 
As with the zeroth-order oracle, we  will omit the dependence on $\lambda$ in the error functions $\epsilon_g(\lambda)$ and $\kappa_g(\lambda)$, except when this dependence must be made explicit. 

\subsubsection{Properties of Deterministic Inexact First-Order Oracles}

We note that \eqref{eq:1_oracle_det} implies (is stronger than)  the non-squared variant, 
\begin{equation} \label{eq:1_oracle_det_nosqr}
	\|g(x) - \nabla \phi(x)\| \leq \epsilon_g+\kappa_g\|\nabla \phi(x)\|.   
\end{equation} 
When $\epsilon_g=0$, condition \eqref{eq:1_oracle_det_nosqr} is well-known in the literature as the {\em norm condition} \cite{carter1991global} and it implies the conditions of  Oracle \ref{def:oracle_det_1} with $\epsilon_g=0$. More generally, however, \eqref{eq:1_oracle_det_nosqr} implies Oracle \ref{def:oracle_det_1} with an increased right-hand side. Specifically, given any $g(x)$ that satisfies 
\begin{equation*}
	\|g(x) - \nabla \phi(x)\| \leq \epsilon^\prime_g+\kappa^\prime_g\|\nabla \phi(x)\|, 
\end{equation*}
it is easy to see that  $g(x)$  satisfies \eqref{eq:1_oracle_det} with $\epsilon_g^2=2{\epsilon_g^\prime}^2$ and  $\kappa_g^2=2{\kappa_g^\prime}^2$. Condition  of Oracle \ref{def:oracle_det_1} better fits the theoretical analysis of the algorithms  and extends naturally to the stochastic oracles setting. We will see that for many examples the condition  of Oracle \ref{def:oracle_det_1} is natural. 

From  \eqref{eq:1_oracle_det_nosqr}, it follows that 
\begin{align} \label{eq:1_normg_det}
	(1-\kappa_g)\| \nabla \phi(x)\| -\epsilon_g\leq \|g(x)\|\leq (1+\kappa_g) \|\nabla \phi(x)\| +\epsilon_g,
\end{align}
which gives a meaningful and nontrivial lower bound on $\|g(x)\|$ when 
\begin{align} \label{eq:1_normg_conditions1}
	\kappa_g < 1 \quad \text{and} \quad \| \nabla \phi (x)\|^2 \geq \frac{\epsilon_g^2}{(1 - \kappa_g)^2}.
\end{align}
Similarly, from  \eqref{eq:1_oracle_det}, it follows that
\begin{align} \label{eq:1_normg_det2}
    \|\nabla\phi(x)\| - \sqrt{\epsilon_g^2 + \kappa_g^2 \|\nabla\phi(x)\|^2} 
    \leq \|g(x)\|
    \leq \|\nabla\phi(x)\| + \sqrt{\epsilon_g^2 + \kappa_g^2 \|\nabla\phi(x)\|^2}, 
\end{align}
which provides a meaningful and nontrivial lower bound on $\|g(x)\|$ when 
\begin{align} \label{eq:1_normg_conditions2}
	\kappa_g < 1 \quad \text{and} \quad \| \nabla \phi (x)\|^2 \geq \frac{\epsilon_g^2}{1 - \kappa_g^2}.
\end{align}

The conditions in  \eqref{eq:1_normg_conditions1} and \eqref{eq:1_normg_conditions2} will play a key role in the convergence analysis of algorithms that utilize Oracle \ref{def:oracle_det_1}. We will demonstrate that the restriction on $\kappa_g$ (that appears in both \eqref{eq:1_normg_conditions1} and \eqref{eq:1_normg_conditions2}) after a suitable scaling, is equivalent to the oracle obtaining a descent direction, a necessary condition for a deterministic algorithm to work. We observe that the lower bound on the squared norm of the gradient is always smaller in \eqref{eq:1_normg_conditions2} than in \eqref{eq:1_normg_conditions1}, i.e., $\frac{\epsilon_g^2}{(1 - \kappa_g)^2} \geq \frac{\epsilon_g^2}{1 - \kappa_g^2}$, 
with the lower bounds only being equal for $\kappa_g = 0$. This is not surprising, as \eqref{eq:1_normg_conditions2} is derived using the stronger squared condition.

The effect of the error term $\epsilon_g$ is different from that of   $\kappa_g<1$, and essentially results in a bound on the best (smallest) value of $\|\nabla \phi(x)\|$ for which an optimization algorithm is guaranteed to make progress.  In fact, by imposing such a lower bound we can essentially remove $\epsilon_g$ from any further consideration. Specifically, we make use of a constant $\tau\in [0,1)$, which we can choose and for which we assume that 
\begin{equation}\label{eq:nablaphi_lower_bound}
    \|\nabla \phi(x)\|^2\geq \frac{\epsilon_g^2}{(1-\kappa_g^2)\tau^2}
\end{equation}
holds. Here and henceforth we will allow $\tau=0$ when $\epsilon_g=0$ with the convention that  the right hand side of \eqref{eq:nablaphi_lower_bound} is then $0$.

This gives a convenient condition on the smallest gradient norm for which we can expect ``good" properties of the gradient approximation $g(x)$. 
From \eqref{eq:nablaphi_lower_bound} it follows that
\begin{equation} \label{eq:1_oracle_det_sqr_noepsg}
	\|g(x) - \nabla \phi(x)\|^2 \leq (\tau^2+(1-\tau^2)\kappa_g^2)\|\nabla \phi(x)\|^2.  
\end{equation} 
Thus, as long as \eqref{eq:nablaphi_lower_bound} is satisfied, $g(x)$ satisfies Oracle~\ref{def:oracle_det_1} with the given $\epsilon_g$ and $\kappa_g$ constants, and also 
\begin{equation}\label{eq:kappa_gt}
	\kappa_{g\tau}^2=\tau^2+(1-\tau^2)\kappa_g^2	
\end{equation}
and no absolute error term. Moreover, $\kappa_g<1$ implies $\kappa_{g\tau}<1$. The conclusion from this is that we can treat the $\epsilon_g$ term in Oracle~\ref{def:oracle_det_1} as a lower bound 
 on $\|\nabla \phi(x)\|$ and ignore it otherwise. Furthermore,  \eqref{eq:1_oracle_det_sqr_noepsg} implies the well-known norm condition \cite{carter1991global}
 \begin{equation}
\label{eq:1_oracle_det_nosqr_noepsg}
	\|g(x) - \nabla \phi(x)\| \leq \kappa_{g\tau} \|\nabla \phi(x)\|.
\end{equation} 

When $\epsilon_g = 0$, then allowing  $\tau=0$ is meant to imply that  \eqref{eq:1_oracle_det_sqr_noepsg} recovers Oracle  \ref{def:oracle_det_1} with the given $\epsilon_g=0$ and the same $\kappa_g$  for arbitrarily small $\|\nabla\phi(x)\|$.

In our analysis we will keep $\epsilon_g$ explicitly present in all the derivations. Bounds such as \eqref{eq:nablaphi_lower_bound} will appear in the analysis explicitly but will vary depending on the assumed properties of $\phi(x)$.

\subsection{Examples of Deterministic Inexact First-Order Oracles}\label{sec.examples_first}

We next present several examples of deterministic inexact first-order oracles that arise in the optimization literature and fit the general framework of Oracle~\ref{def:oracle_det_1}. For each example, we specify the oracle inputs and outputs, characterize the corresponding absolute and relative errors, and discuss the associated computational cost. 

\subsubsection{Forward Finite Difference (FFD)}
Whenever direct computation of the gradient $\nabla \phi(x)$, using calculus, is not possible, the most well-known and straightforward approach to approximate the gradient is forward finite differences (FFD). 
Given access to Oracle~\ref{def:oracle_det_0}, for any $x \in \mathbb{R}^n$, the first-order FFD oracle is implemented as follows,
 \begin{equation} \label{eq:FFD}\tag{FFD}
	g(x) =  \sum_{i=1}^n \frac{f(x+{h}e_i) - f(x)}{{h}} e_i
\end{equation}
where $h>0$  is the finite difference interval (perturbation parameter) and $e_i$ denotes the $i$-th column of the identity matrix. 
By Assumption~\ref{assum:lip_cont}, for any $x \in \mathbb{R}^n$, there exists an upper bound on the absolute error \cite{berahas2022theoretical}, i.e.,  
\begin{equation} \label{eq:FFD error} 
    \|g({x}) - \nabla \phi({x})\| \leq \frac{\sqrt{n} L_1 {h}}{2} + \frac{2\sqrt{n}  \epsilon_f}{h}. 
\end{equation}
Thus, \eqref{eq:FFD} gives an example of Oracle~\ref{def:oracle_det_1} with input $\lambda=h$, the perturbation parameter, and absolute and relative errors given by, 
\begin{equation*}
	\epsilon_g = \frac{\sqrt{n} L_1 {h}}{2} + \frac{2\sqrt{n}  \epsilon_f}{{h}} \quad \text{and} \quad  \kappa_g = 0,
\end{equation*}  
respectively.

In the simplest case where  accuracy $\epsilon_f$ and  cost  $\mathcal{C}_f$ of Oracle~\ref{def:oracle_det_0} are fixed constants, the cost to execute  \eqref{eq:FFD} consists of the $n+1$ calls to Oracle~\ref{def:oracle_det_0} and the subsequent computation of $g$, that is, ${\cal C}_{g} = (n+1)\mathcal{C}_f + \mathcal{O}(n)$. 
Furthermore, if ${h}$ is chosen to satisfy  ${h}\geq 2\sqrt{\epsilon_f / L_1}$, then 
\begin{align*}
	\|g(x)-\nabla\phi(x)\|\leq \sqrt{n}L_1h,
\end{align*}
and hence the FFD approximation satisfies Oracle~\ref{def:oracle_det_1} with $\epsilon_g=\sqrt{n}L_1h$ and $\kappa_g=0$. The smallest value permitted by this condition is $h=2\sqrt{\epsilon_f/L_1}$, which yields $\epsilon_g=2\sqrt{nL_1\epsilon_f}$.

Alternatively, one can consider using $h$ as an input to Oracle~\ref{def:oracle_det_0}. For example, if the error bound $\epsilon_f$ of Oracle~\ref{def:oracle_det_0} can be made to depend on $h$ quadratically, i.e., $\epsilon_f(h)=ch^2$, for some constant $c>0$ with the cost ${\cal C}_{f}(h)$, then in this case, 
\begin{align*}
    \epsilon_g = \sqrt{n} h\left (\frac{L_1}{2}+2c\right ), \quad \kappa_g = 0, \quad\text{ and } \quad {\cal C}_{g}(h)=(n+1){\cal C}_{f}(h) + \mathcal{O}(n). 
\end{align*}

\subsubsection{Central Finite Difference (CFD)}
A more accurate alternative to the FFD method is the central finite difference (CFD) method, where the sample points are chosen antithetically, i.e., for any $x \in \mathbb{R}^n$
\begin{equation} \label{eq:CFD}\tag{CFD} 
    g(x) =  \sum_{i=1}^n \frac{f(x+{h} e_i) - f(x-{h} e_i)}{2{h}} e_i. 
\end{equation}
The first-order oracle obtained by \eqref{eq:CFD} also satisfies the error bound \eqref{eq:FFD error} under Assumption~\ref{assum:lip_cont}, but requires $2n$ calls to the zeroth-order oracle. 
However,  CFD enjoys a better error bound when the Hessian is assumed to be $L_2$-Lipschitz continuous \cite{berahas2022theoretical}, i.e., 
\begin{equation} \label{eq:CFD error}
	\|g({x}) - \nabla \phi({x})\| \leq \frac{\sqrt{n} L_2 {h}^2}{6} + \frac{\sqrt{n} \epsilon_f}{{h}}.
\end{equation}
Thus, under Assumption~\ref{assum:lip_cont2}, the CFD approach gives an example of Oracle~\ref{def:oracle_det_1} with input $\lambda=h$, and absolute and relative errors given by,  
\begin{align*}
    \epsilon_g(h) = \frac{\sqrt{n} L_2 h^2}{6} + \frac{\sqrt{n} \epsilon_{f}}{h}, 
    \quad \text{and} \quad \kappa_g(h)=0, 
\end{align*}
respectively. Here, again, if Oracle~\ref{def:oracle_det_0} is non-adjustable, using  ${h} \geq \sqrt[3]{6\epsilon_f / L_2}$ results in Oracle~\ref{def:oracle_det_1} with 
\begin{align*}
    \epsilon_g = \frac{\sqrt{n} L_2 h^2}{3} \geq \sqrt{n}\sqrt[3]{\frac{4L_2\epsilon_f^2}{3}}, \quad
    \kappa_g=0, \quad \text{and} \quad 
    \mathcal{C}_g = 2n\mathcal{C}_f + \mathcal{O}(n). 
\end{align*}
If the error $\epsilon_f$ of Oracle~\ref{def:oracle_det_0} can be controlled by $h$ as an input, and, specifically, if $\epsilon_f(h)=ch^3$, for some constant $c$ and some cost ${\cal C}_{f}(h)$, then 
\begin{align*}
    \epsilon_g = \sqrt{n} \left (\frac{L_2}{6}+c\right )h^2, \quad 
    \kappa_g = 0, \quad \text{and} \quad 
    {\cal C}_{g}(h) = 2n{\cal C}_{f}(h) + \mathcal{O}(n). 
\end{align*}

\subsubsection{Higher-Order Finite-Difference Approximations}

The FFD and CFD constructions considered above are special cases of more general finite-difference schemes that use additional function evaluations to obtain higher-order approximations. Let $q\geq 1$ denote the order of the approximation, and consider a stencil consisting of $m+1$ nodes $c_j\in\mathbb{R}$ and weights $a_j\in\mathbb{R}$, $j=0,\ldots,m$. The corresponding gradient approximation is
\begin{align}
g(x)
=
\sum_{i=1}^n
\left[
\frac{1}{h}
\sum_{j=0}^{m}
a_j
f\left(x+c_jhe_i\right)
\right]e_i.
\label{eq:higher_order_fd}
\end{align}
A $q$th-order stencil satisfies the moment conditions
\begin{align*}
\sum_{j=0}^{m}a_jc_j = 1, \qquad \text{and} \qquad \sum_{j=0}^{m}a_jc_j^r = 0,
\quad
r=0,2,\ldots,q.
\end{align*}
For example, the three-point one-sided approximation
\begin{align*}
\frac{
-3f(x)
+4f(x+he_i)
-f(x+2he_i)
}{2h}
\end{align*}
has order $q=2$, while the fourth-order centered approximation
\begin{align*}
\frac{
f(x-2he_i)
-8f(x-he_i)
+8f(x+he_i)
-f(x+2he_i)
}{12h}
\end{align*}
has order $q=4$.

Higher-order accuracy requires stronger smoothness and, in general, additional stencil points. For the standard one-sided $q$th-order formulas, $m=q$, so that $q+1$ stencil points are used. Standard centered formulas have even order $q$ and use $q$ nonzero stencil points, corresponding to $m=q-1$. Thus, increasing $q$ improves the truncation error, but also increases the number of oracle evaluations and changes the amplification of errors in the zeroth-order oracle.

Let
\begin{align*}
A_q
=
\sum_{j=0}^{m}|a_j| \qquad \text{and} \qquad
C_q
=
\frac{L_q}{(q+1)!}
\sum_{j=0}^{m}
|a_j||c_j|^{q+1},
\end{align*}
where $\nabla^q\phi$ is assumed to be $L_q$-Lipschitz continuous. The following proposition bounds the error in the q-th order finite difference approximations.

\begin{proposition}
\label{prop:higher_order_fd}
Suppose that $\nabla^q\phi$ is $L_q$-Lipschitz continuous and that the zeroth-order oracle evaluations used in \eqref{eq:higher_order_fd} satisfy Oracle~\ref{def:oracle_det_0} with error bound $\epsilon_f$. Then,
\begin{align*}
\|g(x)-\nabla\phi(x)\|
\leq
\sqrt{n}
\left(
C_qh^q
+
\frac{A_q\epsilon_f}{h}
\right).
\end{align*}
\end{proposition}

The result naturally defines an output of Oracle~\ref{def:oracle_det_1}. In the simplest case where Oracle~\ref{def:oracle_det_0} is not adjustable, we take
\begin{align*}
\lambda
=
\left(
h,
q,
\{(a_j,c_j)\}_{j=0}^{m}
\right),
\end{align*}
and the approximation satisfies Oracle~\ref{def:oracle_det_1} with
\begin{align*}
\epsilon_g(\lambda)
=
\sqrt{n}
\left(
C_qh^q
+
\frac{A_q\epsilon_f}{h}
\right)
\qquad \text{and} \qquad
\kappa_g(\lambda)
=
0.
\end{align*}
Thus, the oracle input determines the finite-difference interval, the approximation order, and the corresponding stencil.

Furthermore, if
$
h
\geq
\left(
\frac{A_q\epsilon_f}{C_q}
\right)^{\frac{1}{q+1}}
$, 
then
\begin{align*}
\epsilon_g(\lambda)
\leq
2\sqrt{n}C_qh^q.
\end{align*}
The smallest admissible value of $h$ yields
\begin{align*}
\epsilon_g
\leq
2\sqrt{n}\,
C_q^{\frac{1}{q+1}}
\left(
A_q\epsilon_f
\right)^{\frac{q}{q+1}}.
\end{align*}
Hence, higher-order stencils can improve the dependence of the gradient approximation error on $\epsilon_f$, provided that the required additional smoothness is available.
Note that it is possible for $C_q=0$ only if $L_q=0$. This means that the finite difference scheme (without noise) can compute $\nabla \phi(x)$ exactly using arbitrarily large $h$. Thus choosing large $h$ reduces the error in the first-order oracle and can make it arbitrarily small, in theory. Hence the inequality $h \geq \left( \frac{A_q\epsilon_f}{C_q} \right)^{\frac{1}{q+1}}$ in the case of $C_q=0$ can be understood as a limiting case for $C_q>0$.

Finally, let $\delta_0=1$ if $c_j=0$ for some $j=0,\ldots,m$, and $\delta_0=0$ otherwise. Assuming the stencil nodes are distinct, the full gradient approximation requires
$
N_f
=
n(m+1)
-
(n-1)\delta_0
$
distinct evaluations of Oracle~\ref{def:oracle_det_0}, and therefore has cost
$
{\cal C}_g(\lambda)
=
N_f{\cal C}_f
+
\mathcal{O}\left(n(m+1)\right).
$
In particular, standard one-sided $q$th-order formulas require $qn+1$ function evaluations, whereas standard centered formulas of even order $q$ require $qn$ function evaluations. Thus, increasing $q$ can improve approximation accuracy, but comes at the cost of stronger smoothness assumptions, additional oracle evaluations, and stencil-dependent noise amplification.

\subsubsection{Block-Coordinate  and Orthogonal Projected Gradient  Directions}\label{sec:det_project}

Deterministic block-coordinate descent methods are a class of methods that, at each iteration, update only a subset of the components of the iterate $x\in\mathbb{R}^n$. Such methods can be effective when partial updates are substantially cheaper than full-dimensional updates \cite{nesterov2012efficiency,wright2015coordinate}. For example, for a dense quadratic objective, computing the full gradient requires $\mathcal{O}(n^2)$ operations, whereas computing a single partial derivative, or updating a cached gradient following a single-coordinate step, requires only $\mathcal{O}(n)$ operations. The coordinates may be selected greedily, as in the Gauss--Southwell rule, which chooses a coordinate corresponding to a largest-magnitude component of the gradient \cite{nutini2015coordinate,Joachims1999}, or according to a fixed cyclic order \cite{bertsekas1997nonlinear,tseng2001convergence,saha2013nonasymptotic}.

While not often presented this way, deterministic block-coordinate descent can be viewed as a gradient descent method with an inexact first-order oracle. Suppose the method updates $s$, for $1 \leq s \leq n$, components of $x \in \mathbb{R}^n$, where the indices of these components belong to the set $\mathcal{S}(x)$, whose cardinality is $s$, and is chosen such that the corresponding elements $[\nabla \phi(x)]_i$ are largest in their absolute value. Such an update can be written as $x - \alpha g(x)$, where $g(x)$ is a vector whose $i$th element is set to $[\nabla \phi(x)]_i$ if $i \in \mathcal{S}(x)$ and is set to 0 otherwise. It is easy to see that
\begin{align}\label{eq:bcd_bound}
	\|g(x)-\nabla \phi(x)\|^2=\sum_{i\notin \mathcal{S}(x)}[\nabla \phi(x)]_i^2\leq \frac{n-s}{n}\|\nabla \phi(x)\|^2.
\end{align}
Thus, $g(x)$ is an output of  Oracle \ref{def:oracle_det_1} with input $\lambda=s$, the number of components that are updated, and the resulting error is 
\begin{align*}
	\epsilon_g(s)=0 \quad  \text{ and } \quad  \kappa_g(s)=\sqrt{\frac{n-s}{n}}.
\end{align*}
The cost ${\cal C}_{g}(s)$ associated with this oracle heavily depends on the form of the objective function. Since it is assumed that the oracle selects the $s$ largest components of $\nabla \phi(x)$, this usually means that the oracle needs to compute the entire gradient. If $\phi(x)$ is quadratic, and only $s$ components of $x$ are updated, then $\nabla \phi(x)$ is computed from its previous value by ${\cal O}(sn)$ operations. Thus, if $s$ is constant and $\phi(x)$ is quadratic, ${\cal C}_{g}(s)$ scales linearly with $n$. This type of coordinate-descent method has limited applications in  practice (aside from quadratic functions) and is presented here as an inexact oracle as an illustration of how block-coordinate directions in general can be viewed in the context of inexact oracles. We will see later that randomized coordinate-descent, which has a much wider applicability, also fits our oracle definition. 

More generally, consider a matrix $S\in \mathbb{R}^{n\times s}$ with orthonormal columns. Then for any $v\in \mathbb{R}^n$, $SS^Tv$ is an orthogonal projection of $v$ onto the subspace induced by the columns of $S$. For example, if $S$ is chosen to be the subset of the columns of the identity matrix, then $SS^Tv$ is a projection onto the subspace of the corresponding coordinates. Thus, $g(x)=SS^T \nabla \phi(x)$ is an output of  Oracle \ref{def:oracle_det_1} with input $\lambda=S$. 
Due to the orthogonality, 
\begin{align}\label{eq:bcd_bound1}
	\|g(x)-\nabla \phi(x)\|^2=\|SS^T \nabla \phi(x)- \nabla \phi(x)\|^2& =\|\nabla \phi(x)\|^2- \|SS^T \nabla \phi(x)\|^2 \nonumber\\ 
	& = (1-\theta^2(x)) \|\nabla \phi(x)\|^2.
\end{align}
where  $\theta(x)=\|SS^T \nabla \phi(x)\|/\|\nabla \phi(x)\| \in [0,1]$.  Let $\theta\leq \theta(x) \in [0,1]$ be a universal lower bound, then  the resulting error terms are  
\begin{align*}
	\epsilon_g(S)=0 \quad  \text{ and } \quad  \kappa_g(S)=\sqrt{1-\theta^2}.
\end{align*}
It is important to note that $\theta$ may be $0$ unless the subspace is chosen with some care. For example, in the case of block-coordinate directions where the selected coordinates correspond to those with the largest $s$ components of $\nabla \phi(x)$,   $\theta^2 \ge s/n$. More generally guaranteeing a lower bound on $\theta^2 $ in a deterministic case  is not simple, but randomized subspace selection ensures it in expectation, as we will see in Section \ref{sec:stochastic}. A randomized scheme can also be used to ensure our deterministic oracle bounds are valid with some sufficiently high probability.

\subsubsection{Block-Coordinate and Projected Gradient  Directions Estimated via Finite Differences}\label{sec:det_project_ffd}
 
Consider a set of coordinates $\mathcal{S}(x)$ and a vector $g\in \mathbb{R}^n$ whose $i$-th component is computed as finite difference approximations of $[\nabla \phi(x)]_i$ for all $i \in \mathcal{S}(x)$ and set to zero otherwise. Specifically with FFD, given $\mathcal{S}(x)$ and access to Oracle~\ref{def:oracle_det_0}, this approximation is given by
  \begin{align*}
	g(x) =  \sum_{i\in \mathcal{S}(x)} \frac{f(x+{h}e_i) - f(x)}{{h}} e_i. 
\end{align*}

Let $\nabla_{\mathcal{S}(x)}\phi(x)$ be the projection of $\nabla \phi(x)$ onto the subspace spanned by components in $\mathcal{S}(x)$. 
From the analysis of the finite difference oracle, it follows that
 \begin{align*}
	\|g(x)-\nabla_{\mathcal{S}(x)}\phi(x)\| \leq \frac{\sqrt{s} L_1 {h}}{2} + \frac{2\sqrt{s}  \epsilon_f}{h}.
 \end{align*}
 This is similar to the bound in \eqref{eq:FFD error} but with $s$ in place of $n$. 
Assuming $\mathcal{S}(x)$ is chosen somehow to satisfy 
$\|\nabla_{\mathcal{S}(x)}\phi(x)-\nabla \phi(x)\|^2\leq \frac{n-s}{n}\|\nabla \phi(x)\|^2$ (discussed above), we can bound the squared error as follows 
\begin{equation} \label{eq.ffd_bcd} \begin{aligned}
	\|g({x}) - \nabla \phi({x})\|^2 
    &= \|g({x}) - \nabla_{\mathcal{S}(x)}\phi(x)\|^2 + \|\nabla_{\mathcal{S}(x)}\phi(x) - \nabla \phi({x})\|^2  \\
    &\leq \left( \frac{\sqrt{s} L_1 {h}}{2} + \frac{2\sqrt{s}  \epsilon_f}{h}\right)^2 + \frac{n-s}{n}\|\nabla \phi(x)\|^2, 
\end{aligned} \end{equation} 
where the equality holds because $g({x}) - \nabla_{\mathcal{S}(x)}\phi(x)$ is in the span of $\mathcal{S}(x)$, while $\nabla_{\mathcal{S}(x)}\phi(x) - \nabla \phi({x})$ is in its null space. 
Thus, this is an example of Oracle~\ref{def:oracle_det_1} with 
\begin{align*}
	\epsilon_g(h,s)= \frac{\sqrt{s} L_1 {h}}{2} + \frac{2\sqrt{s}  \epsilon_f}{h} 
    \quad  \text{ and } \quad  
    \kappa_g(h,s) = \sqrt{\frac{n-s}{n}}.
\end{align*}

We note that one could use CFD approximations to the gradient in the computed components, and this would change $\epsilon_g(h,s)$.

Again, more generally, consider a matrix $S(x)\in {R}^{n\times s}$ with orthonormal columns.
Let us apply standard forward finite difference method to approximate the ``reduced" gradient in the column space of $S(x)$. Since we do not discuss specific ways to choose $S(x)$ except that is has $s$ orthonormal columns, for brevity,  we will use simply $S$, keeping in mind possible dependence on $s$. 

Specifically, given access to Oracle~\ref{def:oracle_det_0}, we compute the $s$-dimensional vector 
\begin{equation}\label{eq:subspace_ffd}
g_s(x)=\sum_{i=1}^s \frac{f(x+h Se_i)-f(x)}{h}e_i 
\end{equation}
where $e_i$ is the $i$th column of the $s\times s$ identity matrix. We then define $g(x)=Sg_s(x)$. By the standard finite difference analysis, the error in the column space of $S$ follows 
\begin{align*}
	\|g(x)-SS^T\nabla {\phi}(x)\| = \|Sg_s(x)-SS^T\nabla {\phi}(x)\| &=\|g_s(x)-S^T\nabla {\phi}(x)\| \\ 
	&\leq  \frac{\sqrt{s}L_1{h}}{2}+ 
\frac{2\sqrt{s}  e_f}{{h} }
\end{align*}
Then, if the error in the null space follows \eqref{eq:bcd_bound1}, we would have 
\begin{align}\label{eq.ffd_randproj}
	\|g({x}) - \nabla \phi({x})\|^2 
    &= \|g({x}) - SS^T\nabla\phi(x)\|^2 + \|SS^T\nabla\phi(x) - \nabla \phi({x})\|^2  \\
    &\leq \left( \frac{\sqrt{s} L_1 {h}}{2} + \frac{2\sqrt{s}  \epsilon_f}{h}\right)^2 + (1-\theta(x)^2) \|\nabla \phi(x)\|^2, 
\end{align}
where $\theta(x)=\|SS^T \nabla \phi(x)\|/\|\nabla \phi(x)\| \in [0,1]$.  Again, let $\theta\leq \theta(x) \in [0,1]$ be a universal lower bound on  $ \theta(x)$. This results in Oracle~\ref{def:oracle_det_1} with 
\begin{align*}
	\epsilon_g(h,s)= \frac{\sqrt{s} L_1 {h}}{2} + \frac{2\sqrt{s}  \epsilon_f}{h} 
    \quad  \text{ and } \quad  
    \kappa_g(h,s)=\sqrt{1-\theta^2}.
\end{align*}

Choosing an appropriate coordinate subset $\mathcal{S}(x)$ or a projection matrix $S$ (without computing $\nabla \phi(x)$) such that $\theta>0$ and thus $\kappa_g(h,s)$ remains suitably bounded away from $1$ is generally not easy. However, we use this oracle for illustration purposes. In later sections we will introduce an oracle that selects $\mathcal{S}(x)$ randomly, where similar derivations will apply, but in probability or in expectation.  

We also note that in either case, $\epsilon_g(h,s)$ does not depend on $n$, while $\kappa_g(h,s)$, generally, does. We will revisit these dependencies when analyzing convergence properties of various algorithms and their dependence on these oracle properties.

\subsubsection{Linear Interpolation}
Linear interpolation (LI) is a generalization of FFD. 
Given an affinely independent interpolation set $\{x, y_1, \ldots, y_n\} \subset \mathbb{R}^n$ and access to Oracle~\ref{def:oracle_det_0}, LI generates a linear approximation model
\begin{equation} \label{eq:linear_model}
    m(y) = f(x) + g(x)^T (y - x), 
\end{equation} 
where $g(x)\in\mathbb{R}^n$, that satisfies {\em interpolation} conditions, i.e.,
\begin{equation} \label{eq:linear_int_cond}
    m(y_i) = f(y_i), \quad i=1, \ldots n
\end{equation} 
The model $m$ is determined by the vector $g$, which can be seen as an approximation of $\nabla\phi(x)$. Specifically, $g$ is obtained by solving the following linear system of equations
\begin{equation} \label{eq:LI_linear_system}
    M g(x) = f_\Delta, \quad \text{ with } \;\; 
    M =\begin{bmatrix} (y_1-x)^T\\ \vdots\\ (y_n-x)^T \end{bmatrix}  \in \mathbb{R}^{n \times n} \;\; \text{ and  } \;\;
   f_\Delta= \begin{bmatrix} f(y_1)-f(x)\\ \vdots\\ f(y_n)-f(x) \end{bmatrix} \in \mathbb{R}^{n}. 
\end{equation} 
When $y_i = x+he_i, i=1,\dots,n$ for some $h>0$, this approach reduces  to FFD. 

To draw further analogy with FFD, we write $y_i = x+hu_i, i=1,\dots,n$ for some linearly independent set of directions $u_i$, for which $\max_{i=1,\ldots n} \|u_i\|= 1$.  We can bound $\|g(x) - \nabla\phi(x)\|$ using the following steps. First, it follows that 
\begin{equation} \label{eq:LI_g_bound1996}
    \| M \big(g(x) - \nabla\phi(x)\big)\| \le \frac{\sqrt{n} L_1 }{2} h^2 + 2\sqrt{n}\epsilon_f. 
\end{equation}
This bound can be found for the case when $\epsilon_f=0$ in \cite{regis2015calculus} and \cite[Chapter XX]{ConnScheVice08c}, and generalizes for $\epsilon_f>0$. 
By letting $\hat{M} = M / h$, we have 
\begin{equation} \label{eq:LI_error_bounds_inexact}
    \|g(x) - \nabla\phi(x)\| \le  \|\hat{M}^{-1}\| \left(\frac{\sqrt{n}L_1h}{2} + \frac{2\sqrt{n}\epsilon_f}{{h}} \right). 
\end{equation}
The quantity  $\|\hat{M}^{-1}\|\geq 1$ characterizes linear independence  of directions $u_i$ and thus contains  geometric information of the interpolation set other than the scaling. For FFD, $\hat{M}^{-1}$ is identity matrix and thus $\|\hat{M}^{-1}\| = 1$.

Thus,  vector $g$  can be seen as an output of  Oracle~\ref{def:oracle_det_1} with 
\begin{equation}\label{eq:LI_kappa_bound}
    \epsilon_g = \|\hat{M}^{-1}\| \left(\frac{\sqrt{n}L_1 h}{2} + \frac{2\sqrt{n}\epsilon_f}{h} \right) \quad
    \text{ and } \quad
    \kappa_g = 0, 
\end{equation}
while the input $\lambda$ consists of the interpolation set $\{y_i\}_{i=1}^n$. 

Considering a non-adjustable Oracle~\ref{def:oracle_det_0}, a naive implementation of the LI oracle requires $n+1$ calls to Oracle~\ref{def:oracle_det_0} and the subsequent solution of the linear system \eqref{eq:LI_linear_system}, implying the cost ${\cal C}_{g} = (n+1){\cal C}_{f} + \mathcal{O}(n^3)$. 
With a generally worse error bound and higher cost, this oracle appears to be strictly inferior to FFD. 
However, since there is no strict requirement for the position of the interpolation points, the cost of LI can be greatly reduced by reusing zeroth-order information that is already computed and available. 
For example, model-based derivative-free optimization methods \cite{MJDPowell_1994,MJDPowell_2002,powell2009bobyqa,ConnScheVice08c}, which are designed to address problems with high $\mathcal{C}_f$ and no directly computable derivatives, determine the optimization direction in each iteration by approximating the objective function $\phi$ locally with interpolation models.  
The interpolation set for each model consists mostly of points from previous sets, and thus only a few calls to the zeroth-order oracle are needed per iteration. 

Thus, we can consider the following somewhat more general oracle. Let $x$, $f(x)$, $\{y_i\}_{i=1}^p$, $\{f(y_i)\}_{i=1}^p$,  be the input to the oracle, for some $p\leq n$. In other words, the oracle receives some sample points around the current iterate $x$ and their respective possibly inexact function values. The oracle would then complete the interpolation set by choosing the remaining $n-p$ points, $\{y_i\}_{i=p+1}^n$, accordingly, call the zeroth-order oracle on them, and solve \eqref{eq:LI_linear_system} to obtain the gradient estimate. In this case the oracle can choose these points with the aim of minimizing $\|\hat M^{-1}\|$. The sampling radius $h$ may be given as an input or inferred based on $\max_{i=1,\ldots n} \|y_i-x\|$ depending on the setting.

\subsubsection{Simplex Gradients} 
Similar to applying FFD in a subspace it makes sense to consider interpolation models that only interpolate 
over a sample set in a subspace, and thus contain fewer than $n+1$ sample points.  
In this  case, one can consider a set  $\{x,y_1,\dots,y_p\}$ of $p+1$ points, with $p \le n$, and solve the problem 
\begin{equation} \label{eq:LI subspace}
	\begin{aligned}
		&\min_{g} &&\|g\|^2\\
		&\text{s.t. } &&M g = f_\Delta, 
	\end{aligned}
	\text{ where } M =  \begin{bmatrix}
		(y_1 - x)^\intercal\\ \vdots\\ (y_p-x)^\intercal, 
	\end{bmatrix} \in \mathbb{R}^{p \times n}
	\text{ and } 
	f_\Delta = \begin{bmatrix}
		f(y_1) - f(x)\\ \vdots\\ f(y_p) - f(x)
	\end{bmatrix}  \in \mathbb{R}^{p}
\end{equation}
to construct a linear interpolation model in the subspace induced  by the linearly independent vectors $\{y_i - x\}_{i=1}^p$. This yields the minimum norm solution to the system of equations.

Alternatively, if the  zeroth-order information is  abundant around the iterate $x$, with the function values readily available at $\{y_i\}_{i=1}^p$ for some $p>n$, one may want to utilize all this information by allowing more than $n+1$ points for interpolation. 
A successful approach in this case is to use higher-order interpolating polynomial (e.g., quadratic), which we will discuss in the section of second order oracles. 
  If one  aims at building only a linear model, then least-square regression is performed in lieu of interpolation to construct the linear model $m(y) = c + g^T(y-x)$ \eqref{eq:linear_model} that minimizes the squared error, i.e., 
  \begin{align*}
  	\min_{c,g} \;\; (m(x)-f(x))^2+\sum_{i=1}^n (m(y_i)-f(y_i))^2,
  \end{align*}
  which can be written compactly as
  \begin{equation} \label{eq:grad_reg}
	\min_{c,g} \left\|\begin{bmatrix} 1 & {\bf 0}_n^T\\  {\bf 1}_p &M \end{bmatrix} \begin{bmatrix} c\\ g \end{bmatrix}- \begin{bmatrix} f(x)\\ f_{Y} \end{bmatrix} \right\|^2,
\end{equation}
where ${\bf 1}_p  \in \mathbb{R}^p$ is the vector of all ones, ${\bf 0}  \in \mathbb{R}^n$ is the vector of all zeros and $f_Y=(f(y_1), \ldots f(y_p))^T$. 

Finally, a gradient approximation can also be obtained by solving the least-square problem 
\begin{equation} \label{eq:alternative linear regression}
	\min_{g} \|M g - f_\Delta\|^2, 
\end{equation}
which is similar to the standard linear regression \eqref{eq:grad_reg} but requires the linear model to interpolate the objective function at $x$ exactly.

The simplex gradient generalizes the gradient $g$ obtained from both \eqref{eq:LI subspace} and \eqref{eq:alternative linear regression} and is even applicable when the vectors $\{y_i-x\}_{i=1}^p$ in \eqref{eq:LI subspace} are linearly dependent. 
Given a set of possibly dependent vectors $\{y_i-x\}_{i=1}^p \subset \mathbb{R}^n$ such that $\max_{i=1,\ldots p} \|y_i-x\|\leq h$, where $p$ can be any positive integer, the simplex gradient, as originally proposed in \cite{custodio2008using}, approximates the gradient as 
\begin{equation}
	g(x) = \frac{1}{h}V\Sigma^{-1}U^T f_\Delta, 
\end{equation}
where $V\Sigma U^T$ is the reduced singular value decomposition of $\hat M = M / h$. 
This concept was later studied in detail in \cite{regis2015calculus} with the following equivalent definition
 \begin{equation}\label{eq:simplex_grad_reg}
    g(x) = M^\dagger f_\Delta,
\end{equation}
where $M^\dagger$ denotes the Moore-Penrose pseudoinverse of the matrix $M$.

Under Assumption~\ref{assum:lip_cont} and an exact zeroth-order oracle (i.e., $\epsilon_f=0$), upper bounds on the approximation error of the simplex gradient are derived in \cite[Proposition 7]{regis2015calculus} and \cite[Proposition 2.1]{hare2024matrix}, and can be summarized as follows,
\begin{equation} \label{eq:simplex_grad_bound} \begin{aligned} 
    \|\proj_{M}\big(g(x) - \nabla\phi(x)\big)\| 
    &= \|g(x) - \proj_{M}\big(\nabla\phi(x)\big)\| \\
    &\le \frac{\sqrt{p}}{2} L_1 \|\hat{M}^\dagger\| {h}, 
\end{aligned} \end{equation}
where $\proj_M(\cdot) = M^\dagger M (\cdot)$. 
Alternatively, if we only have access to a noisy variant of Oracle~\ref{def:oracle_det_0}, the approximation error of the simplex gradient can be bounded as follows. 

\begin{proposition}[Error bound for simplex gradient with Oracle~\ref{def:oracle_det_0}] \label{prop:simplex_grad_noisy}
    Given a sample set $\{x,y_1,\dots,y_p\} \subset B(x,{h})$ and the corresponding $M=M_{p\times n}$ and $\hat{M}=M/{h}$, 
    we have under Assumption~\ref{assum:lip_cont} that 
    \begin{equation} \label{eq:simplex_grad_bound_noisy} \begin{aligned} 
        \|\proj_{M}\big(g(x) - \nabla\phi(x)\big)\| 
        &= \|g(x) - \proj_{M}\big(\nabla\phi(x)\big)\| \\
        &\le \sqrt{p} \|\hat M^\dagger\| \Big(\frac{L_1}{2} h + \frac{2 \epsilon_f}{h} \Big), 
    \end{aligned} \end{equation}
    for $g(x)$ computed by \eqref{eq:simplex_grad_reg} with Oracle~\ref{def:oracle_det_0}, where $\proj_M(\cdot) = M^\dagger M (\cdot)$ is the projection onto row space of $M$. 
\end{proposition}
\begin{proof}
    For any $i = 1,\dots,p$, the absolute value of the $i$th element of the vector $f_\Delta - M \nabla\phi(x)$ is bounded as 
    \begin{align*}
        &|f(y_i) - f(x) - (y_i-x)^\intercal \nabla\phi(x)| \\
        \le &|\phi(y_i) - \phi(x) - (y_i-x)^\intercal \nabla\phi(x)| + |f(y_i) - \phi(y_i)| + |f(x) - \phi(x)| \\
        \le &\frac{L_1}{2} {h}^2 + 2 \epsilon_f. 
    \end{align*}
    Thus, the vector length $\|f_\Delta - M \nabla\phi(x)\| \le \sqrt{p} (L_1 {h}^2/2 + 2\epsilon_f)$ and 
    \begin{align*}
        \|g(x) - M^\dagger M(\nabla\phi(x))\| 
        &= \|M^\dagger f_\Delta - M^\dagger M(\nabla\phi(x))\| 
        \le \|M^\dagger\| \|f_\Delta - M \nabla\phi(x)\| \\
        &\le \sqrt{p} \|M^\dagger\| \Big(\frac{L_1}{2} {h}^2 + 2 \epsilon_f \Big)
        = \sqrt{p} \|\hat{M}^\dagger\| \Big(\frac{L_1}{2} {h} + 2 \frac{\epsilon_f}{h} \Big). 
    \end{align*}
    The final result is obtained using $\|M^\dagger\| = \|\hat M^\dagger\| / h$. 
\end{proof}

Similar to \eqref{eq:bcd_bound1}, let $\theta(x) =\|\proj_{M}\nabla \phi(x)\|/\|\nabla \phi(x)\| \in [0,1]$, then the error in the null space of $M$ follows 
\begin{align*}
  \|\proj_{M}\big(\nabla \phi(x)\big) - \nabla\phi(x)\|^2 
    = (1-\theta(x)^2) \|\nabla\phi(x)\|^2. 
\end{align*}
Thus, the squared error of the simplex gradient is
\begin{align*}
    \|g(x) - \nabla\phi(x)\|^2 
    &= \|g(x) - \proj_{M}\nabla \phi(x)\|^2 + \|\proj_{M}\nabla \phi(x) - \nabla\phi(x)\|^2 \\
    &\le \left[\sqrt{p} \|\hat M^\dagger\| \Big(\frac{L_1h}{2} + \frac{2 \epsilon_f}{h} \Big)\right]^2 + (1-\theta(x)^2) \|\nabla \phi(x)\|^2. 
\end{align*}
Again, letting $\theta$ denote the lower bound on $\theta(x)$ we cast the  simplex gradient  as an output of Oracle~\ref{def:oracle_det_1} with input $\lambda = M$, 
\begin{equation*}
	\epsilon_g(M)= \sqrt{p} \|\hat M^\dagger\| \left(\frac{L_1h}{2} + \frac{2 \epsilon_f}{h} \right) 
    \quad  \text{ and } \quad  
    \kappa_g(M)=\sqrt{1-\theta^2}.
\end{equation*}

We note that while it looks like the error term $\epsilon_g$  scales with $\sqrt{p}$, in the case when $p>n$, under certain geometric conditions on the interpolation set, $\|\hat{M}^\dagger\|$ scales with $\sqrt{n/p}$. 
Loosely speaking, the condition essentially ensures that the rows of $\hat{M}$ can be partitioned into $p/n$ groups of $n$,  all of which are sufficiently linearly independent; see \cite[Chapter 4]{ConnScheVice08c} for more details.  On the other hand, when $p<n$, to ensure $ \kappa_g(M)<1$ requires the space spanned by rows of $M$ to me a sufficiently acute  angle with $\nabla \phi(x)$.

\subsubsection{Gradient Quantization} 
Gradient quantization is a technique widely used in large-scale centralized and decentralized optimization to reduce computational and communication costs. The key idea is to transmit or store lower-precision representations of gradient information, thereby decreasing bandwidth usage and improving scalability without significantly degrading optimization performance. Several quantization mechanisms have been proposed in the literature; one widely adopted strategy is $b$-bit quantization \cite{jacob2018quantization}. If we ignore the computation error and assume the computed gradient is $\nabla\phi(x)$, one needs to first pick a quantization level $2^b$, where $b$ is the number of bits for storing each entry of the gradient; then compute the scale $s = (\max_j \phi'_j(x) - \min_j \phi'_j(x)) / (2^b-1)$, where $\phi'_j(x)$ represents the partial derivative of $\phi$ at $x$ with respect to its $j$th input. 
Assuming that $s>0$,  the gradient is quantized as $q(x)$, a vector of $n$ $b$-bit unsigned integers such that 
\[ q(x)^T e_i = \left\lfloor \frac{\phi'_i(x) - \min_j \phi'_j(x)}{s} \right\rceil, 
\]
where $\lfloor \cdot \rceil$ means rounding to the nearest integer. 
After $q(x)$, $s$ and $\min_j \phi'_j(x)$ are communicated, the gradient can be recovered as $g(x)$ such that 
\[  g(x)^T e_i = s  q(x)^T e_i + \min_j \phi'_j(x). 
\]

Every entry of $g(x)$ is $\min_j \phi'_j(x)$ plus a multiple of $s$, and the maximum deviation is half of the interval length. 
Therefore, the quantization error is bounded as $\|g(x) - \nabla\phi(x)\|_\infty \le s / 2$ and 
\begin{equation} \begin{aligned} 
    \|g(x) - \nabla\phi(x)\| \le \frac{\sqrt{n} s}{2} 
    &= \frac{\sqrt{n}}{2} \frac{\max_j \phi'_j(x) - \min_j \phi'_j(x)}{2^b-1} \\
    &\le \frac{\sqrt{n}}{2} \frac{2\|\nabla\phi(x)\|_\infty}{2^b-1} \\
    &\le \frac{\sqrt{n}}{2^b-1} \|\nabla\phi(x)\|. 
\end{aligned} \end{equation}  
Thus, $b$-bit gradient quantization can be seen as an inexact first-order oracle with input $\lambda = b$, and errors
\begin{equation*}
	\epsilon_g(b) = 0 \quad \text{ and } \quad \kappa_g^2(b) = \frac{n}{\left(2^b-1\right)^2}. 
\end{equation*}
The requirement $\kappa_g<1$ imposes the following lower bound on the quantization level, $b\geq \left\lfloor \log_2(1+\sqrt{n})\right\rfloor+1$. A stochastic variant of this type of quantization was studied in \cite{chen2020statistical}. The stochasticity is introduced  to obtain an unbiased gradient estimator, but as we will see in our analysis, deterministic 
relative error with $\kappa_g<1$ is essentially harmless for optimization algorithms.

\subsubsection{PDE-Constrained Optimization}

In PDE-constrained optimization, gradient computation typically requires solving both a state equation and an adjoint equation \cite{hinze2008optimization,borzi2011computational}. When these systems are solved inexactly, the resulting gradient approximation naturally defines an inexact first-order oracle. Consider, after discretization, an objective of the form
\begin{align*}
\phi(x)
=
F(x,u(x)),
\qquad
R(x,u(x))=0,
\end{align*}
where $u(x)$ denotes the PDE state. Assuming that $R_u(x,u(x))$ is nonsingular, the adjoint $p(x)$ satisfies
\begin{align*}
R_u(x,u(x))^\intercal p(x)
=
\nabla_u F(x,u(x)),
\end{align*}
and the gradient is given by
\begin{align*}
\nabla\phi(x)
=
\nabla_x F(x,u(x))
-
R_x(x,u(x))^\intercal p(x).
\end{align*}

Suppose that approximate state and adjoint solutions $\widetilde u$ and $\widetilde p$ satisfy
\begin{align*}
\|R(x,\widetilde u)\|
\leq \tau_u, \quad \text{and}\quad
\left\|
R_u(x,\widetilde u)^\intercal\widetilde p
-
\nabla_u F(x,\widetilde u)
\right\|
\leq \tau_p.
\end{align*}
The resulting gradient approximation is
\begin{align*}
g(x)
=
\nabla_x F(x,\widetilde u)
-
R_x(x,\widetilde u)^\intercal\widetilde p.
\end{align*}

Under suitable stability and regularity assumptions, suppose that the state error is bounded by $K_u\tau_u$, the inverse adjoint operator is uniformly bounded by $K_p$, and $\|R_x(x,u)\|_2\leq M_x$. Further, assume that the exact adjoint-based gradient expression, evaluated at a given state, is $L_G$-Lipschitz continuous with respect to that state. Then,
\begin{align*}
\|g(x)-\nabla\phi(x)\|
\leq
L_GK_u\tau_u
+
M_xK_p\tau_p.
\end{align*}
Thus, taking $\lambda=(\tau_u,\tau_p)$, the approximation satisfies Oracle~\ref{def:oracle_det_1} with
\begin{align*}
\epsilon_g(\lambda)
=
L_GK_u\tau_u
+
M_xK_p\tau_p \quad \text{and} \quad
\kappa_g(\lambda)
=0.
\end{align*}
The oracle cost consists primarily of the state and adjoint solves and depends on the chosen tolerances:
\begin{align*}
{\cal C}_g(\lambda)
=
{\cal C}_{\rm state}(\tau_u)
+
{\cal C}_{\rm adjoint}(\tau_p)
+
{\cal C}_{\rm assembly}.
\end{align*}
Tighter tolerances improve the gradient approximation but generally increase the computational cost.


\subsubsection{Nested and Bilevel Optimization}

Nested and bilevel optimization problems provide another source of inexact gradient information, since evaluating the outer gradient typically requires solving an inner optimization problem and an associated linear system \cite{ehrhardt2024analyzing,okuno2020bilevel,dempe2020bilevel}. Consider
\begin{align*}
\phi(x)
=
F(x,y^\star(x)),
\qquad
y^\star(x)
=
\underset{y}{\operatorname{argmin}}\,
\ell(x,y),
\end{align*}
where the lower-level objective $\ell(x,\cdot)$ is twice continuously differentiable and $\mu$-strongly convex. By implicit differentiation,
\begin{align*}
\nabla\phi(x)
=
\nabla_x F(x,y^\star)
-
\nabla_{xy}^2\ell(x,y^\star)z^\star,
\end{align*}
where $z^\star$ solves
\begin{align*}
\nabla_{yy}^2\ell(x,y^\star)z^\star
=
\nabla_y F(x,y^\star).
\end{align*}

Suppose that an approximate inner solution $\widetilde y$ and an approximate linear-system solution $\widetilde z$ satisfy
\begin{align*}
\|\nabla_y\ell(x,\widetilde y)\|
\leq \lambda_y \quad \text{and} \quad
\left\|
\nabla_{yy}^2\ell(x,\widetilde y)\widetilde z
-
\nabla_y F(x,\widetilde y)
\right\|
\leq \lambda_z.
\end{align*}
An approximate outer gradient is then
\begin{align*}
g(x)
=
\nabla_x F(x,\widetilde y)
-
\nabla_{xy}^2\ell(x,\widetilde y)\widetilde z.
\end{align*}

Strong convexity gives
\begin{align*}
\|\widetilde y-y^\star(x)\|
\leq
\frac{\lambda_y}{\mu}.
\end{align*}
Suppose additionally that
$\|\nabla_{xy}^2\ell(x,y)\|_2\leq M_{xy}$
and that the exact implicit-gradient expression, evaluated
at a given inner point $y$, is $L_G$-Lipschitz continuous
with respect to $y$. Then,
\begin{align*}
\|g(x)-\nabla\phi(x)\|
\leq
\frac{L_G}{\mu}\lambda_y
+
\frac{M_{xy}}{\mu}\lambda_z.
\end{align*}
Consequently, with $\lambda=(\lambda_y,\lambda_z)$, the
approximation satisfies Oracle~\ref{def:oracle_det_1} with
\begin{align*}
\epsilon_g(\lambda)
=
\frac{L_G\lambda_y+M_{xy}\lambda_z}{\mu} \quad \text{and} \quad 
\kappa_g(\lambda)
=0.
\end{align*}
The oracle cost depends on the accuracy required for the inner optimization problem and the auxiliary linear system:
\begin{align*}
{\cal C}_g(\lambda)
=
{\cal C}_{\rm inner}(\lambda_y)
+
{\cal C}_{\rm linear}(\lambda_z)
+
{\cal C}_{\rm derivatives}.
\end{align*}
Thus, the accuracy and cost of the outer gradient can be controlled through the tolerances of the two inner solves.

\subsection{Connections with Other Deterministic Inexact First-Order Oracles}\label{sec.other_oracles}

We next examine how Oracle~\ref{def:oracle_det_1} relates to other oracle formulations in the literature. In particular, we consider oracles that provide general descent directions and  also  specialized  inexact first-order oracles for convex functions, e.g., \cite{Devolder2013FirstorderMO}.

\subsubsection{General Descent: Angle Condition}\label{sec:angle}

A vast literature of unconstrained optimization explores methods that are based not simply on negative gradient directions but more general directions of descent; see e.g., \cite{NoceWrig06}, \cite{bertsekas1997nonlinear}, and \cite[Equation (3)]{curtis2019exploiting}. Prominent among these are quasi-Newton methods and their limited-memory variants, in which the gradient is premultiplied by a matrix that approximates the inverse Hessian, thereby yielding a step that approximates a Newton direction. The result of this computation,  which we call $\tilde g(x)$, can be seen as an output of an inexact first-order oracle. A key property is that the resulting direction, which is $-{\tilde g(x)}$,   should still be a direction of descent, that is ${\tilde g(x)}$ should make an acute angle with the gradient $\nabla \phi(x)$. The following conditions are common variants in the literature and ensure this property, i.e., that $-{\tilde g(x)}$ is a descent direction, 
\begin{equation} 
\begin{aligned}\label{eq:angle_condition}
	{\tilde g(x)}^T \nabla \phi(x) &\geq \beta_1\|\nabla \phi(x)\|^2\\
	\|{\tilde g(x)}\|^2 &\leq \beta_2^2\|\nabla \phi(x)\|^2,   
\end{aligned}
\end{equation} 
for some   $\beta_2\geq \beta_1>0$. It is easy to verify that \eqref{eq:angle_condition} implies 
\begin{align*}
	\cos({\tilde g(x)},\nabla\phi(x))=\frac{{\tilde g(x)}^T \nabla \phi(x)}{\|\nabla \phi(x)\|\|{\tilde g(x)}\| }\geq \frac{\beta_1}{\beta_2}>0.
\end{align*} 
Equivalently, one can impose \eqref{eq:angle_condition} using a popular alternative, \cite[Equation (3)]{curtis2019exploiting}, 
\begin{equation} 
\begin{aligned}\label{eq:angle_condition_alt}
	{\tilde g(x)}^T \nabla \phi(x) &\geq \tilde{\beta}_1\|\nabla \phi(x)\|\|\tilde g(x)\|\\
	(\tilde{\beta}^\prime_1)^2 \|\nabla \phi(x)\|^2\leq \|{\tilde g(x)}\|^2 &\leq \tilde{\beta}_2^2\|\nabla \phi(x)\|^2,   
\end{aligned}
\end{equation} 
where $\tilde{\beta}_1 \in (0,1]$ and $0< \tilde{\beta}_1^\prime  \leq \tilde{\beta}_2$. It is easy to see that one can derive \eqref{eq:angle_condition} from \eqref{eq:angle_condition_alt}  with constants (${\beta}_1 = \tilde{\beta}_1{\tilde{\beta}_1^\prime}$ and ${\beta}_2 = \tilde{\beta_2}$). The reverse derivation is also fairly straightforward. Thus 
in the remainder of this subsection, we only consider the conditions given in \eqref{eq:angle_condition}.

\begin{figure*}[h!]
    \centering
    \begin{subfigure}[t]{0.24\textwidth}
        \centering
        \includegraphics[width=1\linewidth,trim={10 0 0 0},clip]{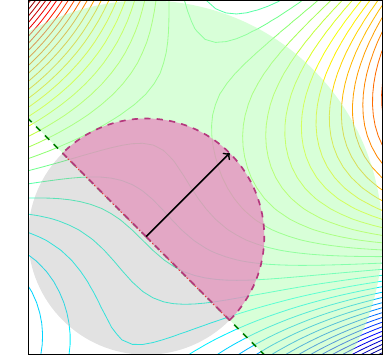}
        \caption{$\beta_1=0$, $\beta_2=1$}
    \end{subfigure}%
    ~ 
     \begin{subfigure}[t]{0.24\textwidth}
        \centering
        \includegraphics[width=1\linewidth,trim={10 0 0 0},clip]{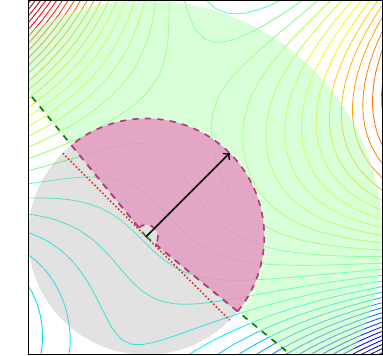}
        \caption{$\beta_1=0.1$, $\beta_2=1$}
    \end{subfigure}%
    ~
     \begin{subfigure}[t]{0.24\textwidth}
        \centering
        \includegraphics[width=1\linewidth,trim={10 0 0 0},clip]{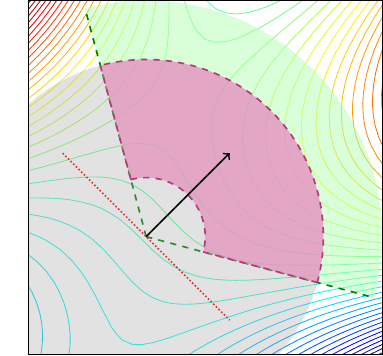}
        \caption{$\beta_1=0.5$, $\beta_2=1.5$}
    \end{subfigure}%
    ~
     \begin{subfigure}[t]{0.24\textwidth}
        \centering
        \includegraphics[width=1\linewidth,trim={10 0 0 0},clip]{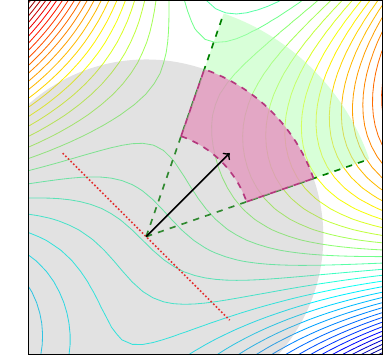}
        \caption{$\beta_1=0.9$, $\beta_2=1.5$}
    \end{subfigure}%
\caption{Illustration of  conditions~\eqref{eq:angle_condition}, ${\tilde g(x)}^T \nabla \phi(x) \geq \beta_1\|\nabla \phi(x)\|^2$ and $\|{\tilde g(x)}\|^2 \leq \beta_2^2\|\nabla \phi(x)\|^2$. 
In the figure, the black arrow represents the true gradient $\nabla \phi(x)$. The red dotted line is orthogonal to the true gradient and separates  all vectors that provide descent directions. 
The grey shaded region in each subplot depicts the region that satisfies $\|\tilde{g}(x)\| \leq \beta_2 \|\nabla \phi(x)\|$. The green region provides both a lower bound and a condition on the angle, and is a result of a combination of the conditions. 
The shaded magenta region in each subplot depicts the set of vectors $\tilde{g}(x)$ that satisfy \eqref{eq:angle_condition} (with $\beta_1=0$ as a limiting case).}
\label{fig:angle}
\end{figure*}

Figure~\ref{fig:angle} provides an illustration of conditions~\eqref{eq:angle_condition} for various values of the parameters $(\beta_1, \beta_2)$, where $\beta_2 \geq \beta_1 > 0$.  The two conditions together bound both the angle between the two vectors as well as provide upper and lower bounds for the approximate search direction $\tilde{g}$. 

\begin{figure*}[h!]
    \centering
    \begin{subfigure}[t]{0.24\textwidth}
        \centering
        \includegraphics[width=1\linewidth,trim={10 0 0 0},clip]{figs/angle_3}
        \caption{Illustration of  conditions~\eqref{eq:angle_condition}: $\beta_1=0.5$, $\beta_2=1.5$.\label{sfig:angle1}}
    \end{subfigure}%
    ~ 
     \begin{subfigure}[t]{0.24\textwidth}
        \centering
        \includegraphics[width=1\linewidth,trim={10 0 0 0},clip]{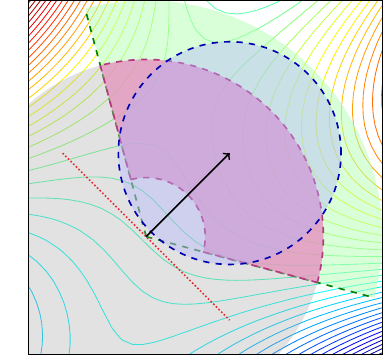}
        \caption{Illustration of  conditions~\eqref{eq:angle_condition}: $\beta_1=0.5$, $\beta_2=1.5$, and condition of Oracle~\ref{def:oracle_det_1}: $\kappa_g = \sqrt{1 - \tfrac{\beta_1^2}{\beta_2^2}} \approx 0.94$ and $\epsilon_g = 0$. \label{sfig:angle2}}
    \end{subfigure}%
    ~
     \begin{subfigure}[t]{0.24\textwidth}
        \centering
        \includegraphics[width=1\linewidth,trim={10 0 0 0},clip]{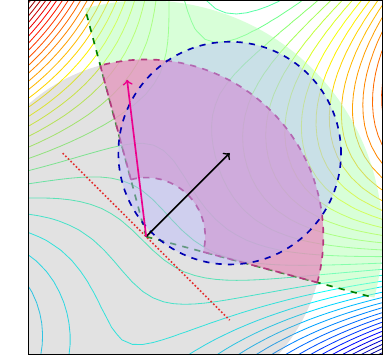}
        \caption{Magenta vector satisfying conditions~\eqref{eq:angle_condition} and not the conditions of Oracle~\ref{def:oracle_det_1}. \label{sfig:angle3_1}}
    \end{subfigure}%
    ~
     \begin{subfigure}[t]{0.24\textwidth}
        \centering
        \includegraphics[width=1\linewidth,trim={10 0 0 0},clip]{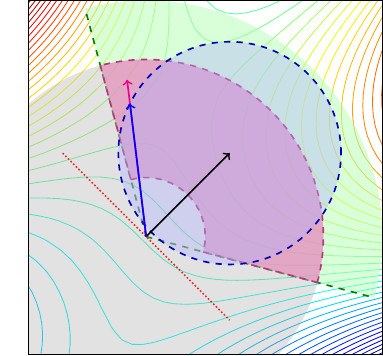}
        \caption{Blue vector is a scaled version of the magenta vector (scaling $g(x)= \frac{\beta_1}{\beta_2^2}{\tilde g(x)}$) satisfying the conditions of Oracle~\ref{def:oracle_det_1}. \label{sfig:angle3}}
    \end{subfigure}%
\caption{Regions satisfying conditions~\eqref{eq:angle_condition} (magenta) and Oracle~\ref{def:oracle_det_1} (blue), and example of scaling a vector (blue) that satisfies \eqref{eq:angle_condition} to satisfy the conditions of Oracle~\ref{def:oracle_det_1} (orange). In the figure, the black arrow represents the true gradient $\nabla \phi(x)$. The red dashed line is orthogonal to the true gradient and separates all vectors that give descent directions. }
\label{fig:first_order_oracle_scaled}
\end{figure*}

We note that  $\beta_2$ can be arbitrarily large and  $\beta_1$ can be arbitrarily small (but $\beta_2 < \beta_1$, there is no region in which the two inequalities can be satisfied simultaneously).  This seemingly implies that condition \eqref{eq:angle_condition} is looser than the conditions of Oracle  \ref{def:oracle_det_1}. In fact, in \cite{bollapragada2018adaptive} this condition is proposed as an alternative to the norm condition (i.e.,   Oracle  \ref{def:oracle_det_1} with $\epsilon_g=0$) since it is easier to satisfy. In Figure~\ref{fig:first_order_oracle_scaled} we illustrate how these conditions relate to each other. Specifically, Figure~\ref{sfig:angle1} denotes the section of annulus centered at the gradient vector that contains the vectors $\tilde g$ satisfying \eqref{eq:angle_condition} (magenta). In Figure~\ref{sfig:angle2}, the ball centered at the gradient vectors is overlaid that contains the vectors $g$ satisfying 
\eqref{eq:1_oracle_det} with $\epsilon_g=0$ and $\kappa_g<1$. This ball has the property that any direction $\tilde g$ satisfying \eqref{eq:angle_condition} can be scaled to satisfy the conditions of Oracle  \ref{def:oracle_det_1}. Figure~\ref{sfig:angle3} illustrates such an example. This means that after a suitable  scaling ${\tilde g(x)}$ satisfies condition \eqref{eq:1_oracle_det}  of Oracle  \ref{def:oracle_det_1} with $\kappa_g<1$ and $\epsilon_g=0$.

We now establish this connection formally by identifying the scaling and the resulting  $\kappa_g<1$.
\begin{proposition}
Suppose $\tilde g(x)$ satisfies \eqref{eq:angle_condition}, and let $g(x)= \frac{\beta_1}{\beta_2^2}{\tilde g(x)}$. Then, $g(x)$ satisfies the conditions of Oracle~\ref{def:oracle_det_1} with $\epsilon_g = 0$ and $\kappa^2_g=1-\frac{\beta_1^2}{\beta_2^2}\in [0, 1)$.
\end{proposition}
\begin{proof} Consider a vector $\tilde g(x)$ satisfying \eqref{eq:angle_condition} and $g(x)= \frac{\beta_1}{\beta_2^2}{\tilde g(x)}$. As a result,
\begin{equation} 
\begin{aligned} \label{eq:angle_condition_scaled}
	g(x)^T \nabla \phi(x) &\geq \frac{\beta_1^2}{\beta_2^2}\|\nabla \phi(x)\|^2\\
	\|g(x)\|^2 &\leq \frac{\beta_1^2}{\beta_2^2}\|\nabla \phi(x)\|^2,   
\end{aligned}
\end{equation} 
from which it follows that,  
\begin{align*}
	\|g(x)-\nabla \phi(x)\|^2 = \|g(x)\|^2 - 2g(x)^T \nabla \phi(x) + \|\nabla \phi(x)\|^2 \leq \left(1+\frac{\beta_1^2}{\beta_2^2}-2\frac{\beta_1^2}{\beta_2^2}\right) \|\nabla \phi(x)\|^2.
\end{align*}
This is the condition of Oracle  \ref{def:oracle_det_1} with constants $\kappa^2_g=1-\frac{\beta_1^2}{\beta_2^2}\in [0, 1)$ and $\epsilon_g = 0$.
\end{proof}

What this derivation shows is that any gradient approximation or search direction that satisfies  \eqref{eq:angle_condition} corresponds to a gradient approximation 
satisfying the conditions of Oracle  \ref{def:oracle_det_1} with constant $\kappa_g\in [0,1)$ and $\epsilon_g=0$. The scaling that is applied can be viewed as the scaling of the step size parameter along $-{\tilde g(x)}$. We will return to this observation when analyzing the corresponding algorithms (Section \ref{sec:fixed_step_det} and \ref{sec:angle_cond_rates}).  We will see then that the analysis 
of a fixed step inexact gradient descent algorithm  that uses gradient estimates $\tilde g(x)$ satisfying condition \eqref{eq:angle_condition} requires the step size parameter to be a factor of $\frac{\beta_1^2}{\beta_2^2}$ smaller than the step size parameter for the algorithm that uses  $g(x)$ satisfying \eqref{eq:1_oracle_det}. In other words, the two algorithms are essentially equivalent. 

\subsubsection{Quasi-Newton Directions}\label{sec:quasi-newton}
\begin{figure*}[h!]
    \centering
    \begin{subfigure}[t]{0.24\textwidth}
        \centering
        \includegraphics[width=1\linewidth,trim={10 0 0 0},clip]{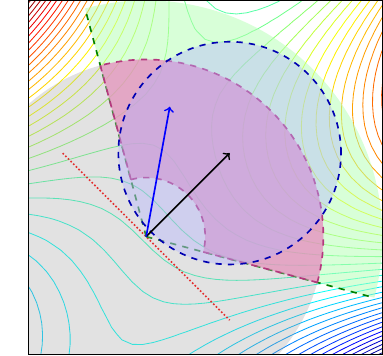}
        \caption{$\epsilon_g=0$, $\kappa_g = 0.94$. Blue vector $g(x)$ satisfies the conditions of Oracle~\ref{def:oracle_det_1}.}
    \end{subfigure}%
    ~ 
     \begin{subfigure}[t]{0.24\textwidth}
        \centering
        \includegraphics[width=1\linewidth,trim={10 0 0 0},clip]{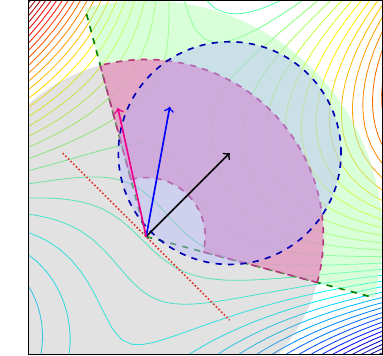}
        \caption{$\beta_1=0.5$, $\beta_2=1.5$. Magenta vector ${\tilde g(x)}=Bg(x)$ satisfies \eqref{eq:angle_condition}.  }
    \end{subfigure}%
    ~
     \begin{subfigure}[t]{0.24\textwidth}
        \centering
        \includegraphics[width=1\linewidth,trim={10 0 0 0},clip]{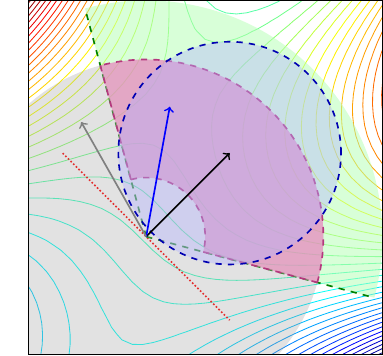}
        \caption{$\beta_1=0.5$, $\beta_2=1.5$. Gray vector ${\tilde g(x)}=Bg(x)$ does not satisfy  \eqref{eq:angle_condition}.}
    \end{subfigure}%
    ~
     \begin{subfigure}[t]{0.24\textwidth}
        \centering
        \includegraphics[width=1\linewidth,trim={10 0 0 0},clip]{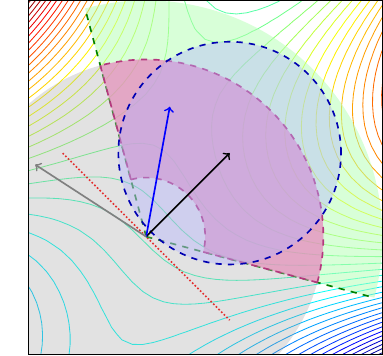}
        \caption{$\beta_1=0.5$, $\beta_2=1.5$. Gray vector ${\tilde g(x)}=Bg(x)$ does not satisfy  \eqref{eq:angle_condition}.}
    \end{subfigure}%
\caption{Example of vector $g(x)$ that satisfies the conditions of Oracle~\ref{def:oracle_det_1} (blue), and instances of vectors $\tilde{g}(x) = Bg(x)$ that satisfy (magenta) and do not satisfy \eqref{eq:angle_condition} and may even not be a descent direction (gray).}
\label{fig:matrix_scaling}
\end{figure*}

Practical performance of first-order optimization algorithms such as gradient descent can be dramatically  improved by applying appropriate rescaling  of the search directions. Generally a symmetric matrix $B$ is computed which is meant to approximate  the inverse of the Hessian matrix and a step is taken along the direction  $-Bg(x)$, where $g(x)$ is the gradient or its approximation. What is important for most algorithms is that $-Bg(x)$ remains a descent direction.  
Thus, what we would like to investigate here are conditions on  $g(x)$ as the output of  Oracle \ref{def:oracle_det_1} with $\kappa_g<1$  and $B$ - a  symmetric positive definite matrix, under which the vector ${\tilde g(x)}=Bg(x)$ satisfies \eqref{eq:angle_condition}. This, in turn, as argued above will establish that after a simple scaling  ${\tilde g(x)}$ can be itself viewed as an output of Oracle \ref{def:oracle_det_1}. Such conditions will depend on the properties of matrix $B$ and constants $\kappa_g$ and $\epsilon_g$. Figure \ref{fig:matrix_scaling} shows such an example.


Let $\lambda_{\min}(B)>0$ and $\lambda_{\max}(B)$ denote the smallest and the largest eigenvalues of $B$, respectively, and $\kappa_B=\lambda_{\max}(B)/\lambda_{\min}(B)$ be the condition number of $B$.

\begin{proposition} \label{prop:angle_condition1}
    Let $\tilde{g}(x) = B g(x)$, where $B$ is some positive definite matrix with the largest and the smallest eigenvalues $\lambda_{\max} \ge \lambda_{\min} > 0$. 
    Then, the following modified variant of conditions  \eqref{eq:angle_condition} holds
 \begin{align}
        \tilde{g}(x)^T \nabla\phi(x) &\ge  \tfrac{1}{2} \lambda_{\min}(B) \| g(x) \|_{2}^2 + \tfrac{1}{2}\left(  \lambda_{\min}(B) - \lambda_{\max}(B)\kappa_g^2 \right)\| \nabla \phi(x) \|_{2}^2 - \tfrac{1}{2}\lambda_{\max}(B)\epsilon_g^2\label{eq:angle_condition_compare_1}\\
        \|\tilde{g}(x)\|^2 &\le \lambda_{\max}(B)^2 \left( (1 + \kappa_g) \|\nabla\phi(x)\| + \epsilon_g\right)^2. \label{eq:angle_condition_compare_2}
    \end{align}
    for any $g(x)$ satisfying condition \eqref{eq:1_oracle_det} of Oracle~\ref{def:oracle_det_1}. Moreover, if $\|\nabla \phi(x)\|^2\geq \frac{\epsilon_g^2}{(1-\kappa_g^2)\tau^2}$ for  some $\tau\in[0,1)$ and $\kappa_B < \frac{\left( (1-\kappa_{g\tau})^2 + 1\right)}{\kappa_{g\tau}^2}$, with $\kappa_{g\tau}$ defined by \eqref{eq:kappa_gt}, the conditions of \eqref{eq:angle_condition} with $0 < \beta_1 \leq \beta_2$ hold, i.e.,
    \begin{align}
        \tilde{g}(x)^T \nabla\phi(x) &\ge \tfrac{1}{2}  \left( \lambda_{\min}(B) \left( (1-\kappa_{g\tau})^2 + 1\right) - \lambda_{\max}(B)\kappa_{g\tau}^2\right) \| \nabla \phi(x) \|_{2}^2\label{eq:angle_condition_compare_3}\\
        \|\tilde{g}(x)\|^2 &\le \lambda_{\max}(B)^2  (1 + \kappa_{g\tau})^2 \|\nabla\phi(x)\|^2. \label{eq:angle_condition_compare_4}
    \end{align}
    Finally, if $\kappa_B \leq (1-\tau_{g\tau})\frac{\left( (1-\kappa_{g\tau})^2 + 1\right)}{\kappa_{g\tau}^2}$ for some  $\tau_{g\tau} < 1 - \frac{\kappa_{g\tau}^2}{ 1+ (1- \kappa_{g\tau})^2}$, then
     \begin{align}
        \tilde{g}(x)^T \nabla\phi(x) &\ge \tfrac{1}{2}  \tau_{g\tau}  \lambda_{\min}(B) \left( (1-\kappa_{g\tau})^2 + 1\right)  \| \nabla \phi(x) \|_{2}^2.\label{eq:angle_1}
    \end{align}
\end{proposition}
\begin{proof}
 Inequality \eqref{eq:angle_condition_compare_2} can be derived using the largest eigenvalue of $B$ and inequality \eqref{eq:1_oracle_det_nosqr}, 
\begin{align*}
        \|\tilde{g}(x)\|^2 &= \|B g(x)\|^2 \\
        	& \le \lambda_{\max}(B)^2 \|g(x)\|^2  \\
        &\le \lambda_{\max}(B)^2 \left( (1 + \kappa_g) \|\nabla\phi(x)\| + \epsilon_g\right)^2. 
    \end{align*}

Next, we derive the inequality in \eqref{eq:angle_condition_compare_1}. Let $\| \cdot\|_{B}$ denote a ${B}$-weighted vector norm. It follows that 
\begin{align*}
	\| g(x) - \nabla \phi(x) \|_{B}^2 &= \| g(x)\|_{B}^2 - 2g(x)^TB\nabla \phi(x) + \| \nabla \phi(x) \|_{B}^2\\
	& = \| g(x)\|_{B}^2 - 2\tilde{g}^T\nabla \phi(x) + \| \nabla \phi(x) \|_{B}^2
\end{align*}
and, by Oracle~\ref{def:oracle_det_1}
\begin{align*}
	\| g(x) - \nabla \phi(x) \|_{B}^2 &\leq \lambda_{\max}(B) \| g(x) - \nabla \phi(x) \|_{2}^2  \leq \lambda_{\max}(B)\epsilon_g^2 + \lambda_{\max}(B)\kappa_g^2\| \nabla \phi(x) \|_{2}^2.
\end{align*}
Thus, 
\begin{align*}
	\| g(x)\|_{B}^2 - 2\tilde{g}^T\nabla \phi(x) + \| \nabla \phi(x) \|_{B}^2 & = \| g(x) - \nabla \phi(x) \|_{B}^2 \\
	& \leq  \lambda_{\max}(B)\epsilon_g^2 + \lambda_{\max}(B)\kappa_g^2\| \nabla \phi(x) \|_{2}^2.
\end{align*}
Re-arranging the above
\begin{align*}
	2\tilde{g}^T\nabla \phi(x) & \geq  \| g(x)\|_{B}^2 +  \| \nabla \phi(x) \|_{B}^2 - \lambda_{\max}(B)\kappa_g^2\| \nabla \phi(x) \|_{2}^2 - \lambda_{\max}(B)\epsilon_g^2\\
	& \geq \lambda_{\min}(B)\| g(x)\|_{2}^2 + \lambda_{\min}(B)\| \nabla \phi(x) \|_{2}^2 - \lambda_{\max}(B)\kappa_g^2\| \nabla \phi(x) \|_{2}^2 - \lambda_{\max}(B)\epsilon_g^2,
\end{align*}
from which \eqref{eq:angle_condition_compare_1} follows.

Moreover, if $\|\nabla \phi(x)\|^2\geq \frac{\epsilon_g^2}{(1-\kappa_g^2)\tau^2}$ for  some $\tau\in(0,1]$, by \eqref{eq:kappa_gt} and \eqref{eq:1_oracle_det_nosqr_noepsg} it follows that

\begin{align*}
	2\tilde{g}^T\nabla \phi(x) & 	\geq \lambda_{\min}(B)\| g(x)\|_{2}^2 + \lambda_{\min}(B)\| \nabla \phi(x) \|_{2}^2 - \lambda_{\max}(B)\kappa_{g\tau}^2\| \nabla \phi(x) \|_{2}^2\\
	& \geq \lambda_{\min}(B) (1-\kappa_{g\tau})^2\| \nabla \phi(x)\|_2^2  + \lambda_{\min}(B)\| \nabla \phi(x) \|_{2}^2 - \lambda_{\max}(B)\kappa_{g\tau}^2\| \nabla \phi(x) \|_{2}^2 \\
	& = \left( \lambda_{\min}(B) \left( (1-\kappa_{g\tau})^2 + 1\right) - \lambda_{\max}(B)\kappa_{g\tau}^2\right) \| \nabla \phi(x) \|_{2}^2 
\end{align*}
When 
\begin{equation}\label{eq:cond_kappaB_1}
	\kappa_B < \frac{\left( (1-\kappa_{g\tau})^2 + 1\right)}{\kappa_{g\tau}^2},
\end{equation} 
the right-hand-side of the above inequality is guaranteed to be positive, and the conditions of \eqref{eq:angle_condition} hold with $0< \beta_1 \leq \beta_2$.

Finally, if $\kappa_B \leq (1-\tau_{g\tau})\frac{\left( (1-\kappa_{g\tau})^2 + 1\right)}{\kappa_{g\tau}^2}$, \eqref{eq:angle_1} follows. This condition on $\kappa_B$ is only meaningful if the right-hand-side is greater than $1$, i.e., the condition number upper bound is larger than 1. As a result, $\tau_{g\tau} < 1 - \frac{\kappa_{g\tau}^2}{ 1+ (1- \kappa_{g\tau})^2}$.
\end{proof}

We make a few comments about this result. The right-hand-side of \eqref{eq:angle_1} is always positive, thus, under the condition 
\begin{equation}\label{eq:cond_kappaB}
    \kappa_B \leq (1-\tau_{g\tau})\frac{\left( (1-\kappa_{g\tau})^2 + 1\right)}{\kappa_{g\tau}^2},
\end{equation}
the vector $\tilde g(x)=Bg(x)$ satisfies \eqref{eq:angle_condition} with
\begin{equation}
\beta_1= \tfrac{1}{2}\tau_{g\tau} \lambda_{\min}(B) \left( (1-\kappa_{g\tau})^2 + 1\right),\quad \text{and} \quad  \beta_2^2= \lambda_{\max}(B)^2  (1 + \kappa_{g\tau})^2,
\end{equation}
where $\kappa_{g\tau}^2 = \tau^2+(1-\tau^2)\kappa_g^2 <1$ and $\tau_{g\tau} < 1 - \frac{\kappa_{g\tau}^2}{ 1+ (1- \kappa_{g\tau})^2}$,  and consequently $\bar g(x)=\frac{\beta_1}{\beta_2^2}Bg(x)$ satisfies conditions of Oracle~\ref{def:oracle_det_1} with $\kappa_g^2=1-\frac{\beta_1^2}{\beta_2^2}$
and $\epsilon_g=0$. The condition, $\kappa_B < \frac{\left( (1-\kappa_{g\tau})^2 + 1\right)}{\kappa_{g\tau}^2}$, holds trivially when $\kappa_B=1$, which is expected since this means that $B$ is a multiple of the  identity and $\tilde g(x)$ is simply a scaling of $g(x)$. When $\kappa_g=\epsilon_g=0$, the condition holds for any positive definite matrix, which is also expected since in this case $\tilde g(x)=B\nabla \phi(x)$ and it is well known that
multiplying the gradient by a symmetric positive definite matrix gives a direction that makes an acute angle with the gradient. 

In summary, we have shown that: 
\begin{itemize}
\item Any direction $\tilde g(x)$ that satisfies \eqref{eq:angle_condition} gives rise to Oracle \ref{def:oracle_det_1} by simply scaling the direction $\tilde g(x)$ appropriately. 
\item  Conversely, given $g(x)$, an output of Oracle \ref{def:oracle_det_1}, with $\kappa_g<1$ and $\epsilon_g>0$, and \eqref{eq:nablaphi_lower_bound}, then $\tilde g(x)=Bg(x)$ satisfies  \eqref{eq:angle_condition} for any positive definite $B$, whose condition number satisfies \eqref{eq:cond_kappaB_1}.
\end{itemize}

\subsubsection{Inexact Oracles for Convex Functions} 
 In \cite{Devolder2013FirstorderMO} a different type of inexact oracle was proposed within the context of first-order algorithms for convex problems. Specifically, $f_{({m}, M)}(x)\in \mathbb{R}$ and $g_{({m}, M)}(x)\in \mathbb{R}^n$ are an output of a $({m}, M)$-oracle if for all $x \in \mathbb{R}^n$ and $y \in \mathbb{R}^n$, 
\begin{equation}\label{eq:tauLoracle}
0\leq \phi(y)-f_{({m}, M)}(x)-g_{({m}, M)}(x)^T(y-x)\leq \frac{M}{2}\|y-x\|^2 +{m}.
\end{equation}
Such function and gradient approximations satisfy the following bounds for any $x$
\begin{align*}
	f_{({m}, M)}(x)\leq \phi(x) &\leq f_{({m}, M)}(x)+{m}\\
	\text{and } \qquad \|\nabla \phi(x)-g_{({m}, M)}(x)\| &\leq \sqrt{2M{m}}.
\end{align*}

On the other hand, the reverse is not true, as is also shown in \cite{Devolder2013FirstorderMO}. In general, Oracles  \ref{def:oracle_det_1} implies \eqref{eq:tauLoracle}  
for $x$ and $y$ from a bounded set. This can be simply derived as follows: let $g(x)$ be an output of Oracles  \ref{def:oracle_det_1} with $\kappa_g=0$. Then,  for any $x,y:\, \|x-y\|\leq D$ from $L_1$-smoothness of $\phi$
\begin{equation}\label{eq:relaxed_ub}
\phi(y)\leq \phi(x)+\nabla \phi(x)^T(y-x)+\frac{L_1}{2}\|y-x\|^2\leq \phi(x)+g(x)^T(y-x)+\epsilon_gD+\frac{L_1}{2}\|y-x\|^2
\end{equation}
and from convexity 
\begin{equation}\label{eq:relaxed_convex}
\phi(y)\geq \phi(x)+\nabla \phi(x)^T(y-x)\geq \phi(x)+g(x)^T(y-x)-\epsilon_gD. 
\end{equation}
It follows that  \eqref{eq:tauLoracle} holds with $M=L_1$, ${m}=2\epsilon_g D$, $g_{({m}, M)}(x)=g(x)$ and  
$f_{({m}, M)}(x)=\phi(x)-\epsilon_gD$. This constructs a restricted $(m,M)$-oracle, not necessarily one valid of all $x,y\in {\mathbb R}^n$.

We note that $\phi(x)$ may not be computable exactly, so $f_{({m}, M)}(x)$ is not computable. But it is not needed for implementation of the fixed step gradient descent algorithm.

\subsection{A Descent Inequality and Complexity Analysis}\label{sec:descent_lemma}
We now analyze algorithms that rely on our deterministic inexact oracles. The analysis is built around a single key inequality, which we term the descent inequality\footnote{In contrast to the standard ‘descent lemma’ in the literature \cite{bertsekas2016nonlinear}.}. This inequality yields classical convergence results for fixed step gradient descent and line-search algorithms in convex, strongly convex, and nonconvex settings, as well as for general smooth functions and functions satisfying Assumptions \ref{assum:KLcond} and \ref{assum:PLcond}. We first state this inequality in abstract terms  and show how the complexity results follow from it. We then demonstrate that the descent inequality (with specified constants) can be ensured by our oracles.

Let us assume we have an algorithm that generates a sequence of points $\{x_k,\, k=0,1, \ldots, \}$ such that 
the following descent inequality holds on every iteration, 
\begin{equation}\label{eq:key_reduction}\tag{DI}
	\phi(x_{k+1}) \leq \phi (x_k) - \alpha \theta \| \nabla \phi (x_k) \|_2^2 + \alpha \omega 
\end{equation}
for some $\alpha>0$,  $\theta >0$ and $\omega\geq0$. From \eqref{eq:key_reduction}, complexity bounds for any such algorithm follow directly in the convex, strongly convex, and general Lipschitz smooth settings. 
 
 In what follows, we again employ a constant $\tau \in (0,1)$ to conveniently define the lower bounds on the best achievable accuracy in each setting. When $\omega=0$, $\tau$ can be chosen arbitrarily close to $0$ and thus removed from the complexity bounds. Below we derive the bounds on $K_{\epsilon_\tau}$ which is the number of iterations 
 within which the desired accuracy is guaranteed. For convenience, we assume that the desired accuracy is not reached at the starting point $x_0$ when deriving the 
 respective results (otherwise they trivially hold). 
 
 The first theorem shows the complexity bound for such an algorithm applied to general (i.e. not necessarily convex) functions  based on the descent inequality \eqref{eq:key_reduction}.

\begin{theorem}\label{thm:nonconvex_bound}
    Suppose Assumptions~\ref{assum:low_bound} and \ref{assum:lip_cont} hold, and \eqref{eq:key_reduction} holds for the sequence $\{x_k,\, k=0,1, \ldots, \}$.  Then, for any $\tau\in (0,1)$
and $\epsilon_\tau \geq  \sqrt{\frac{\omega}{ \tau^2\theta}}$,  $\| \nabla \phi (x_k) \| \leq \epsilon_\tau$ is guaranteed to be satisfied at least once during $K_{\epsilon_\tau}$ iterations 
where 
\begin{align*}
	K_{\epsilon_\tau} = \left\lceil  \frac{1}{\epsilon_\tau^2} \frac{(\phi(x_0) - \hat{\phi})}{ \alpha\theta(1 - \tau^2)} \right\rceil.
\end{align*}
\end{theorem}

\begin{proof} 
By summing \eqref{eq:key_reduction} over  $k=0, \ldots, K_{\epsilon_\tau} - 1$, for some $K_{\epsilon_\tau}>0$,  and rearranging terms it follows that 
      \begin{align*}
      \theta \sum_{k=0}^{K_{\epsilon_\tau}-1} \|\nabla\phi(x_k)\|^2 
        \le \frac{1}{\alpha} \big(\phi(x_0) - \phi(x_{K_{\epsilon_\tau}})\big) + K_{\epsilon_\tau} \omega
        \le \frac{1}{\alpha} \big(\phi(x_0) - \hat{\phi}\big) + K_{\epsilon_\tau} \omega, 
    \end{align*} 
    where the last inequality holds by Assumption \ref{assum:low_bound}. Dividing both sides by $\theta K_{\epsilon_\tau}$, 
      \begin{equation}\label{eq:nonconvex_bound}
        \frac{1}{K_{\epsilon_\tau}} \sum_{k=0}^{K_{\epsilon_\tau}-1} \|\nabla\phi(x_k)\|^2 \le \frac{1}{\theta} \left[ \frac{\phi(x_0) - \hat{\phi}}{\alpha K_{\epsilon_\tau}} + {\omega} \right]. 
    \end{equation}
By  \eqref{eq:nonconvex_bound} and the expression for $K_{\epsilon_\tau}$, it follows that
 \begin{align*}
     \min\{ \|\nabla\phi(x_k)\|^2, k=0, \ldots K_{\epsilon_\tau}-1 \}&\leq  \frac{1}{K_{\epsilon_\tau}} \left[ \frac{(\phi(x_0) - \hat{\phi})}{\theta\alpha } \right] + \frac{\omega}{\theta}\\
     & \leq \epsilon_\tau^2 (1-\tau^2)+\epsilon_\tau^2 \tau^2\leq \epsilon_\tau^2.
    \end{align*}
\end{proof}

We now turn to the case of functions for which Assumption \ref{assum:KLcond} holds, which includes convex functions with bounded level sets. 


\begin{theorem}\label{thm:KL_cond_rate}
Suppose Assumption \ref{assum:lip_cont} holds, Assumption \ref{assum:KLcond} holds with a constant $c=\frac{1}{D}$, for some $D>0$, and \eqref{eq:key_reduction} holds for the sequence $\{x_k,\, k=0,1, \ldots, \}$.
Then, for any $\tau\in (0,1)$
and $\epsilon_\tau \geq \sqrt{\frac{\omega}{\tau^2\theta}}{D}$, $ \phi (x_k) - \phi^* \leq \epsilon_\tau$ is guaranteed to be satisfied at least once after $K_{\epsilon_\tau}$ iterations where
\begin{align*}
	K_{\epsilon_\tau} = \left\lceil \frac{1}{\epsilon_\tau} \frac{ {D}^2}{\alpha \theta (1-\tau^2)}. \right\rceil
\end{align*}
\end{theorem}
\begin{proof} 
By Assumption \ref{assum:KLcond} applied to the iterates $\{x_k,\, k=0,1, \ldots, \}$
\begin{equation}\label{eq:key_reduction3_det}
 	\phi(x_{k}) - \phi^* \leq {D} \| \nabla \phi (x_k) \|.
\end{equation}
Thus, by \eqref{eq:key_reduction} and \eqref{eq:key_reduction3_det} 
\begin{equation}\label{eq:key_reduction4_det}
 	\phi(x_{k+1}) - \phi^* 
	\leq \phi (x_k) - \phi^* 
	-  \frac{\alpha \theta}{{D}^2} \left( \phi (x_k) - \phi^* \right)^2+\alpha\omega.
\end{equation}

If $\phi (x_{k+1}) - \phi^* =0$  or  $\phi (x_k) - \phi^* \leq \epsilon_\tau$ for  some $k=0, \ldots, K_{\epsilon_\tau}-1$ then the statement of the theorem holds. 
Otherwise, for any $k=0, \ldots, K_{\epsilon_\tau}-1$, from   $\phi (x_k) - \phi^* >\epsilon_\tau$ and  from \eqref{eq:key_reduction4_det} it follows that
\begin{align*}
 	\frac{1}{\phi(x_{k+1}) - \phi^*} - \frac{1}{\phi(x_{k}) - \phi^*} &= \frac{\phi(x_{k}) - \phi^* - (\phi(x_{k+1}) - \phi^*)}{(\phi(x_{k}) - \phi^*)(\phi(x_{k+1}) - \phi^*)} \\
	& \geq \frac{\phi(x_{k}) - \phi^* - (\phi(x_{k+1}) - \phi^*)}{(\phi(x_{k}) - \phi^*)^2} \\
	& \geq  \frac{\alpha \theta}{{D}^2}-\frac{\alpha\omega}{(\phi(x_{k}) - \phi^*)^2} \\
	& \geq  \frac{\alpha \theta (1-\tau^2)}{{D}^2}.
\end{align*}
By summing the above inequality over  $k=0, \ldots, K_{\epsilon_\tau}-1$, it follows that 
      \begin{align*}
      		\frac{1}{\phi(x_{K_{\epsilon_\tau}}) - \phi^*} - \frac{1}{\phi(x_{0}) - \phi^*} &\geq \sum_{k=0}^{K_{\epsilon_\tau}-1} \left( \frac{1}{\phi(x_{k+1}) - \phi^*} - \frac{1}{\phi(x_{k}) - \phi^*} \right)  \\
		& \geq \frac{\alpha \theta (1-\tau^2)K_{\epsilon_\tau}}{{D}^2}, 
    \end{align*} 
from which it follows by rearranging the above that
\begin{align}\label{eq:key_reduction5}
	\phi(x_{K_{\epsilon_\tau}}) - \phi^*  \leq \frac{{D}^2}{\alpha \theta(1-\tau^2) K_{\epsilon_\tau}}.
\end{align}
The result follows by setting the right-hand-side equal to $\epsilon_\tau$ and solving for $K_{\epsilon_\tau}$.
\end{proof}

Let us discuss how Theorem \ref{thm:KL_cond_rate} can be applied to convex functions. As shown in Section \ref{sec:assump} a convex function whose  level set ${\cal D}=\{x:\, \phi(x)\leq \phi(x_0)\}$ has a bounded radius $D$ satisfies Assumption \ref{assum:KLcond} on that set. However, Theorem \ref{thm:KL_cond_rate}
 requires that Assumption \ref{assum:KLcond} is satisfied by  the iterate sequence $\{x_k,\, k=0,1, \ldots, K_{\epsilon_{\tau}} \}$. If this sequence lies in ${\cal D}$ then 
 Theorem \ref{thm:KL_cond_rate} directly applies. However, because  \eqref{eq:key_reduction}  does not guarantee that $\phi(x_{k+1})\leq \phi(x_k)$, the sequence does not have to remain in ${\cal D}$.  The argument can be made that 
 descent in the objective function value is guaranteed while the gradient is sufficiently large, i.e., $\| \nabla \phi (x_k) \| >  \sqrt{\frac{\omega}{ \theta}}$. On the other hand, a small gradient 
 can only occur when $ \phi (x_k) - \phi^* \leq \epsilon_\tau$.
 Thus we have the following result.

 \begin{theorem}\label{thm:convex_rate}
Suppose Assumptions~\ref{assum:lip_cont} and ~\ref{assum:convex} hold, and \eqref{eq:key_reduction} holds for the sequence $\{x_k,\, k=0,1, \ldots, \}$.  
Then, for any $\tau\in (0,1)$
and $\epsilon_\tau \geq  \sqrt{\frac{\omega}{ \tau^2\theta}}D$, $ \phi (x_k) - \phi^* \leq \epsilon_\tau$ is guaranteed to be satisfied at least once after $K_{\epsilon_\tau}$ iterations where
\begin{align*}
		K_{\epsilon_\tau} = \left\lceil \frac{1}{\epsilon_\tau} \frac{ {D}^2}{\alpha \theta (1-\tau^2)} \right\rceil.
\end{align*}
\end{theorem}
\begin{proof}
Let us assume that for all iterations $i = 0,\ldots, k - 1$, $\| \nabla \phi (x_i) \| >\frac{\epsilon_\tau}{D}\geq \sqrt{\frac{\omega}{ \tau^2\theta}}$ holds, for some $k\geq 1$. Then by \eqref{eq:key_reduction}, $\phi(x_{i+1})\leq \phi(x_i)$ for all  $i=0, \ldots, k-1$ and hence $\{x_i, \, i=0,\ldots, k \}\subset {\cal D}$ which means that by Assumption~\ref{assum:convex}, Assumption \ref{assum:KLcond} is also satisfied. If $k = K_{\epsilon_\tau}$ then Theorem \ref{thm:KL_cond_rate} applies. 
On the other hand, if $\| \nabla \phi (x_k) \| \leq \frac{\epsilon_\tau}{D}$ occurs for the first time for some $k < \frac{1}{\epsilon_\tau} \frac{ {D}^2}{\alpha \theta (1-\tau^2)}$, then since $x_k\in {\cal D}$, then
$ \phi (x_k) - \phi^* \leq \| \nabla \phi (x_k) \|D\leq \epsilon_\tau$, and, thus, the result again follows. 
\end{proof}

The next result derives a complexity bound for the functions satisfying Assumption \ref{assum:PLcond}, which include strongly convex functions.

\begin{theorem}\label{thm:str_convex_rate} 
Suppose Assumptions~\ref{assum:lip_cont} and ~\ref{assum:PLcond} hold, and \eqref{eq:key_reduction} holds for the sequence $\{x_k,\, k=0,1, \ldots, \}$.  Assume also that $ 2 \alpha \theta \mu \leq 1$. Then, for any $\tau\in (0,1)$
and $\epsilon_\tau \geq \frac{\omega}{2\theta \mu \tau^2}$, 
$ \phi (x_k) - \phi^* \leq \epsilon_\tau$ is guaranteed to be satisfied after $K_{\epsilon_\tau}$ iterations where
\begin{align*}
	K_{\epsilon} = \left\lceil \log_{\frac{1}{1 - 2\alpha \mu \theta(1-\tau^2)}} \left( \frac{\phi (x_0) - \phi^*}{\epsilon_\tau}\right) \right\rceil \approx \left\lceil \frac{1}{ 2\alpha \mu \theta(1-\tau^2)}\log \left( \frac{\phi (x_0) - \phi^*}{\epsilon_\tau}\right) \right\rceil.
\end{align*}
\end{theorem}
\begin{proof} For all iterates $\{x_k,\, k=0,1, \ldots, \}$ for which $ \phi (x_k) - \phi^* > \epsilon_\tau$, by \eqref{eq:key_reduction} and Assumption \ref{assum:PLcond} with $\mu>0$, it follows that
\begin{align*}
	\phi(x_{k+1}) &\leq \phi (x_k) - \alpha \theta  \| \nabla \phi (x_k) \|_2^2  + \alpha \omega \\
	& \leq \phi (x_k) - 2\alpha\theta  \mu \left( \phi (x_k) - \phi^*\right) + \alpha \omega\\ 
	& \leq \phi (x_k) - 2\alpha\theta \mu (1-\tau^2) \left( \phi (x_k) - \phi^*\right). 
\end{align*}
Recursive application of the above yields, 
\begin{align*}
	\phi(x_{k}) - \phi^* &\leq \left( 1 - 2 \alpha \theta \mu (1-\tau^2)\right)^k \left(\phi (x_0) - \phi^* \right).
\end{align*}
The result follows by setting the right-hand-side equal to $\epsilon_\tau>0$ and solving for $k$.
\end{proof}

Theorem~\ref{thm:str_convex_rate} is only meaningful if the step size is positive and sufficiently small, i.e., $0 < \alpha \leq \frac{1}{2\theta \mu (1-\tau^2)}$. This  will be implied automatically by a simpler bound $\alpha\leq \frac{1}{L_1}$ under which \eqref{eq:key_reduction}  will be derived for the gradient descent algorithm.

\subsection{Fixed Step Algorithm based on Oracle \ref{def:oracle_det_1}}\label{sec:fixed_step_det}

We now introduce a basic fixed step inexact gradient descent algorithm that utilizes outputs of Oracle \ref{def:oracle_det_1}, and 
show that descent inequality \eqref{eq:key_reduction} holds, $\theta = \frac{1 - \kappa_g^2}{2}$, and $\omega = \frac{\epsilon_g^2}{2}$ and thus complexity and 
convergence properties can be obtained from Theorems \ref{thm:nonconvex_bound}, \ref{thm:convex_rate}  and \ref{thm:str_convex_rate}. 

\begin{algorithm}[H]
\caption{~\textbf{Inexact Gradient Method with Fixed Step Size}}
\label{alg:first_order_fixed_step}
\begin{algorithmic}
\State \textbf{Inputs:} $x_0 \in \mathbb{R}^n$ (starting point), $\alpha > 0$ (step size), Oracle \ref{def:oracle_det_1} with $\kappa_g<1$
 \State \textbf{1:} \textbf{for $k=0,1,2,\dots$ do}
\State \text{2:}~~~~Compute  $g_k=g(x_k)\approx \nabla \phi(x_k)$ via Oracle \ref{def:oracle_det_1}.  
 \State \text{3:}~~~~Update iterate $x_{k+1}\gets x_k-\alpha g_k$.
\State \textbf{4:} \textbf{end for}
 \end {algorithmic}
\end{algorithm}

We derive the following lemma that guarantees  \eqref{eq:key_reduction} for each iteration of Algorithm \ref{alg:first_order_fixed_step} in terms of constants of Oracle \ref{def:oracle_det_1} ($\epsilon_g$ and $\kappa_g$). 

\begin{lemma}\label{lem.key_lemma_det} Suppose Assumption~\ref{assum:lip_cont} holds. Let the step size parameter in Algorithm~\ref{alg:first_order_fixed_step} satisfy $\alpha \le \frac{1}{L_1}$. Then, for all $k\geq 0$,
\begin{equation} \label{eq.descent_inexact}
	\phi(x_{k+1}) \leq \phi (x_k) - \frac{\alpha}{2} (1-\kappa_g^2) \| \nabla \phi (x_k) \|_2^2 + \alpha\frac{\epsilon_g^2}{2}.
\end{equation}
\end{lemma}
\begin{proof} 
    By Assumption~\ref{assum:lip_cont} it follows that
    \begin{align*}
        \phi(x_{k+1}) &\leq \phi (x_k) - \alpha g_k^T \nabla \phi (x_k) + \frac{\alpha^2 L_1}{2} \| g_k\|_2^2,
    \end{align*}
    and by the definition of Oracle~\ref{def:oracle_det_1}
    \begin{align*}
        -2g_k^T \nabla \phi (x_k) \leq \epsilon^2_g - \left(1-\kappa_g^2\right) \|\nabla \phi(x_k)\|^2 - \|g_k\|^2. 
    \end{align*}
    Combining the above, and by $\alpha \leq 1/L_1$, it follows that
     \begin{align*}
        \phi(x_{k+1}) &\leq \phi (x_k) - \alpha g_k^T \nabla \phi (x_k) + \frac{\alpha^2 L_1}{2} \| g_k\|_2^2\\
       &\leq \phi (x_k) + \frac{\alpha}{2} \left( \epsilon^2_g - \left(1-\kappa_g^2\right) \|\nabla \phi(x_k)\|^2 - \|g_k\|^2\right) + \frac{\alpha^2 L_1}{2} \| g_k\|_2^2\\
       &= \phi (x_k) - \frac{\alpha}{2}\left(1-\kappa_g^2\right) \|\nabla \phi(x_k)\|^2   - \frac{\alpha}{2}\left(1 - \alpha L_1\right) \| g_k\|_2^2 + \frac{\alpha}{2}  \epsilon^2_g\\
       &\leq \phi (x_k) - \frac{\alpha}{2}\left(1-\kappa_g^2\right) \|\nabla \phi(x_k)\|^2 + \frac{\alpha}{2}  \epsilon^2_g,
    \end{align*} 
which completes the proof.
\end{proof}

Lemma~\ref{lem.key_lemma_det} shows that the descent inequality~\eqref{eq:key_reduction} holds whenever $\kappa_g^2<1$. In particular, the objective function decreases whenever $\|\nabla\phi(x_k)\|>\frac{\epsilon_g}{\sqrt{1-\kappa_g^2}}$. Combining Lemma~\ref{lem.key_lemma_det} with Theorems~\ref{thm:nonconvex_bound}, \ref{thm:KL_cond_rate}, \ref{thm:convex_rate}, and \ref{thm:str_convex_rate} yields the following complexity bounds for Algorithm~\ref{alg:first_order_fixed_step} under different assumptions on $\phi$.

\begin{theorem}\label{cor:all_complexity_fixed}
The following complexity bounds follow. 
\begin{itemize}
	\item Suppose $\phi$ is bounded below and $L_1$-smooth (Assumptions \ref{assum:low_bound} and \ref{assum:lip_cont}). Then, Algorithm~\ref{alg:first_order_fixed_step}, applied with a step size $\alpha \leq  \frac{1}{L_1}$, generates an iterate $x_k$ satisfying $\| \nabla \phi (x_k) \| \leq \epsilon_\tau$, with 
$\epsilon_\tau \geq \sqrt{\frac{\epsilon_g^2}{ \tau^2(1-\kappa_g^2)}}$ and $\tau \in (0,1)$, within at most 
\begin{align*}
	K_{\epsilon_\tau} = \left\lceil \frac{1}{\epsilon_\tau^2} \frac{ 2(\phi(x_0) - \hat{\phi})}{\alpha(1-\kappa_g^2)(1 - \tau^2)} \right\rceil
\end{align*}
iterations. Moreover, if $\alpha = \frac{1}{L_1}$, then
\begin{align*}
	K_{\epsilon_\tau} = \left\lceil\frac{1}{\epsilon_\tau^2} \frac{ 2L_1(\phi(x_0) - \hat{\phi})}{(1-\kappa_g^2)(1 - \tau^2)} \right\rceil.
\end{align*}
	\item Additionally, suppose $\phi$ satisfies  Assumption~\ref{assum:KLcond} with $c=\frac{1}{D}$ on the set containing all relevant iterates, e.g., if $\phi$ is convex, with bounded level sets as in  Assumption~\ref{assum:convex}. Then, Algorithm~\ref{alg:first_order_fixed_step}, applied with a step size $\alpha \leq  \frac{1}{L_1}$, generates an iterate satisfying $\phi (x_{k}) -\phi^* \leq \epsilon_\tau$, with 
$\epsilon_\tau \geq \sqrt{\frac{\epsilon_g^2}{ \tau^2(1-\kappa_g^2)}}{D}$ and $\tau \in (0,1)$, within at most 
	\begin{align*}
		K_{\epsilon_\tau} = \left\lceil \frac{1}{\epsilon_\tau} \frac{4D^2}{\alpha(1-\kappa_g^2)(1 - \tau^2)} \right\rceil
	\end{align*}
	iterations. Moreover, if $\alpha = \frac{1}{L_1}$, then
\begin{align*}
	K_{\epsilon_\tau} = \left\lceil \frac{1}{\epsilon_\tau} \frac{4L_1{D}^2}{(1-\kappa_g^2)(1 - \tau^2)} \right\rceil.
\end{align*}
	\item Additionally, suppose $\phi$ satisfies Assumption~\ref{assum:PLcond} on the set containing all relevant iterates, e.g., if $\phi$ is  $\mu$-strongly  convex as in  Assumption~\ref{assum:str.convex}. Then, Algorithm~\ref{alg:first_order_fixed_step}, applied with a step size $\alpha \leq  \frac{1}{L_1}$, generates an iterate satisfying $\phi (x_k) -\phi^* \leq \epsilon_\tau$, with $\epsilon_\tau \geq \frac{\epsilon_g^2}{2(1-\kappa_g^2) \mu \tau^2}$ and $\tau \in (0,1)$, within at most 
	\begin{align*}
	K_{\epsilon_\tau} = \left\lceil \log_{\frac{1}{1 -  \alpha \mu (1-\kappa_g^2)(1-\tau^2)}} \left( \frac{\phi (x_0) - \phi^*}{\epsilon_\tau}\right) \right\rceil
\end{align*}
iterations. Moreover, if $\alpha = \frac{1}{L_1}$, then
\begin{align*}
	K_{\epsilon_\tau} = \left\lceil \log_{\frac{1}{1 - \frac{ (1-\kappa_g^2)(1-\tau^2)\mu}{L_1}}} \left( \frac{\phi (x_0) - \phi^*}{\epsilon_\tau}\right) \right\rceil
	\approx \left\lceil \frac{1}{(1-\kappa_g^2)(1-\tau^2)} \frac{L_1}{\mu}\log \left( \frac{\phi (x_0) - \phi^*}{\epsilon_\tau}\right) \right\rceil.
\end{align*}
\end{itemize}
\end{theorem}

We now show that under the assumption $\|\nabla \phi(x_k)\|^2\geq \frac{\epsilon_g^2}{(1-\kappa_g^2)\tau^2}$ where $\tau \in (0,1)$, the complexity bound can be alternatively derived, by the simple approach proposed earlier; $g(x)$ can be seen as an output of Oracle   \ref{def:oracle_det_1} with $\epsilon_g=0$ and modified relative error parameter $\kappa_{g\tau}^2= \tau^2+(1-\tau^2)\kappa_g^2$. As a result, the descent inequality is
\begin{align}\label{eq.result2}
	\phi (x_k - \alpha g_k) &\leq \phi (x_k) - \frac{\alpha}{2}\left(1-\kappa_g^2\right) \left(1-\tau^2\right) \|\nabla \phi(x_k)\|^2 = \frac{\alpha}{2}\left(1-\kappa_{g\tau}^2\right) \|\nabla \phi(x_k)\|^2.
\end{align}
All of the above derivations 
can then be performed with these new constants 
obtaining exactly the same complexity bounds as we list above. Note that  in the general smooth case $\|\nabla \phi(x_k)\|^2\geq \frac{\epsilon_g^2}{(1-\kappa_g^2)\tau^2}$ holds automatically for all $k<K_{\epsilon_\tau}$ by the definition of $K_{\epsilon_\tau}$. Under the linear KL and PL conditions $\|\nabla \phi(x_k)\|^2\geq \frac{\epsilon_g^2}{(1-\kappa_g^2)\tau^2}$ also holds by applying these conditions and using respective definitions of   $K_{\epsilon_\tau}$.
This shows that this initial simplification does not cause loss in efficiency of the results while significantly simplifying the derivations.

\paragraph{Termination condition} It is common in the optimization literature to state an algorithm as an infinite sequence of iterates and to analyze the behavior of that sequence, as we did above. However, in practice, an algorithm typically requires some termination criterion. Without knowing the target objective function value a typical criterion may be based on the norm of the gradient. While our theory shows that achieving $\| \nabla \phi (x_k) \| \leq \epsilon_\tau$ is possible in the nonconvex case for any $\tau\in(0,1)$, a termination criterion can only be imposed using $\|g(x)\|$, which is computable, rather than $\| \nabla \phi (x_k) \|$. Specifically, we  can impose the following termination condition for Algorithm \ref{alg:first_order_fixed_step}
\begin{align}\label{eq.term_cond}
	\| g(x_k)\| \leq  \epsilon_{g\tau}.
\end{align}
for some chosen $\epsilon_{g\tau} \geq 0$. By \eqref{eq:1_normg_det}, it follows that until the termination condition \eqref{eq.term_cond} is met
\begin{equation*}
	\|\nabla \phi(x_k)\|\geq  \frac{ \epsilon_{g\tau}}{1 + \kappa_g}-\epsilon_g
\end{equation*}
Selecting $\epsilon_{g\tau}$ that satisfies
\begin{equation}\label{eq:epsilongtau_bound}
	\epsilon_{g\tau}\geq \epsilon_g(1+\kappa_g)\left(1+\tfrac{1}{\sqrt{(1-\kappa_g^2)}\tau}\right)
\end{equation}
ensures that until the termination condition is met 
\begin{equation*}
\|\nabla \phi(x_k)\|\geq \frac{\epsilon_g}{\sqrt{(1-\kappa_g^2)}\tau}=\epsilon_{\tau}
\end{equation*}
which is equivalent to \eqref{eq:nablaphi_lower_bound}. This will be used to ensure that until the termination condition is met the algorithm achieves descent. 

On the other hand, again by \eqref{eq:1_normg_det},  once the  termination condition is met we have
\[
\|\nabla \phi(x_k)\|\leq  \frac{ \epsilon_{g\tau}+\epsilon_g}{1 - \kappa_g}
\]
which  is not as tight as $\| \nabla \phi (x_k) \| \leq \epsilon_\tau$. However, this new condition is achieved at the terminal (final) iterate, rather than some iteration prior to the termination. 

We note that ensuring \eqref{eq:epsilongtau_bound} requires knowledge of $\epsilon_g$ and $\kappa_g$ or at least their upper bounds. Such knowledge is available for some specific oracle instances as we discussed in our example section (Section~\ref{sec.examples_first}).

\subsection{Fixed Step Algorithm based on Angle Condition}\label{sec:angle_cond_rates}

Next, we consider another variant of a fixed step inexact gradient descent algorithm that utilizes search directions that satisfy \eqref{eq:angle_condition}, and show that descent inequality \eqref{eq:key_reduction} holds with $\theta = \frac{\beta_1}{2}$ and $\omega = 0$, and thus complexity and 
convergence properties can again be obtained from Theorems \ref{thm:nonconvex_bound}, \ref{thm:convex_rate}  and \ref{thm:str_convex_rate}.

\begin{algorithm}[H]
\caption{~\textbf{Inexact Gradient Method with Fixed Step Size}}
\label{alg:first_order_fixed_step_angle}
\begin{algorithmic}
\State \textbf{Inputs:} $x_0 \in \mathbb{R}^n$ (starting point), $\alpha > 0$ (step size), constants $0<\beta_1\leq\beta_2$
 \State \textbf{1:} \textbf{for $k=0,1,2,\dots$ do}
\State \text{2:}~~~~Compute  $\tilde g_k=\tilde g(x_k)\approx \nabla \phi(x_k)$ such that \eqref{eq:angle_condition} holds.  
 \State \text{3:}~~~~Update iterate $x_{k+1}\gets x_k-\alpha \tilde g_k$.
\State \textbf{4:} \textbf{end for}
 \end {algorithmic}
\end{algorithm}

The following lemma is the version of key descent inequality \eqref{eq:key_reduction} under the angle condition \eqref{eq:angle_condition}. 

\begin{lemma}\label{lem.descent_angle} Suppose Assumption~\ref{assum:lip_cont} holds. Let the step size parameter in Algorithm~\ref{alg:first_order_fixed_step_angle} satisfy 
$\alpha \le \tfrac{\beta_1}{\beta_2^2L_1}$. Then, for all $k\geq 0$,
\begin{equation} \label{eq.descent_inexact}
\phi(x_{k+1}) \leq \phi (x_k) -  \frac{\alpha}{2} \beta_1\| \nabla \phi (x_k) \|_2^2. 
\end{equation}

\end{lemma}
\begin{proof} 
    By Assumption~\ref{assum:lip_cont}  and 
    \eqref{eq:angle_condition} 
     \begin{align*}
        \phi(x_{k+1}) &\leq \phi (x_k) - \alpha \tilde g_k^T \nabla \phi (x_k) + \frac{\alpha^2 L_1}{2} \| \tilde g_k\|_2^2\\
       &\leq \phi (x_k) - \alpha \beta_1 \|\nabla \phi(x_k)\|^2  + \frac{\alpha^2 L_1\beta_2^2}{2} \|\nabla \phi(x_k)\|_2^2\\
       &\leq \phi (x_k) - \frac{\alpha}{2}\left(2\beta_1-\alpha \beta_2^2L_1\right) \|\nabla \phi(x_k)\|^2.
    \end{align*} 
    The step size restriction, $\alpha \le \tfrac{\beta_1}{\beta_2^2L_1}$, completes the proof.
\end{proof}

We conclude that inequality \eqref{eq:key_reduction} holds with $\theta=\frac{\beta_1}{2}$ and $\alpha \le \beta_1/\beta_2^2L_1$.

Now recall, that by defining $g(x)=\frac{\beta_1}{\beta_2^2}\tilde g(x)$ we have shown that $g(x)$ obeys properties of Oracle   \ref{def:oracle_det_1} with
$\kappa^2_g=1-\frac{\beta_1^2}{\beta_2^2}\in [0, 1)$ and $\epsilon_g=0$. For any step $\alpha_2 \le \tfrac{\beta_1}{\beta_2^2L_1}$ used by Algorithm~\ref{alg:first_order_fixed_step_angle}, consider Algorithm~\ref{alg:first_order_fixed_step}  with step $\alpha_1$  such that
  $\alpha_2 \frac{\beta_1}{\beta_2^2}\alpha_1$. We get $\alpha_1 \le \tfrac{1}{L_1}$ and  \eqref{eq:key_reduction}  holds with
  $\theta=\frac{1-\kappa_g^2}{2}=\frac{\beta_1^2}{2\beta_2^2}$. 
Thus, Algorithm~\ref{alg:first_order_fixed_step_angle}, with any permissible step $\alpha_2$ is equivalent to (as it generates the same sequence of iterates) as
Algorithm~\ref{alg:first_order_fixed_step}  with apprpriately scaled  $\alpha_1$.

\subsection{Alternative Analysis of the  Fixed Step Size Inexact Gradient Descent for Convex Functions}\label{sec:telescopic}

Next, we present other results from the literature for the case of convex functions. In particular, we first present an analysis that relies on the classical ``telescopic-sum'' argument instead of descent inequality. We then discuss results in the setting for which the bound $D$ on the level set is not required, and compare the results to those presented in Theorem~\ref{thm:convex_rate}. 

\paragraph{Telescopic Sum Analysis for Convex Functions}

We now provide an alternative analysis of the convergence rate of Algorithm~\ref{alg:first_order_fixed_step} under the assumption that $\phi$ is convex. This derivation relies on the classical ``telescopic-sum'' argument and aligns with standard analyses in the literature. It does not directly rely on the descent inequality \eqref{eq:key_reduction}, and requires the definition of convexity rather than the use of the linear KL condition (Assumption~\ref{assum:KLcond}).  The effect of the gradient error on the convexity property does not appear to benefit from separating the additive and multiplicative components in Oracle~\ref{def:oracle_det_1}.
However, this analysis offers other advantages, as discussed below. In particular, the resulting complexity bound exhibits a more favorable dependence on $D$. The derivation presented here closely follows the analysis in \cite{Devolder2013FirstorderMO}, but is adapted to Oracle~\ref{def:oracle_det_1} rather than the oracle model considered in that work and described in Section~\ref{sec.other_oracles}. Overall, the final complexity result remains the same.

Given $g(x)$, an output of Oracle  \ref{def:oracle_det_1}, similarly to \eqref{eq:relaxed_ub}  
and \eqref{eq:relaxed_convex} derived for the case when $\kappa_g=0$, it can be shown that for any $x,y:\, \|x-y\|\leq D$, by Lipschitz-smoothness (Assumption~\ref{assum:lip_cont}), 
\begin{equation}\label{eq:relaxed_ub2}
\phi(y)\leq  \phi(x)+g(x)^T(y-x)+\epsilon_gD+\kappa_g\|\nabla \phi(x)\|D+\frac{L_1}{2}\|y-x\|^2
\end{equation}
and by convexity (Assumption~\ref{assum:convex}), 
\begin{equation}\label{eq:relaxed_convex2}
\phi(y)\geq  \phi(x)+g(x)^T(y-x)-\epsilon_gD-\kappa_g\|\nabla \phi(x)\|D. 
\end{equation}

The next lemma is the key result in this setting and bounds the optimality gap.

%
%
%
%

\begin{lemma}\label{lem:convex_telescopic}
Suppose Assumptions~\ref{assum:convex} and  \ref{assum:lip_cont} hold. Let the step size parameter in Algorithm~\ref{alg:first_order_fixed_step} satisfy $\alpha \le 1/L_1$. And suppose that $\|x_i-x^*\|\leq D$ for $i=1,\ldots,k$, then
\begin{equation}\label{eq:convex_telescopic_bnd}
 \frac{1}{k}\sum_{i=0}^{k-1}\phi(x_{i+1})-\phi^*\leq \frac{\|x_0-x^*\|^2}{2\alpha k}+\epsilon_g{D}+\frac{\kappa_gD}{k}\sum_{i=0}^{k-1} \|\nabla \phi(x_i)\|.
\end{equation}
%
\end{lemma}
\begin{proof}
 By convexity and Lipschitz
continuity of the gradient, we have
\begin{align*}
\phi(x_{i+1})-\phi^*
&\leq
\nabla\phi(x_i)^T(x_{i+1}-x^*)
+\frac{L_1}{2}\|x_{i+1}-x_i\|^2\\
&\leq
g(x_i)^T(x_{i+1}-x^*)
-(g(x_i)-\nabla\phi(x_i))^T(x_{i+1}-x^*)
+\frac{1}{2\alpha}\|x_{i+1}-x_i\|^2,
\end{align*}
where the second inequality follows from
$\alpha\leq 1/L_1$.

Using $x_{i+1}=x_i-\alpha g(x_i)$, we obtain
\[
2\alpha g(x_i)^T(x_{i+1}-x^*)
+\|x_{i+1}-x_i\|^2
=
\|x_i-x^*\|^2-\|x_{i+1}-x^*\|^2.
\]
Consequently,
\begin{align*}
2\alpha[\phi(x_{i+1})-\phi^*] &\leq
\|x_i-x^*\|^2-\|x_{i+1}-x^*\|^2
+2\alpha (g(x_i)-\nabla\phi(x_i))^T(x^*-x_{i+1})\\
&\leq
\|x_i-x^*\|^2-\|x_{i+1}-x^*\|^2
+2\alpha D\left(
\epsilon_g+\kappa_g\|\nabla\phi(x_i)\|
\right),
\end{align*}
where the last inequality follows from
\ref{def:oracle_det_1} and $\|x_{i+1}-x^*\|\leq D$.

Summing over $i=0,\ldots,k-1$, dropping the
nonpositive term $-\|x_k-x^*\|^2$, and dividing
by $2\alpha k$ completes the proof.
\end{proof}

To interpret this inequality  from a complexity perspective we need to understand the behavior of the term 
$\frac{\kappa_gD}{k}\sum_{i=0}^{k-1} \|\nabla \phi(x_i)\|$.  Through the analysis in Theorem \ref{thm:nonconvex_bound}
we know that $\sum_{i=0}^{k-1} \|\nabla \phi(x_i)\|^2$ is bounded by $C_1+kC_2\epsilon_g$ for some constants $C_1$ and $C_2$. Thus one can expect that while $\|\nabla \phi(x_k)\|$ is large compared to $\epsilon_g$ is decreases  at the rate $\frac{1}{\sqrt{k}}$ and thus so does the term $\frac{\kappa_gD}{k}\sum_{i=0}^{k-1} \|\nabla \phi(x_i)\|$. 
 Consequently, the rate of decrease of  $ \frac{1}{k}\sum_{i=0}^{k-1}\phi(x_{i+1})-\phi^*$ is $\frac{1}{\sqrt{k}}$ which is worse 
 than what we obtained via Theorem \ref{thm:KL_cond_rate}. Once the gradient becomes small, the resulting extra error is proportional to $\kappa_gD\epsilon_g$. 
 Hence the telescopic sum analysis 
  does not seem to be as tolerant to the relative  error term in the oracle model. 
  One could consider  bounding each $ \|\nabla \phi(x_i)\|$ by some constant $G$, giving $\|g(x)-\nabla\phi(x)\|\leq \epsilon_g+\kappa_gG$, which effectively merges the additive and multiplicative errors into a single (potentially large) error term. For many of the examples in Section~\ref{sec.examples_first}, this would lead to a substantial error and, consequently, a very weak complexity bound. Therefore, it appears that the  “telescopic-sum” analysis tends to provide weaker guarantees when $\kappa_g>0$.

Let us focus on the case where   $\kappa_g=0$ and discuss the results further.  
The bound  \eqref{eq:convex_telescopic_bnd} is  on the average optimality gap over $k$ iterations.  
 When the algorithm is guaranteed to provide a decreasing sequence $\{\phi(x_{i})\}$ then
the same bound is also valid for  the last iterate $\phi(x_{k})-\phi^*$. Generally, descent  is not guaranteed  due to the inexact gradient $g_k$, however, 
due to Assumption~\ref{assum:convex}, by the arguments presented in Theorem \ref{thm:convex_rate}, descent is ensured as long as $\phi(x_{i})-\phi^*\geq \frac{\epsilon_g}{\tau}{D}$ for $\tau\in(0,1)$ and $\kappa_g=0$. Thus, we can conclude from \eqref{eq:convex_telescopic_bnd}  that  the algorithm achieves $\phi(x_{k})-\phi^*\leq \epsilon_\tau$, 
for some $\epsilon_\tau \geq \frac{\epsilon_g}{\tau}{D}$ with $\tau\in(0,1)$, after at most 
\begin{equation*}
	K_{\epsilon_\tau}= \frac{1}{\epsilon_\tau}\frac{\|x_0-x^*\|^2}{2\alpha(1-\tau)}
\end{equation*}
iterations.

Next, we compare these results with the result of Theorem \ref{thm:convex_rate}, which can be restated as follows (for $\kappa_g=0$), $\phi (x_k) -\phi^* \leq \epsilon_\tau$, with 
$\epsilon_\tau \geq  \frac{\epsilon_g}{\tau}D$  is achieved after at most
\begin{align*}
	K_{\epsilon_\tau} = \frac{1}{\epsilon_\tau} \frac{{D}^2}{\alpha(1 - \tau^2)}.
\end{align*}

While Lemma \ref{lem:convex_telescopic} establishes the bound for an average function value achieved during $k$ iterations, we know from Theorem \ref{thm:convex_rate} that until certain accuracy it reached, the function value is decreased on each iteration. We can invoke this property and thus obtain the bound 
for optimality gap at the last iterate rather than the average. 
Thus we see that the best achievable optimality gap is the same under both types of analysis. But the resulting complexity bounds differ somewhat. In particular,  the complexity obtained by the telescopic sum analysis under the convexity assumption is smaller by a factor of roughly $(\|x_0-x^*\|/D)^2$. Moreover, in the case of an exact first order oracle, where $\epsilon_g=0$ (and $\kappa_g=0$) the bound $D$ on the radius containing the iterates generated by the algorithm is not required at all, i.e., $D$ can be infinite. 

We note that the quantity $D$ in Theorem \ref{thm:convex_rate} is an upper bound on $\|x_k-x^*\|$ for all iterations $k=0, 1, \ldots, K_{\epsilon_\tau}$. With just slightly more complicated derivations $D^2$ in the bound in Theorem \ref{thm:convex_rate} can be replaced by the average of $\|x_k-x^*\|^2$, which may be larger than $\|x_0-x^*\|$ but also may be smaller. 

In summary, the result of Lemma \ref{lem:convex_telescopic} has an advantage over the result of  Theorem \ref{thm:convex_rate} when $\kappa_g=0$ and 
$\|x_0-x^*\|$ is small compared to the overall bound of the distance of iterates to the optimum and in  particular in the case of exact gradient oracle. In other cases, the analysis of Theorem \ref{thm:convex_rate} based on the descent inequality gives better or comparable results.

\paragraph{Other Results}

We now present some results from the literature, for which the bound $D$ on the level set radius is not required in the analysis, even in the presence of errors in the first-order oracle. However, as we will see, the results are weaker in some sense (than the results of Theorem~\ref{thm:KL_cond_rate}) since the effect of the gradient error accumulates as the algorithm progresses. 
In particular in \cite{Schmidt-Bach-NIPS-2011} proximal gradient and accelerated proximal gradient methods are analyzed  under the condition that the first-order oracle is of the  form of Oracle \ref{def:oracle_det_1} with $\kappa_{g}=0$ and $\epsilon_g$ that can vary with iterations. In the case where there is no proximal term, the basic first-order algorithm in \cite{Schmidt-Bach-NIPS-2011} reduces to Algorithm \ref{alg:first_order_fixed_step} with $\alpha=\frac{1}{L_1}$ and with $g_k$ such that at every iteration $k \geq 0$, the gradient error satisfies $\|g_k-\nabla \phi(x_k)\|\leq \epsilon_{g_k}$. It is then shown that if $\phi(x)$ is convex (without the assumption of the boundedness of the level sets, i.e., Assumption~\ref{assum:convex} with $D=\infty$), then
\begin{equation}
	 \frac{1}{k}\sum_{i=0}^{k-1}\phi(x_{i+1})-\phi^*\leq \frac{L_1}{2k}\left (\|x_0-x^*\|+2\sum_{i=1}^k \frac{{\epsilon_g}_i}{L_1}\right ) ^2
\end{equation}
Note if ${\epsilon_g}_i = \epsilon_g$ the gradient error is constant and not diminishing then this result reduces to  
\begin{equation}
	 \frac{1}{k}\sum_{i=0}^{k-1}\phi(x_{i+1})-\phi^*\leq \frac{L_1}{2k} \|x_0-x^*\|^2 + 2\|x_0-x^*\|{\epsilon_g} + \frac{2k}{L_1}{\epsilon_g}  ^2
\end{equation}
which does not guarantee convergence to an $\epsilon$-approximate solution for any $\epsilon>0$ since the right hand side diverges. The results in  \cite{Schmidt-Bach-NIPS-2011}
imply convergence only if $\epsilon_{g_k}\to 0$ at a certain rate. For example, in the case of strongly convex functions convergence is established when the noise decreases at a linear rate \cite{Schmidt-Bach-NIPS-2011}.

The analysis in \cite{Schmidt-Bach-NIPS-2011} is also extended to an accelerated first-order method; see Algorithm~\ref{alg:first_order_fixed_step_acc}. 
\begin{algorithm}[H]
\caption{~\textbf{Inexact Accelerated Gradient Method with Fixed Step Size}}
\label{alg:first_order_fixed_step_acc}
\begin{algorithmic}
\State \textbf{Inputs:} $x_0 = y_0 \in \mathbb{R}^n$ (starting point), $\alpha > 0$ (step size), $\{\beta_k\}$, Oracle \ref{def:oracle_det_1}. 
 \State \textbf{1:} \textbf{for $k=0,1,2,\dots$ do}
\State \text{2:}~~~~Compute  $g_k = g(y_k)\approx \nabla \phi(y_k)$ via Oracle \ref{def:oracle_det_1}  
 \State \text{3:}~~~~Update iterate $x_{k+1}\gets y_k-\alpha g_k$.
  \State \text{4:}~~~~Update iterate $y_{k+1}\gets x_{k+1}+\beta_k(x_{k+1}-x_{k})$.
\State \textbf{5:} \textbf{end for}
 \end {algorithmic}
\end{algorithm}
The results state that if $\beta_k=\frac{k}{k+3}$, then under Assumption~\ref{assum:convex} with $D=\infty$, 
\begin{equation}
	\phi(x_k)-\phi^*\leq \frac{L_1}{(k+1)^2}\left (\|x_0-x^*\|+2\sum_{i=1}^k  \frac{i{\epsilon_g}_i}{L_1}\right ) ^2.
\end{equation}
The consequence of this is that convergence is achieved as long as $i{\epsilon_g}_i$ diminishes to zero at a certain rate. It is not clear if one can establish convergence to an
$\epsilon$-optimal solution for an accelerated gradient method with a constant step size and under fixed (non-diminishing) error in the first-order oracle.

Even under the assumption of $D$-bounded level sets, for an accelerated first order algorithm (similar to Algorithm \ref{alg:first_order_fixed_step_acc}) analyzed in 
  \cite{Devolder2013FirstorderMO}, the error in the convergence bound accumulates as follows 
\begin{equation}
	\phi(x_k)-\phi^*\leq \frac{4L_1}{(k+1)(k+2)} \|x_0-x^*\|^2+\frac{k+1}{3} 2\epsilon_gD,
\end{equation}
and thus establishing convergence is not straightforward in the setting with non-diminishing noise in the first-order oracle.

\subsection{Line Search Algorithm based on Oracles \ref{def:oracle_det_0} and \ref{def:oracle_det_1}}\label{sec:line_search}

We now consider first-order methods that employ a line search procedure to adaptively select the step size parameter $\alpha$. This approach offers two main advantages over the fixed step size strategy considered in Algorithm \ref{alg:first_order_fixed_step}. First, selecting $\alpha\leq \tfrac{1}{L_1}$ requires knowledge of the global Lipschitz constant $L_1$. Second, this choice can be overly conservative, as the bound $\alpha\leq \tfrac{1}{L_1}$ is sufficient to guarantee descent but is not necessary. In particular, larger step sizes may be appropriate when the local Lipschitz constant of $\nabla\phi$ near the iterates is significantly smaller than the global constant $L_1$.

A classical Armijo line search \cite[Section 3.1]{NW} selects the step size parameter $\alpha$ by checking whether a sufficient decrease condition is satisfied. Unlike the fixed step size method considered above, this approach requires access to function value information. We therefore consider a line search method based on the inexact zeroth- and first-order oracles (Oracles \ref{def:oracle_det_0} and \ref{def:oracle_det_1}). At each iteration, the algorithm uses $-g(x_k)$ as a search direction and selects a step size $\alpha_k$ by checking the following Armijo-like condition,
\begin{align}\label{eq:mod_LS}
    f(x_k-\alpha_k g(x_k)) \leq f(x_k) - c_1 \alpha_k \|g(x_k)\|^2.
\end{align}
We refer to this as an ``Armijo-like'' condition because the classical Armijo condition (see, e.g., \cite{armijo1966minimization}) involves the true objective function $\phi$ and the directional derivative $g(x_k)^T\nabla\phi(x_k)$ in place of $\|g(x_k)\|^2$. Since the exact function values and gradients are unavailable, condition \eqref{eq:mod_LS} is formulated entirely in terms of the available oracle outputs.

We consider a simple, yet practical, backtracking procedure in which each iteration begins with an initial trial step size $\alpha_{\max}$, which is reduced, if necessary, until the Armijo-like condition \eqref{eq:mod_LS} is satisfied. For simplicity of presentation and analysis, we assume that $\alpha_{\max}>\frac{2}{L_1}$. This assumption is not essential, and smaller values of $\alpha_{\max}$ can also be accommodated, albeit at the expense of more complicated expressions in the complexity bounds. Although $L_1$ is generally unknown, choosing $\alpha_{\max}$ significantly larger than $\frac{2}{L_1}$ incurs only a logarithmic overhead in the number of backtracking steps. Alternatively, $\alpha_{\max}$ could be adjusted dynamically throughout the algorithm, potentially reducing the number of trial evaluations. However, such strategies require a more involved analysis and may introduce additional logarithmic factors in the complexity bounds. We therefore restrict our attention to a fixed $\alpha_{\max}$ satisfying the above assumption. The complete procedure is presented in Algorithm~\ref{alg:first_order_linesearch}.

\begin{algorithm}[H]
\caption{~\textbf{Inexact gradient method with adaptive step size (line-search)}}
\label{alg:first_order_linesearch}
\begin{algorithmic}
\State \textbf{Inputs:} $x_0 \in \mathbb{R}^n$, $\alpha_{\max} > \frac{2}{L_1}$,  $c_1\in (0,1),\  \gamma\in (0,1)$, Oracle \ref{def:oracle_det_0}   and Oracle \ref{def:oracle_det_1} 
with $\kappa_g<1$. 
 \State \text{1:} \textbf{ for $k=0,1,2,\cdots$ do}
\State \text{2:}~~~~Compute  $g_k = g(x_k)\approx \nabla \phi(x_k)$ via Oracle \ref{def:oracle_det_1}
\State \text{3:}~~~~Set  $\alpha_k = \alpha_{\max}$  
 \State \text{4:}~~~~\textbf{while} $f(x_k-\alpha_k g_k) > f(x_k) - c_1 \alpha_k \| g_k\|_2^2$
 \State \text{5:}~~~~~~~~Update $\alpha_k = \gamma \alpha_k$
 \State \text{6:}~~~~\textbf{end while.}
 \State \text{7:}~~~~Update iterate $x_{k+1}=x_k-\alpha_k g(x_k)$. 
\State \text{8:} \textbf{ end for.}
 \end {algorithmic}
\end{algorithm}

 We now turn to the analysis of Algorithm~\ref{alg:first_order_linesearch}. We begin by establishing an analogue of Lemma~\ref{lem.key_lemma_det}, which bounds the change in the inexact (observable) objective function values induced by a trial step. This result will play a key role in analyzing the Armijo-like sufficient decrease condition and establishing the convergence guarantees of the algorithm.
 
 \begin{lemma}\label{lem.key_lemma_det_LS_3}
Suppose Assumption~\ref{assum:lip_cont} holds. Then, for all $k\geq 0$,
\begin{align*}
	f (x_k - \alpha_k g_k) &\leq f (x_k) - \frac{\alpha_k}{2}\left(1-\kappa_g^2\right) \|\nabla \phi(x_k)\|^2   - \frac{\alpha_k}{2}\left(1 - \alpha_k L_1\right) \| g_k\|_2^2 + \frac{\alpha_k}{2}  \epsilon^2_g + 2 \epsilon_f, 
\end{align*}
Moreover, suppose that $\kappa_g<1$ and, for some $\tau\in(0,1)$, $\|\nabla \phi(x_k)\|^2\geq \frac{1}{(1-\kappa_g^2)\tau^2}\left( \epsilon_g^2 + \frac{4\epsilon_f}{\alpha_k}\right) $
Then,
\begin{align}\label{eq.result3}
	f (x_k - \alpha_k g_k) &\leq f (x_k) - \frac{\alpha_k}{2}\left(1-\kappa_g^2\right) \left(1-\tau^2\right) \|\nabla \phi(x_k)\|^2 - \frac{\alpha_k}{2}\left(1 - \alpha_k L_1\right) \| g_k\|_2^2.
\end{align}
\end{lemma}

\begin{proof}
By Assumption~\ref{assum:lip_cont} and following the arguments in the proof of Lemma~\ref{lem.key_lemma_det}, we have
\begin{align*}
    \phi(x_k-\alpha_k g_k)
    &\leq \phi(x_k)
    -\alpha_k g_k^T\nabla\phi(x_k)
    +\frac{\alpha_k^2 L_1}{2}\|g_k\|^2\\
    &\leq \phi(x_k)
    -\frac{\alpha_k}{2}(1-\kappa_g^2)
    \|\nabla\phi(x_k)\|^2-\frac{\alpha_k}{2}(1-\alpha_k L_1)\|g_k\|^2
    +\frac{\alpha_k}{2}\epsilon_g^2,
\end{align*}
where the last inequality follows from Oracle~\ref{def:oracle_det_1}. The first result then follows from Oracle~\ref{def:oracle_det_0}, since $|f(x_k - \alpha_k g_k) -\phi(x_k - \alpha_k g_k)|\leq \epsilon_f$ and 
$|f(x_k ) -\phi(x_k )|\leq \epsilon_f$.

For the second result, the assumed lower bound on
$\|\nabla\phi(x_k)\|^2$ implies that
\begin{align*}
    \frac{\alpha_k}{2}\epsilon_g^2+2\epsilon_f
    \leq
    \frac{\alpha_k}{2}(1-\kappa_g^2)\tau^2
    \|\nabla\phi(x_k)\|^2.
\end{align*}
Substituting this bound into the first inequality yields
\eqref{eq.result3}.
\end{proof}

We now use Lemma~\ref{lem.key_lemma_det_LS_3} to establish conditions under which the Armijo-like condition \eqref{eq:mod_LS} is satisfied. In particular, we show that when the true gradient norm is sufficiently large relative to the oracle errors, a sufficiently small step size $\alpha_k$ satisfies the Armijo-like condition, and hence the trial step is accepted, yielding $x_{k+1}=x_k-\alpha_k g_k$. We also establish a bound on the corresponding change in the true objective function.

\begin{lemma}\label{lem.key_lemma_det_LS}
Suppose Assumption~\ref{assum:lip_cont} holds. Let $\tau\in(0,1)$, define $\kappa_{g\tau}^2 = \tau^2+\left(1-\tau^2\right)\kappa_g^2$ where $\kappa_g<1$, assume that $0<c_1<\frac{1}{1+\kappa_{g\tau}}$, and define $\bar{\alpha} = \frac{2 \left( 1 - c_1 (1+ \kappa_{g\tau})\right)}{L_1(1+ \kappa_{g\tau})} $. For any iterate $x_k$ for which $\|\nabla \phi(x_k)\|^2\geq \frac{1}{(1-\kappa_g^2)\tau^2}\left( \epsilon_g^2+\frac{4\epsilon_f}{\gamma \bar{\alpha} } \right)$, the Armijo-like condition \eqref{eq:mod_LS} is satisfied for every trial step size
$\alpha\in[\gamma\bar{\alpha},\bar{\alpha}]$.
Consequently, the backtracking procedure accepts a step size
$\alpha_k>\gamma\bar{\alpha}$, and
\begin{align}\label{eq:key_descent_epsillonf}
\phi(x_{k+1})
\leq \phi(x_k)
-c_1\gamma\bar{\alpha}
(1-\kappa_{g\tau})^2
\|\nabla\phi(x_k)\|^2
+2\epsilon_f.
\end{align}
Moreover, if $\|\nabla \phi(x_k)\|^2\geq \max\left\{ \frac{1}{(1-\kappa_g^2)\tau^2}\left( \epsilon_g^2+\frac{4\epsilon_f}{\gamma \bar{\alpha} } \right), \frac{1}{(1-\kappa_{g\tau})^2\tau^2}\left( \frac{2\epsilon_f}{c_1 \gamma \bar{\alpha} } \right) \right\}$, then
\begin{align}\label{eq:key_descent}
\phi(x_{k+1})
\leq \phi(x_k)
-c_1\gamma\bar{\alpha}
(1-\kappa_{g\tau})^2
(1-\tau^2)
\|\nabla\phi(x_k)\|^2.
\end{align}
\end{lemma}

\begin{proof}
By the gradient lower bound it follows that $\|\nabla \phi(x_k)\|^2\geq \frac{\epsilon_g^2}{(1-\kappa_g^2)\tau^2}$, and by Oracle~\ref{def:oracle_det_1},
\begin{align}\label{eq.ls_case2}
	\|g_k-\nabla\phi(x_k)\|^2 &\leq
\kappa_g^2\|\nabla\phi(x_k)\|^2
+\epsilon_g^2 \nonumber\\
&\leq
\left(\tau^2+(1-\tau^2)\kappa_g^2\right)
\|\nabla\phi(x_k)\|^2 =
\kappa_{g\tau}^2
\|\nabla\phi(x_k)\|^2.
\end{align}
Since $\kappa_g<1$ and $\tau\in(0,1)$, we have $\kappa_{g\tau}<1$, and hence
\begin{align}\label{eq.ls_case2b}
	(1-\kappa_{g\tau})\|\nabla\phi(x_k)\|	\leq \|g_k\| \leq (1+\kappa_{g\tau})\|\nabla\phi(x_k)\|.
\end{align}

Now consider any trial step size $\alpha\in[\gamma\bar{\alpha},\bar{\alpha}]$. Since $\alpha\geq\gamma\bar{\alpha}$, the gradient lower bound implies
\begin{align*}
	\|\nabla\phi(x_k)\|^2 \geq \frac{1}{(1-\kappa_g^2)\tau^2} \left( \epsilon_g^2+\frac{4\epsilon_f}{\alpha} \right).
\end{align*}
Thus, by Lemma~\ref{lem.key_lemma_det_LS_3}, the bound \eqref{eq.result3} applies to this trial step. Using \eqref{eq.ls_case2b} and the identity $(1-\kappa_g^2)(1-\tau^2) =1-\kappa_{g\tau}^2$, we obtain
\begin{align*}
	f(x_k-\alpha g_k) &\leq f(x_k) -\frac{\alpha}{2}(1-\kappa_{g\tau}^2) \|\nabla\phi(x_k)\|^2-\frac{\alpha}{2}(1-\alpha L_1)\|g_k\|^2\\
	&\leq f(x_k) -\frac{\alpha}{2} \frac{1-\kappa_{g\tau}^2} {(1+\kappa_{g\tau})^2} \|g_k\|^2-\frac{\alpha}{2}(1-\alpha L_1)\|g_k\|^2\\
	&= f(x_k) -\alpha \left( \frac{1}{1+\kappa_{g\tau}} -\frac{\alpha L_1}{2} \right)\|g_k\|^2\\
	&\leq f(x_k)-c_1\alpha\|g_k\|^2,
\end{align*}
where the last inequality follows from $\alpha\leq\bar{\alpha}$ and the definition of $\bar{\alpha}$. Hence, the Armijo-like condition \eqref{eq:mod_LS} is satisfied for every
$\alpha\in[\gamma\bar{\alpha},\bar{\alpha}]$.

Next, since $0<\bar{\alpha}<\frac{2}{L_1}
<\alpha_{\max}$, 
the backtracking procedure starts with a trial step size larger than $\bar{\alpha}$. If a step size larger than $\bar{\alpha}$ is accepted, then clearly $\alpha_k>\gamma\bar{\alpha}$. Otherwise, the first trial step satisfying $\alpha\leq\bar{\alpha}$ necessarily belongs to $(\gamma\bar{\alpha},\bar{\alpha}]$, since consecutive trial step sizes differ by the factor $\gamma\in(0,1)$. As shown above, this step size satisfies the Armijo-like condition and is therefore accepted. Consequently, in either case, $\alpha_k>\gamma\bar{\alpha}$.

By Oracle~\ref{def:oracle_det_0},
the Armijo-like condition \eqref{eq:mod_LS},
and \eqref{eq.ls_case2b}, we obtain
\begin{align}\label{eq.descent_inexact1}
\phi(x_{k+1})
&\leq f(x_{k+1})+\epsilon_f \nonumber\\
&\leq f(x_k)-c_1\alpha_k\|g_k\|^2
+\epsilon_f \nonumber\\
&\leq \phi(x_k)-c_1\alpha_k\|g_k\|^2
+2\epsilon_f \nonumber\\
&\leq \phi(x_k)
-c_1\gamma\bar{\alpha}
(1-\kappa_{g\tau})^2
\|\nabla\phi(x_k)\|^2
+2\epsilon_f,
\end{align}
which establishes \eqref{eq:key_descent_epsillonf}.

Finally, under the stronger gradient lower bound, we also have
\[
2\epsilon_f
\leq
c_1\gamma\bar{\alpha}
(1-\kappa_{g\tau})^2\tau^2
\|\nabla\phi(x_k)\|^2.
\]
Substituting this bound into
\eqref{eq.descent_inexact1} yields
\begin{align*}
	\phi(x_{k+1}) \leq \phi(x_k) -c_1\gamma\bar{\alpha} (1-\kappa_{g\tau})^2 (1-\tau^2) \|\nabla\phi(x_k)\|^2,
\end{align*}
which establishes \eqref{eq:key_descent} and completes the proof.
\end{proof}

Combining the descent guarantees of Lemma~\ref{lem.key_lemma_det_LS} with the general results of Theorems~\ref{thm:nonconvex_bound}, \ref{thm:KL_cond_rate}, \ref{thm:convex_rate}, and \ref{thm:str_convex_rate}, we obtain the following complexity bounds for Algorithm~\ref{alg:first_order_linesearch} under different assumptions on the objective function $\phi$.


\begin{corollary}\label{cor:all_complexity_ls}

Let $\kappa_{g\tau}^2 = \tau^2+\left(1-\tau^2\right)\kappa_g^2$ and $ \bar{\alpha} = \frac{2 \left( 1 - c_1 (1+ \kappa_{g\tau})\right)}{L_1(1+ \kappa_{g\tau})} $, and assume $c_1 < \frac{1}{1 + \kappa_{g\tau}}$. The following complexity bounds follow.
\begin{itemize}
	\item Suppose $\phi$ is bounded below and $L_1$-smooth (Assumptions \ref{assum:low_bound} and \ref{assum:lip_cont}). Then, Algorithm~\ref{alg:first_order_linesearch}, generates an iterate $x_k$ satisfying $\| \nabla \phi (x_k) \| \leq \epsilon_\tau$, 
	
	with $\epsilon_\tau \geq \sqrt{\max\left\{ \frac{1}{(1-\kappa_g^2)\tau^2}\left( \epsilon_g^2+\frac{4\epsilon_f}{\gamma \bar{\alpha} } \right), \frac{1}{(1-\kappa_{g\tau})^2\tau^2}\left( \frac{2\epsilon_f}{c_1 \gamma \bar{\alpha} } \right) \right\}}$ and $\tau \in (0,1)$, within at most 
\begin{align*}
	K_{\epsilon_\tau} = \left\lceil \frac{1}{\epsilon_\tau^2} \frac{ (\phi(x_0) - \hat{\phi})}{c_1\gamma \bar{\alpha}(1-\kappa_{g\tau})^2(1 - \tau^2)} \right\rceil
\end{align*}
iterations.
	\item Suppose $\phi$ satisfies  Assumption~\ref{assum:KLcond} with $c=\frac{1}{D}$, e.g., if $\phi$ is convex, with bounded level sets as in  Assumption~\ref{assum:convex}. Then, Algorithm~\ref{alg:first_order_linesearch} generates an iterate satisfying $\phi (x_{k}) -\phi^* \leq \epsilon_\tau$, with 
$\epsilon_\tau \geq \sqrt{\max\left\{ \frac{1}{(1-\kappa_g^2)\tau^2}\left( \epsilon_g^2+\frac{4\epsilon_f}{\gamma \bar{\alpha} } \right), \frac{1}{(1-\kappa_{g\tau})^2\tau^2}\left( \frac{2\epsilon_f}{c_1 \gamma \bar{\alpha} } \right) \right\}} {D}$ and $\tau \in (0,1)$, within at most 
\begin{align*}
	K_{\epsilon_\tau} = \left\lceil \frac{1}{\epsilon_\tau} \frac{D^2}{c_1\gamma\bar{\alpha}(1-\kappa_{g\tau})^2(1 - \tau^2)} \right\rceil
\end{align*}
iterations. 
	\item Suppose $\phi$ satisfies Assumption~\ref{assum:PLcond}, e.g., if $\phi$ is  $\mu$-strongly  convex as in  Assumption~\ref{assum:str.convex}. Then, Algorithm~\ref{alg:first_order_linesearch}, generates an iterate satisfying $\phi (x_k) -\phi^* \leq \epsilon_\tau$, with $\epsilon_\tau \geq \max\left\{ \frac{1}{(1-\kappa_g^2)\tau^2 \mu}\left( \epsilon_g^2+\frac{4\epsilon_f}{\gamma \bar{\alpha} } \right), \frac{1}{(1-\kappa_{g\tau})^2\tau^2\mu}\left( \frac{2\epsilon_f}{c_1 \gamma \bar{\alpha} } \right) \right\}$ and $\tau \in (0,1)$, within at most
\begin{align*}
	K_{\epsilon_\tau} = \left\lceil \log_{\frac{1}{1 - c_1\gamma \bar{ \alpha} \mu (1-\kappa_{g\tau})^2(1-\tau^2)}} \left( \frac{\phi (x_0) - \phi^*}{\epsilon_\tau}\right) \right\rceil
\end{align*} 
iterations.
\end{itemize}
\end{corollary}

Corollary~\ref{cor:all_complexity_ls} is the line-search counterpart of Theorem~\ref{cor:all_complexity_fixed} for first-order methods with deterministic inexact oracles. The two results exhibit several similarities, as well as important differences. In both cases, the algorithms achieve convergence rates of the same order as their exact-oracle counterparts, up to accuracy thresholds determined by the oracle errors. However, the attainable accuracy differs between the two methods. In particular, the line-search procedure requires access to a zeroth-order oracle, and errors in function evaluations can lead to a larger neighborhood of stationarity. This can be viewed as a cost of adaptively selecting the step size. When exact function evaluations are available, this additional source of error disappears, and the two methods achieve accuracy neighborhoods of the same order, determined solely by the first-order oracle parameters.

Another important distinction concerns the computational cost per iteration. While the fixed step size method requires only one call to the first-order oracle, each iteration of Algorithm~\ref{alg:first_order_linesearch} requires one call to Oracle~\ref{def:oracle_det_1} and potentially multiple calls to Oracle~\ref{def:oracle_det_0} to determine an acceptable step size. In particular, under the gradient-norm condition of Lemma~\ref{lem.key_lemma_det_LS}, the number of zeroth-order oracle evaluations at trial points is bounded above by $1+\left\lceil \log_{\frac{1}{\gamma}}\left(\frac{\alpha_{\max}}{\bar{\alpha}}\right)\right\rceil$. Thus, the additional cost associated with backtracking grows only logarithmically with the ratio $\alpha_{\max}/\bar{\alpha}$.

%
%

Several works \cite{berahas2018derivative,berahas2021global} relax the Armijo condition to account for errors in the objective function evaluations by adding a constant ($r\geq 0$) to its right-hand side, i.e., 
\begin{align}\label{eq:mod_LS_relax}
	f(x_k-\alpha_k g(x_k)) \leq f(x_k) - c_1 \alpha_k \| g(x_k)\|^2 + r.
\end{align} 
This yields analogous convergence and complexity guarantees, with the main difference appearing in the convergence neighborhood. In particular, when ($r=2\epsilon_f$), the second term in the lower bound on ($\epsilon_\tau$) vanishes; for general $r\geq 0$, it is replaced by
\begin{align*}
	\epsilon_\tau \geq \sqrt{\max\left\{ \frac{1}{(1-\kappa_g^2)\tau^2}\left( \epsilon_g^2+\frac{\max\{4\epsilon_f-2r,0\}}{\gamma \bar{\alpha} } \right), \frac{1}{(1-\kappa_{g\tau})^2\tau^2}\left( \frac{2\epsilon_f+r}{c_1 \gamma \bar{\alpha} } \right) \right\}}
\end{align*} 
see \cite{berahas2026combining} for more details. 

Our main conclusion here is that the relaxing  the Armijo condition does not in fact provide noticeable theoretical advantages. This is in contrast with the use of inexact oracles 
in various randomized line-search schemes where relaxing Armijo condition is essential \cite{berahas2021global}. 

A similar observation will hold for the deterministic trust region methods which we study in the next section.

\subsection{Analysis of a trust region method based on inexact oracles}\label{sec:tr}
Trust-region (TR) methods constitute a well-established alternative to line-search methods, as they also regulate the step taken between iterates based on the reduction achieved in the objective function. Rather than selecting a step length along a prescribed search direction, TR methods compute a step by approximately minimizing a local model of the objective within a neighborhood of the current iterate.

The central ingredient of a TR method is a model $m(\cdot)$ of the objective function $\phi(\cdot)$ whose accuracy is trusted within a region around a given point $x \in \mathbb{R}^n$. This region, referred to as the trust region, restricts the step $s$ to lie within a prescribed radius. The model is typically chosen to be a quadratic polynomial of the form
\begin{equation} \label{eq:tr_quad_model}
    m(x+s) = \phi(x) + g(x)^T s + \frac{1}{2} s^T H(x) s,
\end{equation}
where $g(x)$ and $H(x)$ denote either the true gradient and Hessian of $\phi$ at $x$, or approximations thereof, depending on the oracle employed. Although the expression includes the exact function value $\phi(x)$, it is not necessary to compute it when forming the model; rather, it appears for notational convenience and is primarily used in the analysis.

Trust-region methods were originally developed in the context of second-order algorithms that rely on second-order oracles. However, first-order variants play an important role, particularly when establishing convergence guarantees to (approximate) first-order stationary points. For this reason, we begin by discussing first-order TR methods and defer the second-order case until after introducing second-order oracles and their properties.

By a first-order trust-region method, we mean a variant in which convergence guarantees are established solely with respect to first-order stationarity, without requiring $H(x)$ to approximate the true Hessian $\nabla^2 \phi(x)$ with any prescribed accuracy. In particular, the matrix $H(x)$ may be an arbitrary symmetric matrix, including the choice $H(x)=0$, although such a choice is rarely competitive in practice. The role of $H(x)$ in this setting is primarily to shape the model and influence practical performance, rather than to ensure second-order accuracy. We defer explicit assumptions on $H(x)$ to below, where we specify the mild conditions under which first-order convergence can be established.

We consider a first-order trust region method where the model~\eqref{eq:tr_quad_model} is constructed using the output of the deterministic first-order oracle defined as Oracle~\ref{def:oracle_det_1}. In particular, our analysis relies on a weaker gradient approximation bound satisfying
\begin{equation}\label{eq:epskappa_g}
	\|g(x) - \nabla \phi(x)\| \leq \epsilon_g + \kappa_g \|\nabla \phi(x)\|,
\end{equation}
which is a consequence of Oracle~\ref{def:oracle_det_1}. The deterministic zeroth-order Oracle~\ref{def:oracle_det_0} is employed by the trust-region algorithm to assess the agreement between the model and the objective function, and thus to determine whether a trial step should be accepted and how the trust-region radius should be updated. It may also be utilized in the construction of  Oracle~\ref{def:oracle_det_1} but we do not explicitly consider it, instead treating  Oracle~\ref{def:oracle_det_1} as defined by its general properties, regardless of  its implementation.

We now formally present the trust-region method (Algorithm~\ref{alg:TR1}) together with the oracle assumptions under which  iteration complexity can be established for computing approximate first-order stationary points. In the subsequent subsection, we relax these assumptions and investigate how weaker oracle conditions influence the resulting complexity guarantees.


\begin{algorithm}[H] 
    \caption{~\textbf{Trust region method based on general  Oracles \ref{def:oracle_det_1} and \ref{def:oracle_det_0}}}
    \label{alg:TR1}
    \begin{algorithmic}
    \State \textbf{Inputs:} $x_0 \in \mathbb{R}^n$, TR radius $\delta_0 > 0$, $\eta_1\in(0,1)$, $\eta_2 > 0$, and $\gamma\in(0,1)$, Oracle \ref{def:oracle_det_0}   and Oracle \ref{def:oracle_det_1} 
     \State \text{1:} \textbf{ for $k=0,1,2,\cdots$ do}
\State \text{2:}~~~~Compute $g_k=g(x_k)$ via Oracle  \ref{def:oracle_det_1}, and set $H_k=H(x_k)$. Define model $m_k(x_k+s)=\phi(x_k)+g_k^Ts+\frac{1}{2}s^TH_ks$. 
 \State \text{3:}~~~~\textbf{if} $\|g_k\| \ge \eta_2 \delta_k$
  \State \text{4:}~~~~~~~~Compute a trial step $x_k^+ \approx \arg\min_y \{m_k(y):~y \in B(x_k,\delta_k)\}$.
  \State \text{5:}~~~~~~~~Compute $f(x_k)$ and $f(x_k^+)$ via  Oracle \ref{def:oracle_det_0} and calculate the ratio $\rho_k$ as 
        \begin{align*} 
        		\rho_k = \frac{f(x_k) - f(x_k^+)}{m_k(x_k) - m_k(x_k^+)}.
	\end{align*}
	 \State \text{6:}~~~~\textbf{end if.}
\State \text{7:}~~~~Update the iterate and the TR radius as 
        \begin{align*}  
        	(x_{k+1}, \delta_{k+1}) \gets \left\{ \begin{aligned} 
            &(x_k^+, \gamma^{-1} \delta_k) &&\text{if } \|g_k\| \ge \eta_2 \delta_k \text{ and } \rho_k \ge \eta_1  \\
            &(x_k, \gamma \delta_k) &&\text{otherwise. }
        \end{aligned} \right. 
	\end{align*}

 \State \text{8:} \textbf{ end for.}
    \end{algorithmic}
\end{algorithm}

\begin{assumption}[\textbf{Trust Region Method Assumptions}] \label{assum:TR1}
We make the following standard \cite{Conn-Gould-Toint-00} assumptions:
    \begin{enumerate}
           \item The trust-region subproblem is solved sufficiently accurately in each iteration $k$ such that $x_k^+$ provides at least a fraction of Cauchy decrease, i.e., for some constant $0 < {c_1} < \frac{1}{2}$ 
          \begin{equation}\label{eq:fcd_cond}
            m_k(x_k) - m_k(x_k^+) \ge {c_1} \|g_k\| \min\bigg\{\frac{\|g_k\|}{\|H_k\|}, \delta_k\bigg\}. 
          \end{equation}
              \item There exists a constant $\kappa_{bhm} > 0$ such that, for all $x_k$ generated by Algorithm~\ref{alg:TR1},  
         $ \|H_k\| \le \kappa_{bhm}$. 
    \end{enumerate}
\end{assumption}

The above  assumptions are standard for trust region methods and are easily enforced \cite{Conn-Gould-Toint-00}. 

\begin{lemma}\label{lem:TRD0D1_model_error}
    Suppose Assumptions~\ref{assum:lip_cont} and \ref{assum:TR1} hold. For any $s\in B(x_k,\delta_k)$  
    \begin{equation}
        |m_k(x_k+s) - \phi(x_k+s)| \le \frac{L_1 + \kappa_{bhm}}{2} \delta_k^2 + \big(\epsilon_{g} + \kappa_{g}\|\nabla\phi(x_k)\|\big) \delta_k.
    \end{equation}
\end{lemma}
\begin{proof} By Assumption \ref{assum:lip_cont} and \ref{assum:TR1}, and \eqref{eq:epskappa_g} 
    \begin{align*}
        &|m_k(x_k+s) - \phi(x_k+s)| \\
        &= |\phi(x_k) + g(x_k)^Ts + s^T H_k s / 2 - \phi(x_k+s)| \\
        &\le |\phi(x_k) + \nabla\phi(x_k)^Ts - \phi(x_k+s)| + |(g(x_k) - \nabla\phi(x_k))^T s| + |s^T H_k s / 2| \\
        &\le \frac{L_1}{2} \|s\|^2 + \big(\epsilon_{g} + \kappa_{g}\|\nabla\phi(x_k)\|\big)\|s\| + \kappa_{bhm} \|s\|^2 / 2 \\ 
        &\le \frac{L_1 + \kappa_{bhm}}{2} \delta_k^2 + \big(\epsilon_{g} + \kappa_{g}\|\nabla\phi(x_k)\|\big) \delta_k.
    \end{align*}
\end{proof}

Note that, unlike the gradient-descent-type algorithms analyzed in Sections~\ref{sec:fixed_step_det} and~\ref{sec:line_search}, the trust-region method may not produce an iterate update at every iteration. Specifically, we call iteration $k$ \emph{successful} if $x_{k+1} = x_k^+$ and the trust-region radius is increased (i.e., $\delta_{k+1} = \gamma^{-1} \delta_k$), and \emph{unsuccessful} otherwise. Even after an unsuccessful iteration, we assume that a new $g_k = g(x_k)$ (and possibly $H_k = H(x_k)$) is computed, despite the fact that the iterate $x_k$ itself remains unchanged. This is motivated by the reduction of the trust-region radius $\delta_k$ on unsuccessful iterations. A smaller trust-region radius may warrant a more accurate gradient approximation, if permitted by the oracle, that is, a smaller value of $\epsilon_g$ in~\eqref{eq:epskappa_g}. Our analysis will show that the results hold provided $\delta_k \ge \frac{\epsilon_g}{\kappa_{\delta g}}$, for some constant $\kappa_{\delta g} > 0$. Consequently, as $\delta_k$ is reduced, it may be desirable to decrease $\epsilon_g$ accordingly, when the oracle allows adaptive accuracy. We return to this relationship after presenting the complexity analysis. For now, we assume that $\epsilon_g$ is fixed.

The following two lemmas  together are, in some sense, the counterparts of the descent inequality for the trust region methods. They state that when the trust region radius (which controls the step sizes) is sufficiently small, the iteration is successful and ensures sufficient progress. The following lemma provides a sufficient condition for an iteration to be successful. 
\begin{lemma}\label{lem:TRD0D1_success}
    Under Assumptions~\ref{assum:lip_cont} and \ref{assum:TR1}, if  the following conditions
    \begin{subequations} \label{eq:TR1_success} \begin{align}
        & \frac{1-\kappa_g}{\kappa_g}> \frac{1}{{c_2}}  \label{eq:TR1_success1}\\
        &\| \nabla\phi(x_k)\|\geq \frac{1+{c_2}}{{c_2}(1-\kappa_g)-\kappa_g}\epsilon_g \label{eq:TR1_success2}\\
        & \sqrt{\frac{4\epsilon_f}{L_1+\kappa_{bhm}}}  \le  \delta_k \label{eq:TR1_success3}\\
        & \delta_k\leq C_1 \left[ (1-\kappa_g)\|\nabla \phi(x_k)\| - \epsilon_g \right] \label{eq:TR1_success4}
    \end{align} \end{subequations}
    hold with 
    \[
        C_1=  \min\left \{\frac{1}{\eta_2}, \frac{{c_2}}{L_1+\kappa_{bhm}}\right \}
        \text{ and } 
        {c_2} = \frac{(1-\eta_1){c_1}}{2}<1, 
    \]
    then $\rho_k \ge \eta_1$, $\|g_k\| \ge \eta_2 \delta_k$, and $x_{k+1} = x_k^+$. 
\end{lemma}

\begin{proof}
    First, notice that $C_1\leq \frac{1}{\max\{\eta_2, \kappa_{bhm}\}}$ holds by the definition of $C_1$, and \eqref{eq:TR1_success4} implies  
    \[  \delta_k \le 
        \frac{(1-\kappa_g) \|\nabla\phi(x_k)\| - \epsilon_{g}}{\max\{\kappa_{bhm},\eta_2\}}, 
    \]
    which further implies that 
    \begin{align*}
        \|g(x_k)\| - \max\{\kappa_{bhm},\eta_2\} \delta_k 
        \ge (1 - \kappa_{g}) \|\nabla\phi(x_k)\| - \epsilon_{g} - \max\{\kappa_{bhm},\eta_2\} \delta_k 
        \ge 0. 
    \end{align*}
    This means $\min \left\{\frac{\|g(x_k)\|}{\|H_k\|}, \delta_k\right\} = \delta_k$ and $\|g(x_k)\| \ge \eta_2 \delta_k$. 
    To show $\rho_k \ge \eta_1$, we start by applying Lemma~\ref{lem:TRD0D1_model_error} to obtain 
    \[ \begin{aligned}
        \rho_k &= \frac{m(x_k) - m(x_k^+) + (f(x_k) - m(x_k)) - (f(x_k^+) - m(x_k^+))}{m(x_k) - m(x_k^+)} \\
        &=\frac{m(x_k) - m(x_k^+) + (f(x_k) - \phi(x_k)) - (f(x_k^+) - \phi(x_k^+))+ (\phi(x_k^+) - m(x_k^+))}{m(x_k) - m(x_k^+)} \\
        &\ge 1 - \frac{2\epsilon_{f} + (L_1 + \kappa_{bhm}) \delta_k^2/2 + (\epsilon_{g} + \kappa_g\|\nabla\phi(x_k)\|) \delta_k}{m(x_k) - m(x_k^+)}, 
    \end{aligned} \] 
    where the denominator $m(x_k) - m(x_k^+) \ge {c_1} \|g_k\| \delta_k  \ge {c_1} [(1-\kappa_g) \|\nabla\phi(x_k)\| - \epsilon_{g}]\delta_k$. 
    Thus we have 
    \[  \rho_k \geq 1- \frac{2\epsilon_{f} + (L_1 + \kappa_{bhm}) \delta_k^2/2 }{{c_1} [(1-\kappa_g)\|\nabla\phi(x_k)\|  - \epsilon_{g}]\delta_k}-  \frac{\epsilon_{g} + \kappa_g\|\nabla\phi(x_k)\|}{{c_1} [(1-\kappa_g)\|\nabla\phi(x_k)\|  - \epsilon_{g}]}. 
    \]   

    To show that $\rho_k \geq \eta_1$, we derive two inequalities:
    \begin{align}
        \frac{2\epsilon_{f} + (L_1 + \kappa_{bhm}) \delta_k^2/2 }{{c_1} [(1-\kappa_g)\|\nabla\phi(x_k)\|  - \epsilon_{g}]\delta_k} &\leq \frac{1-\eta_1}{2} \label{eq:suffdec1} \\
        \text{ and } \frac{\epsilon_{g} + \kappa_g\|\nabla\phi(x_k)\|}{{c_1} [(1-\kappa_g)\|\nabla\phi(x_k)\|  - \epsilon_{g}]} &\leq \frac{1-\eta_1}{2}.  \label{eq:suffdec2} 
    \end{align} 

    The bound \eqref{eq:suffdec1} is equivalent to 
    \begin{equation*} 
        \frac{2\epsilon_{f}}{ \delta_k} + \frac{(L_1 + \kappa_{bhm}) \delta_k }{2 }\leq {c_2} \left[ (1-\kappa_g)\|\nabla\phi(x_k)\| -\epsilon_g \right]. 
    \end{equation*}
    And thus  by \eqref{eq:TR1_success3},  \eqref{eq:suffdec1} is implied by
    \begin{equation*} 
        \delta_k\leq \frac{{c_2}}{L_1 + \kappa_{bhm} } \left[ (1-\kappa_g)\|\nabla\phi(x_k)\| -\epsilon_g \right]
    \end{equation*}
    which follows from \eqref{eq:TR1_success4}. 

    Bound \eqref{eq:suffdec2} is equivalent to 
    \begin{equation*} 
        \epsilon_{g} + \kappa_g\|\nabla\phi(x_k)\|\leq {c_2} \left[ (1-\kappa_g)\|\nabla\phi(x_k)\| -\epsilon_g \right], 
    \end{equation*}
    which is equivalent to 
    \begin{equation*} 
        (1+{c_2})\epsilon_g \leq ({c_2} (1-\kappa_g)-\kappa_g )\|\nabla\phi(x_k)\|
    \end{equation*}
    which is implied by \eqref{eq:TR1_success1} and \eqref{eq:TR1_success2}.
\end{proof}

The next lemma provides a lower bound on the function decrease in a successful step. 
\begin{lemma}[guaranteed progress on successful iterations] \label{lem:TRD0D1_progress}
    Let Assumptions~\ref{assum:lip_cont} and \ref{assum:TR1} hold. 
    If iteration $k$ is successful, and 
    \begin{equation}\label{eq:TRProgess_lowbnd} \begin{aligned} 
        \delta_k \ge \sqrt{\frac{4\epsilon_f}{C_2}},  
        \text{ where } C_2 =\eta_1\eta_2{c_1} \min\{1, \frac{\eta_2}{\kappa_{bhm}}\}, 
    \end{aligned} \end{equation}
    the following inequality holds  
    \begin{equation} \label{eq:TR1_Progess} \begin{aligned} 
        &\phi(x_k) - \phi(x_{k+1}) \geq  C_2\delta_k^2-2\epsilon_f\geq \frac{C_2}{2}\delta_k^2; 
    \end{aligned} \end{equation}
    otherwise, we have $x_{k+1} = x_k$ and $\phi(x_k) - \phi(x_{k+1}) = 0$.  
\end{lemma}
\begin{proof}
    When the trial point is accepted, we have $\rho_k \ge \eta_1$ and $\|g_k\| \ge \eta_2 \delta_k$, which guarantees
    \begin{align*}
        \phi(x_k) - \phi(x_{k+1}) 
        \ge f(x_k) - f(x_k^+) - 2\epsilon_f 
        \ge \eta_1 \big(m(x_k) - m(x_k^+)\big) - 2\epsilon_f, 
    \end{align*}
    where 
    \begin{align*}
        m(x_k) - m(x_k^+)
        \ge {c_1} \|g_k\| \min\bigg\{\frac{\|g_k\|}{\|H_k\|}, \delta_k\bigg\} 
        \ge {c_1} \eta_2\delta_k \min\bigg\{\frac{\eta_2\delta_k}{\kappa_{bhm}}, \delta_k\bigg\}. 
    \end{align*}
\end{proof}

In previous sections, we used  $K_\epsilon$ to denote  the number of iterations of the corresponding algorithm that ensured that this algorithm would reach $\epsilon$-accuracy, defined differently depending on properties of $\phi(x)$. 
Bounding $K_\epsilon$ gives us complexity of the corresponding algorithm. Here we will focus solely on nonconvex case, since this is the most common setting when trust region methods are employed. 
Specifically, for Algorithm~\ref{alg:TR1}, let 
\[
    K_\epsilon =\min \{ k\in {\mathbb N},\, \|\nabla \phi(x_k)\|\leq \epsilon \}
\]
i.e., the first iteration when Algorithm~\ref{alg:TR1} reaches $\epsilon$-stationary point.  
We would like to bound $K_\epsilon$ as a function of $\epsilon$. To achieve this, we first bound the number of successful iterations that can occur while $\|\nabla \phi(x)\| > \epsilon$ and then bound the number of unsuccessful iterations that can occur in the meantime. 
We define index set 
\[ {\cal S}_\epsilon := \{k\in\{0,\dots,K_\epsilon-1\}:~ \text{iteration $k$ is successful}\}.
\]  
(i.e. the index set of successful iterations for which $\|\nabla \phi(x_k)\| > \epsilon$). 

Now let us discuss the lower bound on $\epsilon$. Due to the error terms $\epsilon_{g}$ in Oracles~\ref{def:oracle_det_1} and $\epsilon_{f}$ in Oracle~\ref{def:oracle_det_0}, the algorithm is not capable of finding a solution with arbitrarily small $\epsilon$, and it is only possible to prove that Algorithm~\ref{alg:TR1}  converges to a near $\epsilon$-stationary point with an $\epsilon$ sufficiently bounded away from zero. 
We establish a lower bound on $\epsilon$ using the two previous lemmas. 
Due to the mechanism of the TR radius update, a TR radius $\delta_k$ satisfying \eqref{eq:TR1_success3}, \eqref{eq:TR1_success4}, and \eqref{eq:TRProgess_lowbnd} will be found by Algorithm~\ref{alg:TR1} as long as $\delta_0 \ge \sqrt{\frac{4\epsilon_f}{\min\{L_1+\kappa_{bhm},C_2\}}}$ and 
\[
    \max\left \{ \sqrt{\frac{4\epsilon_f}{C_2}}, \sqrt{\frac{4\epsilon_f}{L_1+\kappa_{bhm}}}\right \} 
    \le \gamma C_1\left( (1-\kappa_g)\|\nabla \phi(x_k)\|-\epsilon_g \right). 
\]
Combined with \eqref{eq:TR1_success2}, we can conclude that when $\delta_0 \ge \sqrt{\frac{4\epsilon_f}{\min\{L_1+\kappa_{bhm},C_2\}}}$ and $\|\nabla \phi(x_k)\| > \epsilon$ for 
\begin{equation}\label{eq:epsilon_lowbndtr}
    \epsilon \geq  \max\left \{
        \frac{1}{1-\kappa_g} \left (\sqrt{\frac{4\epsilon_f}{ \gamma^2 C_1^2\min\{ C_2,L_1+\kappa_{bhm}\} }} +\epsilon_g\right ), 
        \frac{1+{c_2}}{{c_2}(1-\kappa_g)-\kappa_g}\epsilon_g 
    \right \}, 
\end{equation}
the step is successful and leads to a decrease in function value satisfying \eqref{eq:TR1_Progess}.

 Now, let us consider any $\epsilon$ that satisfies \eqref{eq:epsilon_lowbndtr}.    Define $\tau=1-\kappa_g-\frac{\epsilon_g}{\epsilon}$. 
 It is simple to show that  $\tau\in(0,1]$ and 
 \begin{equation}\label{eq:tau_epsilon}
 \tau\epsilon \leq  (1-\kappa_g)\epsilon-\epsilon_g. 
 \end{equation}

The following result establishes upper bound on the total number of successful iterations. 
\begin{lemma}\label{lem:tr_bound_success}
    Let Assumptions~\ref{assum:lip_cont} and \ref{assum:TR1} hold. Let  $\kappa_g$ satisfy \eqref{eq:TR1_success1}. 
    For Algorithm~\ref{alg:TR1} and any given $\epsilon$ satisfying \eqref{eq:epsilon_lowbndtr}, suppose the initial trust-region radius $\delta_0 \ge \gamma C_1 \epsilon$ and $\tau\in (0,1]$ that satisfies \eqref{eq:tau_epsilon}. Then, 
    \begin{equation} \begin{aligned} 
        |{\cal S}_\epsilon | &\le \frac{2(\phi(x_0) - \phi^\star)}{C_2 (\gamma C_1 \tau \epsilon)^2}.
    \end{aligned} \end{equation}
\end{lemma}
\begin{proof} 
    According to Lemma~\ref{lem:TRD0D1_success} and the discussion above, any iteration $k \in \{0,\dots, K_\epsilon - 1\}$ must be successful when $\delta_k in [ \gamma C_1\tau \epsilon, C_1\tau \epsilon]$. 
    Thus, given $\delta_0 \ge \gamma C_1 \epsilon$ and $\epsilon$ satisfying \eqref{eq:epsilon_lowbndtr}, we must have $\delta_k \ge \gamma C_1 \tau \epsilon$ for all $k \in \{0,\dots, K_\epsilon - 1\}$. 
    From that and  Lemma~\ref{lem:TRD0D1_progress} 
    \[
     \phi(x_k) - \phi(x_{k+1}) \geq \frac{1}{2} C_2\delta_k^2\geq \frac{1}{2} C_2(\gamma C_1\tau \epsilon)^2
     \]
      if $k\in{\cal S}_\epsilon$ and $\phi(x_k) - \phi(x_{k+1}) =0$ otherwise. 
Using this we have 
    \begin{align*} 
        \phi (x_0) - \phi^\star &\ge \sum_{k=0}^{K_\epsilon-1} \phi(x_k) - \phi(x_{k+1})  >\frac{1}{2}  |{\cal S}_\epsilon | C_2 (\gamma C_1 \tau \epsilon)^2
            \end{align*}
            and the result follows. 
\end{proof}

    The total complexity bound is now derived as follows. 
    
\begin{theorem}\label{thm:tr_complexity_basic}
  Let Assumptions~\ref{assum:lip_cont} and \ref{assum:TR1} hold. Let  $\kappa_g$ satisfy \eqref{eq:TR1_success1}. 
    For Algorithm~\ref{alg:TR1} and any given $\epsilon$ satisfying \eqref{eq:epsilon_lowbndtr} and $\tau\in (0,1]$ that satisfies \eqref{eq:tau_epsilon}, 
   suppose the initial trust-region radius $\delta_0 \ge  \gamma C_1 \epsilon$.
      \begin{equation} \begin{aligned} 
        K_\epsilon&\le 2  |{\cal S}_\epsilon | + \log_{\gamma^{-1}} \frac{\delta_0}{C_1\tau\epsilon}  \\  
        &\le \max\left \{\eta_2, \frac{2(L_1+\kappa_{bhm})}{(1-\eta_1){c_1}}\right \}^2\frac{4(\phi(x_0) - \phi^\star)}{\eta_1\eta_2{c_1} \gamma^2\tau^2\min\{1, \eta_2/\kappa_{bhm}\}}\frac{1}{\epsilon^2} \\
        &\quad + \log_{\gamma^{-1}}
         \left( \max\left \{ \eta_2, \frac{2(L_1+\kappa_{bhm})}{(1-\eta_1){c_1}}\right \}
         \frac{\delta_0}{\tau\epsilon} \right)
    \end{aligned} \end{equation}
\end{theorem}
\begin{proof} 
    Assume that $|{\cal S}_\epsilon | \le \ceil*{\Big( K_\epsilon -  \log_\gamma \frac{C_1\tau \epsilon}{\delta_0} \Big) / 2} - 1 < \Big( K_\epsilon -  \log_\gamma \frac{C_1\tau\epsilon}{\delta_0} \Big) / 2$, then because $\delta_k$ is multiplied by $\gamma^{-1}$ for each $k\in {\cal S}_\epsilon$ and by $\gamma$ for each $k\notin {\cal S}_\epsilon$
    \[  \delta_{K_\epsilon-1} = \gamma \delta_{K_\epsilon} 
        = \gamma \gamma^{-|{\cal S}_\epsilon |} \gamma^{K_\epsilon - |{\cal S}_\epsilon |} \delta_0  < \gamma C_1\tau \epsilon, 
    \] 
    leading to contradiction. 
    The rest of the bound follows from Lemma \ref{lem:tr_bound_success}. 
\end{proof}

{\bf Discussion.} The complexity result in Theorem \ref{thm:tr_complexity_basic} as well as the bound on the best $\epsilon$ \eqref{eq:epsilon_lowbndtr} involve a lot of constants. Let us consider different cases and reduce the complexity bound in order to make it easier to comprehend. 
In summary, the bounds we have are as follows:
\begin{align*}
    &C_1 = \min\left \{\frac{1}{\eta_2}, \frac{c}{L_1+\kappa_{bhm}}\right \}, 
    \quad C_2 =\eta_1\eta_2{c_1} \min\{1, \frac{\eta_2}{\kappa_{bhm}}\}, \\
    &{c_2} = \frac{(1-\eta_1){c_1}}{2}, 
    \quad \frac{1-\kappa_g}{\kappa_g}> \frac{1}{{c_2}}, \\
    &\epsilon \geq \max\left \{\frac{1}{1-\kappa_g} \sqrt{\frac{4\epsilon_f}{ \gamma^2 C_1^2\min\{ C_2,L_1+\kappa_{bhm}\} }} +\epsilon_g, \frac{1+{c_2}}{{c_2}(1-\kappa_g)-\kappa_g}\epsilon_g \right \}, 
\end{align*}
\begin{equation} \begin{aligned} 
    K_\epsilon &\le \frac{4(\phi(x_0) - \phi^\star)}{C_2 (\gamma C_1\tau  \epsilon)^2}  
    + \log_{\gamma^{-1}} \left( \frac{\delta_0}{C_1\tau\epsilon} \right). 
\end{aligned} \end{equation}

Note that ${c_1} \in (0,1/2)$ and $\eta_1 \in (0,1)$, where neither of which needs to be very close to $0$ or $1$. Thus, ${c_2}$ is a constant in $(0,1/4)$ and not very close to $0$, and having ${c_2}$ in denominator or numerator will have a moderate effect.  
The bounds then heavily depend on the choice of the constants $\eta_2$ and $\kappa_{bhm}$, especially relative to one another and to $L_1$. Let us recall that $L_1$ is the Lipschitz constant of $\nabla \phi(x)$ and may be unknown; $\kappa_{bhm}$ is an  upper bound on the Lipschitz constant of $\nabla m_k(x)$ and can be chosen; and $\eta_2$ is the smallest ratio between $\|g_k\|$ and $\delta_k$ that is acceptable on successful steps and is also chosen. 

If we are able to set $\eta_2\sim L_1\sim \kappa_{bhm}$ then $C_1 \sim \frac{1}{L_1}$, $C_2 \sim L_1$ and 
\begin{equation}\label{eq:epsilon_bnd_bigO}
    \epsilon\geq  {\Omega} \left( \max \left\{ \frac{\sqrt{{L_1\epsilon_f }}}{1-\kappa_g} +\epsilon_g, \frac{1+{c_2}}{{c_2}(1-\kappa_g)-\kappa_g}\epsilon_g \right\} \right)
\end{equation}
and 
\begin{equation*} 
    K_\epsilon \le  \mathcal{O} \left( \frac{L_1(\phi(x_0) - \phi^\star)}{\epsilon^2} 
    + \log_{\gamma^{-1}} \Big( \frac{L_1\delta_0}{\epsilon}\Big) \right). 
\end{equation*}
This complexity is essentially the same as that of gradient descent methods such as line search. The neighborhood of convergence (i.e. best achievable $\epsilon$)  is similar to that for line search, and has the same type of dependencies on $\epsilon_g$ and $\epsilon_f$. 

If $\kappa_{bhm}$ is chosen significantly larger than $L_1$ and $\eta_2\sim\kappa_{bhm}$, then we can simply replace $L_1$ with $\kappa_{bhm}$ as all the above derivations remain valid, with $\kappa_{bhm}$ upper bounding the Lipschitz constant of the gradient of the models and  the function $\phi(x)$. 
If $\eta_2$ is chosen significantly larger than $\max\{L_1, \kappa_{bhm}\}$ then we can assume $L_1=\kappa_{bhm}=\eta_2$ and all above derivations and bounds apply.  
Finally, if $\eta_2$ is chosen significantly smaller than  $\kappa_{bhm}$, and assuming $L_1\sim\kappa_{bhm}$,   $C_1={\cal O}(\frac{1}{\max\{L_1, \kappa_{bhm}\}})$, $C_2={\cal O}(\frac{\eta_2^2}{\max\{L_1, \kappa_{bhm}\}})$  and the \eqref{eq:epsilon_bnd_bigO} is still valid, while the complexity bound increases roughly by a factor of $\left (\frac{\max\{L_1, \kappa_{bhm}\}}{\eta_2}\right )^2$. 
The best achievable accuracy in this case  also increases by a factor of  $\frac{\max\{L_1, \kappa_{bhm}\}}{\eta_2}$
The overall order of complexity can be written as 
\begin{equation*} 
    {\cal O} \left( \frac{\max\{L_1, \kappa_{bhm},\eta_2\}^3}{\eta_2^2} \cdot \frac{(\phi(x_0) - \phi^\star)}{\epsilon^2} \right). 
\end{equation*}

The most problematic assumption of this analysis is  $\frac{1-\kappa_g}{\kappa_g}> \frac{1}{c_2}$, which does not allow $\kappa_g$ to be close to $1$, which rules out multiple common oracles. In fact, the conditions on $c_2$ imply that $\kappa_g<\frac{1}{5}$. Thus our results so far are mostly of interest for  oracles with $\kappa_g=0$.  In this case the achievable  accuracy becomes 
\begin{equation}\label{eq:epsilon_bnd_nokappag}
    \epsilon \geq {\Omega} \left ( \max \left \{  \sqrt{{L_1\epsilon_f }}+\epsilon_g, \frac{1+{c_2}}{{c_2}}\epsilon_g \right \}\right ). 
\end{equation}
 
If $\kappa_g=0$ and  $\epsilon_f$ and $\epsilon_g$ can be adjusted dynamically, i.e. the accuracy of the oracle can be set on each iteration depending on the input parameter $\lambda=\delta_k$, then Algorithm~\ref{alg:TR1} can achieve any desired accuracy. 
We perform this analysis in the next section and then comment on how it expands to certain oracles with a general $\kappa_g\in (0,1)$ but under certain  additional assumptions.

\subsubsection{Analysis of the trust region method based on inexact oracles with dynamic accuracy and $\kappa_g=0$}
As before we consider Algorithm~\ref{alg:TR1} based on Oracles~\ref{def:oracle_det_0}  and  \ref{def:oracle_det_1}, 
with $\delta_k$ serving as input to the oracles at each iteration. Specifically, we have the following assumption. 
\begin{assumption} \label{assum:TR2}
    Oracle~\ref{def:oracle_det_1} takes $\delta$ as input and satisfies $\epsilon_g= \kappa_{\delta g}\delta$ and $\kappa_g=0$.  Oracle~\ref{def:oracle_det_0}, takes $\delta$ as input and satisfies $\epsilon_f=\kappa_{\delta f}\delta^2$. 
\end{assumption}

We will consider generalizing this condition to account for $\kappa_g\neq 0$ at the end of this section. 

With the addition of Assumption~\ref{assum:TR2}, Lemma~\ref{lem:TRD0D1_model_error} can be refined as follows. 
\begin{lemma}[model error] \label{lem:TR_fullin_model_error}
    Under Assumptions~\ref{assum:lip_cont}, \ref{assum:TR1} and \ref{assum:TR2}, for any $s\in B(x_k,\delta_k)$  
    \[
|m(x_k+s)-\phi(x_k+s)|\leq \kappa_{em}\delta^2
\]
with $\kappa_{em} = \tfrac{1}{2}(L_1 +  \kappa_{bhm} + 2\kappa_{\delta g})$. 
\end{lemma}
\begin{proof}
    This new bound can be obtained by inserting $\epsilon_g = \kappa_{\delta g}\delta$ and $\kappa_g = 0$ into the result of Lemma~\ref{lem:TRD0D1_model_error}. 
\end{proof}

With these oracle properties, the result of Lemma \ref{lem:TRD0D1_success} can be simplified and strengthened. 
\begin{lemma}\label{lem:TR1_success}
    Under Assumptions~\ref{assum:lip_cont}, \ref{assum:TR1} and \ref{assum:TR2}, if 
        \begin{equation} \label{eq:TR1kappa_success}
   \delta_k \le C_1 \|\nabla \phi(x_k)\|, 
    \end{equation}
    where
    \begin{equation} \label{eq:C1}
  C_1 = \bigg(\max\left\{\eta_2,\ \kappa_{bhm},\ \frac{2\kappa_{\delta f}+\kappa_{em}}{(1-\eta_1){c_1}}\right\}+\kappa_{\delta g}\bigg)^{-1}, 
    \end{equation}
    then $\rho_k \ge \eta_1$, $\|g_k\| \ge \eta_2 \delta_k$, and $x_{k+1} = x_k^+$. 
\end{lemma}

\begin{proof}
    By Assumption \ref{assum:TR2},  
    \begin{equation}\label{eq:TR1_success_aux}
        \|\nabla \phi(x_k)\| \le \|g_k\| + \kappa_{\delta g}\delta_k. 
    \end{equation}
    Combining this with \eqref{eq:TR1kappa_success} leads to 
     \begin{align*}
        (\max\{\kappa_{bhm},\eta_2\} + \kappa_{\delta g}) \delta_k &\le \|\nabla \phi(x_k)\| \le \|g_k\| + \kappa_{\delta g}\delta_k \\
        \max\{\kappa_{bhm},\eta_2\} \delta_k &\le \|g_k\|.
    \end{align*} 
    This not only proves that $\|g_k\| \ge \eta_2\delta_k$, but also indicates $m_k(x_k) - m_k(x_k^+) \ge {c_1} \|g_k\| \delta_k$ in \eqref{eq:fcd_cond}. 
    Then, using Lemma~\ref{lem:TR_fullin_model_error}, the ratio
    \[ \begin{aligned}
        \rho_k &= \frac{m(x_k) - m(x_k^+) + (f(x_k) - m(x_k)) - (f(x_k^+) - m(x_k^+))}{m(x_k) - m(x_k^+)} \\
        &=\frac{m(x_k) - m(x_k^+) + (f(x_k) - \phi(x_k)) - (f(x_k^+) - \phi(x_k^+))+ (\phi(x_k^+) - m(x_k^+))}{m(x_k) - m(x_k^+)} \\
        &\ge 1 - \frac{(2\kappa_{\delta f}+\kappa_{em}) \delta_k^2}{m(x_k) - m(x_k^+)} 
        \ge 1 - \frac{(2\kappa_{\delta f} +\kappa_{em})\delta_k^2}{{c_1} \|g_k\| \delta_k } \\
        &\ge 1 - \frac{(2\kappa_{\delta f} +\kappa_{em})\delta_k}{{c_1} (\|\nabla \phi(x_k)\| - \kappa_{\delta g}\delta_k)} \ge \eta_1, 
    \end{aligned} \]
    where the last step is true because $\|\nabla \phi(x_k)\| \ge \left(\frac{(2\kappa_{\delta f}+\kappa_{em})}{(1-\eta_1){c_1}} + \kappa_{\delta g}\right ) \delta_k$ is assumed by  \eqref{eq:TR1kappa_success}. 
\end{proof}

Lemma \ref{lem:TRD0D1_progress} also simplifies.

\begin{lemma}[guaranteed progress in each step] \label{lem:TR1_progress}
    Let Assumptions~\ref{assum:lip_cont},  \ref{assum:TR1} and \ref{assum:TR2} hold. 
    For each $k$, if iteration $k$ is successful, then we have  
    \begin{equation}
        \phi(x_k) - \phi(x_{k+1}) \ge C_2\delta_k^2 \text{, where } C_2 = \eta_1\eta_2{c_1} \min\{1, \frac{\eta_2}{\kappa_{bhm}}\}-2\kappa_{\delta f}; 
    \end{equation}
    otherwise, we have $x_{k+1} = x_k$ and $\phi(x_k) - \phi(x_{k+1}) = 0$.  
\end{lemma}
\begin{proof}
    When the trial point is accepted, we have $\rho_k \ge \eta_1$ and $\|g_k\| \ge \eta_2 \delta_k$, which guarantees 
    (due to our assumption on Oracle  \ref{def:oracle_det_0} and $\delta_k$)
    \begin{align*}
        \phi(x_k) - \phi(x_{k+1}) 
        &\geq  f(x_k) - f(x_k^+) - |f(x_k) - \phi(x_k) |-|f(x_k^+) - \phi(x_k^+) |\\
        &\ge \eta_1 \big(m(x_k) - m(x_k^+)\big)-2\kappa_{\delta f}\delta_k^2    \ge \eta_1 {c_1} \|g_k\| \min\bigg\{\frac{\|g_k\|}{\|H_k\|}, \delta_k\bigg\}  -2\kappa_{\delta f}\delta_k^2\\
        &\ge \eta_1  \eta_2{c_1} \delta_k \min\bigg\{\frac{\eta_2\delta_k}{\kappa_{bhm}}, \delta_k\bigg\} -2\kappa_{\delta f}\delta_k^2
        = C_2 \delta_k^2. 
    \end{align*}
\end{proof}

This lemma is only useful when $C_2 > 0$, which implies a condition on $\eta_2$ with respect to the error term $\kappa_{\delta f}$ of the zeroth-order oracle. This condition essentially dictates that $\eta_2$ should not be very small with respect to $\kappa_{\delta f}$ (and $\kappa_{bhm}$). 

We want to note here that instead of implying the lower bound $\|g_k\|\geq \eta_2\delta_k$ in the analysis, one can utilize other conditions to impose a bound on $\|g_k\|$. However, we prefer to keep the condition on $\eta_2$. It simplifies the analysis and has a reasonable interpretation: since $\delta_k$ controls the accuracy of Oracles~\ref{def:oracle_det_1} and \ref{def:oracle_det_0} called by the algorithm, it is reasonable to ensure that it is not too large compared to the evaluated gradient. 
Finally, this condition was first found useful in the case of trust region methods based on probabilistic oracles \cite{Bandeiraetal2014,ChenMenickellyScheinberg2014}.

The complexity result is then derived identically to the general case, but with modified values of constants $C_1$ and $C_2$, $\tau=1$ and using the expression for $\kappa_{em}$. The final result is as follows.

\begin{theorem}\label{thm:tr_complexity_fullin}
    Let Assumptions~\ref{assum:lip_cont}, \ref{assum:TR1}  and \ref{assum:TR2} hold.  For  any given $\epsilon $, suppose the initial trust-region radius $\delta_0 \ge \gamma C_1 \epsilon$. Assume that $\eta_2$ is chosen so that
 $C_2 = \eta_1\eta_2{c_1} \min\{1, \frac{\eta_2}{\kappa_{bhm}}\}-2\kappa_{\delta f} >0$.  
   \begin{equation*} \begin{aligned} 
        K_\epsilon &\le 2  |{\cal S}_\epsilon | + \log_\gamma \frac{C_1\epsilon}{\delta_0}  
        \leq \frac{4(\phi(x_0) - \phi^\star)}{C_2 (\gamma C_1 \epsilon)^2} + \log_\gamma \frac{C_1\epsilon}{\delta_0} \\  
        &\le \frac{4\Big(\max\left\{\eta_2,\ \kappa_{bhm},\ \frac{4\kappa_{\delta f}+L_1 +  \kappa_{bhm} + 2\kappa_{\delta g}}{2(1-\eta_1){c_1}}\right\}+\kappa_{\delta g}\Big)^2}{\gamma^2 (\eta_1\eta_2{c_1} \min\{1,\frac{\eta_2}{\kappa_{bhm}}\}-2\kappa_{\delta f})} \cdot \frac{\phi(x_0) - \phi^\star}{\epsilon^2} \\
        &\quad + \log_{\gamma^{-1}} \frac{\delta_0}{\Big(\max\left\{\eta_2,\ \kappa_{bhm},\ \frac{4\kappa_{\delta f}+L_1 +  \kappa_{bhm} + 2\kappa_{\delta g}}{2(1-\eta_1){c_1}}\right\}+\kappa_{\delta g}\Big)\epsilon}. 
    \end{aligned} \end{equation*}
\end{theorem}

The advantage of this new result is that it shows the method achieving any accuracy $\epsilon$ under the assumption 
that Oracles  \ref{def:oracle_det_1} and \ref{def:oracle_det_0}  can deliver appropriate accuracy. It is easy to see that if the best possible accuracy of the oracles is bounded from below, then an appropriate lower bound on $\epsilon$ will need to be imposed to guarantee that $\delta_k$ satisfying Assumption \ref{assum:TR2} and condition $\delta_k\leq C_1\|\nabla \phi(x)\|$ exists. 
In several oracle settings, such as in the finite difference gradient approximation, $\kappa_{\delta g}$ depends on $L_1$ but also on the dimension $n$ and it often dominates other constants in the complexity bound. We will show examples of complexity bounds for some specific oracles cases below. 

\subsubsection{Extension to the case $\kappa_g\in (0,1)$.}

Now let us consider the case when $\kappa_{g}\in (0,1)$. The  changes in the analysis are in Lemma 
\ref{lem:TR1_success} where  \eqref{eq:TR1_success_aux} becomes 
     \begin{align*}
     (1-\kappa_g)\|\nabla \phi(x_k)\| \le \|g_k\| + \kappa_{\delta g}\delta_k 
    \end{align*} 
    and in the fact that Lemma \ref{lem:TR_fullin_model_error} no longer holds. Note, however, that we only need the bound of Lemma~\ref{lem:TR_fullin_model_error} to hold at $x_k^+$ instead of 
any $x_k+s$ with $\|s\|\leq \delta$. We can impose the following assumption on $x_k^+$ computed by the algorithm. 
\begin{assumption} \label{assum:TR3}
There exists a constant $\kappa_{em}$ such that $|m_k(x_k^+)-\phi(x_k^+)|\leq \kappa_{em}\delta_k^2$ holds at every iteration $k$. 
\end{assumption}

Under this assumption Lemma \ref{lem:TR1_success} holds with 
\begin{equation*} 
    C_1 = (1-\kappa_g)\bigg(\max\left\{\eta_2,\ \kappa_{bhm},\ \frac{2\kappa_{\delta f}+\kappa_{em}}{(1-\eta_1){c_1}}\right\}+\kappa_{\delta g}\bigg)^{-1}, 
\end{equation*}
and all other derivations are carried out without change. 

Thus, the final complexity bound in Theorem~\ref{thm:tr_complexity_fullin} applies under the additional Assumption~\ref{assum:TR3}, for any $\epsilon>0$ with a modified value of $C_1$, which means the complexity is multiplied by 
$(1-\kappa_g)^{-2}$. 

In some special cases, for example when $\kappa_g>0$ happens due to projection of the gradient onto a subspace, 
the result can be improved by replacing factor $1-\kappa_g$ by $\sqrt{1-\kappa_g^2}$ and the extra factor in complexity bound becomes $(1-\kappa_g^2)^{-1}$.  Specifically, in \cite{ChaudhryScheinberg2026ICM,ChaudhryScheinbergSun2026} a trust region method whose iterations are computed in randomly generated subspaces is analyzed under Oracle~\ref{def:oracle_det_0} and a particular 
somewhat stronger version of Oracle~\ref{def:oracle_det_1} stated by the assumption below. 
\begin{assumption} \label{assum:TR4}
    Assume that Oracle  \ref{def:oracle_det_1}, with input $\lambda=\delta$ computes $g(x)$ such that the error $g(x)-\nabla \phi(x)$ is decomposable in two vectors $g(x)-d(x)$ and $d(x)-\nabla \phi(x)$ 
    so that $g(x)\perp d(x)-\nabla \phi(x)$ and $d(x)\perp d(x)-\nabla \phi(x)$, moreover 
    \[
        \|d(x)-\nabla \phi(x)\|^2\leq \kappa_g^2\|\nabla \phi(x)\|^2
    \]
    and
    \[
        \|g(x)-d(x)\|^2\leq \epsilon_g^2=\kappa^2_{\delta g}\delta^2
    \]
\end{assumption}

Assumption  \ref{assum:TR3} is also satisfied in \cite{ChaudhryScheinberg2026ICM,ChaudhryScheinbergSun2026} and complexity results similar to Theorem~\ref{thm:tr_complexity_fullin}  are obtained. However, since the algorithm employs random subspaces and analysis needs to be extended to bound expected complexity   using techniques for probabilistic oracles, which combine the analysis we present here, with stochastic processes and is   beyond the scope of this article.

\newpage

\section{Stochastic Oracles with Bounded Second Moments}\label{sec:stochastic}

We now turn our exposition to the case of stochastic zeroth- and first-order oracles, that is, oracles whose output involves realizations of one or more random variables. This randomness can be ``external'' when the oracle  computation involves some random data coming from a possibly unknown distribution, or it can be ``internal'' when this computation employs some form of randomization. Both types of oracles are ubiquitous in modern optimization literature and we will not attempt to give a  full review of existing work in this article. Instead we will try to extract some key properties and their effect on the basic optimization schemes, such as fixed step size gradient descent. We aim to provide one reasonably concise and unifying framework that will cover some particularly popular cases of both oracles and analysis and by doing so we hope to provide a clear route for analyzing other possible stochastic oracles. 

We will limit our discussions to oracles whose error is bounded in expectation in some form. While not all stochastic oracles fall under this category, many commonly used stochastic oracles satisfy such bounds, as illustrated by the examples that follow. 
For zeroth-order oracles, we consider potentially biased oracles and present two different models of bounded variance. For first-order oracles, we present three related but distinct types, characterized by their error bounds. These types differ in how they bound the bias and variance of the estimation error (jointly or separately), but in all cases, the bounds consist of two components: an absolute error term and a relative error term. As in the case of deterministic oracles we will observe that the effect of the relative error term is quite different from that of the absolute error terms. 

Analogously to the deterministic case, we will base most of our analysis on a descent inequality, first stated in the abstract as a version of \eqref{eq:key_reduction}, but with more terms and holding  in expectation. We will then show that this descent inequality 
holds for our examples of oracles and derive complexity bounds for a fixed step size gradient descent method based on a  variety of cases. We will show that many existing complexity results can be recovered using this simple approach. However, we will observe some new results as well. Firstly, in the case of convex functions the telescopic sum analysis (when applicable) gives better complexity bounds than the analysis based on the linear KL condition (Assumption \ref{assum:KLcond}). This is consistent with some past results, but has not been clearly exposed in the prior literature. On the other hand, analysis under the  linear KL condition is more robust to noise. Secondly, when analyzing gradient descent based on first-order oracles that are constructed via randomized finite differences and utilize certain stochastic zeroth-order oracle, we discover some differences in how complexity depends on the dimension depending on what random directions are used. Finally, we clearly observe the difference in complexity of all finite difference gradients based on stochastic zeroth-order oracles, depending on the type of the zeroth-order oracle. The phenomenon was noted in simulation optimization under the term  ``common-random-numbers'' however, and simple and clear connection to the improved complexity has not been presented to the best of our knowledge. 

We will not analyze line search and trust region methods for stochastic oracles, as the analysis requires more involved techniques with respect to stochastic processes, however, examples of such analysis can be found in \cite{BlanchetCartisMenickellyScheinberg2018,paquette2018stochastic,jin2021high}.

 \subsection{Stochastic zeroth-order oracles}

Stochastic zeroth-order oracles are ubiquitous in modern applications, including simulation, machine learning, artificial intelligence, reinforcement learning, and many others.  
In optimization algorithms, such oracles can be used directly to evaluate optimization progress. They can also be used to construct first-order oracle estimates, for example through finite-differences or other methods discussed in Section~\ref{sec.examples_first}. In this article  we will mainly focus on the latter use.

 It is important to distinguish between two fundamentally different settings of stochastic zeroth-order oracles. 
In the first setting, the oracle is queried at a point $x$ and returns a random estimate of the objective value $\phi(x)$, generated from a realization of an internal random variable. If an estimate at another point $x'$ is needed, a separate oracle call is made using an independent realization of the randomness. 
We will refer to it as a ``single'' stochastic zeroth-order oracle.
In the second setting, we will consider  a stochastic oracle that,
with one realization of the  random variable, creates an access to a random function that approximates $\phi(\cdot)$ and that can be measured at multiple points. This function is assumed to have some desirable properties. We will refer to this oracle as ``multi-point'' stochastic zeroth-order oracle. 
The accuracy properties of the stochastic zeroth-order oracles can be quite different under these two settings and have a profound effect on the properties of the algorithm discussed later in this section. 
Below we will introduce two stochastic zeroth-order oracle definitions that capture these two settings. 
 Here and throughout  $f(\cdot,\xi^{(0)})$ is viewed as a random function and 
 expectations, such as $\E_{\xi^{(0)}} [f(\cdot,\xi^{(0)})]$ are understood pointwise, that is,
 \[
 \E_{\xi^{(0)}} [f(\cdot,\xi^{(0)})](x) = \E_{\xi^{(0)}} [f(x,\xi^{(0)})]. 
 \]
 
 \begin{customoracle}{S-0.1}[Single-Point Stochastic Zeroth-Order Oracle] \label{def:oracle_stoch_0.1}
    Given a point $x\in \mathbb{R}^n$ and an input $\lambda$, the oracle returns a random value $f(x,\xi^{(0)})$ that satisfies 
        \begin{equation} \label{eq:0_oracle_bias}
        (\E_{\xi^{(0)}} [f(x,\xi^{(0)})] - \phi(x))^2 \leq \epsilon_f^2(\lambda), 
    \end{equation} 
    and for $e(x,\xi^{(0)})=f(x,\xi^{(0)}) - \E_{\xi^{(0)}} [f(x,\xi^{(0)})]$
 \begin{equation} \label{eq:0_oracle_var}
        \E_{\xi^{(0)}} [e(x,\xi^{(0)})^2] \leq \sigma_f^2(\lambda), 
    \end{equation} 
    where  $\epsilon_f^2(\lambda)$ and $\sigma_{f}^2(\lambda)$  are some functions of the input $\lambda$ and $x$ that generally may depend on the distribution of $\xi^{(0)}$ and   properties of $\phi$. 
\end{customoracle}
This oracle subsumes the deterministic Oracle \ref{def:oracle_det_0}
 as a special case where $\sigma_f(\lambda)=0$, but also includes many stochastic settings, examples of which we will list below.

 \begin{customoracle}{S-0.2}[Multi-Point Stochastic Zeroth-Order Oracle] \label{def:oracle_stoch_0.2}
    Given a point $x\in \mathbb{R}^n$,  a finite query set ${Y}\subset   \mathbb{R}^n$   and an input parameter $\lambda$, the oracle samples a realization $\xi^{(0)}$ of the underlying randomness and returns the values $f(x, \xi^{(0)})$ and $\{ f(y,  \xi^{(0)}): \, y\in { Y}\}$ where the random function $f(\cdot, \xi^{(0)})$ satisfies 
    \begin{equation} \label{eq:0_oracle_biases}
        \left (\E_{\xi^{(0)}} [f(y,\xi^{(0)})] - \phi(y)\right )^2 \leq \epsilon_f^2(\lambda), \quad \forall y\in {Y}\cup \{x\}
            \end{equation} 
    and for $e(\cdot,\xi^{(0)})=f(\cdot,\xi^{(0)}) - \E_{\xi^{(0)}} [f(\cdot,\xi^{(0)})]$
    \begin{align} \label{eq:0_oracle_lipsch}
       & \E_{\xi^{(0)}} \left [\left ( e(y,\xi^{(0)})- e(x, \xi^{(0)})\right)^2\right ]  \leq \left(L_c(\lambda)^2+L_v(\lambda)^2\|\nabla \phi(x)\|^2\right) \|y-x\|^2, \quad \forall y\in {Y}
    \end{align} 
    where  $\epsilon_f^2(\lambda)$, $L_c(\lambda)$ and $L_v(\lambda)$ are some functions of the input $\lambda$ and $x$ that generally may depend on the distribution of $\xi^{(0)}$.

\end{customoracle}

One important instance of this oracle is obtained under the {\em common-random-number} paradigm in simulation optimization.  The term "common random numbers" refers to a simulation process which can be repeated for multiple inputs with the same sequence of internal random variables (often generated from  the same seed).   This setting is particularly useful in simulation of queueing systems. Details can be found in \cite{Glynn1987,GlassermanYao1992, Fu1994, Kleinman1999}.
We note that the exact properties of the resulting oracle differ somewhat in the literature. Here we impose assumptions that are, to the best of our knowledge, more general than those appearing in the existing literature while still retaining the improved bounds associated with common random numbers.

In essence the multi-point oracle captures the setting where  estimates of $\phi(\cdot)$ around a given point $x$ are generated by sampling a realization of the random  function $f(\cdot,\xi^{(0)})$ which has some nice properties. For example, a locally Lipschitz continuous function (around $x$) with the square of the local Lipschitz constant bounded in expectation by $L_c^2+L_v^2\|\nabla \phi(x)\|^2$ satisfies the oracle properties. 
At the cost of somewhat more cumbersome notation, this assumption is a relaxed  version of the common assumption in  the simulation optimization literature where  a global bound is assumed on the expectation of the   Lipschitz constant. 
Like in the case of the deterministic first-order oracle (and later with the stochastic one), it is useful to consider the term 
$L_v^2\|\nabla \phi(x)\|^2$ separately because it has a less harmful effect on convergence properties of the  resulting algorithms than the $L_c^2$ term. On the other hand, even simple functions like quadratics are not globally Lipschitz continuous, thus accounting for the dependence of the oracle bound  on $\|\nabla \phi(x)\|^2$ is essential. 

The important difference between the two oracle definitions is that the first one essentially bounds the variance of the function noise, while the second only bounds the rate of change of that noise - the noise itself may be quite large. This difference is especially consequential when the oracle is used in finite difference schemes. 

We note that other stochastic zeroth-order oracle cases are also possible, for example, where the random variables used in consecutive oracle calls are drawn from a time series or exhibit other forms of dependence. Many more complex oracles than those presented here exist  in various applications  and are considered in the literature. However, here we limit our consideration to Oracles \ref{def:oracle_stoch_0.1} and \ref{def:oracle_stoch_0.2} as they are the most common and yet give a rich variety of cases to analyze. 

 \subsection{Examples of stochastic zeroth-order oracles}
 
 We now present several examples of stochastic zeroth-order oracles and examine how their errors can be characterized within the models introduced above.
 
 \subsubsection{A single-point  oracle with Gaussian noise}
 The simplest example of a stochastic zeroth-order oracle is
 \[
 f(x,\xi^{(0)})=\phi(x)+\xi^{(0)},\quad \xi^{(0)} \sim {\cal N}(\epsilon_f, \sigma_f^2)
 \]
 This clearly satisfies conditions of Oracle \ref{def:oracle_stoch_0.1}. While it is not clear how often such oracles arise in practice it is one of the 
 most commonly used models  when testing algorithms and software as well as theoretical derivations. 
 One important property of Gaussian noise is that it is light-tailed, however,  our zeroth-order oracle definition in this section only requires finite variance. In this respect 
 Gaussian noise oracle does not have an advantage over other oracle examples we  consider below.

 \subsubsection{ Single-point oracles in expectation minimization}
 Let us consider one of the most common applications of stochastic oracles --- the expected risk minimization, or simply expectation minimization. 
In this case, the objective function is defined as  
\begin{equation} \label{eq:ERM} \tag{EM}
	\phi(x)=\E_{\xi\sim\Xi}[\ell(x,\xi)]
\end{equation}
where $x$ is the vector of model parameters over which optimization is performed, $\xi$ is a data sample drawn from the distribution $\Xi$ that typically consists of a feature vector $z$ and a label $y$, and $\ell(x,\xi)$ is the loss incurred by the model parameterized by $x$ on the sample $\xi$. This setting includes the empirical risk minimization (ERM) as a special case when ${\Xi}$ is a uniform distribution over a finite set of samples. 
In this setting, a natural stochastic zeroth-order oracle returns the random loss value $\ell(x,\xi)$ for a sampled realization of $\xi$. This oracle is unbiased, since $\E_{\xi\sim\Xi}[\ell(x,\xi)]=\phi(x)$. A common accuracy assumption is that the stochastic loss has bounded variance, namely
\[
    \E_{\xi\sim\Xi}\left[(\ell(x,\xi)-\phi(x))^2\right] \leq V(x),
\]
where $V(x)$ may depend on the point $x$.
For example, if $\ell(x,\xi)$ is the zero-one loss or a sigmoid loss, then $(\phi(x)-\ell(x,\xi))^2 \leq 1$ deterministically for all $x$. On the other hand, if $\ell(x,\xi)$ is a least squares loss then $V(x)$ depends on the data distribution, but is  finite under suitable fourth-moment assumptions.

For \eqref{eq:ERM}, the zeroth-order oracle can be implemented by applying sample averaging over a set of i.i.d. samples.
Given a positive integer $\lambda$ as input, the oracle generates $\xi_1,\ldots,\xi_\lambda$ as i.i.d.\ samples from $\Xi$. We then define
\[
    \xi^{(0)}:=(\xi_1,\ldots,\xi_\lambda),
\]
and set
\[
    f(x,\xi^{(0)})=\frac1\lambda\sum_{j=1}^\lambda \ell(x,\xi_j).
\]

From the basic properties of sample averages, we have 
\begin{equation} 
    \E_{\xi^{(0)}}[f(x, \xi^{(0)})] = \phi(x) \text{ and } 
	\E_{\xi^{(0)}}\big[\big(f(x,\xi^{(0)})-\phi(x)\big)^2\big] \leq \frac{V(x)}{\lambda}.
\end{equation} 
To treat this as an example of Oracle~\ref{def:oracle_stoch_0.1},  we need to have a uniform bound on $V(x)$. 
This usually applies if $x$ is restricted to a bounded domain. In conclusion, if $V(x)\leq V$, then an example of Oracle~\ref{def:oracle_stoch_0.1} with $\lambda$ as the input and $\epsilon_f = 0$ and $\sigma_f^2 = V / \lambda$. 
The cost of this oracle  depends linearly on the  size of the minibatch $\lambda$, which, in principle, can be chosen to depend on  the current state of whatever algorithm utilizes this oracle.

\subsubsection{Multi-point oracles in expectation minimization with globally Lipschitz loss}
Now, let us consider the \eqref{eq:ERM} problem with the additional assumption  that for every $\xi\in \Xi$, 
the function $\ell(x,\xi)$ is Lipschitz continuous with constant $L_\ell$ with respect to $x$, i.e.,
\[
    |\ell(x,\xi) - \ell(y,\xi)|\le L_\ell \|x-y\|,\qquad \forall x,y\in \mathbb{R}^n.
\]
This assumption is satisfied by many commonly used loss functions such as logistic loss and hinge loss, when data vectors are uniformly bounded.
By Jensen's inequality, the objective function $\phi$ is also Lipschitz continuous with constant $L_\ell$: 
\[ \begin{aligned}
    |\phi(x) - \phi(y)| &= \left|\E_{\xi\sim\Xi} \bigl[\ell(x,\xi) - \ell(y,\xi)\bigr]\right| \\
    &\le \E_{\xi\sim\Xi} \left[ \bigl|\ell(x,\xi) - \ell(y,\xi)\bigr| \right]
    \le L_\ell \|x-y\|,\qquad \forall x,y\in \mathbb{R}^n.
\end{aligned} \]

Fix a positive integer $\lambda$, and let $\xi_1,\ldots,\xi_\lambda$ be i.i.d.\ samples from $\Xi$. Define
\[
\xi^{(0)}:=(\xi_1,\ldots,\xi_\lambda)
\]
and a function
\[
f(x,\xi^{(0)}):=\frac1\lambda\sum_{j=1}^\lambda \ell(x,\xi_j),
\]
We observe that the oracle uses common random numbers, since the same realization $\xi^{(0)}$ is used to define the function of  $x$, that can be queried at multiple points and which satisfies
\begin{equation}
 \E_{\xi^{(0)}}[f(x,\xi^{(0)})]=\phi(x)
\end{equation}
so that \eqref{eq:0_oracle_biases} holds with
\[
\epsilon_f(\lambda)=0.
\]
Now,  define $e(x,\xi^{(0)}) := f(x,\xi^{(0)}) - \E_{\xi^{(0)}}[f(x,\xi^{(0)})] = f(x,\xi^{(0)})-\phi(x)$. 
Then for every realization of $\xi^{(0)}$ and $x,y\in \mathbb{R}^n$, using independence of samples in $\xi^{(0)}$
\begin{align*}
    \E \left [(e(x,\xi^{(0)}) - e(y,\xi^{(0)}))^2\right ]&=
        \frac1\lambda\sum_{j=1}^\lambda  \E \left [(e(x,\xi_j) - e(y,\xi_j))^2\right ]\\
       &\le \frac1\lambda\sum_{j=1}^\lambda  \E \left [(\ell(x,\xi_j) - \ell(y,\xi_j))^2\right ]
    \le \frac{L_\ell^2 }{\lambda}\|x-y\|^2.
\end{align*}

and therefore \eqref{eq:0_oracle_lipsch} holds and  the conditions of Oracle \ref{def:oracle_stoch_0.2} are satisfied with
\[\epsilon_f=0, \qquad
    L_c^2= \frac{L_\ell^2 }{\lambda},
    \qquad
    L_v=0.
\]

The cost of the oracle depends linearly on $\lambda$ and the number of points in  $\{x\cup Y\}$.

\subsubsection{Multi-point oracles for least-squares loss}\label{sec:leastsquares}
The assumption that each $\ell(x,\xi)$ is Lipschitz with a global Lipschitz constant is quite strong and is violated by even very classical loss functions, such as least-squares loss. Fortunately in such cases it is often possible to satisfy conditions \eqref{eq:0_oracle_lipsch} with $L_v>0$. 

To illustrate this, we consider a specific case of \eqref{eq:ERM} where
\[
    \ell(x,\xi)=\frac12(a(\xi)^Tx-b(\xi))^2,  
\]
with the distribution of $(a(\xi), b(\xi))$ satisfying the following properties: 
\begin{itemize}
    \item
        there exists a constant $M_1 \ge 0$ such that 
        \begin{equation} \label{assum:M1}
            \E_{\xi\sim\Xi}\bigl[\bigl(a(\xi)a(\xi)^T\bigr)^2\bigr]
            \preceq (M_1+1) \bigl(\E_{\xi\sim\Xi}\bigl[a(\xi)a(\xi)^T\bigr]\bigr)^2; 
        \end{equation}
    \item
        there exists a constant $M_2 > 0$ such that 
        \begin{equation} \label{assum:M2}
            \E_{\xi\sim\Xi}\bigl[a(\xi)a(\xi)^T\bigr] \preceq M_2 I. 
        \end{equation}
\end{itemize}

The meaning of the first property will be discussed later. 
Now, as before, given a positive integer $\lambda$ as input, the oracle generates $\xi_1,\ldots,\xi_\lambda$ as i.i.d.\ samples from $\Xi$. We then define
\[
    \xi^{(0)}:=(\xi_1,\ldots,\xi_\lambda),
\]
and set
\[
    f(x,\xi^{(0)})=\frac1\lambda\sum_{j=1}^\lambda \ell(x,\xi_j).
\]
This is again an unbiased estimator since $\E_{\xi^{(0)}}[f(x,\xi^{(0)})]=\phi(x)$ and thus \eqref{eq:0_oracle_biases} holds with $\epsilon_f(\lambda)=0$.

Let $x^\star$ be the global minimizer of $\phi$, and let 
\begin{equation*} 
    \nabla f(x^\star,\xi^{(0)}) = \frac1\lambda\sum_{j=1}^\lambda \nabla\ell(x^\star,\xi_j) 
\end{equation*}
be the sample gradient at this optimum. Since
$
    \nabla\phi(x^\star)=0,
$
it follows that
\[
    \E_{\xi^{(0)}}[\nabla f(x^\star,\xi^{(0)})] = 0.
\]
Define
\[
    H_{\xi^{(0)}}:=\frac1\lambda\sum_{j=1}^\lambda a(\xi_j)a(\xi_j)^T
    \quad \text{ and } \quad 
    H:=\E_{\xi\sim\Xi}[a(\xi)a(\xi)^T]. 
\]
Since $f(x,\xi ^{(0)})$ is a quadratic function of $x$ with Hessian $H_{\xi ^{(0)}}$, we have 
\begin{equation*} 
    f(x,\xi ^{(0)}) = f(x^\star,\xi ^{(0)}) + \nabla f(x^\star,\xi^{(0)})^T(x-x^\star)
    + \frac12 (x-x^\star)^T H_{\xi ^{(0)}} (x-x^\star).
\end{equation*}

Therefore, given $x,y\in{\mathbb R}^n$, 
\begin{align*}
    f(y,\xi ^{(0)}) - f(x,\xi ^{(0)})
    &= \nabla f(x^\star,\xi^{(0)})^T (y-x) 
    + (x-x^\star)^T H_{\xi ^{(0)}} (y-x) \\
    &\quad + \frac12 (y-x)^T H_{\xi ^{(0)}} (y-x).
\end{align*}

Taking expectations and using
$
    \E_{\xi^{(0)}}[\nabla f(x^\star,\xi^{(0)})] = 0 
$
gives
\begin{align*}
    \phi(y)-\phi(x)
    &= (x-x^\star)^T H (y-x) 
    + \frac12 (y-x)^T H (y-x).
\end{align*}
Hence, 
\begin{align*}
    e(y,\xi^{(0)}) - e(x,\xi^{(0)}) 
    &= \Big(f(y,\xi ^{(0)}) - f(x,\xi ^{(0)})\Big) - \Big(\phi(y) - \phi(x)\Big) \\
    &= \nabla f(x^\star,\xi^{(0)})^T(y-x) 
    + (x-x^\star)^T (H_{\xi ^{(0)}}-H) (y-x) \\
    &\quad + \frac12 (y-x)^T (H_{\xi ^{(0)}}-H) (y-x).
\end{align*}
Then, using $(a+b+c)^2 \le 3a^2+3b^2+3c^2$, 
\begin{align*}
    &\E_{\xi^{(0)}}\left[\left(
        e(y,\xi^{(0)}) - e(x,\xi^{(0)}) 
    \right)^2\right] \\
    &\le 3\E_{\xi^{(0)}}\left[\left(
        \nabla f(x^\star,\xi^{(0)})^T(y-x) 
    \right)^2\right] 
    + 3\E_{\xi^{(0)}}\left[\left(
        (x - x^\star)^T (H_{\xi^{(0)}} - H) (y - x)
    \right)^2\right] \\
    &\quad + 3\E_{\xi^{(0)}}\left[\left(
        \frac{1}{2} (y-x)^T (H_{\xi^{(0)}} - H) (y - x)
    \right)^2\right] \\
    &\le 3\E_{\xi^{(0)}}\left[
        \|\nabla f(x^\star,\xi^{(0)})\|^2
    \right] \|y - x\|^2 
    + 3\E_{\xi^{(0)}} \left[
        \left\|(H_{\xi^{(0)}} - H) (x - x^\star)\right\|^2
    \right] \|y - x\|^2 \\
    &\quad + \frac{3}{4} \left[
        \left\|(H_{\xi^{(0)}} - H) (y - x)\right\|^2
    \right] \|y - x\|^2. 
\end{align*}
Note that $\E_{\xi^{(0)}}[ H_{\xi^{(0)}} ]= H$ and that our assumption indicates 
\[ \begin{aligned} 
    \E_{\xi^{(0)}} [H_{\xi^{(0)}}^2] 
    &= \E_{\xi^{(0)}} \left[ \Big( \frac1\lambda\sum_{j=1}^\lambda a(\xi_j)a(\xi_j)^T \Big)^2 \right] \\ 
    &= \frac{1}{\lambda^2} \left( (\lambda^2 - \lambda) H^2 + \sum_{j=1}^\lambda \E_{\xi_j\sim\Xi} \bigl[\bigl(a(\xi_j)a(\xi_j)^T\bigr)^2\bigr] \right) 
    \preceq \left( 1 + \frac{M_1}{\lambda} \right) H^2. 
\end{aligned} \] 
These lead to 
\[ \begin{aligned}
    \E_{\xi^{(0)}} \left[ (H_{\xi^{(0)}} - H)^2 \right]
    = \E_{\xi^{(0)}} \left[ H_{\xi^{(0)}}^2 \right] - H^2 
    \preceq \frac{M_1}{\lambda} H^2 
\end{aligned} \]
and, for any vector $v\in\mathbb{R}^n$, 
\begin{equation} \label{eq:least_square M1 implication}
    \E_{\xi^{(0)}} \left[ \left\|(H_{\xi^{(0)}} - H) v\right\|^2 \right] 
    = v^T \E_{\xi^{(0)}} \left[ (H_{\xi^{(0)}} - H)^2 \right] v 
    \le v^T \Big(\frac{M_1}{\lambda} H^2\Big) v 
    = \frac{M_1}{\lambda} \|H v\|^2. 
\end{equation}
Thus, 
\begin{align*}
    &\E_{\xi^{(0)}}\left[\left(
        e(x_i,\xi^{(0)}) - e(x_0,\xi^{(0)}) 
    \right)^2\right] \\ 
    &\le \left( 3\E_{\xi^{(0)}}\left[\|\nabla f(x^\star,\xi^{(0)})\|^2\right] + \frac{3M_1}{\lambda} \left\|H(x_0 - x^\star)\right\|^2 + \frac{3M_1}{4\lambda} \|H(x_i - x_0)\|^2\right) \|x_i - x_0\|^2 \\
    &\le \left( 3\E_{\xi^{(0)}}\left[\|\nabla f(x^\star,\xi^{(0)})\|^2\right] + \frac{3M_1}{\lambda} \left\|\nabla\phi(x_0)\right\|^2 + \frac{3M_1M_2^2}{4\lambda} \|x_i - x_0\|^2\right) \|x_i - x_0\|^2. 
\end{align*}

Since $\xi^{(0)}$ is a set of i.i.d. random variables, 
\[
    \E_{\xi^{(0)}}\left[\|\nabla f(x^\star,\xi^{(0)})\|^2\right] 
    = \frac{1}{\lambda} \E_{\xi\sim\Xi}\left[\|\nabla \ell(x^\star,\xi)\|^2\right]. 
\]
Let $h>0$ be a constant such that   $\|y - x\| \le h$, for $x$ and $y\in { Y}$ - input to Oracle ~\ref{def:oracle_stoch_0.2}.  This bound $h$ is not a property of the  Oracle ~\ref{def:oracle_stoch_0.2}, but rather of how it is used. As we will see in later examples the bound $h$ is usually a small  number. 

Given such $h$, the conditions of Oracle ~\ref{def:oracle_stoch_0.2} hold with
\[
    \epsilon_f(\lambda)=0,
    \qquad
    L_c(\lambda)^2 = \frac{3}{\lambda} \E_{\xi\sim\Xi}\left[\|\nabla \ell(x^\star,\xi)\|^2\right] + \frac{3M_1M_2^2h^2}{4\lambda},
    \qquad
    L_v(\lambda)^2 = \frac{3M_1}{\lambda}.
\]

Generally a finite bound on  $\E_{\xi\sim\Xi}\left[\|\nabla \ell(x^\star,\xi)\|^2\right]$ is assumed and can be guaranteed via suitable bounds on the 
moments of $(a,b)$.  When there exists $x^\star \in \mathbb{R}^n$ such that $\phi(x^\star) = 0$, the nonnegative loss $\ell(x^\star,\xi) = 0$ almost surely, and the prediction matches the observation 
\begin{equation} \label{eq:interp}
    a(\xi)^T x^\star = b(\xi) \quad \text{almost surely. } 
\end{equation}
This situation is usually referred to in the literature as {\em interpolation}. 
The closer we are to satisfying interpolation condition the smaller  the value of $\E_{\xi\sim \Xi}\left[\|\nabla \ell(x^\star,\xi )\|^2\right]$. If $h$ is small and interpolation nearly holds, we observe that $L_c$ constant can be quite small, especially compared to $L_v$. Since their effect on algorithmic complexity is different (as we will show later) it is beneficial to track these constants separately.

We now discuss  property 
\[
    \E_{\xi\sim\Xi}\bigl[\bigl(a(\xi)a(\xi)^T\bigr)^2\bigr]
    \preceq (M_1+1) \bigl(\E_{\xi\sim\Xi}\bigl[a(\xi)a(\xi)^T\bigr]\bigr)^2 
\]
which offers great convenience in our derivations but lacks a clear meaning without further explanation. 
We discuss its applicability here to provide some intuition. 
First, Jensen's inequality implies $\E_{\xi\sim\Xi}\bigl[\bigl(a(\xi)a(\xi)^T\bigr)^2\bigr] \succeq \bigl(\E_{\xi\sim\Xi}\bigl[a(\xi)a(\xi)^T\bigr]\bigr)^2$, therefore this property only makes sense when $M_1 \ge 0$. 
When $a(\xi) \sim \mathcal{N}(0,I)$, we have $\E_{\xi\sim\Xi}\bigl[\bigl(a(\xi)a(\xi)^T\bigr)^2\bigr] = (n+2) I$ and $\bigl(\E_{\xi\sim\Xi}\bigl[a(\xi)a(\xi)^T\bigr]\bigr)^2 = I^2 = I$, thus the inequality holds with $M_1 = n+1$. 
Alternatively, when $1 \le \|a(\xi)\| \le R$ for all $\xi$ for some $R>0$, and 
the smallest positive eigenvalue of $\E_{\xi\sim\Xi}\bigl[a(\xi)a(\xi)^T\bigr]$ is at least $\mu$, then 
\[
    \E_{\xi\sim\Xi}\bigl[\bigl(a(\xi)a(\xi)^T\bigr)^2\bigr]
    \preceq R^2 \E_{\xi\sim\Xi}\bigl[a(\xi)a(\xi)^T\bigr]
    \preceq \frac{R^2}{\mu} \bigl(\E_{\xi\sim\Xi}\bigl[a(\xi)a(\xi)^T\bigr]\bigr)^2. 
\]

\subsubsection{Multi-point oracles in expectation minimization with globally Lipschitz-smooth loss}
Let us now consider a general case of a  loss function with  Lipschitz continuous gradients.  
Specifically, under assumption  that 
\[
    \|\nabla\ell(x,\xi) - \nabla\ell(y,\xi)\| \le L_{\nabla\ell} \|x-y\|,\qquad \forall x,y\in \mathbb{R}^n, 
\]
we have 
\[ \begin{aligned}
    \|\nabla\phi(x) - \nabla\phi(y)\| &= \left\|\E_{\xi\sim\Xi} \bigl[\nabla\ell(x,\xi) - \nabla\ell(y,\xi)\bigr]\right\| \\
    &\le \E_{\xi\sim\Xi} \left[ \bigl\|\nabla\ell(x,\xi) - \nabla\ell(y,\xi)\bigr\| \right]
    \le L_{\nabla\ell} \|x-y\|,\qquad \forall x,y\in \mathbb{R}^n.
\end{aligned} \]

Again, given a positive integer $\lambda$ as input, we generate $\xi_1,\ldots,\xi_\lambda$ as i.i.d.\ samples from $\Xi$, define
\[
    \xi^{(0)}:=(\xi_1,\ldots,\xi_\lambda),
\]
and set the oracle output to 
\[
    f(x,\xi^{(0)})=\frac1\lambda\sum_{j=1}^\lambda \ell(x,\xi_j).
\]

Then, since the random variables
$e(y,\xi_j)-e(x,\xi_j)$ are independent and have
zero mean, we have
\begin{align*}
\E_{\xi^{(0)}}\left[
    \left(e(y,\xi^{(0)})-e(x,\xi^{(0)})\right)^2
\right] &= \E_{\xi^{(0)}}\left[
    \left(\frac1\lambda\sum_{j=1}^\lambda
    (e(y,\xi_j)-e(x,\xi_j))\right)^2
\right] \\ &= \frac1\lambda\E_{\xi\sim\Xi}\left[
    (e(y,\xi)-e(x,\xi))^2
\right].
\end{align*}
Moreover,
\begin{align*}
|e(y,\xi)-e(x,\xi)|
&\le
\bigl|\ell(y,\xi)-\ell(x,\xi)
-\nabla\ell(x,\xi)^T(y-x)\bigr|+
\bigl|\phi(y)-\phi(x)
-\nabla\phi(x)^T(y-x)\bigr| \\ &+
\bigl|(\nabla\ell(x,\xi)-\nabla\phi(x))^T(y-x)\bigr|\le L_{\nabla\ell}\|y-x\|^2
+\|\nabla\ell(x,\xi)-\nabla\phi(x)\|\|y-x\|.
\end{align*}

Now, assume that the variance of the gradient is bounded as
\[
\E_{\xi\sim\Xi}\left[
\|\nabla\ell(x,\xi)-\nabla\phi(x)\|^2
\right]
\le M_c^2+M_v^2\|\nabla\phi(x)\|^2
\]
for all $x\in\mathbb{R}^n$
(which is a common assumption, see e.g.
\cite{bottou2018optimization}).
Then, using $(a+b)^2\le 2a^2+2b^2$, we have
\begin{align*}
\E_{\xi^{(0)}}\left[
(e(y,\xi^{(0)})-e(x,\xi^{(0)}))^2
\right] \le \frac{2L_{\nabla\ell}^2}{\lambda}
\|y-x\|^4
+\frac{2M_c^2}{\lambda}\|y-x\|^2 +
\frac{2M_v^2}{\lambda}
\|\nabla\phi(x)\|^2\|y-x\|^2.
\end{align*}

As before, let $h>0$ be a constant such that
$\|y-x\|\le h$ for all $y\in Y$, where $x$
and $Y$ are inputs to
Oracle~\ref{def:oracle_stoch_0.2}.
Then the conditions of
Oracle~\ref{def:oracle_stoch_0.2} hold with
\[
\epsilon_f(\lambda)=0,\qquad
L_c^2(\lambda)=
\frac{2L_{\nabla\ell}^2h^2+2M_c^2}{\lambda},
\qquad
L_v^2(\lambda)=\frac{2M_v^2}{\lambda}.
\]
The cost of the oracle depends linearly on $\lambda$ and the number of points in  $\{x\}\cup Y$.

\subsection{Stochastic first-order oracles definitions} 
We now extend the conditions on deterministic first-order oracles defined in Oracle \ref{def:oracle_det_1} to the stochastic setting, where oracle outputs may be biased and have bounded variance. 
There are several such extensions, each serving a somewhat different purpose. We start with the simplest extension, which replaces the deterministic bound of Oracle \ref{def:oracle_det_1} with the bound on the expectation.

\begin{customoracle}{S-1.1}[Stochastic First-Order Oracle, Case 1] \label{def:oracle_stoch_1.1}
    Given a point $x\in \mathbb{R}^n$ and an input $\lambda$, the oracle generates random output $g(x,\xi^{(1)})$ that satisfies
    \begin{equation} \label{eq:oracle_stoch_1.1}
        \E_{\xi^{(1)}} \left[ \|g(x,\xi^{(1)}) - \nabla \phi(x)\|^2 \right] \leq \epsilon_g^2(\lambda)+\kappa_g^2(\lambda)\|\nabla \phi(x)\|^2, 
    \end{equation} 
    where  $\epsilon_g^2(\lambda)$, $\kappa_g^2(\lambda)$ are some nonnegative real-valued functions of $\lambda$ that generally depend on the global properties of $\phi$ and the distribution of $\xi^{(1)}$. 
\end{customoracle}
This definition clearly includes the deterministic Oracle \ref{def:oracle_det_1} as a special case, so $g(x,\xi^{(1)})$ can be a biased estimator with the square bias that is as large as $\E_{\xi^{(1)}} \left[ \|g(x,\xi^{(1)}) - \nabla\phi(x)\|^2 \right]$. 
On the other hand, if $g(x,\xi^{(1)})$ happens to be an unbiased estimator, i.e. $\E_{\xi^{(1)}}[g(x,\xi^{(1)}] = \nabla \phi(x)$, then the left-hand side of \eqref{eq:oracle_stoch_1.1} is the variance. 
The idea behind using Oracle \ref{def:oracle_stoch_1.1} is to avoid having to bound the bias and variance of the estimator separately, since in some cases it is easier to bound them together and apply this bound in the analysis.

Next we define a more refined condition on the stochastic first-order oracle that accounts for the bias and the variance separately. 
\begin{customoracle}{S-1.2}[Stochastic First-Order Oracle, Case 2] \label{def:oracle_stoch_1.2}
    Given point $x\in \mathbb{R}^n$ and an input $\lambda$, the oracle generates random output $g(x,\xi^{(1)})$ that satisfies
    \begin{equation} \label{eq:oracle_stoch_1.2_bias}
        \|\E_{\xi^{(1)}} [g(x,\xi^{(1)})] - \nabla \phi(x)\|^2 \leq \epsilon_g^2(\lambda)+\kappa_g^2(\lambda)\|\nabla \phi(x)\|^2, 
    \end{equation} 
    \begin{equation} \label{eq:oracle_stoch_1.2_var}
        \E_{\xi^{(1)}} [\| g(x,\xi^{(1)}) - E_{\xi^{(1)}} [g(x,\xi^{(1)})]\|^2] \leq \sigma_c^2(\lambda)+\sigma_v^2(\lambda)\|\nabla \phi(x)\|^2, 
    \end{equation} 
    where  $\epsilon_g^2(\lambda)$, $\kappa_g^2(\lambda)$, $\sigma_{c}^2(\lambda)$ and $\sigma_{v}^2(\lambda)$ are some nonnegative  real-valued functions of $\lambda$ that generally depend on the global properties of $\phi$ and the distribution of $\xi^{(1)}$. 
\end{customoracle}

We observe that under Oracle~\ref{def:oracle_stoch_1.2} we recover Oracle~\ref{def:oracle_stoch_1.1} since
\begin{align*}
    \E [ \|g(x,\xi^{(1)}) - \nabla \phi(x)\|^2]
    &= \E[ \|g(x,\xi^{(1)}\|^2] - 2\E[ g(x,\xi^{(1)})]^T\nabla \phi(x)+\|\nabla \phi(x)\|^2\\
    &= \|\E [g(x,\xi^{(1)})] - \nabla \phi(x)\|^2 + \E [\| g(x,\xi^{(1)}) - \E [g(x,\xi^{(1)})]\|^2] \\
    &\le \epsilon_g^2(\lambda)+\sigma_{c}^2(\lambda)+(\kappa_g^2(\lambda)+\sigma_{v}^2(\lambda))\|\nabla \phi(x)\|^2.
\end{align*}

Finally, we propose a third version which combines Case 1 and Case 2. Specifically, we consider the case when the random variable $\xi^{(1)}$ is itself a pair of independent random variables: $\xi^{(1)}=(\zeta_1, \zeta_2)$. 

\begin{customoracle}{S-1.3}[Stochastic First-Order Oracle, Case 3] \label{def:oracle_stoch_1.3}
    Given point $x\in \mathbb{R}^n$ and an input $\lambda$, the oracle generates random output $g(x,\zeta_1, \zeta_2)$ that satisfies
    \begin{equation} \label{eq:oracle_stoch_1.3_bias}
        \E_{\zeta_1}   \|\E_{\zeta_2} [g(x,\zeta_1, \zeta_2)] - \nabla \phi(x)\|^2 \leq \epsilon_g^2(\lambda)+\kappa_g^2(\lambda)\|\nabla \phi(x)\|^2, 
    \end{equation} 
    \begin{equation} \label{eq:oracle_stoch_1.3_var}
        \E_{\zeta_1,\zeta_2} \left[\| g(x,\zeta_1, \zeta_2) - \E_{\zeta_2} [g(x,\zeta_1, \zeta_2)]\|^2 \right] \leq \sigma_c^2(\lambda)+\sigma_v^2(\lambda)\|\nabla \phi(x)\|^2, 
    \end{equation} 
    where  $\epsilon_g^2(\lambda)$, $\kappa_g^2(\lambda)$, $\sigma_{c}^2(\lambda)$ and $\sigma_{v}^2(\lambda)$ are some nonnegative real-valued functions of $\lambda$ that generally depend on the global properties of $\phi$ and the distributions of $(\zeta_1, \zeta_2)$. 
\end{customoracle}

\subsection{Examples of stochastic first-order oracles} \label{sec:stoch_first_order_ex}
We now present several examples of the stochastic first-order oracles covering cases of Oracles \ref{def:oracle_stoch_1.1}, \ref{def:oracle_stoch_1.2} and \ref{def:oracle_stoch_1.3}.

\subsubsection{Minibatch stochastic gradient in EM and ERM}\label{sec:erm}
For expected/empirical  risk minimization problems \eqref{eq:ERM}, the first-order oracles can also be implemented by applying sample averaging over a set of i.i.d. samples. As before, given a positive integer $\lambda$ as input, the oracle generates $\xi_1,\ldots,\xi_\lambda$ as i.i.d.\ samples from $\Xi$. We then define
\[
    \xi^{(1)}:=(\xi_1,\ldots,\xi_\lambda),
\]
and set
\begin{equation} \label{eq:erm_oracle1}
    g(x, \xi^{(1)}) = \frac{1}{\lambda}\sum_{j=1}^\lambda \nabla_x \ell(x,\xi_j).
\end{equation}
Under the common assumption that 
$
\nabla\phi(x) = \E_{\xi \sim \Xi}[\nabla_x \ell(x,\xi)], 
$
the output of the oracle $g(x,\xi^{(1)})$ is an unbiased estimator of $\nabla\phi(x)$. 
Meanwhile, it is also common to assume (e.g., see \cite{bottou2018optimization}) 
\begin{equation}\label{BCN}
	\E_{\xi\sim\Xi}\|\nabla_x \ell(x,\xi)-\nabla\phi(x)\|^2\leq M_c^2 + M_v^2\|\nabla\phi(x)\|^2. 
\end{equation} 
Then from the basic properties of sample averages, we have 
\begin{equation}\label{BCN2}
	\E_{\xi^{(1)}}\|g(x,\xi^{(1)})-\nabla\phi(x)\|^2\leq \frac{M_c^2}{\lambda}+\frac{M_v^2}{\lambda}\|\nabla\phi(x)\|^2, 
\end{equation} 
giving an example of Oracle~\ref{def:oracle_stoch_1.2} with $\lambda$ as the input, $\epsilon_g = 0$, $\kappa_g = 0$, $\sigma_c^2 = M_c^2 / \lambda$, and $\sigma_v^2 = M_v^2 / \lambda$. 
The cost of this oracle is also linearly dependent on the size of the minibatch $\lambda$ and can be chosen to depend on $x$ or on the current state of the algorithm. 

\subsubsection{Stochastic Gradient Oracle for Least Squares Loss}\label{strong_growth_reg}
Let us consider the specific case of \eqref{eq:ERM} for linear least squares, as discussed  in Section~\ref{sec:leastsquares} and derive the particular properties of the first-order oracle for this case. 

Define the sample gradient 
\[
    g(x,\xi^{(1)} ) = \frac{1}{\lambda} \sum_{j=1}^\lambda \nabla\ell(x,\xi_j) 
\]
as our stochastic oracle, where $\xi^{(1)} = (\xi_1, ...,\xi_\lambda)$ are i.i.d. samples from $\Xi$, making $g(x,\xi^{(1)})$ an unbiased estimator of $\nabla\phi(x)$. 
Following the derivation in Section~\ref{sec:leastsquares}, we have 
\begin{align*}
    g(x,\xi^{(1)}) &= \nabla f(x^\star,\xi^{(0)}) + H_{\xi^{(0)}} (x-x^\star) \\ 
    \text{ and }
    \nabla\phi(x) &= \nabla\phi(x^\star) + H(x-x^\star) = H(x-x^\star), 
\end{align*}
and, by the i.i.d. property of the components of $\xi^{(1)}$, $\|u+v\|^2 \le 2\|u\|^2 + 2\|v\|^2$, and \eqref{eq:least_square M1 implication}, we have 
\[ \begin{aligned}
    &\E_{\xi^{(1)}} [\|g(x,\xi^{(1)}) - \nabla\phi(x)\|^2] 
    = \E_{\xi^{(1)}} [\|\nabla f(x^\star,\xi^{(1)}) + (H_{\xi^{(1)}} - H) (x-x^\star)\|^2] \\
    &\le 2\E_{\xi^{(1)}} [\|\nabla f(x^\star,\xi^{(1)}\|] + 2\E_{\xi^{(1)}} [\|(H_{\xi^{(1)}} - H) (x-x^\star)\|^2] \\
    &\le \frac{2}{\lambda} \E_{\xi\sim\Xi}\left[\|\nabla \ell(x^\star,\xi)\|^2\right] + \frac{2M_1}{\lambda} \|\nabla\phi(x)\|^2. 
\end{aligned} \]
Thus, this is an example of Oracle~\ref{def:oracle_stoch_1.2} with 
\[
    \epsilon_g = 0, \quad 
    \kappa_g = 0, \quad
    \sigma_c^2 = \frac{2}{\lambda} \E_{\xi\sim\Xi}\left[\|\nabla \ell(x^\star,\xi)\|^2\right], \quad 
    \text{ and } \sigma_v^2 = \frac{2M_1}{\lambda}. 
\]

We note that $\sigma_c^2=0$ under the interpolation condition \ref{eq:interp}. In this case the variance 
of the stochastic oracle is bounded only by $\sigma_v^2\|\nabla \phi(x)\|^2$. This is known as the {\em strong growth condition} \cite{schmidt2013fast, vaswani2019fast} and 
implies faster convergence rates of stochastic gradient descent. We will see exactly how this faster rate is obtained in the case of simple fixed step gradient descent when we analyze this method under Oracle~\ref{def:oracle_stoch_1.2}. 

\subsubsection{Other stochastic gradients satisfying strong growth condition}\label{strong_growth_class}
Another notable example of the strong growth condition, aside from the linear least squares,   can be found in \cite{vaswani2019fast}. 
The general conditions in the context of \eqref{eq:ERM} is defined as \eqref{BCN} with $M_c = 0$ and $M_v > 0$: 
\begin{equation} \label{eq:SGC}
    \E_{\xi\sim\Xi} \|\nabla_x \ell(x,\xi) - \nabla\phi(x)\|^2 \le M_v^2 \|\nabla\phi(x)\|^2 
\end{equation}
and is always assumed together with $\E_{\xi\sim\Xi} \nabla_x \ell(x,\xi) = \nabla\phi(x)$. 
In other words, it assumes the stochastic gradient $\nabla_x \ell(x,\xi)$ to be an instance of Oracle~\ref{def:oracle_stoch_1.2} with $\epsilon_g = \kappa_g = \sigma_c^2 = 0$ and $\sigma_v^2 = M_v^2$. 

Condition \eqref{eq:SGC} bounds the variance and the second moment of the stochastic gradient $\nabla l(x,\xi)$ with the squared norm of the true gradient. 
It requires any stationary point of $\phi$ to be also stationary for every individual loss function, in other words $\nabla \ell(x,\xi) = 0$ for all $\xi$ when the model parameters $x$ are optimal. 
Though this may appear to be strong, it can hold true when the model is able to interpolate or fit the labelled training data completely. 
A concrete example was provided in \cite{vaswani2019fast}, and a slightly more generalized version is presented below using our notation. 

\begin{proposition}
  Suppose the random variable $\xi = (y,z)$ represents a label-feature pair that follows a distribution $\Xi$ with $c$ equally likely possible outcomes, where the label $y \in \{-1,1\}$, and the feature vector $z$ is normalized so that all $\|z\| \le {M}$ for some ${M} > 0$. 
  Let $\ell$ be the squared-hinge loss $\ell(x,\xi) = \ell(x,(y,z)) = (1-y \cdot z^T x)_+^2$, where $(\cdot)_+ = \max (0, \cdot)^2$. 
  If the data is linearly separable with margin ${r} \le {M}$, i.e., there exists a classifier $x^\star$ such that $\|x^\star\| = 1$ and $y \cdot d^T x^\star \ge {r}$ for all outcomes, then the objective function $\phi(x) = \mathbb{E}_{\xi \sim \Xi} \ell(x,\xi)$ satisfies \eqref{eq:SGC} with $M_v = c {M}^2 / {r}^2 - 1$. 
\end{proposition}

\begin{proof}
    Let $a = y \cdot z$. We have 
    \[ \begin{aligned}
        \|\nabla\phi(x)\|^2 &= \|\mathbb{E}[-2 (1 - a^T x)_+ \cdot a]\|^2
        \ge 4\big(\mathbb{E}[(1 - a^T x)_+ \cdot a]^T x^\star \big)^2 \\
        &\ge 4 {r}^2 \big( \mathbb{E}[(1 - a^T x)_+] \big)^2
        \ge \frac{4{r}^2}{c} \mathbb{E}[(1 - a^T x)_+^2], 
    \end{aligned} \]
    where the last step is true because, for any random variable $A$ taking nonnegative values $a_1,\dots,a_c$ with equal possibilities, $(\mathbb{E} A)^2 = \sum_{i=1}^c \sum_{j=1}^c \frac{1}{c} \frac{1}{c} a_i a_j \ge \frac{1}{c}\sum_{i=1}^c \frac{1}{c} a_i^2 = \frac{1}{c} \mathbb{E} A^2$. 
    Meanwhile, since $\|a\| = \|z\| \le {M}$ for all outcomes, 
    \[ \begin{aligned}
        \E_{\xi \sim\Xi} \|\nabla_x \ell(x,\xi)\|^2 = \E \|-2 (1 - a^T x)_+ \cdot a\|^2 \le 4{M}^2 \mathbb{E}[(1 - a^T x)_+^2]. 
    \end{aligned}\] 
    Combine the two inequalities to obtain 
    \[  \E_{\xi \sim\Xi} \|\nabla_x \ell(x,\xi)\|^2 
        \le \frac{c{M}^2}{{r}^2} \|\nabla\phi(x)\|^2. 
    \]
    Finally, since $\nabla\phi(x) = \E_{\xi \sim \Xi} \nabla_x \ell(x,\xi)$, we have $\E_{\xi \sim\Xi} \|\nabla_x \ell(x,\xi) - \nabla\phi(x)\|^2 = \E_{\xi \sim\Xi} \|\nabla_x \ell(x,\xi)\|^2  - \|\nabla\phi(x)\|^2$. 
\end{proof}

This proposition gives a concrete example of the averaged gradient \eqref{eq:erm_oracle1} serving as  Oracle~\ref{def:oracle_stoch_1.2} with $\lambda$ as input, $\epsilon_g = 0$, $\kappa_g = 0$, $\sigma_c^2 = 0$, and $\sigma_v^2 = (c{M}^2/r^2-1) / \lambda$.

\subsubsection{Robust Stochastic Gradient in Byzantine Machine Setting} \label{sec:Byzantine}
There is a large literature that studies stochastic gradient (and other algorithms) under the assumption that some of the data samples coming from distribution $\Xi$ are corrupted. There is a rich variety of settings depending on the assumed  nature of corruption. 
Here, as a particular example of our stochastic oracle,  we discuss a setting known as {\em Byzantine machine} corruption and an
 established approach in this setting (e.g. see \cite{chen2017distributed,Yin2018byzantine}) of computing gradients via robust mean estimation. While prior works do not 
 analyze the resulting gradient estimate in the way we do it here, we show that properties of Oracle~\ref{def:oracle_stoch_1.1} can be established and thus gradient descent can be analyzed in a fairly straightforward manner.

We again consider expected risk minimization problems \eqref{eq:ERM}, 
where we assume that the data is stored on $m$ machines, so that the algorithm can sample from each machine's distribution.
Let us assume now that some machines are corrupted (``Byzantine'') so that the data on these machines can come from an arbitrary (possibly adversarially chosen) distribution. We also assume that the data on all  uncorrupted machines  obeys distribution $\mathcal{D}$. Which machines are corrupted is unknown. 
  
Following the robust mean  estimation approach used in  \cite{chen2017distributed,Yin2018byzantine}  we first compute $m$ gradient  estimators by sample averaging using data from each individual machine,  
\begin{equation}\label{eq:m_machine_grad}
   \tilde g(x,{\xi}_i) = \frac{1}{\lambda} \sum_{j=1}^{\lambda} \nabla\ell(x,\xi_{ij}) 
\end{equation}
where each $\xi_i=(\xi_{i1}, ...,\xi_{i\lambda})$ is a set of $\lambda$ i.i.d samples following $\Xi$ from the $i$-th machine.  The sizes of the sample sets do not have to be the same for all machines, but under the assumption that all uncorrupted machines have the same distribution of the data, this is the most common approach. If the variance of the data differs on different machines, it could be useful to select different minibatch sizes so that the variances of all uncorrupted gradient estimators are as close to each other as possible. 

Using the $m$ ``weak'' estimators $\tilde g(x, \xi_i)$, $i=1,\ldots, m$, we can compute the final ``strong" estimator via the geometric median approach. Specifically, $y^*$ is an approximate geometric median of $\tilde g(x, \xi_i)$, $i=1,\ldots, m$, if it satisfies
\begin{equation}\label{eq:geom_med_gamma}
    \sum_{i=1}^m\|y^*-\tilde g(x, \xi_i)\|\leq (1+\gamma)\min_y \sum_{i=1}^m\|y-\tilde g(x, \xi_i)\|. 
\end{equation}

We now have the following stochastic first order oracle: 
\begin{itemize}
    \item Input: $x$, $\lambda$, $\gamma$. 
    \item For each machine $i=1, \ldots, m$, sample $\lambda$ data points uniformly (with replacement) and compute $\tilde g(x, {\xi}_i)$ using \eqref{eq:m_machine_grad}.
    \item Compute $y^*$ to satisfy  \eqref{eq:geom_med_gamma} and set $g_{\rm med}(x,\xi) = y^*$, where $\xi$ denotes the set of random variables $(\xi_1, \xi_2, \ldots, \xi_m)$. 
\end{itemize}

To analyze the behavior of the output of this oracle, let us state the following deterministic results. 
\begin{theorem}{\bf Exact case \cite{Minsker_2015}}\label{thm:geom_med_bound}
    Let $g_{\rm med}(x, {\xi})$ be the exact geometric median, i.e. $\gamma=0$ in   \eqref{eq:geom_med_gamma}. For any $\theta \in (0,\frac{1}{2})$ and any $r \in \mathbb{R}$, if $\sum_{i=1}^{m} \mathbbm{1}\{\|g(x, {\xi_i})-\nabla \phi(x)\|\leq r\}\geq (1-\theta)m$, then 
    \[ 
        \|g_{\rm med}(x, {\xi})-\nabla \phi(x)\|\leq  \frac{1-\theta}{\sqrt{1-2\theta}} r. 
    \]
\end{theorem}

\begin{theorem}{\bf Inexact case \cite{chen2017distributed}}
    Let $g_{med}(x, {\xi})$ be a $\gamma$-approximate geometric median, satisfying \eqref{eq:geom_med_gamma}, then
    for any $\theta \in (0,\frac{1}{2})$ and any $r \in \mathbb{R}$, if $\sum_{i=1}^{m} \mathbbm{1}\{\|g(x, {\xi_i})-\nabla \phi(x)\|\leq r\}\geq (1-\theta)m$, then
    \[ 
        \|g_{\rm med}(x, {\xi})-\nabla \phi(x)\|\leq  \frac{2(1-\theta)}{1-2\theta} r +\gamma \frac{\max_{i}\|g(x, {\xi}_i)-\nabla \phi(x)\|}{ 1-2\theta}. 
    \]
\end{theorem}

Let us consider the case when $\gamma=0$ (or is negligibly small). 
Let ${\cal M}$ be the set of non-corrupted machines. As long as $|{\cal M}|\geq (1-\theta)m$ with $\theta <\frac{1}{2}$, then 
\[
    \sum_{i=1}^{m} \mathbbm{1}\left\{ \|g(x, {\xi_i})-\nabla \phi(x)\|\leq  \max_{i\in  {\cal M}}\{\|g(x, \xi_i )-\nabla \phi(x)\|\} \right\}\geq (1-\theta)m. 
\] 
By setting $r(\xi)=  \max_{i\in {\cal M}}\{\|g(x, \xi_i )-\nabla \phi(x)\|\}$, we have
    \[ 
        \|g_{\rm med}(x, {\xi})-\nabla \phi(x)\|\leq  \frac{1-\theta}{\sqrt{1-2\theta}} r(\xi) 
    \]
holds for all realization of $\xi$.   
Thus, we can conclude that 
\[ \begin{aligned}
    \E_{\xi} [\|g_{\rm med}(x, \xi )-\nabla \phi(x)\|^2]
    &\le \E_{\xi} \Big[\Big(\frac{1-\theta}{\sqrt{1-2\theta}} r(\xi) \Big)^2 \Big] \\ 
    &= \frac{(1-\theta)^2}{1-2\theta} \E_{\xi, i\in{\cal M}} \Big[ \max_{i\in  {\cal M}}\{\|g(x, \xi_i )-\nabla \phi(x)\|^2\} \Big]. 
\end{aligned} \]
 To bound the right-hand side, one needs assumptions on the distributions of 
 $g(x, \xi_i )$ for $i\in  {\cal M}$. Let us assume here that all $g(x,\xi_i)$,   $i\in  {\cal M}$, are independent and identically distributed and 
  satisfy \eqref{BCN2}, i.e. 
 \begin{equation*}
	\E_{\xi_i}[\|g(x,\xi_i)-\nabla\phi(x)\|^2]
    \leq \frac{M_c^2}{\lambda}+\frac{M_v^2}{\lambda}\|\nabla\phi(x)\|^2,  
    \quad \forall i \in {\cal M}.
\end{equation*} 
This assumption is valid if all uncorrupted data satisfy \eqref{BCN}. 
Without further assumptions on the distribution of $\cal M$ or the tails of the distribution of $\|g(x,\xi_i)-\nabla\phi(x)\|$, we can apply basic inequality 
\[
    \E_{\xi_i, i\in{\cal M}}\Big[ \max_{i\in  {\cal M}}\{\|g(x, \xi_i )-\nabla \phi(x)\|^2\} \Big]
    \leq \frac{\cal M}{\lambda} (M_c^2 + M_v^2 \|\nabla\phi(x)\|^2) 
\]
to obtain the following bound 
\[
    \E_{\xi}[\|g_{\rm med}(x, \xi )-\nabla \phi(x)\|^2] 
    \leq \frac{(1-\theta)^2 {\cal M}}{(1-2\theta)\ \lambda} \left( M_c^2+M_v^2\|\nabla\phi(x)\|^2\right). 
\]

In other words, $g_{\rm med}(x, \xi )$ can be viewed as Oracle~\ref{def:oracle_stoch_1.1} with 
\[
    \epsilon_g = \sqrt{\frac{(1-\theta)^2 {\cal M}}{(1-2\theta)\ \lambda}} M_c
    \text{ and } 
    \kappa_g =\sqrt{\frac{(1-\theta)^2 {\cal M}}{(1-2\theta)\ \lambda}} M_v. 
\] 
While this is a rather crude bound which essentially translates the variance of the uncorrupted gradients to the bias of the robust estimator, it is still possible to have a convergent gradient descent algorithm based on such estimator, provided that $\lambda$ is chosen large enough, which is the premise analyzed in \cite{chen2017distributed,Yin2018byzantine}).  

The case $\gamma>0$ will result in an extra bias term which will not shrink with $\lambda$, but can be made small by solving the geometric median problem more precisely. The cost of this calculation (which contributes to the cost of the oracle) depends on the method being used. For various approaches to geometric median problem and their complexity we refer the reader to \cite{BeckSabach2015} and \cite{CohenLeeMillerPachockiSidford2016}.

We note that there is at least one other  robust estimator which can be computed exactly and whose properties can be found for example in \cite{Hsu2016}. It satisfies a bound similar to that in Theorem \ref{thm:geom_med_bound} but 
with $\theta>\frac{1}{2}$ and the term $\frac{1-\theta}{\sqrt{1-2\theta}} $ replaced by $3$.

\subsubsection{Gradient Projection onto  Random Subspaces}\label{sec:gradproj}
Consider the vector
\begin{equation} \label{eq:subspace_grad}
	g(x,{\cal U}) =  \sum_{u_i\in\mathcal{U}} (\nabla \phi(x)^Tu_i)u_i, 
\end{equation}
where ${\cal U}=\{u_1, \ldots, u_s\} \subset \mathbb{R}^n$ is a set of $s$ random vectors following a distribution that satisfies Assumption~\ref{ass:haar}. 
Under this assumption, each vector $(\nabla \phi(x)^Tu_i)u_i$ in \eqref{eq:subspace_grad} is a projection of the gradient $\nabla\phi(x)$ onto the one-dimensional subspace defined by $u_i$. Since ${\cal U}$ is orthonormal, the sum $g(x,{\cal U})$ in \eqref{eq:subspace_grad} is the projection of $\nabla \phi(x)$ onto $\text{span}(\mathcal{U})$. 
\begin{assumption}\label{ass:haar}
    The set of vectors ${\cal U}$ is orthonormal ($s \le n$, $\|u_i\|=1$, and $u_i^T u_j=0$ when $i \neq j$) and $\E[\sum_{i=1}^s u_i u_i^T] = \frac{s}{n} I$. 
\end{assumption}

Two notable distributions that satisfy Assumption~\ref{ass:haar} are 
\begin{itemize}
    \item {\bf Block Coordinate} $\mathcal{C}$: ${\cal U}$ is a uniformly distributed random subset (without repetition) of $\{e_1, e_2, \ldots e_n\}$; 
    \item {\bf Haar} $\mathcal{H}$: The vectors $u_1, \ldots u_s$ are columns of matrix $S\in\mathbb{R}^{n\times s}$ distributed uniformly over $\{S\in\mathbb{R}^{n\times s}:~ S^TS = I_{s}\}$.  
\end{itemize}
When $\mathcal{U} \sim \mathcal{C}$, the gradient descent based on $g(x,{\cal U})$ is the well-known randomized block coordinate descent algorithm due to Nesterov \cite{nesterov2012efficiency}. It is very efficient in many large-scale applications, especially with certain sparsity structures, where  the partial derivatives $\nabla \phi(x)^Tu_i$ are much cheaper to compute than the entire $\nabla \phi(x)$.  
The case of ${\cal H}$ is a generalization of this method to an arbitrary subspace.
The form \eqref{eq:subspace_grad} is useful when computing directional derivatives is cheaper (in some sense) that computing the entire gradient. Such methods have been recently utilized in  Machine Learning to save memory use compared to the full backpropagation \cite{Baydin2022Gradients, Shukla2023Randomized}.  Moreover, this form will  prove very useful when combined with finite difference approximations, which we will discuss later. 

When $s=1$ and $\mathcal{U} = \{u\}$, any distribution with $\|u\|=1$ and ${\mathbb E} [uu^T] = \frac{1}{n} I_{n}$ would satisfy Assumption~\ref{ass:haar}. 
One interesting  example of this  is $u$ whose elements independently take  values $1/\sqrt{n}$ or $-1/\sqrt{n}$,  with probability $\frac{1}{2}$.  
This is a scaled version of the Rademacher vector used in conjunction with finite difference in SPSA \cite{spall1992multivariate} for gradient estimation. 

The projection $g(x,\mathcal{U})$ can be treated as a stochastic first-order oracle that takes the sample size $s$ and the distribution of ${\cal U}$ as the input $\lambda$. 
To see how accurately it approximates $\nabla\phi(x)$, we derive the constants in the conditions of Oracle~\ref{def:oracle_stoch_1.1} and \ref{def:oracle_stoch_1.2}. 
Let $S\in\mathbb{R}^{n\times s}$ be the matrix with columns $u_1, \ldots, u_s$. 
We observe that 
\begin{equation}\label{eq:S Haar}
    g(x,\mathcal{U})=SS^T \nabla\phi(x), \quad 
    S^TS = I_s, \quad \text{ and }  
    \E\left[SS^T\right] = \frac{s}{n} I_{n}. 
\end{equation}

\paragraph{As Oracle~\ref{def:oracle_stoch_1.1}}
Using \eqref{eq:S Haar}, we have 
\[ 
\begin{aligned} 
    &\E\left[\|g(x,\mathcal{U}) - \nabla\phi(x)\|^2\right] 
    = \E\left[\| (SS^T -I)\nabla\phi(x) \|^2\right]\\ 
    &= \nabla\phi(x)^T\E\left [ SS^TSS^T - 2SS^T + I\right ] \nabla\phi(x)
    =\frac{n-s}{n} \|\nabla\phi(x)\|^2. 
\end{aligned} \]
Thus, the conditions of Oracle~\ref{def:oracle_stoch_1.1} are satisfied with $\epsilon_g = 0$ and $\kappa_g^2 = \frac{n-s}{n}$. 

\paragraph{As Oracle~\ref{def:oracle_stoch_1.2}}
We observe that $\E [g(x,\mathcal{U})] = \frac{s}{n} \nabla\phi(x)$. This leads to the squared bias
\begin{equation*}
    \|\E[g(x,\mathcal{U})] - \nabla\phi(x)\|^2 
    = \left\|\frac{s}{n} \nabla\phi(x) - \nabla\phi(x)\right\|^2 
    = \frac{(n-s)^2}{n^2} \|\nabla\phi(x)\|^2. 
\end{equation*}
The variance is 
\[ \begin{aligned} 
    \E\left[\big\|g(x,\mathcal{U}) - \E[g(x,\mathcal{U})] \big\|^2\right]
    &= \E\left[\|g(x,\mathcal{U}) - \nabla\phi(x)\|^2\right] - \|\E[g(x,\mathcal{U})] - \nabla\phi(x)\|^2 \\
    &= \Big(\frac{n-s}{n} - \frac{(n-s)^2}{n^2}\Big) \|\nabla\phi(x)\|^2 
    = \frac{(n-s)s}{n^2} \|\nabla\phi(x)\|^2. 
\end{aligned} \]
Thus, the conditions of Oracle~\ref{def:oracle_stoch_1.2} are satisfied with  $\epsilon_g = 0$,  $\kappa_g^2 = \frac{(n-s)^2}{n^2}$, $\sigma_c = 0$, and $\sigma_v^2 = \frac{(n-s)s}{n^2}$. 

\paragraph{Scaled Gradient Projection onto  Random Subspaces}
The above derivation show that $g(x,{\cal U})$ computed by \eqref{eq:subspace_grad} is a biased gradient estimator. 
For this exact reason, a scaled version 
 \begin{equation} \label{eq:subspace_grad_scale}
	g(x,{\cal U}) = \frac{n}{s} \sum_{u_i\in\mathcal{U}} (\nabla \phi(x)^Tu_i)u_i, 
\end{equation}
which gives an unbiased stochastic gradient $\E [g(x,\mathcal{U})] = \nabla\phi(x)$, is often used instead. 
Its variance is scaled by $n^2/s^2$, giving 
\begin{equation} \label{eq:subspace_grad_S1.2}
    \E\left[\|g(x,\mathcal{U}) -  \E[g(x,\mathcal{U})]\|^2\right]
    = \frac{n-s}{s} \|\nabla\phi(x)\|^2.
\end{equation}

Thus, it satisfies the conditions of Oracle~\ref{def:oracle_stoch_1.1} with $\epsilon_g = 0$ and $\kappa_g^2=(n-s)/s$; and the conditions of Oracle~\ref{def:oracle_stoch_1.2} with $\epsilon_g = \kappa_g = \sigma_c = 0$, and  $\sigma_v^2= (n-s)/s$. 
We will see later that viewing this case of $g(x,\mathcal{U})$ as Oracle~\ref{def:oracle_stoch_1.1} is nonsensical for $s<\frac{n}{2}$, since $\kappa_g>1$ would cause analysis of appropriate algorithms to fail. 
Viewed as Oracle~\ref{def:oracle_stoch_1.2}, the scaled unbiased estimator \eqref{eq:subspace_grad_scale} will fit  the analysis and the resulting algorithm and its complexity  bound will be the same as for the unscaled case  viewed as Oracle~\ref{def:oracle_stoch_1.1}, as our complexity analysis will show in later sections. Thus scaling of the projected gradient (and removal of the bias) does not seem to have theoretical advantage over the unscaled biased estimator \eqref{eq:subspace_grad} due to the benign nature of the bias in the latter.

\subsubsection{Gradient Approximation via Random Directional Derivatives with Independent Directions}\label{sec:rand_dirder}
Let us again consider the vector in \eqref{eq:subspace_grad}
\begin{equation*} 
	g(x,{\cal U}) =  \sum_{u_i\in\mathcal{U}} (\nabla \phi(x)^Tu_i)u_i, 
\end{equation*}
but with ${\cal U}=\{u_1, \ldots, u_s\}$ being a set of random vectors that  satisfies the following assumption rather than Assumption~\ref{ass:haar}. 
\begin{assumption}\label{ass:normal}
    The vectors in ${\cal U}$ are independent random vectors, not necessarily identically distributed.
    The distribution followed by each $u\in\mathcal{U}$ satisfies the following properties.
    \begin{enumerate}
        \item Given any signed permutation matrix (a matrix having exactly one nonzero entry in each row and each column, with each nonzero entry being either $+1$ or $-1$) $Q\in\mathbb{R}^{n\times n}$, the random vector $Qu$ has the same distribution as $u$. 

        \item $\mathbb{E}[\|u\|^2] = 1$ and $\mathbb{E}[\|u\|^{2(k+1)}] \le 1 + {\cal O}(n^{-1})$ for every fixed $k\in\mathbb{N}_+$. 
    \end{enumerate}
\end{assumption}

Instead of being orthonormal, the vectors in $\mathcal{U}$ are now independent and follow distributions invariant to signed permutation. 
Examples of such distributions include:
\begin{itemize}
    \item {\bf Gaussian} $\mathcal{N}$: the multivariate normal distribution
    $\mathcal{N}(0,\frac{1}{n}I)$;

    \item {\bf Uniform on sphere} $\mathcal{S}$: the uniform distribution on the unit sphere
    $\{x\in\mathbb{R}^n: \|x\|=1\}$;

    \item {\bf Uniform over coordinates}: the uniform distribution over
    $\{\pm e_i\}_{i=1}^n$;

    \item {\bf Scaled Rademacher}: the elements in $u$ are i.i.d. random variables that are equal to either $-1/\sqrt{n}$ or $1/\sqrt{n}$, each with probability $1/2$.
\end{itemize}

For the latter three distributions, $\|u\| = 1$ holds deterministically, thus the ${\cal O}(n^{-1})$ term in Assumption~\ref{ass:normal} simply disappears. 
When $u \sim \mathcal{N}$, the vector length is random, and it holds that 
\[ 
    \mathbb{E}[\|u\|^{2(k+1)}] = \frac{(n+2)(n+4) \cdots (n+2k)}{n^k}
    \leq  1 + \frac{k(k+1)}{n} + \mathcal{O}(n^{-2}). 
\] 
Apart from being invariant to signed permutation, the distributions $\mathcal{N}$ and $\mathcal{S}$ are also invariant to orthogonal transformations, i.e., for any orthogonal matrix $Q \in \mathbb{R}^{n \times n}$, the random vector $Qu$ has the same distribution as $u$. 
This property is also referred to as rotational invariance.
It is stronger than signed-permutation invariance and is not satisfied by the other two distributions, which have finite supports. 
More properties that only hold under $\mathcal{N}$ and $\mathcal{S}$ but not under Assumption~\ref{ass:normal} in general will be discussed  in Section~\ref{sec:smoothing}. 

We first note that under Assumption~\ref{ass:normal} the vector $g(x,\mathcal{U})$ is no longer a projection of $\nabla\phi(x)$, but it can still be treated as a stochastic first-order oracle that requires  sample size $s$ and the distribution of ${\cal U}$ as the input $\lambda$. 
As in the case of random projection, this oracle is of interest when computing directional derivatives is significantly cheaper  than computing the entire gradient and it is specifically  very useful when  $\nabla \phi(x)^Tu_i$ is replaced by a finite difference approximation. 

We also note that Assumption~\ref{ass:normal} itself is very general, as it allows vectors in {\cal U} to be sampled from different distributions. Yet, we can still derive useful Oracle properties under this assumption.
We again analyze $g(x,\mathcal{U})$ as Oracle~\ref{def:oracle_stoch_1.1} and \ref{def:oracle_stoch_1.2}, starting with the following lemma regarding the essential properties held by each $u \in {\cal U}$ under Assumption~\ref{ass:normal}. 

\begin{lemma} \label{lem:normal}
    When $\mathcal{U}=\{u_1,\ldots,u_s\}$ satisfies Assumption~\ref{ass:normal}, it holds for every $u\in\mathcal{U}$ and every fixed $k\in\mathbb{N}_+$ that 
    \begin{equation*}
        \mathbb{E}[u] = 0, \quad
        \mathbb{E}[uu^T] = \frac{1}{n}I_n, \quad 
        \text{and } \mathbb{E}\left[\|u\|^k uu^T\right]
        \preceq \frac{1 + \mathcal{O}(n^{-1})}{n} I_n.
    \end{equation*}
\end{lemma}

\begin{proof}
Consider the signed permutation matrix $-I$. 
Since the distribution is invariant to signed permutation, $-u$ have the same distribution as $u$. 
Hence $\mathbb{E}[u] = \mathbb{E}[-u]$, implying $\mathbb{E}[u] = 0$. 

Next, for any $i\in\{1,\ldots,n\}$, let $Q$ be the signed permutation matrix that flips the sign of the $j$th coordinate and leaves all other coordinates unchanged. 
Then, all off-diagonal elements in row $i$ and column $i$ of $\mathbb{E}[uu^T]$ and $\mathbb{E}[(Qu)(Qu)^T] = Q\mathbb{E}[uu^T]Q^T$ are negatives of each other. 
However, by signed-permutation invariance, $\mathbb{E}[uu^T] = Q\mathbb{E}[uu^T]Q^T$. 
Hence all off-diagonal elements of $\mathbb{E}[uu^T]$ are zero. 

For any $i,j\in\{1,\ldots,n\}$ such that $i \neq j$, redefine $Q$ as the signed permutation matrix that swaps the $i$th coordinate with the $j$th coordinate and leaves all other coordinates unchanged. 
Then the $i$th diagonal element of $\mathbb{E}[uu^T]$ is equal to the $j$th diagonal element of $\mathbb{E}[(Qu)(Qu)^T]$. 
However, $\mathbb{E}[uu^T] = \mathbb{E}[(Qu)(Qu)^T]$, so all diagonal elements of $\mathbb{E}[uu^T]$ are equal. 

The above results establishe that $\mathbb{E}[uu^T]$ is a scalar multiple of the identity matrix $I_n$. 
Since its trace is 1 due to $\mathbb{E}[\|u\|^2]=1$, we determine that $\mathbb{E}[uu^T] = \frac{1}{n}I_n$.

Since $\|u\|$ is invariant under signed permutations, the same argument applies to $\|u\|^k uu^T$. 
Its expectation is also a scalar multiple of the identity matrix $I_n$ and is bounded as  
\begin{equation*}
    \mathbb{E}\left[\|u\|^k uu^T\right]
    = \frac{\mathbb{E}[\|u\|^{k+2}]}{n}I_n
    \preceq \frac{1 + \mathcal{O}(n^{-1})}{n} I_n.
\end{equation*}
\end{proof}

\paragraph{As Oracle~\ref{def:oracle_stoch_1.1}}
By the independence of $\mathcal{U}$ and Lemma~\ref{lem:normal}, 
\begin{equation} \label{eq:e_proj} \begin{aligned}
    &\E\left[\big\|g(x,\mathcal{U}) - \nabla\phi(x)\big\|^2\right] 
    = \nabla\phi(x)^T \E\bigg[\Big(\sum_{i=1}^s u_i u_i^T - I\Big)^2\bigg] \nabla\phi(x) \\
    &= \nabla\phi(x)^T \E\bigg[\Big(\sum_{i=1}^s \|u_i\|^2 u_i u_i^T\Big) + \Big(\sum_{i=1}^s \sum_{j=1, j \neq i}^s\!\! u_i u_i^T u_j u_j^T\Big) - 2\Big(\sum_{i=1}^s u_i u_i^T\Big) + I\bigg] \nabla\phi(x) \\
    &= \bigg(s \cdot \frac{1 + \mathcal{O}(n^{-1})}{n} + (s^2-s) \cdot \frac{1}{n^2} - 2s \cdot \frac{1}{n} + 1\bigg) \|\nabla\phi(x)\|^2 \\
    &= \bigg(\Big(\frac{n-s}{n}\Big)^2 + \frac{(1+\mathcal{O}(n^{-1}))s}{n}\bigg) \| \nabla\phi(x)\|^2. 
\end{aligned} \end{equation} 
Thus the conditions of Oracle~\ref{def:oracle_stoch_1.1} are satisfied with $\epsilon_g=0$ and 
\[
    \kappa_g^2 = \Big(\frac{n-s}{n}\Big)^2 + \frac{(1+\mathcal{O}(n^{-1}))s}{n}. 
\]

\paragraph{As Oracle~\ref{def:oracle_stoch_1.2}}
The expectation of $g(x,\mathcal{U})$ is 
\[
    \E[g(x,\mathcal{U})] 
    = \E\left[\sum_{i=1}^s [u_i u_i^T] \nabla\phi(x)\right]
    = \frac{s}{n} \nabla\phi(x). 
\]
Then the bias is 
 \begin{equation*}
    \|\E[g(x,\mathcal{U})] - \nabla\phi(x)\|^2 
    = \left\|\frac{s}{n} \nabla\phi(x) - \nabla\phi(x)\right\|^2 
    = \Big(\frac{n-s}{n}\Big)^2 \|\nabla\phi(x)\|^2. 
\end{equation*}
The variance is 
\[ \begin{aligned} 
    \E\left[\big\|g(x,\mathcal{U}) - \E[g(x,\mathcal{U})] \big\|^2\right]
    &= \E\left[\|g(x,\mathcal{U}) - \nabla\phi(x)\|^2\right] - \|\E[g(x,\mathcal{U})] - \nabla\phi(x)\|^2 \\ 
    &= \frac{(1+\mathcal{O}(n^{-1}))s}{n} \|\nabla\phi(x)\|^2
\end{aligned} \]
Thus the conditions of Oracle~\ref{def:oracle_stoch_1.2} are satisfied with $\epsilon_g = 0, \kappa_g^2 = \left(\frac{n - s}{n}\right)^2, \sigma_c^2 = 0$, and 
\[
    \sigma_v^2 = \frac{(1 + \mathcal{O}(n^{-1})) s}{n}. 
\]

\paragraph{Scaled/Unbiased Gradient Approximation via Random Directional Derivatives with Independent Directions}
Similar to the example of gradient projection onto random subspaces, the derivation above shows that $g(x,{\cal U})$ computed by \eqref{eq:subspace_grad} is also a biased gradient estimator. 
It can also be scaled by $\frac{n}{s}$ to \eqref{eq:subspace_grad_scale}, so that $\E [g(x,\mathcal{U})] = \nabla\phi(x)$. 
The scaling also changes the variance by a factor of $\frac{n^2}{s^2}$. 
Thus, \eqref{eq:subspace_grad_scale} with the distribution of $\mathcal{U}$ satisfying Assumption~\ref{ass:normal} gives an example of Oracle~\ref{def:oracle_stoch_1.2} with the same input as before, $\epsilon_g =\kappa_g=\sigma_c= 0$, and 
\[
    \sigma_v^2 =  \frac{(1 + \mathcal{O}(n^{-1})) n}{s}. 
\]

\subsubsection{Randomized Forward Finite Differences with Orthonormal Directions}\label{sec:RFFD_orth}
We now turn to the finite difference version of \eqref{eq:subspace_grad}, which replaces each directional derivative $(\nabla\phi(x)^T u_i) u_i$ with its finite difference approximation. 
Given a zeroth-order oracle $f$, we define 
\begin{equation} \label{eq:RFFD} \tag{RFFD}
	g(x,\mathcal{U},\xi) = \sum_{i=1}^s \frac{f(x+ h u_i, \xi_i) - f(x,\xi_0)}{h} u_i, 
\end{equation}
where ${\cal U} = \{u_1, \dots, u_s\} \subseteq \mathbb{R}^n$ is a set of random vectors satisfying Assumption~\ref{ass:haar}. 
Here, we assume $f$ is either the single-point stochastic zeroth-order oracle \ref{def:oracle_stoch_0.1} or its multi-point variant \ref{def:oracle_stoch_0.2}. 
With Oracle~\ref{def:oracle_stoch_0.1}, $\xi=\{\xi_0, \xi_1, \dots, \xi_s\}$ is the set of independent random variables invoked by the zeroth-order oracle when estimating the function values at $x$ and $\{x+ h u_i\}_{i=1}^s$. 
With Oracle~\ref{def:oracle_stoch_0.2}, $\xi$ is instead the  common random number for the multi-point oracle, i.e., $\xi_i = \xi$ for all $i=0,1,\dots,s$. 
In the upcoming analysis, we assume ${\cal U}$ and $\xi$ are independent from each other. 

Since ${\cal U}$ is orthonormal under Assumption~\ref{ass:haar}, formula \eqref{eq:RFFD} is estimating the projection of $\nabla\phi(x)$ onto the random subspace $\text{span}({\cal U})$ using $s+1$ (approximate) function values. 
Estimating random gradient projections is a practical and popular tool for reducing computational cost of finite difference gradient estimation. 
We have discussed applying finite difference in a subspace in Sections \ref{sec:det_project} and \ref{sec:det_project_ffd}. Ensuring bounds of Oracle \ref{def:oracle_det_1} with any fixed $\kappa_g<1$  essentially requires that the projection of the  gradient onto the subspace is not too small --- something that is not always easy to guarantee if a subspace is chosen deterministically but, as we have shown in Section~\ref{sec:gradproj}, holds in expectation for random subspaces. 

The vector $g(x,\mathcal{U},\xi)$ defined in \eqref{eq:RFFD} is again a stochastic first-order oracle that requires the sample size $s$ and the distribution of ${\cal U}$ as the input $\lambda$. 
We analyze properties of  $g(x,\mathcal{U},\xi)$  in terms of Oracle~\ref{def:oracle_stoch_1.3}. 
We start by decomposing $g(x,\mathcal{U},\xi)$ as 
\begin{equation} \label{eq:RFF decomposition}
   g(x,\mathcal{U},\xi) = \nabla\phi(x) + \epsilon_{\rm proj}(x,\mathcal{U}) + \epsilon_{\rm ffd}(x,\mathcal{U}) + \epsilon_{\rm noise}(x,\mathcal{U},\xi), 
\end{equation}
where, 
\begin{align*}
	\epsilon_{\rm proj}(x,\mathcal{U}) &= \sum_{i=1}^s (\nabla\phi(x)^T u_i) u_i  - \nabla\phi(x)  \\
	\epsilon_{\rm ffd}(x,\mathcal{U}) &= \sum_{i=1}^s \bigg[\frac{\phi(x+ {h} u_i) - \phi(x) - \nabla\phi(x)^T ({h} u_i)}{h} + \\ 
    &\qquad \frac{\E_{\xi }[f(x+ {h} u_i, \xi_i)]-\phi(x+{h}u_i) - (\E_{\xi }[f(x, \xi_0)]-\phi(x))}{ {h}} \bigg] u_i  \\
    \text{and } \epsilon_{\rm noise}(x,\mathcal{U},\xi) &= \sum_{i=1}^s \frac{\big(f(x+ {h} u_i, \xi_i)-\E_{\xi}[f(x+ {h} u_i, \xi_i)]\big) - \big(f(x, \xi_0)-\E_{\xi}[f(x,\xi_0)]\big)}{h} u_i. 
\end{align*} 

The first term $\epsilon_{\rm proj}$ is the difference between the gradient and its projection onto ${\rm span}({\cal U})$, which has been analyzed in Section \ref{sec:gradproj}. Moreover, we recall that this vector is orthogonal to ${\rm span}({\cal U})$. The term $\epsilon_{\rm ffd}$ is the ``deterministic'' error caused by applying finite difference to estimate the projected gradient. It is not really a deterministic quantity, since it depends on ${\cal U}$, but it does not depend on $\xi$ and can be deterministically bounded under Assumptions~\ref{assum:lip_cont}. 
Finally, the $\epsilon_{\rm noise}$ term is the stochastic error induced by using stochastic zeroth-order oracle in the finite difference scheme. 
Both $\epsilon_{\rm ffd}$ and $\epsilon_{\rm noise}$ reside in ${\rm span}({\cal U})$. 
We now derive bounds on these three error terms. 

\paragraph{Bounding $\epsilon_{\rm proj} $.}
Under Assumption~\ref{ass:haar}, it has been proved previously when developing conditions for Oracle~\ref{def:oracle_stoch_1.1} in Section~\ref{sec:gradproj} that 
\begin{equation} \label{eq:RFFD_proj_error}
    \E_{\cal U}\left[\|\epsilon_{\rm proj}(x,\mathcal{U})\|^2\right] = \frac{n-s}{n}\|\nabla \phi(x)\|^2, 
\end{equation}

\paragraph{Bounding $\epsilon_{\rm ffd} $.}
Let $\epsilon_1(x,{h} u_i)=\phi(x+ {h} u_i) - \phi(x) - \nabla\phi(x)^T ({h} u_i)$ be the error caused by the nonlinearity of $\phi$, which follows $|\epsilon_1(x,{h} u_i)| \le L_1 {h}^2 / 2$ under Assumption~\ref{assum:lip_cont}.
Let $\epsilon_2(x,{h} u_i) = \E_{\xi }[f(x+ {h} u_i, \xi_i)]-\phi(x+{h}u_i) - (\E_{\xi }[f(x, \xi_0)]-\phi(x))$ be the bias of the zeroth-order oracle, which is bounded as $|\epsilon_2(x,{h} u_i)| \le 2\epsilon_f$ using either Oracles~\ref{def:oracle_stoch_0.1} or~\ref{def:oracle_stoch_0.2} as the zeroth-order oracle. 
Since $\mathcal{U}$ is orthonormal, we have 
\begin{equation} \label{eq:ortho eps_ffd bound} \begin{aligned}
    \E_{\cal U}\left[\|\epsilon_{\rm ffd}(x,\mathcal{U})\|^2\right]
    &= \E_{\cal U}\left[\bigg\|\sum_{i=1}^s\Big(\frac{\epsilon_1(x,{h} u_i)}{h}+ \frac{\epsilon_2(x,{h} u_i)}{h}\Big) u_i\bigg\|^2\right] \\
    &= \sum_{i=1}^s\E_{\cal U}\left[\Big(\frac{\epsilon_1(x,{h} u_i)}{h}+ \frac{\epsilon_2(x,{h} u_i)}{h} \Big)^2\right] 
    \le s \Big(\frac{L_1 {h}}{2} + \frac{2\epsilon_f}{h}\Big)^2. 
\end{aligned} \end{equation}

\paragraph{Bounding $\epsilon_{\rm noise}$}
Let us denote the function noise term as $\epsilon_3(x,\xi)=f(x,\xi)-\E_{\xi}[f(x,\xi)]$. Since $\mathcal{U}$ is orthonormal, we obtain 
\begin{equation} \begin{aligned}
    \E_{{\cal U},\xi} \left[ \|\epsilon_{\rm noise}(x,\mathcal{U},\xi)\|^2 \right] 
    &= \E_{{\cal U},\xi} \left[\bigg\|\sum_{i=1}^s \Big(\frac{\epsilon_3(x+ {h} u_i, \xi_i) - \epsilon_3(x, \xi_0)}{h} \Big) u_i \bigg\|^2\right] \\
    &= \frac{1}{{h}^2} \sum_{i=1}^s \E_{{\cal U},\xi} \left[  \left(\epsilon_3(x+ {h} u_i, \xi_i) - \epsilon_3(x, \xi_0) \right)^2 \right].  
\end{aligned} \end{equation}
Under Oracle \ref{def:oracle_stoch_0.1}, since  $\xi_i$, $i=0, 1, \ldots, s$, are independent, then each term in the summation is deterministically bounded (regardless of $u_i$) by $2\sigma_f^2$, which results in 
\[
    \E_{{\cal U},\xi} \left[ \|\epsilon_{\rm noise}(x,\mathcal{U},\xi)\|^2 \right] 
    \leq \frac{2s\sigma_f^2}{{h}^2}. 
\]
Under Oracle~\ref{def:oracle_stoch_0.2}, we have $\xi_0 = \xi_1 = \cdots = \xi_s$, and each term  of the sum is bounded by 
\[
  (L_c^2+L_v^2\|\nabla \phi(x)\|^2){h}^2\|u_i\|^2, 
\]
which gives (recalling that $\|u_i\|=1$)
\[  \E_{{\cal U},\xi} \left[ \|\epsilon_{\rm noise}(x,\mathcal{U},\xi)\|^2 \right] 
    \leq s(L_c^2+L_v^2\|\nabla\phi(x)\|^2). 
\]

\paragraph{As Oracle~\ref{def:oracle_stoch_1.3}} 
Since $\E_{\xi}\left[\epsilon_{\rm noise}(x,\mathcal{U},\xi)\right] = 0$ and that $\epsilon_{\rm proj}(x,\mathcal{U})$ is perpendicular to $\epsilon_{\rm ffd}(x,\mathcal{U})$, we have 
\[ \begin{aligned}
    \E_{\cal U}\left[\|\E_{\xi}[g(x,\mathcal{U},\xi)] - \nabla\phi(x)\|^2\right] 
    &= \E_{\cal U}\left[\|\epsilon_{\rm proj}(x,\mathcal{U}) + \epsilon_{\rm ffd}(x,\mathcal{U})\|^2\right] \\ 
    &= \E_{\cal U}\left[\|\epsilon_{\rm proj}(x,\mathcal{U})\|^2\right] + \E_{\cal U}\left[\|\epsilon_{\rm ffd}(x,\mathcal{U})\|^2\right] \\
    &\le \frac{n-s}{n} \|\nabla \phi(x)\|^2 + s \Big(\frac{L_1 {h}}{2} + \frac{2\epsilon_f}{h}\Big)^2. 
\end{aligned} \] 
Meanwhile, 
\[ \begin{aligned}
    \E_{\cal U,\xi}\left[\|g(x,\mathcal{U},\xi) - \E_{\xi}[g(x,\mathcal{U},\xi)]\|^2\right] 
    &= \E_{\cal U,\xi}\left[\|\epsilon_{\rm noise}(x,\mathcal{U},\xi)\|^2\right]. 
\end{aligned} \] 
Thus the conditions of Oracle~\ref{def:oracle_stoch_1.3} are satisfied with 
\[ 
    \epsilon_g^2 = s \big(\frac{L_1 {h}}{2} + \frac{2\epsilon_f}{h}\big)^2, \quad \text{ and }
    \kappa_g^2 = \frac{n-s}{n},
\] 
while $\sigma_c^2 = 2s\sigma_f^2 / {h}^2$, $\sigma_v = 0$ when using Oracle~\ref{def:oracle_stoch_0.1}, or $\sigma_c^2 = sL_c^2$, $\sigma_v = sL_v^2$ when using Oracle~\ref{def:oracle_stoch_0.2}. 

\subsubsection{Randomized Forward Finite Differences with Random Independent Directions}\label{sec:RFFD}
Let us now analyze \eqref{eq:RFFD} under Assumption~\ref{ass:normal} instead of Assumption~\ref{ass:haar}. 
In this case, the projection error $\epsilon_{\rm proj}(x,\mathcal{U})$  in the decomposition \eqref{eq:RFF decomposition} no longer fully lies in the orthogonal complement of span$(\mathcal{U})$. 
We further decompose it into two mutually orthogonal components
\begin{equation*} \begin{aligned}
    \epsilon_{\rm proj}(x,\mathcal{U}) = \epsilon_{\rm proj \perp}(x,\mathcal{U}) + \epsilon_{\rm proj \parallel}(x,\mathcal{U}), 
\end{aligned} \end{equation*}
where $\epsilon_{\rm proj \perp}(x,\mathcal{U})$ is the part contained in the orthogonal complement of span$(\mathcal{U})$, and $\epsilon_{\rm proj \parallel}(x,\mathcal{U})$ is the part that resides in span$(\mathcal{U})$. 
To define $\epsilon_{\rm proj \perp}(x,\mathcal{U})$ and $\epsilon_{\rm proj \parallel}$ rigorously, we introduce the basis matrix $S$ and the expected dimension $\bar s$ in the following lemma. 

\begin{lemma} \label{lem:sbar}
    Define
    \begin{equation*}
        \bar{s}
        = \mathbb{E}\big[\dim(\operatorname{span}(\mathcal{U}))\big] 
    \end{equation*} 
    and let $S$ be any matrix whose columns form an orthonormal basis ($S^TS = I$)
    of $\operatorname{span}(\mathcal{U})$. Then
    \begin{equation*}
        \bar s \le s \quad \text{ and } \quad 
        \mathbb{E}[SS^T] = \frac{\bar{s}}{n}I_n.
    \end{equation*}
\end{lemma}

\begin{proof}
    Because $\dim(\operatorname{span}(\mathcal{U})) \le s$ always holds, its expectation $\bar s$ must also follow $\bar s \le s$. 
    With any signed permutation matrix $Q$, $QS$ is an orthonormal basis of $\text{span}(\{Qu_i\}_{i=1}^s)$ and a projection onto  $\text{span}(\{Qu_i\}_{i=1}^s)$, which is $QSS^TQ^T$
     has the same distribution as the projection only $\mathcal{U}$. 
    Therefore, $\mathbb{E}[SS^T] = \mathbb{E}[QS(QS)^T] = Q\mathbb{E}[SS^T]Q^T$. 
    The same argument in Lemma~\ref{lem:normal} again implies $\mathbb{E}[SS^T]$ is also a scalar multiple of $I_n$ --- sign-flip invariance implies that $\mathbb{E}[SS^T]$ is diagonal, while permutation invariance implies that all its diagonal entries are equal. 
    Then, 
    \[
        \mathbb{E}[SS^T] 
        = \frac{\text{tr}(\mathbb{E}[SS^T])}{n} I_n 
        = \frac{\mathbb{E}[\text{tr}(S^TS)]}{n} I_n 
        = \frac{\mathbb{E}[\dim(\operatorname{span}(\mathcal{U}))]}{n} I_n 
        = \frac{\bar s}{n} I_n. 
    \]
\end{proof}

When $\cal U$ is a set of orthonormal vectors satisfying Assumption~\ref{ass:haar}, we  have $\bar s = s$, and $S$ can be the $n \times s$ matrix with columns $u_1, \dots, u_s$. 
On the other hand, when $\cal U$ is a set of independent vectors under Assumption~\ref{ass:normal}, the values of $\bar s$ and $S$ depend on the specific distribution. 
Assume vectors in $\cal U$ are also identically distributed.\footnote{Recall that Assumption~\ref{ass:normal} only imposes independence.} 
If $\cal U \sim \cal N$ or $\cal S$, the dimension of $\operatorname{span}(\cal U)$ is $\min\{s, n\}$ almost surely, so $\bar s = \min\{s, n\}$. 
The matrix $S$, though only used in the analysis, can be computed using the Gram-Schmidt process on $\cal U$, and the columns of $S$ would follow the Haar distribution if the Gram-Schmidt process is done in a specific way \cite{kozak2023zeroth}. 
When $\cal U$ follows the uniform distribution over $\{\pm e_i\}_{i=1}^n$, the span of $\cal U$ is the subspace consisting of the coordinates contained in $\cal U$, and $S$ is its basis. 
Since the probability of each coordinate not getting picked by any of $u_1,\dots,u_s$ is $(1 - \frac{1}{n})^s$, we have $\bar s = n[1 - (1 - \frac{1}{n})^s]\approx s$ when $s\ll n$. 
When $\cal U$ are i.i.d. random vectors following the scaled Rademacher distribution, $\bar s \approx \min\{s,n\}$, but we are not aware of any closed-form expression of $\bar s$ or $S$. 
In fact, even the probability of a square i.i.d. Rademacher matrix $R$ of order $n$ being singular, i.e. $\mathbf{P}[\operatorname{rank}(R) < n]$, is still an open question. 
It is only recently proved in \cite{hunter2026random} that there exists an absolute constant $c > 0$ such that $\mathbf{P}[\operatorname{rank}(R) \le n-k] \le \exp(-cnk)$.  For what follows we will use $\bar s$ in the intermediate computation but our upper bounds do not depend on it. 

With the matrix $S$ defined above, we formally define 
\[ \begin{aligned}
    \epsilon_{\rm proj \perp}(x,\mathcal{U}) &= SS^T \nabla\phi(x) - \nabla\phi(x) \\
    \text{ and } \epsilon_{\rm proj \parallel}(x,\mathcal{U}) &= \bigg(\sum_{i=1}^s (\nabla\phi(x)^T u_i) u_i \bigg) - SS^T \nabla\phi(x). 
\end{aligned} \] 
We again analyze \eqref{eq:RFFD} as Oracle~\ref{def:oracle_stoch_1.3}, starting with bounding the four errors $\epsilon_{\rm proj \perp}(x,\mathcal{U})$, $\epsilon_{\rm proj \parallel}(x,\mathcal{U})$, $\epsilon_{\rm ffd}(x,\mathcal{U})$, and $\epsilon_{\rm noise}(x,\mathcal{U},\xi)$. 

\paragraph{Bounding $\epsilon_{\rm proj \perp}$ and $\epsilon_{\rm proj \parallel}$.}
Using Lemma~\ref{lem:normal}, we have under Assumption~\ref{ass:normal} that 
\begin{equation} \label{eq:RFFD_proj_error_perp} \begin{aligned}
    \E_{\cal U} \left[\|\epsilon_{\rm proj \perp}(x,\mathcal{U}) \|^2\right] 
    &= \nabla\phi(x)^T \E_{\cal U} \left[SS^TSS^T\!\! - 2SS^T\!\! + I\right] \nabla\phi(x) \\
    &= \frac{n - \bar s}{n} \|\nabla \phi(x)\|^2. 
\end{aligned} \end{equation}

Now using the derivations for  $\epsilon_{\rm proj}(x,\mathcal{U})$ in Section \ref{sec:rand_dirder}  
and specifically the value of $\E \left[\|\epsilon_{\rm proj}(x,\mathcal{U}) \|^2\right]$ established in \eqref{eq:e_proj}, by the orthogonality between $\epsilon_{\rm proj \perp}(x,\mathcal{U})$ and $\epsilon_{\rm proj \parallel}(x,\mathcal{U})$, we have 
\begin{equation} \label{eq:RFFD_proj_error_parallel} \begin{aligned}
    \E \left[\|\epsilon_{\rm proj \parallel}(x,\mathcal{U}) \|^2\right] 
    &= \E \left[\|\epsilon_{\rm proj}(x,\mathcal{U}) \|^2\right] - \E \left[\|\epsilon_{\rm proj \perp}(x,\mathcal{U}) \|^2\right] \\
    &= \bigg(\Big(\frac{n-s}{n}\Big)^2 + \frac{(1+\mathcal{O}(n^{-1}))s}{n}\bigg) \| \nabla\phi(x)\|^2 - \frac{n - \bar s}{n} \|\nabla \phi(x)\|^2 \\
    &= \bigg(\frac{s^2}{n^2} + \frac{\mathcal{O}(n^{-1})s}{n} - \frac{s - \bar{s}}{n}\bigg) \|\nabla \phi(x)\|^2. 
\end{aligned} \end{equation}

\paragraph{Bounding $\epsilon_{\rm ffd}$.}
We now turn to bounding $\epsilon_{\rm ffd}$. 
As before, let $\epsilon_1(x,{h} u_i)=\phi(x+ {h} u_i) - \phi(x) - \nabla\phi(x)^T ({h} u_i)$ and $\epsilon_2(x,{h} u_i)=\E_{\xi }[f(x+ {h} u_i, \xi_i)]-\phi(x+{h}u_i) - (\E_{\xi }[f(x, \xi)]-\phi(x))$, and we have 
\begin{align*}
    \|\epsilon_{\rm ffd}(x,\mathcal{U})\|^2 = \left \|\sum_{i=1}^s\left (\frac{\epsilon_1(x,{h} u_i)}{h} + \frac{\epsilon_2(x,{h} u_i)}{ {h}} \right ) u_i\right \|^2 
    = \frac{1}{h^2} \sum_{i=1}^s \left( \epsilon_1(x,{h} u_i) + \epsilon_2(x,{h} u_i) \right)^2 \|u_i\|^2& \\ 
    + \frac{2}{h^2} \sum_{i=1}^s \sum_{j=i+1}^s \left( \epsilon_1(x,{h} u_i) + \epsilon_2(x,{h} u_i) \right) \left(\epsilon_1(x,{h} u_j) + \epsilon_2(x,{h} u_j) \right) u_i^T u_j.& 
\end{align*}

When $\phi$ satisfies Assumption~\ref{assum:lip_cont} and its value is obtained through either Oracle~\ref{def:oracle_stoch_0.1} or Oracle~\ref{def:oracle_stoch_0.2}, we have $|\epsilon_1(x,{h} u_i)| \le L_1 {h}^2 \|u_i\|^2 / 2$ and $|\epsilon_2(x,{h} u_i)| \le 2\epsilon_f$. 
Then, for all $i = 1,\dots,s$, 
\begin{equation*} \begin{aligned}
    \mathbb{E}_{\cal U} \left[\left( \epsilon_1(x,{h} u_i) + \epsilon_2(x,{h} u_i) \right)^2 \|u_i\|^2\right] 
    &\le \mathbb{E}_{\cal U} \left[\left( L_1 h^2 \|u_i\|^2/2 + 2\epsilon_f \right)^2 \|u_i\|^2\right] \\
    &= \Big(\frac{(1+\mathcal{O}(n^{-1})) L_1 h^2}{2} + 2\epsilon_f\Big)^2. 
\end{aligned} \end{equation*}
Under Assumption~\ref{ass:normal}, by the independence among ${\cal U}$, each cross term ($i \neq j$) is equal to 
\begin{equation*} \begin{aligned}
    &\mathbb{E}_{\cal U} \left[\left( \epsilon_1(x,{h} u_i) + \epsilon_2(x,{h} u_i) \right) \left(\epsilon_1(x,{h} u_j) + \epsilon_2(x,{h} u_j) \right) u_i^T u_j\right] \\
    &=\mathbb{E}_{\cal U} \left[\left( \epsilon_1(x,{h} u_i) + \epsilon_2(x,{h} u_i) \right)\right ]^T \mathbb{E}_{\cal U} \left[\left(\epsilon_1(x,{h} u_j) + \epsilon_2(x,{h} u_j) \right) u_i^T u_j\right] \\
    &\leq \left\|\mathbb{E}_{\cal U} \left[\left( \epsilon_1(x,{h} u_i) + \epsilon_2(x,{h} u_i) \right) u_i\right] \right\| \left\|\mathbb{E}_{\cal U} \left[\left( \epsilon_1(x,{h} u_j) + \epsilon_2(x,{h} u_j) \right) u_j\right] \right\|, 
\end{aligned} \end{equation*}
where, by the Cauchy–Schwarz inequality, we have for all $i = 1, \dots, s$ that 
\begin{equation*} \begin{aligned}
    &\left\|\mathbb{E}_{\cal U} \left[\left( \epsilon_1(x,{h} u_i) + \epsilon_2(x,{h} u_i) \right) u_i\right] \right\| \\ 
    &= \max_{v\in\mathbb{R}^n:~\|v\| = 1} \left\{\mathbb{E}_{\cal U} \left[\left( \epsilon_1(x,{h} u_i) + \epsilon_2(x,{h} u_i) \right) u_i^T v\right]\right\} \\
    &\le \max_{v\in\mathbb{R}^n:~\|v\| = 1} \left\{\sqrt{\mathbb{E}_{\cal U} \left[\left( \epsilon_1(x,{h} u_i) + \epsilon_2(x,{h} u_i) \right)^2\right] \mathbb{E}_{\cal U} \left[(u_i^T v)^2\right]}\right\} \\
    &\le \sqrt{\mathbb{E}_{\cal U} \left[\left( L_1 h^2 \|u_i\|^2 / 2 + 2\epsilon_f \right)^2\right]} \max_{v\in\mathbb{R}^n:~\|v\| = 1} \left\{\sqrt{v^T \mathbb{E}_{\cal U} \left[u_i u_i^T\right] v} \right\} \\
    &\leq \Big(\frac{(1+\mathcal{O}(n^{-1})) L_1 h^2}{2} + 2\epsilon_f\Big) \cdot \frac{1}{\sqrt{n}}. 
\end{aligned} \end{equation*}

Thus, when $\phi$ satisfies Assumption~\ref{assum:lip_cont}, the distribution of $\mathcal{U}$ satisfies Assumption~\ref{ass:normal}, and the function values are obtained through either Oracle~\ref{def:oracle_stoch_0.1} or Oracle~\ref{def:oracle_stoch_0.2}, 
\begin{equation*} \begin{aligned}
    \mathbb{E}_{\cal U} \left[\|\epsilon_{\rm ffd}(x,\mathcal{U})\|^2\right] 
    &\le \frac{1}{h^2} \Big(\frac{(1+\mathcal{O}(n^{-1})) L_1 h^2}{2} + 2\epsilon_f\Big)^2 \cdot \Big( s + \frac{s^2-s}{n} \Big) \\
    &\le \frac{n+s-1}{n} \cdot s\, \Big(\frac{(1+\mathcal{O}(n^{-1})) L_1 h}{2} + \frac{2\epsilon_f}{h} \Big)^2. 
\end{aligned} \end{equation*}

\paragraph{Bounding $\epsilon_{\rm noise}$.}
Moving on to the stochastic noise error, let $\epsilon_3(x,\xi) = f(x, \xi)-\E_{\xi}[f(x, \xi)]$. 
When using Oracle~\ref{def:oracle_stoch_0.1}, it follows that $\E[\epsilon_3(x,\xi)] = 0$ and $\E[\epsilon_3(x,\xi)^2] \le \sigma_f^2$. 
Subsequently, we have $\E_\xi [\epsilon_{\rm noise}(x,\mathcal{U},\xi)] = 0$ and, since $\mathcal{U}, \xi_0, \xi_1, \dots, \xi_s$ are independent from each other, that 
\begin{align*}
    \E_{\mathcal{U},\xi} &\left[ \|\epsilon_{\rm noise}(x,\mathcal{U},\xi)\|^2 \right] 
    = \frac{1}{{h}^2} \mathbb{E}_{\mathcal{U},\xi} \left[ \Big\|\sum_{i=1}^s \left(\epsilon_3(x+{h} u_i,\xi_i) -  \epsilon_3(x,\xi_0)\right)  u_i \Big\|^2 \right] \\
	&= \frac{1}{{h}^2} \mathbb{E}_{\mathcal{U},\xi} \left[\Big\|\Big(\sum_{i=1}^s \epsilon_3(x+{h} u_i,\xi_i) u_i\Big) - \Big(\epsilon_3(x,\xi_0) \sum_{i=1}^s u_i\Big) \Big\|^2 \right] \\
    &= \frac{1}{{h}^2} \mathbb{E}_{\mathcal{U},\xi} \left[ \Big(\sum_{i=1}^s \|\epsilon_3(x+{h} u_i,\xi_i) u_i\|^2\Big) + \| \epsilon_3(x,\xi_0) \sum_{i=1}^s u_i \|^2 \right] \\
    &\le \frac{\sigma_f^2}{{h}^2} \mathbb{E}_{\cal U} \left[ \Big(\sum_{i=1}^s \|u_i\|^2\Big) + \Big\|\sum_{i=1}^s u_i \Big\|^2 \right] 
    = \frac{\sigma_f^2}{{h}^2} \mathbb{E}_{\cal U} \left[ 2 \sum_{i=1}^s \|u_i\|^2 \right] 
     = \frac{2s \sigma_f^2}{{h}^2}, 
\end{align*}
where the second-to-last equality holds be Lemma~\ref{lem:normal}. 

When using Oracle~\ref{def:oracle_stoch_0.2}, we still have $\E_\xi [\epsilon_{\rm noise}(x,\mathcal{U},\xi)] = 0$, but the $\epsilon_3$ terms in 
\begin{align*}
    \E_{\mathcal{U},\xi} \left[ \|\epsilon_{\rm noise}(x,\mathcal{U},\xi)\|^2 \right] 
    = \frac{1}{{h}^2} \mathbb{E}_{\mathcal{U},\xi} \left[ \Big\|\sum_{i=1}^s \left(\epsilon_3(x+{h} u_i,\xi) -  \epsilon_3(x,\xi)\right)  u_i \Big\|^2 \right] 
\end{align*}
are no longer independent due to the common random number $\xi$. 
We bound the pure quadratic terms (involving $u_i^T u_i$) and the cross terms (involving $u_i^T u_j$, $i \neq j$) separately. 
For $i = 1,\dots, s$, we have 
\[ \begin{aligned}
    &\frac{1}{{h}^2} \mathbb{E}_{\mathcal{U},\xi} \left[ \Big\|\left(\epsilon_3(x+{h} u_i,\xi) -  \epsilon_3(x,\xi)\right)  u_i \Big\|^2 \right] \\
    &= \frac{1}{{h}^2} \mathbb{E}_{\mathcal{U}} \left[ \mathbb{E}_{\xi} \left[ \left(\epsilon_3(x+{h} u_i,\xi) -  \epsilon_3(x,\xi)\right)^2 \right] \|u_i\|^2 \right] \\
    &\le (L_c^2+L_v^2\|\nabla\phi(x)\|^2) \mathbb{E}_{\mathcal{U}} \left[ \|u_i\|^4 \right] 
    = (L_c^2+L_v^2\|\nabla\phi(x)\|^2) (1 + \mathcal{O}(n^{-1})). 
\end{aligned} \]
By the independence among $\{u_i\}$ and the Cauchy–Schwarz inequality, the cross terms satisfy 
\[ \begin{aligned}
    &\frac{1}{{h}^2} \mathbb{E}_{\mathcal{U},\xi} \left[ \Big(\left(\epsilon_3(x+{h} u_i,\xi) -  \epsilon_3(x,\xi)\right)  u_i \Big)^T \Big(\left(\epsilon_3(x+{h} u_j,\xi) -  \epsilon_3(x,\xi)\right)  u_j \Big) \right] \\
    &= \frac{1}{{h}^2} \mathbb{E}_{\xi} \left[ \max_{v\in\mathbb{R}^n:~ \|v\|=1} \Big\{ \Big(\mathbb{E}_{\mathcal{U}} \left[ (\epsilon_3(x+{h} u_i,\xi) -  \epsilon_3(x,\xi))  u_i^T v\right] \Big)^2 \Big\} ~\middle|~ \xi \right] \\
    &\le \frac{1}{{h}^2} \mathbb{E}_{\mathcal{U},\xi} \left[ (\epsilon_3(x+{h} u_i,\xi) -  \epsilon_3(x,\xi))^2\right] \max_{v\in\mathbb{R}^n:~ \|v\|=1} \Big\{ \mathbb{E}_{\mathcal{U}} \left[ (u_i^T v)^2 \right] \Big\} \\
    &\le \frac{1}{{h}^2} \mathbb{E}_\mathcal{U}\left[ (L_c^2+L_v^2\|\nabla\phi(x)\|^2) \|h u_i\|^2 \right] \cdot \frac{1}{n} 
    = \frac{1}{n} (L_c^2+L_v^2\|\nabla\phi(x)\|^2). 
\end{aligned} \] 

Finally, 
\begin{align*}
    \E_{\mathcal{U},\xi} \left[ \|\epsilon_{\rm noise}(x,\mathcal{U},\xi)\|^2 \right] 
    &\le \Big(s \cdot (1 + \mathcal{O}(n^{-1})) + (s^2 - s) \cdot \frac{1}{n} \Big) (L_c^2+L_v^2\|\nabla\phi(x)\|^2)  \\
    &= \Big(\frac{n+s-1}{n} + \mathcal{O}(n^{-1})\Big) \cdot s\, (L_c^2+L_v^2\|\nabla\phi(x)\|^2). 
\end{align*}

\paragraph{As Oracle~\ref{def:oracle_stoch_1.3}.} 
To cast \eqref{eq:RFFD} as instances of Oracle~\ref{def:oracle_stoch_1.3}, we combine the bounds derived above for 
\begin{equation*} \begin{aligned}
    &\E_{\cal U} \left[ \big\| \E_{\xi }[ g(x,\mathcal{U},\xi)] - \nabla\phi(x)\big\|^2 \right] 
    \begin{aligned}[t] 
        &= \E_{\cal U} \Big[ \|\epsilon_{\rm proj\perp}\|^2 + \|\epsilon_{\rm proj\parallel} + \epsilon_{\rm ffd}\|^2 \Big] \\
        &\le \E_{\cal U} \Big[\|\epsilon_{\rm proj\perp}\|^2 + 2\|\epsilon_{\rm proj\parallel}\|^2 + 2\|\epsilon_{\rm ffd} \|^2\Big]
    \end{aligned} \\
    \text{ and } &\E_{\mathcal{U},\xi} \left[ \big\|\E_{\xi} \left[g(x,\mathcal{U},\xi)\right] - g(x,\mathcal{U},\xi) \big\|^2\right]
    = \E \left[ \|\epsilon_{\rm noise}(x,\mathcal{U},\xi)\|^2 \right]. 
\end{aligned} \end{equation*}
The input  $\lambda$ includes the sample size $s$ and the distribution of $\mathcal{U}$. 
Under Assumptions~\ref{assum:lip_cont} and \ref{ass:normal}, the constants in the bounds are 
\begin{align*}
    \epsilon_g^2 &= \frac{2(n+s-1)}{n} \cdot s\, \Big(\frac{(1+\mathcal{O}(n^{-1}))L_1 h}{2}+ \frac{2\epsilon_f}{h} \Big)^2, \\
    \kappa_g^2 &= 1 - \frac{s}{n} + \frac{2s^2}{n^2} + \frac{\mathcal{O}(n^{-1})s}{n},  
\end{align*}
while $\sigma_c^2 = 2s\sigma_f^2 / {h}^2$, $\sigma_v = 0$ when using Oracle~\ref{def:oracle_stoch_0.1}, or 
\begin{equation*}
    \sigma_c^2 = \Big(\frac{n+s-1}{n} + \mathcal{O}(n^{-1})\Big) \cdot sL_c^2, \quad 
    \sigma_v^2 = \Big(\frac{n+s-1}{n} + \mathcal{O}(n^{-1})\Big) \cdot sL_v^2 
\end{equation*}
when using Oracle~\ref{def:oracle_stoch_0.2}. 

\subsubsection{Scaled Randomized Forward Finite Differences}\label{sec:SRFFD}
Similar to approximating the gradient via directional derivatives, the randomized finite differences \eqref{eq:RFFD} can also be scaled by a factor of $n/s$, yielding
\begin{equation} \label{eq:SRFFD} \tag{SRFFD}
	g(x,\mathcal{U},\xi) = \frac{n}{s} \sum_{i=1}^s \frac{f(x+ h u_i, \xi_i) - f(x,\xi_0)}{h} u_i. 
\end{equation} 
This scaling is a popular choice in the literature because it is based on the unbiased estimator \eqref{eq:subspace_grad_scale}. 
However, as we show later, this choice may lead to worse bounds and, consequently, worse theoretical complexity results than the unscaled estimator \eqref{eq:RFFD}.

In this case, treating the scheme \eqref{eq:SRFFD} as a case of Oracle~\ref{def:oracle_stoch_1.3} will not be meaningful since $\kappa_g$ will no longer remain in $(0,1)$. 
Thus, we analyze it in terms of Oracle~\ref{def:oracle_stoch_1.2}. 
The vector $g(x,\mathcal{U},\xi)$ can be decomposed as 
\begin{equation} \label{eq:FFD decomposition} \begin{aligned}
	g(x,\mathcal{U},\xi) = \nabla\phi(x) + \epsilon_{\rm proj \perp}(x,\mathcal{U}) + \epsilon_{\rm proj \parallel}(x,\mathcal{U})  + \epsilon_{\rm ffd}(x,\mathcal{U}) + \epsilon_{\rm noise}(x,\mathcal{U},\xi), 
\end{aligned} \end{equation}
where
\begin{align*}
    \epsilon_{\rm proj \perp}(x,\mathcal{U}) 
    &= SS^T \nabla\phi(x) - \nabla\phi(x), \\
    \epsilon_{\rm proj \parallel}(x,\mathcal{U}) 
    &= \Big(\frac{n}{s} \sum_{i=1}^s (\nabla\phi(x)^T u_i) u_i\Big) - SS^T \nabla\phi(x), \\
	\epsilon_{\rm ffd}(x,\mathcal{U}) 
    &= \frac{n}{s} \sum_{i=1}^s \bigg[\frac{\phi(x+ {h} u_i) - \phi(x) - \nabla\phi(x)^T ({h} u_i)}{h} + \\ 
    &\hspace{5em} \frac{\E_{\xi }[f(x+ {h} u_i, \xi_i)]-\phi(x+{h}u_i) - (\E_{\xi }[f(x, \xi_0)]-\phi(x))}{ {h}}\bigg] u_i, \\
    \text{and } \epsilon_{\rm noise}(x,\mathcal{U},\xi) 
    &= \frac{n}{s} \sum_{i=1}^s \frac{f(x+ {h} u_i, \xi_i)-\E_{\xi}[f(x+ {h} u_i, \xi_i)] - f(x, \xi_0) + \E_{\xi}[f(x,\xi_0)]}{h} u_i.  
\end{align*} 
Additionally, we still denote the sum of $\epsilon_{\rm proj \perp}$ and $\epsilon_{\rm proj \parallel}$ as $\epsilon_{\rm proj}$. 
It holds that $\E_{{\cal U}}[\epsilon_{\rm proj}(x,\mathcal{U})] = 0$ when $\mathcal{U}$ satisfies either Assumption~\ref{ass:haar} or Assumption~\ref{ass:normal}, because $\E[u_i u_i^T] = \frac{1}{n} I$. 
In this decomposition, the first component $\epsilon_{\rm proj \perp}(x,\mathcal{U})$ is in the orthogonal complement of $\text{span}(\mathcal{U})$, while the other three are in $\text{span}(\mathcal{U})$. 
They are bounded as follows. 

\paragraph{Bounding $\epsilon_{\rm proj \perp}$ and $\epsilon_{\rm proj \parallel}$.}
Equations \eqref{eq:RFFD_proj_error} and \eqref{eq:RFFD_proj_error_perp}, both stating 
\[
    \E_{\cal U} \left[\|\epsilon_{\rm proj \perp}(x,\mathcal{U}) \|^2\right] = \frac{n-\bar{s}}{n} \|\nabla \phi(x)\|^2, 
\]
still hold. 
Under Assumption~\ref{ass:haar}, we have $\bar{s} = s$, and combining \eqref{eq:RFFD_proj_error} with \eqref{eq:subspace_grad_S1.2} yields that 
\begin{equation*} \begin{aligned}
    \E_{{\cal U}} \left[\|\epsilon_{\rm proj \parallel}(x,\mathcal{U}) \|^2\right] 
    &= \E_{{\cal U}} \left[\|\epsilon_{\rm proj}(x,\mathcal{U}) \|^2\right] - \E_{{\cal U}}\left[\|\epsilon_{\rm proj \perp}(x,\mathcal{U}) \|^2\right] \\
    &= \frac{n-s}{s} \|\nabla\phi(x)\|^2 - \frac{n-s}{n} \|\nabla \phi(x)\|^2
    = \frac{(n-s)^2}{ns} \|\nabla \phi(x)\|^2. 
\end{aligned} \end{equation*}
Alternatively under Assumption~\ref{ass:normal}, using the independence of $\mathcal{U}$ and Lemma~\ref{lem:normal}, we have 
\begin{equation*} \begin{aligned}
    &\E_{{\cal U}} \left[\|\epsilon_{\rm proj \parallel}(x,\mathcal{U}) \|^2\right] 
    = \E_{{\cal U}} \left[\Big\|\Big(\frac{n}{s} \sum_{i=1}^s (\nabla\phi(x)^T u_i) u_i\Big) - SS^T \nabla\phi(x)\Big\|^2\right] \\
    &= \nabla\phi(x)^T \E_{{\cal U}} \left[\frac{n^2}{s^2} \Big(\sum_{i=1}^s \|u_i\|^2 u_iu_i^T\Big) - \frac{2n}{s} \Big(\sum_{i=1}^s u_i u_i^T\Big)  + SS^TSS^T\right] \nabla\phi(x) \\ 
    &= \Big(\frac{(1+\mathcal{O}(n^{-1})) n}{s} - 2 + \frac{\bar s}{n}\Big) \|\nabla\phi(x)\|^2, 
\end{aligned} \end{equation*}
where the result also equals to $\frac{(n-s)^2}{ns}$ assuming $s \le n$ and the $\mathcal{O}(1/n)$ term is zero. 

\paragraph{Bounding $\epsilon_{\rm ffd}$.}
Other than being scaled by $\frac{n}{s}$, the error $\epsilon_{\rm ffd}$ is identical to its unscaled counterpart in Sections~\ref{sec:RFFD_orth} and \ref{sec:RFFD}, so the bounds on $\E_{{\cal U}} \left[ \|\epsilon_{\rm ffd}(x,\mathcal{U})\|^2\right]$ are simply what we have derived in those two sections scaled by $\frac{n^2}{s^2}$. 

However, we also need to bound $\|\E_{\cal U}[\epsilon_{\rm ffd}(x,\mathcal{U})]\|^2$, so again let $\epsilon_1(x,{h} u_i)=\phi(x+ {h} u_i) - \phi(x) - \nabla\phi(x)^T ({h} u_i)$ and assume Assumption~\ref{assum:lip_cont}, and let $\epsilon_2(x,{h} u_i)=\E_{\xi}[f(x+ {h} u_i, \xi_i)]-\phi(x+{h}u_i) - (\E_{\xi }[f(x, \xi)]-\phi(x))$ and suppose the function values come from either Oracle~\ref{def:oracle_stoch_0.1} or \ref{def:oracle_stoch_0.2}. 
Under Assumption~\ref{ass:normal}, we have 
\begin{equation} \label{eq:||E eps_ffd||} \begin{aligned}
    &\E_{\cal U}[\|\epsilon_{\rm ffd}(x,\mathcal{U})\|]^2
    = \left\|\E_{\cal U}\left[\frac{n}{s} \sum_{i=1}^s\left(\frac{\epsilon_1(x,{h} u_i)}{h}+ \frac{\epsilon_2(x,{h} u_i)}{ {h}} \right) u_i\right] \right\|^2 \\
    &= \frac{n^2}{s^2h^2} \max_{v\in\mathbb{R}^n:~ \|v\|=1} \left\{\bigg(\sum_{i=1}^s \E_{\cal U}\left[\left(\epsilon_1(x,h u_i) + \epsilon_2(x,{h} u_i)\right) u_i^T v\right]\bigg)^2\right\} \\
    &\le \frac{n^2}{s^2h^2} \max_{v\in\mathbb{R}^n:~ \|v\|=1} \left\{\bigg(\sum_{i=1}^s \sqrt{\E_{\cal U}\left[\left(\epsilon_1(x,h u_i) + \epsilon_2(x,h u_i)\right)^2\right]\E_{\cal U}\left[(u_i^T v)^2\right]}\bigg)^2\right\} \\
    &\le \frac{n^2}{s^2h^2} \Big(\frac{(1+\mathcal{O}(n^{-1})) L_1h^2}{2} + 2\epsilon_f\Big)^2 \max_{v\in\mathbb{R}^n:~ \|v\|=1} \left\{\bigg(\sum_{i=1}^s \sqrt{v^T \E_{\cal U}\left[u_i u_i^T\right] v}\bigg)^2\right\} \\
    &= \frac{n^2}{s^2h^2} \Big(\frac{(1+\mathcal{O}(n^{-1})) L_1h^2}{2} + 2\epsilon_f\Big)^2 \left(\frac{s}{\sqrt{n}}\right)^2 
    = n \Big(\frac{(1+\mathcal{O}(n^{-1})) L_1 {h}}{2} + \frac{2\epsilon_f}{h}\Big)^2, 
\end{aligned} \end{equation}
where the first inequality holds by the Cauchy-Schwarz inequality. 
This derivation also works under Assumption~\ref{ass:haar}, except the second last equality becomes inequality, which holds by the concavity of the square root function and Jensen’s inequality: 
\[ \begin{aligned}
    &\max_{v\in\mathbb{R}^n:~ \|v\|=1} \left\{\bigg(\sum_{i=1}^s \sqrt{v^T \E_{\cal U}\left[u_i u_i^T\right] v}\bigg)^2\right\} \\
    &= \max_{v\in\mathbb{R}^n:~ \|v\|=1} \left\{\bigg(s \cdot \frac{1}{s} \sum_{i=1}^s \sqrt{v^T \E_{\cal U}\left[u_i u_i^T\right] v}\bigg)^2\right\} \\
    &\le \max_{v\in\mathbb{R}^n:~ \|v\|=1} \left\{\Bigg(s \cdot \sqrt{\frac{1}{s} \sum_{i=1}^s v^T \E_{\cal U}\left[u_i u_i^T\right] v}\Bigg)^2\right\} 
    = \frac{s^2}{n}. 
\end{aligned} \]
The resulting bound is therefore identical to the one under Assumption~\ref{ass:normal} except without the $\mathcal{O}(n^{-1})$ term: 
\[
    \|\E_{\cal U}[\epsilon_{\rm ffd}(x,\mathcal{U})]\|^2 
    \le n \Big(\frac{L_1 {h}}{2} + \frac{2\epsilon_f}{h}\Big)^2. 
\]

\paragraph{Bounding $\epsilon_{\rm noise}$.}
Similar to $\epsilon_{\rm ffd}$, other than being scaled by $\frac{n}{s}$, the stochastic noise term $\epsilon_{\rm noise}$ is identical to its unscaled counterpart in Sections~\ref{sec:RFFD_orth} and \ref{sec:RFFD}, so the bounds on $\E_{{\cal U},\xi} \left[ \|\epsilon_{\rm noise}(x,\mathcal{U},\xi)\|^2\right]$ are simply what we have derived in those two sections scaled by $\frac{n^2}{s^2}$. 

\paragraph{As Oracle~\ref{def:oracle_stoch_1.2}.} 
Since $\E_{\cal U} [\epsilon_{\rm proj}(x,\mathcal{U})] = 0$ and $\E_{\cal \xi} [\epsilon_{\rm noise}(x,\mathcal{U},\xi)] = 0$, the only bias comes from $\epsilon_{\rm ffd}$, and $\E_{{\cal U}, \xi}[g(x,\mathcal{U},\xi)] = \E_{{\cal U}, \xi}[\epsilon_{\rm ffd}(x,\mathcal{U})]$. 
The squared bias is therefore 
\begin{equation*} \begin{aligned}
    \|\E_{\cal U,\xi} [g(x,{\cal U},\xi)] - \nabla\phi(x)\|^2 
    = \|\E_{\cal U} [\epsilon_{\rm ffd}(x,\mathcal{U})]\|^2. 
\end{aligned} \end{equation*}
Since only $\epsilon_{\rm proj\perp}(x,\mathcal{U})$ is in the orthogonal complement of $\text{span}(\cal U)$, and $\E_{\cal \xi} [\epsilon_{\rm noise}(x,\mathcal{U},\xi)] = 0$ with $\xi$ being independent of $\cal U$, the variance can be bounded as 
\begin{equation*} \begin{aligned} 
    &\E_{{\cal U}, \xi}\left[\|g(x,\mathcal{U},\xi) - \E_{{\cal U}, \xi}[g(x,\mathcal{U},\xi)]\|^2\right] 
    \le \E_{{\cal U}, \xi}\left[\|g(x,\mathcal{U},\xi)\|^2\right] \\ 
    &= \E_{{\cal U}, \xi}\left[\|\epsilon_{\rm proj\perp} + \epsilon_{\rm proj\parallel} + \epsilon_{\rm ffd} + \epsilon_{\rm noise}\|^2\right] \\ 
    &= \E_{\cal U}\left[\|\epsilon_{\rm proj\perp}\|^2\right] + \E_{{\cal U}, \xi}\left[\|\epsilon_{\rm noise}\|^2\right] + \E_{{\cal U}}\left[\|\epsilon_{\rm proj\parallel} + \epsilon_{\rm ffd}\|^2\right] \\
    &\le \E_{\cal U}\left[\|\epsilon_{\rm proj\perp}\|^2\right] + \E_{{\cal U}, \xi}\left[\|\epsilon_{\rm noise}\|^2\right] + 2\E_{{\cal U}}\left[\|\epsilon_{\rm proj\parallel}\|^2\right] + 2\E_{{\cal U}}\left[\|\epsilon_{\rm ffd}\|^2\right]. 
\end{aligned} \end{equation*}

Combining the results above, we conclude that under Assumption~\ref{ass:normal} the vector $g(x,\mathcal{U},\xi)$ defined in \eqref{eq:SRFFD} satisfies the conditions of Oracle~\ref{def:oracle_stoch_1.2} with 
\[ \begin{aligned}
    \epsilon_g^2 &= n \Big(\frac{(1+\mathcal{O}(n^{-1})) L_1 {h}}{2} + \frac{2\epsilon_f}{h}\Big)^2, \quad 
    \kappa_g^2 = 0, \\
    \sigma_c^2 &= \frac{2n^2\sigma_f^2}{s h^2} + \frac{n+s-1}{n} \cdot \frac{n^2}{s} \Big(\frac{(1+\mathcal{O}(n^{-1})) L_1 h}{2} + \frac{2\epsilon_f}{h} \Big)^2, \quad 
    \sigma_v^2 = C 
\end{aligned} \]  
when using Oracle~\ref{def:oracle_stoch_0.1}, where 
\[ 
    C := \frac{n-\bar{s}}{n} + 2 \Big(\frac{(1+\mathcal{O}(n^{-1})) n}{s} - 2 + \frac{\bar s}{n}\Big), 
\]
which simplifies to 
\[
    C = \frac{(2n-s)(n-s)}{ns}
\]
when $\bar{s} = s$ and the $\mathcal{O}(n^{-1})$ term is zero. 
When using Oracle~\ref{def:oracle_stoch_0.2}, $\sigma_c$ and $\sigma_v$ change to 
\[ \begin{aligned}
    \sigma_c^2 &= \Big(\frac{n+s-1}{n} + \mathcal{O}(n^{-1})\Big) \cdot \frac{n^2 L_c^2}{s} + \frac{n+s-1}{n} \cdot \frac{n^2}{s} \Big(\frac{(1+\mathcal{O}(n^{-1})) L_1 h}{2} + \frac{2\epsilon_f}{h} \Big)^2, \\
    \sigma_v^2 &= \Big(\frac{n+s-1}{n} + \mathcal{O}(n^{-1})\Big) \cdot \frac{n^2 L_v^2}{s}+ C. 
\end{aligned} \]  

Alternatively when $\cal U$ are orthonormal vectors satisfying Assumption~\ref{ass:haar}, the conditions of Oracle~\ref{def:oracle_stoch_1.2} are satisfied with 
\[ \begin{aligned}
    \epsilon_g^2 = n \Big(\frac{L_1 {h}}{2} + \frac{2\epsilon_f}{h}\Big)^2, \quad 
    \kappa_g^2 = 0, 
\end{aligned} \]
while 
\[
    \sigma_c^2 = \frac{2n^2\sigma_f^2}{s h^2} + \frac{n^2}{s} \Big(\frac{L_1 {h}}{2} + \frac{2\epsilon_f}{h}\Big)^2, \quad 
    \sigma_v^2 = \frac{(2n-s)(n-s)}{ns}
\]
when using Oracle~\ref{def:oracle_stoch_0.1}, or 
\[
    \sigma_c^2 = \frac{n^2 L_c^2}{s} + \frac{n^2}{s} \Big(\frac{L_1 {h}}{2} + \frac{2\epsilon_f}{h}\Big)^2, \quad 
    \sigma_v^2 = \frac{n^2 L_v^2}{s} + \frac{(2n-s)(n-s)}{ns}
\]
when using Oracle~\ref{def:oracle_stoch_0.2}.

\subsubsection{Randomized Central Finite Differences with Orthogonal Directions}\label{sec:RCFD}
The randomized forward finite difference formula \eqref{eq:RFFD} can be replaced by its antithetic variant 
\begin{equation} \label{eq:RCFD} \tag{RCFD}
	g(x,\mathcal{U},\xi) = \sum_{i=1}^s \frac{f(x+ {h} u_i, \xi_{i+}) - f(x- {h} u_i,\xi_{i-})}{2{h}} u_i, 
\end{equation} 
with $\xi= \{\xi_{i+}, \xi_{i-}\}_{i=1}^s$ when using the single-point oracle \ref{def:oracle_stoch_0.1}; or with $\xi = \xi_{i+} = \xi_{i-}$ for all $i$ when using the multi-point oracle \ref{def:oracle_stoch_0.2}. 
This CFD (Central Finite Difference) scheme is very popular in the literature, since it has some theoretical and numerical advantages over the FFD scheme (similarly to their deterministic counterparts). 
CFD can be analyzed under all the same cases that we considered for the FFD, and the techniques are similar. 
Here we will present the analysis only under Assumption~\ref{ass:haar}. 

As in Section~\ref{sec:RFFD_orth}, we write 
\begin{equation} \label{eq:CFD decomposition}
   g(x,\mathcal{U}) = \nabla\phi(x) + \epsilon_{\rm proj}(x,\mathcal{U}) + \epsilon_{\rm cfd}(x,\mathcal{U}) + \epsilon_{\rm noise}(x,\mathcal{U}), 
\end{equation}
where $\epsilon_{\rm proj}(x,\mathcal{U})$ remains the same as in \eqref{eq:RFF decomposition}, 
but the (central) finite difference error and the error caused by noise are now 
\begin{align*}
    \epsilon_{\rm cfd}(x,\mathcal{U}) = \sum_{i=1}^s &\left(
        \frac{\phi(x+ {h} u_i) - \phi(x-{h} u_i) - 2\nabla\phi(x)^T ({h} u_i)}{2h}\right . \\
        & \left .+ \frac{ 
            \E_{\xi}[f(x+ {h} u_i, \xi_{i+})] - \phi(x+{h}u_i) - \E_{\xi}[f(x- {h} u_i, \xi_{i-})] + \phi(x- {h} u_i)
            }{2h} 
    \right) u_i, 
\end{align*} 
and, using $\epsilon_3(x,\xi)$ to denote $f(x, \xi)-\E_{\xi}[f(x, \xi)]$, 
\[
    \epsilon_{\rm noise}(x,\mathcal{U},\xi) = \sum_{i=1}^s \frac{
        \epsilon_3(x+ {h} u_i, \xi_{i+}) - \epsilon_3(x - {h} u_i, \xi_{i-}) 
    }{2{h}} u_i.  
\]
As in the forward finite difference case, it holds that $\E_{\xi }[\epsilon_{\rm noise}(x,\mathcal{U},\xi)] = 0$, and we need to derive bounds for $\epsilon_{\rm cfd}(x,\mathcal{U})$ and $\epsilon_{\rm noise}(x,\mathcal{U},\xi)$. 

\paragraph{Bounding $\epsilon_{\rm cfd}$} 
We first note, with $-\mathcal{U} = \{-u_i\}_{i=1}^s$, we have $\epsilon_{\rm cfd}(x,\mathcal{U}) = [\epsilon_{\rm ffd}(x,\mathcal{U}) + \epsilon_{\rm ffd}(x,-\mathcal{U})] / 2$ and that $-\cal U$ satisfies Assumption~\ref{ass:haar} if $\cal U$ does. 
Thus, the bound \eqref{eq:ortho eps_ffd bound}, which holds for $\epsilon_{\rm ffd}(x,\mathcal{U})$, also holds for $\epsilon_{\rm cfd}(x,\mathcal{U})$. 
The advantage of the central finite difference schemes manifests itself under a stronger smoothness assumption on $\phi(x)$ than assumed for the FFD schemes. 
Specifically, $\epsilon_{\rm cfd}(x,\mathcal{U})$ also satisfies a different bound under Assumption~\ref{assum:lip_cont2}, which we derive below. 
Let $\epsilon_1(x,{h} u_i)=\phi(x+ {h} u_i) - \phi(x) - \nabla\phi(x)^T ({h} u_i) - \frac{1}{2}  ({h} u_i)^T \nabla^2 \phi(x)  ({h} u_i)$, which follows $|\epsilon_1(x,{h} u_i)| \le L_2 {h}^3 / 6$ under Assumption~\ref{assum:lip_cont2}; and let $\epsilon_2(x,{h} u_i)=\E_{\xi }[f(x+ {h} u_i, \xi_{i+})]-\phi(x+{h}u_i) - \E_{\xi }[f(x- {h} u_i, \xi_{i-})] + \phi(x- {h} u_i)$, which follows $|\epsilon_2(x,{h} u_i)| \le 2\epsilon_f$ when using Oracles~\ref{def:oracle_stoch_0.1} or ~\ref{def:oracle_stoch_0.2}. 
Then, under Assumption~\ref{ass:haar}, we have 
\begin{align*}
    \E_{\cal U} &\left[ \|\epsilon_{\rm cfd}(x,\mathcal{U})\|^2 \right]
    = \E_{\cal U}\left[ \sum_{i=1}^s\left\|\left(\frac{\epsilon_1(x,{h} u_i) - \epsilon_1(x,- {h} u_i)}{2{h}}+ \frac{\epsilon_2(x,{h} u_i)}{2{h}} \right) u_i \right\|^2 \right] \\
    &=\sum_{i=1}^s \E_{\cal U}\left[ \left(\frac{\epsilon_1(x,{h} u_i) - \epsilon_1(x,- {h} u_i)}{2{h}}+ \frac{\epsilon_2(x,{h} u_i)}{2{h}} \right)^2 \right]
    \le s \Big( \frac{L_2 {h}^2}{6} + \frac{\epsilon_f}{h} \Big)^2. 
\end{align*}

\paragraph{Bounding $\epsilon_{\rm noise}$} 
Again let $\epsilon_3(x,\xi)=f(x,\xi)-\E_{\xi}[f(x,\xi)]$. 
When Oracle~\ref{def:oracle_stoch_0.1} is used, we can use the same logic for bounding $\epsilon_{\rm noise}$ in Section~\ref{sec:RFFD_orth} to obtain  
\[
    \E_{\xi} \left[ 
        \big\| \epsilon_{\rm noise}(x,\mathcal{U},\xi) \big\|^2
    \right] 
    \leq \frac{s\sigma_f^2}{2 {h}^2}. 
\]
This bound is a quarter of what is for forward finite difference in Section~\ref{sec:RFFD_orth}, because the definition of $\epsilon_{\rm noise}(x,\mathcal{U})$ for central finite difference has an extra 2 in the denominator. 
Alternatively, if Oracle~\ref{def:oracle_stoch_0.2} is used, we have that 
\[ \begin{aligned}
    \E_{\xi} &\left[ \big\| 
        \epsilon_{\rm noise}(x,\mathcal{U},\xi) 
    \big\|^2 \right] 
    = \E_{\xi} \bigg[ \Big\| 
        \sum_{i=1}^s\frac{ \epsilon_3(x+ h u_i, \xi) - \epsilon_3(x - h u_i, \xi) }{2 h} u_i 
    \Big\|^2 \bigg] \\
    &= \frac{1}{4h^2} \sum_{i=1}^s \E_{\xi} \left[ \big( 
        \epsilon_3(x+ h u_i, \xi) - \epsilon_3(x, \xi) + \epsilon_3(x, \xi) - \epsilon_3(x - h u_i, \xi) 
    \big)^2 \right] \\
    &\le \frac{1}{4h^2} \sum_{i=1}^s \E_{\xi} \left[ 
        2\big( \epsilon_3(x+ h u_i, \xi) - \epsilon_3(x, \xi) \big)^2 
        + 2\big( \epsilon_3(x, \xi) - \epsilon_3(x - h u_i, \xi) \big)^2
    \right] \\
    &\le sL_c^2+sL_v^2\|\nabla \phi(x)\|^2. 
\end{aligned} \] 

We note that $f(x)$ is not used in the \eqref{eq:RCFD} scheme, thus it does not need to be computed, but we consider $x$ as the input to Oracle~\ref{def:oracle_stoch_0.2}
for the sake of centering the upper bound around $x$ and getting the dependence on $\|\nabla \phi(x)\|^2$ rather than $\|\nabla \phi(\cdot)\|^2$ 
measured at a different point around $x$, which will slightly complicate the resulting bounds. 

\paragraph{As Oracle~\ref{def:oracle_stoch_1.3}.} 
Following the same derivation for the forward finite difference, we can conclude that, under Assumptions~\ref{assum:lip_cont}, \ref{assum:lip_cont2}, and \ref{ass:haar}, the central finite difference scheme \eqref{eq:RCFD} satisfies the conditions of Oracle~\ref{def:oracle_stoch_1.3} with 
\[ 
    \epsilon_g^2 =  s \cdot \min \left\{\Big( \frac{L_1 {h}}{2} + \frac{2\epsilon_f}{h} \Big)^2, \Big( \frac{L_2 {h}^2}{6} + \frac{\epsilon_f}{h} \Big)^2\right\}
    \quad \text{ and } \quad 
    \kappa_g^2 = \frac{n-s}{n}, 
\] 
while $\sigma_c^2 = s\sigma_f^2 / (2{h}^2)$, $\sigma_v = 0$ when using Oracle~\ref{def:oracle_stoch_0.1}, or $\sigma_c^2 = sL_c^2$, $\sigma_v = sL_v^2$ when using Oracle~\ref{def:oracle_stoch_0.2}. 

\subsubsection{Scaled Randomized Central Finite Differences with Orthogonal Directions} \label{sec:SRCFD}
The central finite difference gradient estimators \eqref{eq:RCFD} are also often used with the scaling factor $n/s$: 
\begin{equation} \label{eq:SRCFD} \tag{SRCFD}
	g(x,\mathcal{U},\xi) = \frac{n}{s} \sum_{i=1}^s \frac{f(x+ {h} u_i, \xi_{i+}) - f(x- {h} u_i,\xi_{i-})}{2{h}} u_i, 
\end{equation} 
which can be decomposed as 
\begin{equation*} 
   g(x,\mathcal{U}) = \nabla\phi(x) + \epsilon_{\rm proj \perp}(x,\mathcal{U}) + \epsilon_{\rm proj \parallel}(x,\mathcal{U}) + \epsilon_{\rm cfd}(x,\mathcal{U}) + \epsilon_{\rm noise}(x,\mathcal{U}), 
\end{equation*}
where $\epsilon_{\rm proj \perp}(x,\mathcal{U})$ and $\epsilon_{\rm proj \parallel}(x,\mathcal{U})$ are the same as in \eqref{eq:FFD decomposition}, while $\epsilon_{\rm cfd}(x,\mathcal{U})$ and $\epsilon_{\rm noise}(x,\mathcal{U})$ are the same as in \eqref{eq:CFD decomposition} except with the extra coefficient $n/s$. 
Most bounds on the error terms in this decomposition are readily available from previous sections, except that under Assumption~\ref{assum:lip_cont2}
\[
    \|\E_{\cal U}[\epsilon_{\rm cfd}(x,\mathcal{U})]\|^2 
    \le n \Big(\frac{L_2 h^2}{6} + \frac{\epsilon_f}{h}\Big)^2, 
\]
which can be proved analogously to \eqref{eq:||E eps_ffd||}. 
Do not forget the bounds on $\epsilon_{\rm ffd}$ also hold for $\epsilon_{\rm cfd}$ because $\epsilon_{\rm cfd}(x,\mathcal{U}) = [\epsilon_{\rm ffd}(x,\mathcal{U}) + \epsilon_{\rm ffd}(x,-\mathcal{U})] / 2$ and that $-\cal U$ satisfies Assumption~\ref{ass:haar} if $\cal U$ does. 
For example, we can use \eqref{eq:||E eps_ffd||} to show 
\[ \begin{aligned} 
    \|&\E_{\cal U}[\epsilon_{\rm cfd}(x,\mathcal{U})]\|^2 
    = \frac{1}{4} \|\E_{\cal U}[\epsilon_{\rm ffd}(x,\mathcal{U}) + \epsilon_{\rm ffd}(x,-\mathcal{U})]\|^2 \\
    &\le \frac{1}{4} \Big( 2\|\E_{\cal U}[\epsilon_{\rm ffd}(x,\mathcal{U})]\|^2 + 2\|\E_{\cal U}[\epsilon_{\rm ffd}(x,-\mathcal{U})]\|^2 \Big) 
    \le n \Big(\frac{L_1 {h}}{2} + \frac{2\epsilon_f}{h}\Big)^2 
\end{aligned} \]
when $\phi$ satisfies Assumption~\ref{assum:lip_cont}. 

In summary, \eqref{eq:SRCFD} satisfies the conditions of Oracle~\ref{def:oracle_stoch_1.2} with 
\[ \begin{aligned}
    \epsilon_g^2 = n \cdot \min \left\{\Big( \frac{L_1 {h}}{2} + \frac{2\epsilon_f}{h} \Big)^2, \Big( \frac{L_2 {h}^2}{6} + \frac{\epsilon_f}{h} \Big)^2\right\}, \quad 
    \kappa_g^2 = 0, 
\end{aligned} \]
while 
\[
    \sigma_c^2 = \frac{n^2\sigma_f^2}{2s h^2} + \frac{n^2}{s} \min \left\{\Big( \frac{L_1 {h}}{2} + \frac{2\epsilon_f}{h} \Big)^2, \Big( \frac{L_2 {h}^2}{6} + \frac{\epsilon_f}{h} \Big)^2\right\}, \quad 
    \sigma_v^2 = \frac{(2n-s)(n-s)}{ns}
\]
when using Oracle~\ref{def:oracle_stoch_0.1}, or 
\[
    \sigma_c^2 = \frac{n^2 L_c^2}{s} + \frac{n^2}{s} \min \left\{\Big( \frac{L_1 {h}}{2} + \frac{2\epsilon_f}{h} \Big)^2, \Big( \frac{L_2 {h}^2}{6} + \frac{\epsilon_f}{h} \Big)^2\right\}, \quad 
    \sigma_v^2 = \frac{n^2 L_v^2}{s} + \frac{(2n-s)(n-s)}{ns}
\]
when using Oracle~\ref{def:oracle_stoch_0.2}.

\subsubsection{View of randomized finite differences via smoothing}\label{sec:smoothing}
\paragraph{Function Smoothing and Gradient Approximation}
When $\cal U$ are i.i.d. random vectors following $\cal N$ or $\cal S$, the randomized finite difference schemes \eqref{eq:SRFFD} and \eqref{eq:SRCFD}  are often referred to as  smoothing schemes. This is because the expectation of the gradient estimators computed by \eqref{eq:SRFFD} and \eqref{eq:SRCFD}  form a gradient of a smoothing (convolution) of the following function
\[
    f(x) = \E_{\xi^{(0)}} [f(x,\xi^{(0)})],
\]
which is the expected output of a zeroth-order oracle of $\phi$ at $x$. We note that is this zeroth order oracle is unbiased, then $f(x)=\phi(x)$, but otherwise it can be a different function, possibly nonsmooth and even discontinuous.  

The expected value of \eqref{eq:SRFFD} is equal to 
\[ \begin{aligned} 
    \E_{\cal U,\xi} [g(x,{\cal U},\xi)]
    &= \E_{\cal U,\xi} \Big[\frac{n}{s} \sum_{i=1}^s \frac{f(x+ h u_i, \xi_i) - f(x,\xi_0)}{h} u_i\Big] \\
    &= \E_{\cal U} \Big[\frac{n}{s} \sum_{i=1}^s \frac{f(x+ h u_i) - f(x)}{h} u_i\Big] \\
    &= n \E_{u} \Big[ \frac{f(x+hu) - f(x)}{h} u \Big]
    = \frac{n}{h} \E_{u} [f(x+hu) u]. 
\end{aligned} \]
Similarly, since $-u$ follows the same distribution as $u$, the expected value of \eqref{eq:SRCFD} is equal to 
\[ \begin{aligned} 
    \E_{\cal U,\xi} [g(x,{\cal U},\xi)]
    &= \E_{\cal U,\xi} \Big[\frac{n}{s} \sum_{i=1}^s \frac{f(x+ {h} u_i, \xi_{i+}) - f(x- {h} u_i,\xi_{i-})}{2{h}} u_i \Big] \\
    &= \E_{\cal U} \Big[\frac{n}{s} \sum_{i=1}^s \frac{f(x+ h u_i) - f(x-hu_i)}{2h} u_i\Big] \\
    &= n \E_{u} \Big[ \frac{f(x+hu) - f(x-hu)}{2h} u \Big]
    = \frac{n}{h} \E_{u} [f(x+hu) u]. 
\end{aligned} \]
Now consider the convolution function 
\[
    f_{{\cal B},h}(x) := \E_{u\sim \cal B} [f(x+hu)], 
\]
where $\cal B$ means the uniform distribution over the unit ball $\{x\in\mathbb{R}^n: \|x\| \le 1\}$, as opposed to $\cal S$, which is the uniform distribution over the unit sphere $\{x\in\mathbb{R}^n: \|x\| = 1\}$. 
This function can be seen as a smoothed version of $f$, and its gradient is shown in \cite{flaxman2005online} to be equal to 
\[
    \nabla f_{{\cal B},h}(x) = \frac{n}{h} \E_{u \sim \cal S} [f(x+hu) u] 
\] 
under some mild conditions on $f$. 
This means the expectation $\E_{\cal U,\xi} [g(x,{\cal U},\xi)]$ of \eqref{eq:SRFFD} and \eqref{eq:SRCFD} are both the gradient of the smoothed function $f_{{\cal B},h}(x)$, and \eqref{eq:SRFFD} and \eqref{eq:SRCFD} themselves can be seen as sample average approximations of $\nabla f_{{\cal B},h}(x)$. 

Similarly, it is shown in \cite{nesterov2017random} that 
\[
    \nabla f_{{\cal N}(0,I),h}(x)
    = \E_{u \sim \mathcal{N}(0,I)} \Big[ \frac{f(x+hu) - f(x)}{h} u \Big], 
\] 
where 
\[
    f_{{\cal N}(0,I),h}(x) = \E_{u \sim \mathcal{N}(0,I)} [f(x+hu)]. 
\]
The Gaussian distribution here is standard. 
It is not the same as the one we use, which is ${\cal N}(0, \frac{1}{n} I)$ and simply denoted as $\cal N$. 
With our Gaussian distribution, let 
\[
    f_{{\cal N},h}(x) := \E_{u\sim \cal N} [f(x+hu)], 
\]
then the above equation changes to 
\[
    \nabla f_{{\cal N},h}(x)
    = n \E_{u \sim \mathcal{N}} \Big[ \frac{f(x+hu) - f(x)}{h} u \Big], 
\] 
again showing \eqref{eq:SRFFD}and \eqref{eq:SRCFD} are sample average approximations of the gradient of a smoothed version of $f$. 

\paragraph{Bias of the Gradient of the Smoothed Function}
Under both distributions $\cal N$ and $\cal S$, the difference between the gradient of the smoothed function $\nabla f_{\cdot,h}(x)$ and $\nabla\phi(x)$ is equal to the expectation of the finite difference error term used in the previous sections, which in turn is equal to the bias of the first-order oracle. 
That is 
\begin{equation*} \begin{aligned}
    \|\nabla f_{\cdot,h}(x) - \nabla\phi(x)\|
    = \|\E_{\cal U} [\epsilon_{\rm ffd}(x,\mathcal{U})]\|
    = \|\E_{\cal U,\xi} [g(x,{\cal U},\xi)] - \nabla\phi(x)\| 
\end{aligned} \end{equation*}
for \eqref{eq:SRFFD}, where $\epsilon_{\rm ffd}(x,\mathcal{U})$ is defined in \eqref{eq:FFD decomposition}, and 
\begin{equation*} \begin{aligned}
    \|\nabla f_{\cdot,h}(x) - \nabla\phi(x)\|
    = \|\E_{\cal U} [\epsilon_{\rm cfd}(x,\mathcal{U})]\|
    = \|\E_{\cal U,\xi} [g(x,{\cal U},\xi)] - \nabla\phi(x)\| 
\end{aligned} \end{equation*}
for \eqref{eq:SRCFD}, where $\epsilon_{\rm cfd}(x,\mathcal{U})$ is defined in \eqref{eq:CFD decomposition}. 

Recall that by  $\cal N$ we denote  Gaussian distribution ${\cal N}(0, \frac{1}{n}I)$. In instead a standard  Gaussian distribution ${\cal N(0, I)}$ is used and the zeroth-order oracle is exact, it is shown in \cite{nesterov2017random} that 
\begin{equation} \label{eq:Nesterov bias} 
    \|\nabla f_{{\cal N}(0,I),h}(x) - \nabla\phi(x)\| \le \frac{L_1 h}{2} (n+3)^{3/2} 
\end{equation}
when $\phi$ follows Assumption~\ref{assum:lip_cont} and 
\begin{equation}\label{eq:Nesterov bias 2} 
    \|\nabla f_{{\cal N}(0,I),h}(x) - \nabla\phi(x)\| \le \frac{L_2 h^2}{6} (n+4)^2 
\end{equation}
when $\phi$ follows Assumption~\ref{assum:lip_cont2}. 
Improved bounds can also be found in \cite{berahas2022theoretical}, where $f(x)$ is assumed to be the output of an inexact deterministic zeroth-order oracle that satisfies $|f(x) - \phi(x)| \le \epsilon_f$ for any $x\in\mathbb{R}^n$. 
This exactly fits  our Oracle~\ref{def:oracle_det_0}, but   also can be viewed as a special case of  either stochastic Oracles~\ref{def:oracle_stoch_0.1} or \ref{def:oracle_stoch_0.2}. 
It is shown in \cite{berahas2022theoretical} that under Assumption~\ref{assum:lip_cont} 
\begin{align} 
	\|\nabla f_{{\cal N}(0,I),h}(x) - \nabla\phi(x)\| &\le \sqrt{n} L_1 h +  \frac{\sqrt{n} \epsilon_f}{h}, \label{eq:GSG_bound1} \\
	\|\nabla f_{{\cal B},h}(x) - \nabla \phi(x)\| &\le L_1h  + \frac{n \epsilon_f}{h}, \label{eq:BSG_bound1}
\end{align} 
and under Assumption~\ref{assum:lip_cont2} 
\begin{align} 
	\|\nabla f_{{\cal N}(0,I),h}(x) - \nabla\phi(x)\| &\le nL_2h^2 + \frac{\sqrt{n}\epsilon_f}{h}, \label{eq:GSG_bound2} \\
	\|\nabla f_{{\cal B},h}(x) - \nabla \phi(x)\| &\le L_2h^2  + \frac{n \epsilon_f}{h}. \label{eq:BSG_bound2}
\end{align}
A bound on the bias of \eqref{eq:SRFFD} 
\[
    \|\E_{\cal U,\xi} [g(x,{\cal U},\xi)] - \nabla\phi(x)\| 
    \le \frac{\sqrt{n}(1 + {\cal O}(n^{-1})) L_1 {h}}{2} + \frac{2\sqrt{n} \epsilon_f}{h}
\]
is derived in Section~\ref{sec:SRFFD} (as $\epsilon_f$ of Oracle~\ref{def:oracle_stoch_1.2}) for the case when $\cal U$ satisfies Assumption~\ref{ass:normal}. 
Notice $\E_{\cal U,\xi} [g(x,{\cal U},\xi)] = \nabla f_{{\cal N},h}(x)$ when $\cal U \sim \cal N$, and $\E_{\cal U,\xi} [g(x,{\cal U},\xi)] = \nabla f_{{\cal B},h}(x)$ when $\cal U \sim \cal S$. 
Using the same derivation, we can conclude that 
\begin{align}
    \|\nabla f_{{\cal N},h}(x) - \nabla\phi(x)\| 
    &\le \frac{\sqrt{n} L_1 {h}}{2} + \frac{2\sqrt{n} \epsilon_f}{h}, \label{eq:our bias normal} \\
    \text{ and } \|\E_{\cal U,\xi} [g(x,{\cal S},\xi)] - \nabla\phi(x)\| 
    &\le \frac{\sqrt{n} L_1 {h}}{2} + \frac{2\sqrt{n} \epsilon_f}{h}. \label{eq:our bias sphere}
\end{align} 

These three sets of bounds on the difference between the gradient of the smoothed function and $\nabla\phi(x)$ (\eqref{eq:Nesterov bias} and \eqref{eq:Nesterov bias 2} from \cite{nesterov2017random}, \eqref{eq:GSG_bound1} - \eqref{eq:BSG_bound2} from \cite{berahas2022theoretical}, \eqref{eq:our bias normal} and \eqref{eq:our bias sphere} from Section~\ref{sec:SRFFD} of this paper) can all be improved.

\paragraph{Updated Bounds under Assumption~\ref{assum:lip_cont}}
First assume $\cal U \sim \cal N$. 
Starting with the triangular inequality, we have 
\begin{equation} \label{eq:smoothing triangular ineq} \begin{aligned}
    \|&\nabla f_{{\cal N},h}(x)  - \nabla\phi(x)\|
    = \left\| \frac{n}{h} \E_{u\sim \mathcal{N}} \Big[f(x+hu)  u\Big] - \nabla\phi(x)\right\| \\
    &\le \left\| \frac{n}{h} \E_{u\sim \mathcal{N}} \Big[\phi(x+hu)  u\Big] - \nabla\phi(x)\right\| + \left\| \frac{n}{h} \E_{u\sim \mathcal{N}} \Big[\big(f(x+hu) - \phi(x+hu)\big) u\Big]\right\|. 
\end{aligned} \end{equation}
The first term 
\[ \begin{aligned}
    \left\| \frac{n}{h} \E_{u\sim \mathcal{N}} \Big[\phi(x+hu)  u\Big] - \nabla\phi(x)\right\| 
    = \left\| \E_{u\sim \mathcal{N}} [\nabla\phi(x+hu)] - \nabla\phi(x)\right\|& \\
    \le \E_{u\sim \mathcal{N}} \left[\| \nabla\phi(x+hu) - \nabla\phi(x)\|\right]
    \le L_1 h \E_{u\sim \mathcal{N}} \left[\|u\|\right]
    \le L_1 h&,
\end{aligned} \]
where the first equality holds because, using the smoothing argument, 
\[  
    \frac{n}{h} \E_{u \sim \cal N} [\phi(x+hu) u] 
    = n \E_{u \sim \mathcal{N}} \Big[ \frac{\phi(x+hu) - \phi(x)}{h} u \Big]
    = \nabla \E_{u \sim \cal N} [\phi(x+hu)] 
\] 
and differentiation and integration can be exchanged when $\phi$ satisfies Assumption~\ref{assum:lip_cont}: 
\[
    \nabla \bigl(\E_{u\sim \mathcal{N}} [\phi(x+hu)] \bigr)
    = \E_{u\sim \mathcal{N}} [\nabla\phi(x+hu)]. 
\]
Using the Cauchy-Schwarz inequality, the second term of \eqref{eq:smoothing triangular ineq} satisfies 
\[ \begin{aligned}
    &\left\| \frac{n}{h} \E_{u\sim \mathcal{N}} \Big[\big(f(x+hu) - \phi(x+hu)\big) u\Big]\right\| \\
    &= \frac{n}{h} \max_{v\in\mathbb{R}^n:~ \|v\|=1} \Big\{\E_{u\sim \mathcal{N}} \Big[\big(f(x+hu) - \phi(x+hu)\big) u^Tv\Big]\Big\} \\
    &\le \frac{n}{h} \max_{v\in\mathbb{R}^n:~ \|v\|=1} \Big\{ \sqrt{\E_{u\sim \mathcal{N}} \Big[\big(f(x+hu) - \phi(x+hu)\big)^2\Big] \E_{u\sim \mathcal{N}} [(u^Tv)^2]} \Big\} \\
    &= \frac{n}{h} \sqrt{\E_{u\sim \mathcal{N}} \Big[\big(f(x+hu) - \phi(x+hu)\big)^2\Big]} \max_{v\in\mathbb{R}^n:~ \|v\|=1} \Big\{ \sqrt{v^T \E_{u\sim \mathcal{N}} [uu^T] v} \Big\} \\
    &\le \frac{\sqrt{n}\epsilon_f}{h}. 
\end{aligned} \]

Thus, when $\phi$ satisfies Assumption~\ref{assum:lip_cont} and the zeroth-order oracle satisfies $|f(x) - \phi(x)| \le \epsilon_f$ for all $x\in\mathbb{R}^n$ (holds for Oracles~\ref{def:oracle_det_0}, \ref{def:oracle_stoch_0.1}, and \ref{def:oracle_stoch_0.2}), the bias of \eqref{eq:SRFFD} under $\mathcal{U} \sim \mathcal{N}$ has a better bound 
\[ \begin{aligned}
    \|\E_{{\cal U}, \xi}[g(x,\mathcal{U},\xi)] - \nabla\phi(x)\|
    = \|\nabla f_{{\cal N},h}(x)  - \nabla\phi(x)\|
    \le L_1 h + \frac{\sqrt{n}\epsilon_f}{h}. 
\end{aligned} \] 
With an almost identical proof, one can show the same bound also holds for the case when $\mathcal{U} \sim \mathcal{S}$: 
\[ \begin{aligned}
    \|\E_{{\cal U}, \xi}[g(x,\mathcal{U},\xi)] - \nabla\phi(x)\|
    = \|\nabla f_{{\cal B},h}(x)  - \nabla\phi(x)\|
    \le L_1 h + \frac{\sqrt{n}\epsilon_f}{h}. 
\end{aligned} \] 
Thus, the constant $\epsilon_g$ for Oracle~\ref{def:oracle_stoch_1.2} in Section~\ref{sec:SRFFD} can be updated 
\[
    \text{ from } \quad 
    \frac{\sqrt{n} L_1 h}{2} + \frac{2\sqrt{n}\epsilon_f}{h}
    \quad \text{ to } \quad 
    L_1 h + \frac{\sqrt{n}\epsilon_f}{h} 
\]
when $\cal U \sim \cal N$ or $\cal S$. 

For standard Gaussian ${\cal N}(0, \frac{1}{n}I)$, slightly modifying  the proof one can also show 
\begin{equation*}
	\|\nabla f_{{\cal N}(0,I),h}(x) - \nabla\phi(x)\| \le \sqrt{n} L_1 h +  \frac{\epsilon_f}{h}. 
\end{equation*}

\paragraph{Updated Bounds under Assumption~\ref{assum:lip_cont2}}
If $\phi$ satisfies Assumption~\ref{assum:lip_cont2}, the first term in \eqref{eq:smoothing triangular ineq} instead follows 
\[ \begin{aligned}
    &\left\| \frac{n}{h} \E_{u\sim \mathcal{N}} \Big[\phi(x+hu)  u\Big] - \nabla\phi(x)\right\| 
    = \left\| \E_{u\sim \mathcal{N}} [\nabla\phi(x+hu)] - \nabla\phi(x)\right\| \\
    &= \left\| \frac{1}{2}\E_{u\sim \mathcal{N}} [\nabla\phi(x+hu) + \nabla\phi(x-hu)] - \nabla\phi(x)\right\| \\
    &= \frac{1}{2} \left\| \E_{u\sim \mathcal{N}} [\nabla\phi(x+hu) - \nabla\phi(x) - \nabla^2\phi(x) (hu)] + \E_{u\sim \mathcal{N}} [\nabla\phi(x-hu) - \nabla\phi(x)- \nabla^2\phi(x) (-hu)]\right\| \\
    &\le \frac{1}{2} \E_{u\sim \mathcal{N}} [\|\nabla\phi(x+hu) - \nabla\phi(x) - \nabla^2\phi(x) (hu)\|] + \frac{1}{2} \E_{u\sim \mathcal{N}} [\|\nabla\phi(x-hu) - \nabla\phi(x)- \nabla^2\phi(x) (-hu)\|] \\
    &\le \frac{1}{2} \cdot \E_{u\sim \mathcal{N}} \Big[\frac{L_2}{2} \|hu\|^2\Big] + \frac{1}{2} \cdot \E_{u\sim \mathcal{N}} \Big[\frac{L_2}{2} \|hu\|^2\Big]
    = \frac{L_2h^2}{2}. 
\end{aligned} \]
Thus, under Assumption~\ref{assum:lip_cont2}, we also have 
\[ \begin{aligned}
    \|\E_{{\cal U}, \xi}[g(x,\mathcal{U},\xi)] - \nabla\phi(x)\|
    = \|\nabla f_{{\cal N},h}(x)  - \nabla\phi(x)\|
    \le \frac{L_2h^2}{2} + \frac{\sqrt{n}\epsilon_f}{h}. 
\end{aligned} \] 
Again, with an almost identical proof, we have 
\[ \begin{aligned}
    \|\E_{{\cal U}, \xi}[g(x,\mathcal{U},\xi)] - \nabla\phi(x)\|
    = \|\nabla f_{{\cal B},h}(x)  - \nabla\phi(x)\|
    \le \frac{L_2h^2}{2} + \frac{\sqrt{n}\epsilon_f}{h}. 
\end{aligned} \] 
When $\phi$ satisfies both Assumptions~\ref{assum:lip_cont} and \ref{assum:lip_cont2}, the constant $\epsilon_g$ for Oracle~\ref{def:oracle_stoch_1.2} in Section~\ref{sec:SRCFD} can be updated from
\[
    \min \left\{\frac{\sqrt{n}L_1 {h}}{2} + \frac{2\sqrt{n}\epsilon_f}{h}, \frac{\sqrt{n}L_2 {h}^2}{6} + \frac{\sqrt{n} \epsilon_f}{h}\right\}
\]
to 
\[
    \min \left\{L_1 {h} + \frac{\sqrt{n}\epsilon_f}{h}, \frac{L_2 {h}^2}{2} + \frac{\sqrt{n} \epsilon_f}{h}\right\}
\]

Moreover, for the standard Gaussian distribution ${\cal N}(0, \frac{1}{n}I)$, we can also use a similar derivation to obtain 
\begin{equation*}
	\|\nabla f_{{\cal N}(0,I),h}(x) - \nabla\phi(x)\| \le \frac{nL_2h^2}{2} +  \frac{\epsilon_f}{h}. 
\end{equation*}

\subsection{Other stochastic first-order oracle conditions} 
There is a large variety of conditions that have been proposed  in the literature for stochastic first-order oracles. 
They often include the assumption that the oracle is unbiased and impose various bounds on the second moment of the gradient estimator. 

For example one of the earliest and simplest assumptions
 from \cite{ghadimi2013stochastic}, bounds the variance of the unbiased estimator  by a constant which is a special case of
 our  Oracle~\ref{def:oracle_stoch_1.2} with $\kappa_g=\epsilon_g=\sigma_v=0$. 

The previously mentioned strong growth condition studied in \cite{vaswani2019fast} can also be seen as a special case of 
our Oracle~\ref{def:oracle_stoch_1.2} with $\kappa_g=\epsilon_g=\sigma_c=0$.  
Our Oracle~\ref{def:oracle_stoch_1.2}  also generalizes conditions in  \cite{bottou2018optimization} where they essentially impose similar bounds on the variance, allowing $\sigma_v$ and $\sigma_c$ to be positive, but also allow for relative bias which is imposed via the angle condition discussed in Section \ref{sec:angle} and  which as we have shown can be reduced to our oracle conditions.  These are examples of conditions that are special cases of Oracle~\ref{def:oracle_stoch_1.2}.  

\paragraph{The ABC Conditions} 
There are conditions that do not fit into our oracle definitions. Most notably a more general condition on the 
 second  moment for unbiased gradient estimators can be found in \cite{khaled2022better}. The conditions is stated as follows
\begin{assumption} \label{assum:ABC}
    The stochastic gradient estimator $g(x)$ is unbiased, and its second moment satisfies 
    \begin{equation} \label{eq:ABC}
        \mathbb{E}[ \|g(x)\|^2 ] \le 2A (\phi(x) - \hat\phi) + B \|\nabla\phi(x)\|^2 + C
    \end{equation}
    for some $A,B,C \ge 0$ and all $x \in \mathbb{R}^n$.
\end{assumption}

In the case of unbiased oracles, our Oracle~\ref{def:oracle_stoch_1.2}  is a special case of \eqref{eq:ABC} with $A=0$, $B=1+\sigma_v^2$ and $C=\sigma_c^2$. 

The inclusion of the extra term when $A>0$ is meant to subsume some other cases of oracle conditions. In most cases $\phi(x)$ is assumed to be convex and  $\hat\phi$ is a global minimum value $\phi^*$.  For example,
the {\em weak growth condition} from \cite{vaswani2019fast} is one example of Assumption~\ref{assum:ABC} with $B=C=0$  and $A>0$.   
A case of Assumption~\ref{assum:ABC} with $B=0$ was used in \cite{gower2021stochastic}  to establish the convergence rate of  stochastic subgradient methods when applied to convex non-Lipschitz functions. 
A condition  for convex functions  called {\it convex expected smoothness}  is studied in \cite{khaled2022better} and   \cite{gower2021stochastic} and is stated as 
\begin{equation}
    \mathbb{E}[ \|g(x) - g(x^\star)\|^2 ] \le 2A (\phi(x) - \phi(x^\star)), 
\end{equation}
where $x^\star$ is a global optimizer of $\phi$, and both stochastic gradients share the same source of randomness through random sampling. 

Assumption~\ref{assum:ABC} has been shown to hold with $A>0$ in several scenarios.  
For instance, the problem of solving a consistent linear system was reformulated in \cite{richtarik2020stochastic} as a stochastic optimization whose objective function follows \eqref{eq:ABC} with equality and $(A,B,C) = (1,0,0)$. 
Multiple subsampling strategies for computing the stochastic gradient of a finite sum problem have been shown to satisfy the ABC assumption in \cite{khaled2022better}, and the corresponding values for $(A,B,C)$ were provided.  Some of the analysis in \cite{khaled2022better} for when $A>0$ recovers the same results (or special cases thereof) that we will present later. The results with $A>0$ are not as general,  do not seem to provide significant improvement  compared to results and also do not lead to as general analysis as what we aim for in this article.  Moreover, the cases of oracles that satisfy Assumption~\ref{assum:ABC} with $A>0$ but that cannot be stated in terms of Oracle~\ref{def:oracle_stoch_1.2} are relatively rare and specialised  \cite{khaled2022better}. 
Clearly for functions satisfying Assumption \ref{assum:PLcond} globally the ABC condition~\ref{eq:ABC}  implies conditions of Oracle~\ref{def:oracle_stoch_1.2} since the function value gap can be bounded by the squared gradient norm. For all those reasons we do not consider this extended oracle condition.

A biased version of the ABC assumption was proposed in \cite{demidovich2023guide} where the conditions on the bias are 
based on the angle condition, similar to  \cite{bottou2018optimization}  but include an absolute error term. 

Other related conditions but for an unbiased estimator were proposed in  \cite{fatkhullin2022sharp}, and a bound on the expected optimality gap  was derived under Assumptions~\ref{assum:lip_cont} and \ref{assum:KLcond}.

\paragraph{Stochastic quasi-Newton directions}
A stochastic analogue of the angle condition discussed in Section~\ref{sec:angle} was used in \cite{bottou2018optimization}. It is particularly useful when 
 the stochastic gradient is unbiased and the matrix defining the search direction is uniformly positive definite. In particular, consider
\begin{align*}
    d=-B g,
    \qquad
    \beta_1 I\preceq B\preceq\beta_2 I,
\end{align*}
where $g$ is an unbiased instance of Oracle~\ref{def:oracle_stoch_1.1}. If $B$ is independent of the current gradient sample, conditional on the iteration history, then
\begin{align*}
    \E\left[
        \nabla\phi(x)^Td
        \mid x
    \right]
    \leq
    -\beta_1\|\nabla\phi(x)\|^2.
\end{align*}
Thus, uniformly positive-definite curvature approximations preserve a descent property in expectation. This framework arises naturally in stochastic quasi-Newton methods, where positive-definite approximations are constructed using stochastic gradients or subsampled curvature information \cite{byrd2016stochastic,berahas2016multibatch,curtis2016selfcorrecting}.

Our analysis in Sections \ref{sec:angle} and Section \ref{sec:quasi-newton} that connect the general descent directions to our oracles can be extended to the stochastic case analogously.

\subsection{Expected Descent Inequality and Expected Complexity Analysis}
We now turn to deriving expected complexity bounds for a general algorithm based only on an expected descent inequality. 
In other words we extend results of Section \ref{sec:descent_lemma} to the stochastic case. We observe that there are two types of results that can be derived in the stochastic case: the expected rate of convergence and expected complexity. While the former is more common in the stochastic optimization literature, we find the latter more instructive and better aligned with the deterministic case, especially in cases when oracles are biased. We will present both types of results in this section but will focus more on the expected complexity in the later sections when we particularize the results to the various oracle cases.

Consider an iterative stochastic algorithm that, given a starting point $x_0$, generates a sequence of stochastic iterates $X_1, X_2, \dots$ in $\mathbb{R}^n$. 
For any $k\in\mathbb{N}$, let random variable $\xi_k$ represent all randomness in iteration $k\in\mathbb{N}$ and $\sigma$-algebra $\mathcal{F}_k = \sigma(\xi_0,\xi_1,\dots,\xi_{k-1})$ be the filtration representing all randomness generated in the optimization process prior to iteration $k$. 
We assume the expected descent inequality 
\begin{equation}\label{eq:key_stoch_reduction}\tag{EDI}
	\E[\phi(X_{k+1})|\, \mathcal{F}_k] \leq \phi (X_k) - \alpha \theta \| \nabla \phi (X_k) \|_2^2 + \alpha \omega +\alpha^2\nu, 
\end{equation}
holds at every iteration $k\in\mathbb{N}$ for some $\alpha>0$, $\theta >0 $, $\omega\geq 0$, and $\nu\geq 0$. 
This is a stochastic analogue of the descent result \eqref{eq:key_reduction} in Section~\ref{sec:descent_lemma} with an extra term $\alpha^2\nu$. 
We will show that Oracles~\ref{def:oracle_stoch_1.1},  \ref{def:oracle_stoch_1.2} and \ref{def:oracle_stoch_1.3} all give rise to algorithms that satisfy \eqref{eq:key_stoch_reduction}, but first we derive the expected rates and expected complexity bounds that follow from this inequality.

\subsubsection{Expected rate results}\label{sec:expected_rates}
Expected convergence rate results bound the rate of decrease of some measure of distance to optimality or criticality as a function of the number of iterations. 
Two such results are presented below. 
One regards the expectation of the squared norm of the gradient at $X_k$ averaged over all iterations, the other regards the expected function value at the last iteration. 
They are prevalent in the literature \cite{wright2022optimization}. 

\paragraph{Expected rate for general smooth functions.}
\begin{theorem} \label{thm:nonconvex_exp_rate}
    Let Assumption~\ref{assum:low_bound} hold. 
    If the expected descent inequality \eqref{eq:key_stoch_reduction} holds for the sequence $\{X_k\}$ on each $k$, then for any integer $K>0$  
    \begin{equation}\label{eq:nonconvex_exp_rate}
        \frac{1}{K} \sum_{k=0}^{K-1} \E[\|\nabla\phi(X_k)\|^2 ]\le \frac{1}{\theta} \left ( \frac{\phi(x_0) - \hat{\phi}}{\alpha K} + {\omega} +\alpha\nu \right ). 
    \end{equation}
\end{theorem}

\begin{proof} 
By summing \eqref{eq:key_stoch_reduction} over $k=0, \dots, K-1$, taking total expectation, and rearranging terms, we obtain 
    \begin{align*}
        \alpha \theta \sum_{k=0}^{K-1} \E\left[\|\nabla\phi(X_k)\|^2 \right]
        &\le \phi(x_0) - \E[\phi(X_K) ] + K \alpha\omega + K\alpha^2\nu \\
        &\le \phi(x_0) - \hat{\phi} + K \alpha\omega + K\alpha^2\nu. 
    \end{align*} 
    Dividing both sides by $\alpha\theta K$, we obtain the result. 
 \end{proof}

Let us consider the implication of this theorem on the type of convergence we might expect. 
The result provides the rate at which the average expected squared gradient norm diminishes. It is clear that this bound cannot guarantee that the average expected squared norm decreases below  $\frac{\omega}{ \theta}$. This  observation is essentially the same as in Theorem \ref{thm:nonconvex_bound} where the goal is to reach an iterate whose gradient  norm is not larger than $\epsilon_\tau \geq  \sqrt{\frac{\omega}{ \tau^2\theta}}$, for some $\tau\in(0,1)$. Thus, here we also can aim at 
achieving average expected squared norm for some $\epsilon_\tau^2 \geq  \frac{\omega}{ \tau^2\theta}$. Given such a goal 
we can select $\alpha \leq \min \left \{\frac{\theta \epsilon_\tau^2(1-\tau^2)}{2\nu}, \bar \alpha\right \}$ and $K\geq \frac{2(\phi(x_0) - \hat{\phi})}{(1-\tau^2)\theta \epsilon_\tau^2\alpha}$, where $\bar \alpha$ is some upper bound on $\alpha$ imposed 
when $\nu$ is small or zero. 

Plugging this into \eqref{eq:nonconvex_exp_rate} we get 
           \begin{equation*}
        \frac{1}{K} \sum_{k=0}^{K-1} E[\|\nabla\phi(X_k)\|^2 ]\le  \frac{\epsilon_\tau^2(1-\tau^2)}{2}+ \frac{\omega}{\theta} + \frac{\epsilon_\tau^2(1-\tau^2)}{2}\leq  \epsilon_\tau^2(1-\tau^2) +\epsilon_\tau^2\tau^2\leq \epsilon_\tau^2. 
\end{equation*}        

Finally, it can be concluded that the average expected squared gradient norm can be driven below $\epsilon_\tau^2 \geq  \frac{\omega}{ \tau^2\theta}$ for any $\tau\in (0,1)$ with 
\[
    K = {\cal O}\left( \max\left\{\frac{1}{\epsilon_\tau^2\bar \alpha},\frac{\nu}{\epsilon_\tau^4 }\right\} \right).
\]
 
In the literature, $\alpha$ is often selected as $\alpha=\frac{1}{\sqrt{K}}$  in which case \eqref{eq:nonconvex_exp_rate} becomes 
\[
    \frac{1}{K} \sum_{k=0}^{K-1} E[\|\nabla\phi(X_k)\|^2 ]\le \frac{1}{\theta} \left ( \frac{(\phi(x_0) - \hat{\phi})+\nu}{\sqrt{ K}} + {\omega}  \right ). 
\]
However, this simplification conceals the possibility of a faster convergence for small $\nu$, thus we prefer to select $\alpha$ based on the target $\epsilon_\tau$. 
   
\paragraph{Expected rate for functions satisfying the P{\L} condition.}

If the P{\L} condition (Assumption ~\ref{assum:PLcond}) holds, then we have the following convergence rate on the expected function value at the last iterate. 
\begin{theorem} \label{thm:PL_exp_rate} 
    Let Assumptions~\ref{assum:low_bound} and \ref{assum:PLcond} hold. 
    If the descent inequality \eqref{eq:key_stoch_reduction} holds for the sequence $\{X_k\}$ on each $k$, and $\alpha\theta\mu\leq \frac{1}{2}$, then for any integer $K$
    \begin{equation}\label{eq:strconvex_exp_rate}
	    \E\left[ \phi(X_{K}) \right] - \phi^* \leq \left( 1 - \theta \alpha \mu \right)^K \left(\phi (x_0) - \phi^*\right) + \frac{\omega+\alpha \nu}{2\theta \mu}. 
	\end{equation}
   \end{theorem}
\begin{proof} 
    By \eqref{eq:key_stoch_reduction} and Assumption \ref{assum:PLcond}, it follows that
    \begin{align*}
        \E[\phi(X_{k+1}) | \mathcal{F}_k] - \phi^* 
        &\leq \phi (X_k) - \phi^* - \alpha \theta \|\nabla \phi(X_k)\|_2^2  + \alpha \omega + \alpha^2\nu \\
        & \leq \phi (X_k) - \phi^* - 2\alpha\theta \mu \left(\phi (X_k) - \phi^*\right) + \alpha \omega + \alpha^2\nu \\ 
        & \leq (1-2\alpha\theta\mu)\left(\phi (X_k) - \phi^*\right)  + \alpha \omega+\alpha^2\nu. 
    \end{align*}
    Subtracting $\frac{\omega+\alpha \nu}{2\theta \mu}$ from both sides, we get 
    \begin{align*}
        \E[\phi(X_{k+1})| \mathcal{F}_k] - \phi^* -\frac{\omega+\alpha \nu}{2\theta \mu}
        & \leq (1-2\alpha\theta\mu) \left(\phi (X_k)-\phi^* -\frac{\omega+\alpha \nu}{2\theta \mu}\right). 
    \end{align*}
    Recursive application of the above for $k=0, \ldots K-1$ yields 
    \begin{equation*}
	    \E\left[ \phi(X_{K}) \right] - \phi^* -\frac{\omega+\alpha \nu}{2\theta \mu} \leq \left( 1 - 2\theta \alpha \mu \right)^K \left(\phi (x_0) - \phi^* -\frac{\omega+\alpha \nu}{2\theta \mu} \right), 
	\end{equation*}
    which relaxes to \eqref{eq:strconvex_exp_rate}. 
\end{proof}

Here again the best expected optimality gap $\E\left[ \phi(X_{K}) \right] - \phi^* \leq \epsilon_\tau$ that can be achieved  when $\epsilon_\tau \geq \frac{\omega}{2\mu\theta\tau^2}$ for some $\tau\in (0,1)$. 
Selecting $\alpha = \min \{\frac{\epsilon_\tau \theta\mu(1-\tau^2)}{\nu}, \bar \alpha\}$,  where $\bar \alpha$ is again some upper bound on $\alpha$ imposed  when $\nu$ is small or zero, and setting 
\[
    K 
    \ge \frac{1}{\theta \alpha \mu} \log\frac{2(\phi(x_{0}) - \phi^* )}{\epsilon_\tau(1-\tau^2)}\geq \log_{1-\theta \alpha \mu}\frac{\epsilon_\tau(1-\tau^2)}{2(\phi(x_{0}) - \phi^* )}
\]
ensures that 
\begin{align*}
    \E\left[\phi(X_{K})\right] - \phi^* 
    &\leq \left( 1 - \theta \alpha \mu \right)^K \left(\phi (x_0) - \phi^*  \right) +\frac{\omega+\alpha \nu}{2\theta \mu} \\& \leq \frac{\epsilon_\tau(1-\tau^2)}{2}+\epsilon_\tau \tau^2+\frac{\epsilon_\tau(1-\tau^2)}{2}\leq \epsilon_\tau, 
\end{align*}
where the second inequality in the lower bound for $K$ is true because ${-\log (1-\theta\alpha\mu)^{-1}} \leq \frac{1}{\theta\alpha \mu}$, which itself follows from the inequality $-\log (1-t) \ge t $ by setting $t = \theta\alpha\mu$. 
	
In summary, we can reduce the expected optimality gap below any $\epsilon_\tau \geq \frac{\omega}{2\mu\theta\tau^2}$ for some $\tau\in (0,1)$ with 
\[
    K = {\cal O}\left( \max \left\{1, \frac{\nu}{\epsilon_\tau}\right\} \log\left (\frac{1}{\epsilon_\tau }\right )\right)
\]
iterations. 

To the best of our knowledge it is not known how to derive an expected convergence rate under Assumption~\ref{assum:KLcond} or for convex $\phi(x)$ using only \eqref{eq:key_stoch_reduction}, without Assumption~\ref{assum:PLcond} or assuming $\phi$ is strongly convex. As we show later, the expected convergence rate for {\em convex } $\phi(x)$ can be derived using the telescopic sum analysis by explicitly using convexity and unbiased gradient estimator assumption. 

\subsubsection{Expected complexity bounds.}
We will now turn to deriving the expected complexity bounds for any algorithm satisfying \eqref{eq:key_stoch_reduction}. Analogously to the deterministic case, we consider the three cases where $\phi(x)$ satisfies Assumption~\ref{assum:lip_cont}, only, or satisfies Assumption~\ref{assum:KLcond} or Assumption~\ref{assum:PLcond} in addition. 

While a complexity guarantee can be immediately obtained after proving a convergence rate in the deterministic case, as shown in Section~\ref{sec:descent_lemma}, an expected convergence rate does not immediately lead to an expected complexity in the stochastic case. 
The expected-complexity bound is an upper bound on certain {\em stopping time} of the stochastic optimization process. 
Given a stochastic process $\{Y_k\}_{k\in\mathbb{N}}$, where $\mathbb{N} = \{1, 2, \dots\}$, we define its natural filtration by
\[
    \mathcal{F}_0 \text{ as the trivial } \sigma\text{-algebra}, 
    \qquad 
    \mathcal{F}_k=\sigma(Y_1,\dots,Y_k), \quad k\in\mathbb{N}.
\]
An $\mathbb{N}\cup \{\infty\}$-valued random variable $N$ is called a stopping time with respect to $\{\mathcal{F}_k\}_{k\in\mathbb{N}}$ if $\{N\le k\}\in \mathcal{F}_k$ for every $k\in\mathbb{N}$. 
To bound its expectation, we use the following lemma. 

\begin{theorem}\label{thm:walds_inequal}
    Let $(\Omega, \mathcal{F}, (\mathcal{F}_i)_{i\geq 0}, \mathbb{P})$ be a filtered probability space. 
    Let $(Y_{k})_{k\in\mathbb{N}}$ be a sequence of random variables such that $Y_k\in L^1$ is $\mathcal{F}_k$-measurable. 
    Let $N$ be a stopping time with respect to $(\mathcal{F}_k)_{k\in\mathbb{N}}$. 
    Suppose there is a constant $\mathcal{Y}>0$ such that 
    \[ 
    \E[Y_{k}\mid \mathcal{F}_{k-1}] \geq \mathcal{Y} \text{ on } \{N\geq k\}
    \]
    for every $k\geq 1$. 
        For any positive integer $m$, let $N \wedge m = \min(N, m)$ be the truncated stopping time. 
    Define
    \[
      S_m=  \sum_{k=1}^{N \wedge m} Y_k.
    \]
Suppose there exists a finite constant $Z$ such that 
\[
S_m\leq Z\quad \text{a.s. for all } m\geq 1
\]
Then 
    \[
\E\left(  N\right)  \mathcal{Y}\leq    Z
    \]
    and thus $N<\infty$ almost surely. 
\end{theorem}
\begin{proof}

    \[
        S_m = \sum_{k=1}^m Y_k \mathbf{1}_{\{N \ge k\}}.
    \]
    Since $N$ is a stopping time, $\{N \ge k\} = \{N > k-1\} \in \mathcal{F}_{k-1}$, the random indicator $\mathbf{1}_{\{N \ge k\}}$ is $\mathcal{F}_{k-1}$-measurable. 
    Thus, taking expectation of the above equality yields 
    \[
      \mathbb{E} \left[ S_m\right ]
        = \sum_{k=1}^m \mathbb{E} \left[ Y_k  \mathbf{1}_{\{N \ge k\}} \right]
        = \sum_{k=1}^m \mathbb{E} \left[ \mathbf{1}_{\{N \ge k\}} \mathbb{E}[Y_k \mid \mathcal{F}_{k-1}] \right].
    \]
    Since $\mathbb{E} [Y_k \mid \mathcal{F}_{k-1}] \geq \mathcal{Y}$ on $\{N\geq k \}$, we have 
    \[
        \mathbb{E} \left[ \mathbf{1}_{\{N \ge k\}} \mathbb{E}[Y_k \mid \mathcal{F}_{k-1}] \right]
        \geq \mathcal{Y} \mathbb{E}[ \mathbf{1}_{\{N \ge k\}} ], 
    \]
    which leads to  
    \[
        \mathbb{E} \left[ \sum_{k=1}^{N \wedge m} Y_k \right] 
        \geq \mathcal{Y} \sum_{k=1}^m \mathbb{P}(N \geq k)
        = \mathcal{Y}  \mathbb{E}[N \wedge m].
    \]
   Since $S_m\leq Z$ then  $\mathcal{Y}  \mathbb{E}[N \wedge m]\leq Z$. 
   Taking  the limit as $m \to \infty$ by the monotone convergence theorem
   \[ 
    \mathcal{Y}  \mathbb{E}[N]\leq Z. 
    \]
    which completes the proof. 
%
\end{proof}

Now we show how the results for deterministic complexity can be adapted to the stochastic case. 
We consider the stochastic process $\{X_k\}$ and let ${\cal F}_{k}$ be the corresponding filtration. 
For convenience we will assume that whatever accuracy we aim to reach with the algorithm, it is not attained at $x_0$.

\paragraph{Expected complexity for general smooth functions.}
The first theorem is the analogue of Theorem~\ref{thm:nonconvex_bound} for general smooth functions.  
\begin{theorem}\label{thm:nonconvex_exp_bound}
    Suppose Assumption~\ref{assum:low_bound} holds for the objective function $\phi$, and inequality \eqref{eq:key_stoch_reduction} holds for the stochastic process $\{X_k\}_{k\in\mathbb{N}}$.  
    For any $\tau\in (0,1)$ and $\epsilon_\tau \geq \sqrt{\frac{\omega+\alpha\nu}{ \tau^2\theta}}$, define the stopping time $K_{\epsilon_\tau}=\min\{k \in \mathbb{N}: \| \nabla \phi (X_k) \| \leq \epsilon_\tau\}$.
    Then, we have the bound 
    \begin{align*}
        \E[K_{\epsilon_\tau}] \leq \frac{\phi(x_0) - \hat{\phi}}{ \alpha\theta(1 - \tau^2) \epsilon_\tau^2}. 
    \end{align*}
\end{theorem}

\begin{proof} 
    Without loss of generality, assume $K_{\epsilon_\tau}\geq 1$. Define $Y_k=\phi(X_{k-1})-\phi(X_{k})$. For $k=1, \ldots, K_{\epsilon_\tau}$, from \eqref{eq:key_stoch_reduction} we have 
    \[ \begin{aligned}
        \E [Y_k \mid {\cal F}_{k-1}] 
        &\ge \alpha\theta\|\nabla\phi(X_{k-1})\|^2 - \alpha \omega -\alpha^2\nu \\
        &> \alpha\theta \epsilon_\tau^2 - \alpha \omega -\alpha^2\nu 
        \geq \alpha\theta\epsilon_\tau^2(1-\tau^2), 
    \end{aligned} \]
    where the last inequality follows from the definition of $\epsilon_\tau$. 
    Defining  $N=K_{\epsilon_\tau}$ and ${\cal Y}= \alpha\theta\epsilon_\tau^2(1-\tau^2)$ in  Theorem~\ref{thm:walds_inequal} 
    we observe that $S_m= \phi(x_0) - \phi(X_{K_{\epsilon_\tau}\wedge m}) \leq Z$ with $Z=\phi(x_0) - \hat{\phi}$. Thus by Theorem~\ref{thm:walds_inequal} 
    \[
        \E[K_{\epsilon_\tau}] \alpha\theta\epsilon_\tau^2(1-\tau^2)
        \leq  \phi(x_0) - \hat{\phi}
    \]
    which proves the statement of the theorem. 
\end{proof}

Let us again note the connection between $\epsilon_\tau$ and $\alpha$. 
The solution accuracy $\epsilon_\tau$ is bounded below by $\sqrt{\frac{\omega}{\theta}}$ regardless of $\alpha$, but any $\epsilon_\tau> \sqrt{\frac{\omega}{\tau^2 \theta}}$ is achievable by setting $\alpha $
sufficiently small, and in particular $\alpha \sim \frac{\epsilon_{\tau}^2}{\nu}$. 
Thus, 
\[
    \E[K_{\epsilon_\tau}] \leq {\cal O}\left(\max\left\{ \frac{1}{\epsilon_\tau^2} , \frac{\nu}{\epsilon_\tau^4} \right\}\right). 
\]

\paragraph{Expected complexity for  functions satisfying generalized K{\L} condition.}

We now turn, as in the deterministic case, to functions for which Assumption \ref{assum:KLcond} holds, which includes convex functions with bounded level sets. We have the following analogue of Theorem~\ref{thm:KL_cond_rate}.

\begin{theorem}\label{thm:KL_cond_stoch_complex}
    Suppose Assumption~\ref{assum:KLcond} holds with a constant $c=\frac{1}{D}$, for some $D>0$, and inequality \eqref{eq:key_stoch_reduction} holds for the stochastic process $\{X_k\}_{k\in\mathbb{N}}$.
    For any $\tau\in (0,1)$ and $\epsilon_\tau\geq \sqrt{\frac{\omega+\alpha\nu}{\tau^2\theta}}D$, define the stopping time $K_{\epsilon_\tau}=\min\{k\in\mathbb{N}: \phi (X_k) - \phi^*  \leq \epsilon_\tau\}$.
    Then we have the bound 
    \begin{align*}
        \E\left[K_{\epsilon_\tau}\right] \leq \frac{2D^2}{\alpha \theta (1-\tau^2)\epsilon_\tau}
    \end{align*}
\end{theorem}

\begin{proof}
Let
\[
\Delta_k=\phi(X_k)-\phi^*.
\]

Define the function $\Psi_{\epsilon_\tau}:\mathbb{R}+\to\mathbb{R}$ by
\[
\Psi_{\epsilon_\tau}(t)=\left \{
\begin{aligned}
&\frac{1}{t}, & t\geq\epsilon_\tau,\\
&\frac{2}{\epsilon_\tau}-\frac{t}{\epsilon_\tau^2}, 
&0\leq t<\epsilon_\tau.
\end{aligned} \right .
\]
The function $\Psi_{\epsilon_\tau}$ is convex, decreasing, and satisfies
\[
\Psi_{\epsilon_\tau}(t)\leq \frac{2}{\epsilon_\tau}
\qquad\text{for all }t\geq0.
\]

For any $k\leq K_{\epsilon_\tau}$, we have $\Delta_{k-1}>\epsilon_\tau$. Using Assumption~\ref{assum:KLcond}, it follows that
\[
\|\nabla\phi(X_{k-1})\|
\geq \frac{1}{D}\Delta_{k-1}.
\]
Hence, \eqref{eq:key_stoch_reduction} implies
\[
\mathbb{E}[\Delta_k\mid\mathcal{F}_{k-1}]
\leq
\Delta_{k-1}
-\frac{\alpha\theta}{D^2}\Delta_{k-1}^2
+\alpha\omega+\alpha^2\nu.
\]

Using $\Delta_{k-1}>\epsilon_\tau$ we have 
\[
\Psi_{\epsilon_\tau}^\prime(\Delta_{k-1})
=-\frac{1}{\Delta_{k-1}^2}.
\]

By convexity of $\Psi_{\epsilon_\tau}$ 
\[
\Psi_{\epsilon_\tau}(\Delta_k)
\geq
\Psi_{\epsilon_\tau}(\Delta_{k-1})
-\frac{1}{\Delta_{k-1}^2}
(\Delta_k-\Delta_{k-1}).
\]
Taking conditional expectation with respect to $\mathcal{F}_{k-1}$, we obtain
\[
\begin{aligned}
\mathbb{E}[\Psi_{\epsilon_\tau}(\Delta_k)\mid\mathcal{F}_{k-1}]
-\Psi_{\epsilon_\tau}(\Delta_{k-1})
&\geq
-\frac{1}{\Delta_{k-1}^2}
\left(
\mathbb{E}[\Delta_k\mid\mathcal{F}_{k-1}]
-\Delta_{k-1}
\right)\\
&\geq
\frac{\alpha\theta}{D^2}
-\frac{\alpha\omega+\alpha^2\nu}{\Delta_{k-1}^2}\\
&\geq
\frac{\alpha\theta}{D^2}
-\frac{\alpha\omega+\alpha^2\nu}{\epsilon_\tau^2}.
\end{aligned}
\]
Since $\epsilon_\tau$ is chosen so that
\[
\epsilon_\tau
\geq
D\sqrt{\frac{\omega+\alpha\nu}{\theta\tau^2}},
\]
we have 
\[
\frac{\alpha\omega+\alpha^2\nu}{\epsilon_\tau^2}
\leq
\frac{\alpha\theta\tau^2}{D^2}.
\]
Consequently,
\[
\mathbb{E}[\Psi_{\epsilon_\tau}(\Delta_k)\mid\mathcal{F}_{k-1}]
-\Psi_{\epsilon_\tau}(\Delta_{k-1})
\geq
\frac{\alpha\theta(1-\tau^2)}{D^2}.
\]

Let us define   ${\cal Y}=  \frac{\alpha \theta (1-\tau^2)}{D^2}$ and $Y_k=\Psi_{\epsilon_\tau}(\Delta_k)-\Psi_{\epsilon_\tau}(\Delta_{k-1})$
 then  in Theorem \ref{thm:walds_inequal} 
  for any $m\geq1$,
\[
\begin{aligned}
S_m
&=
\sum_{k=1}^{K_{\epsilon_\tau}\wedge m}Y_k\\
&=
\Psi_{\epsilon_\tau}
\left(\Delta_{K_{\epsilon_\tau}\wedge m}\right)
-\Psi_{\epsilon_\tau}(\Delta_0)\
&\leq
\frac{2}{\epsilon_\tau}\\
\end{aligned}
\]
Taking $Z=\frac{2}{\epsilon_\tau}$ and $N=K_{\epsilon_\tau}$
we apply Theorem \ref{thm:walds_inequal} to conclude that  $K_{\epsilon_\tau}<\infty$ almost surely and
 \[
 \E[K_{\epsilon_\tau}]\leq \frac{2D^2}{\alpha \theta (1-\tau^2)\epsilon_\tau}
\]
\end{proof}

We note again that any $\epsilon_\tau> \sqrt{\frac{\omega}{ \tau^2\theta}}D$ is achievable by setting $\alpha $
sufficiently small, and in particular $\alpha \sim \frac{\epsilon_{\tau}^2}{\nu}$. 
Thus the result of Theorem \ref{thm:KL_cond_stoch_complex} gives 
\[
	E[K_{\epsilon_\tau}] \leq  {\cal O}\left(\frac{1}{\epsilon_\tau} \max\left\{1, \frac{\nu}{\epsilon_\tau^2}\right\}\right). 
\] 

We note that Although convexity and boundedness of the initial level set imply Assumption~\ref{assum:KLcond}  on that level set, stochastic iterates need not remain in it, since function decrease is not guaranteed at every iteration.
This issue is often remedied by modifying the algorithm so that it remains in some bounded set. As we will show later, the convex case is more commonly analyzed under different conditions using the telescopic sum approach and convergence rate analysis. It is important to note that the complexity that corresponds to the convergence rate obtained through that approach is better than what we get from Theorem \ref{thm:KL_cond_stoch_complex}. This is likely due to the fact that convexity is a stronger assumption than Assumption~\ref{assum:KLcond} but it also may be the result of the analysis. As in the deterministic case, the telescopic sum analysis is not as easily amenable to bias, as the analysis based on  \eqref{eq:key_stoch_reduction}. 

\paragraph{Expected complexity for functions satisfying P{\L} condition.} 

Finally, we present the expected complexity for the case when $\phi(x)$ satisfies the P{\L} condition, that is Assumption~\ref{assum:PLcond}, which includes strongly convex functions. 
\begin{theorem}\label{thm:PL_cond_stoch_complex}
Suppose Assumption~\ref{assum:PLcond} holds, and inequality \eqref{eq:key_stoch_reduction} holds for the sequence $\{X_k\}_{k\in\mathbb{N}}$.  For any $\tau\in (0,1)$ and $\epsilon_\tau \geq \frac{\omega + \alpha\nu}{2\theta \mu \tau^2}$, define the stopping time $K_{\epsilon_\tau}  =\min\{k\in\mathbb{N}: \phi (X_k) - \phi^*  \leq \epsilon_\tau\}$.
Then, we have 
\begin{align*}
    \E[K_{\epsilon_\tau}]
    \le  \frac{1}{2\alpha\theta \mu(1-\tau^2)}\left (\log\left(\frac{\phi (x_0) - \phi^*}{\epsilon_\tau} \right)+1\right ). 
\end{align*}
\end{theorem}

\begin{proof}
The steps of this proof closely follow those of Theorem \ref{thm:KL_cond_stoch_complex}. We again define 
\[
\Delta_k=\phi(X_k)-\phi^*,
\]
and this time we make a use of  the function $\Psi_{\epsilon_\tau}:\mathbb{R}+\to\mathbb{R}$ defined by
\[
\Psi_{\epsilon_\tau}(t)=\left \{
\begin{aligned}
&-\log t, & t\geq\epsilon_\tau,\\
&-\log\epsilon_\tau+1-\frac{t}{\epsilon_\tau}, 
&0\leq t<\epsilon_\tau.
\end{aligned} \right .
\]
The function $\Psi_{\epsilon_\tau}$ is again convex, decreasing, and satisfies
\[
\Psi_{\epsilon_\tau}(t)\leq -\log\epsilon_\tau+1
\qquad\text{for all }t\geq0.
\]

For any $k\leq K_{\epsilon_\tau}$, we have $\Delta_{k-1}>\epsilon_\tau$. By \eqref{eq:key_stoch_reduction} and Assumption~\ref{assum:PLcond},
\[
\begin{aligned}
\mathbb{E}[\Delta_k\mid\mathcal{F}_{k-1}]
&\leq
\Delta_{k-1}
-\alpha\theta|\nabla\phi(X_{k-1})|^2
+\alpha\omega+\alpha^2\nu\\
&\leq
\Delta_{k-1}
-2\alpha\theta\mu\Delta_{k-1}
+\alpha\omega+\alpha^2\nu.
\end{aligned}
\]
Because $\Delta_{k-1}>\epsilon_\tau$ 
\[
\Psi_{\epsilon_\tau}^\prime(\Delta_{k-1})
=-\frac{1}{\Delta_{k-1}},
\]
and by convexity of $\Psi_{\epsilon_\tau}$
\[
\Psi_{\epsilon_\tau}(\Delta_k)
\geq
\Psi_{\epsilon_\tau}(\Delta_{k-1})
-\frac{1}{\Delta_{k-1}}
(\Delta_k-\Delta_{k-1}).
\]
Taking conditional expectation with respect to $\mathcal{F}_{k-1}$, we obtain
\[
\begin{aligned}
\mathbb{E}[\Psi{\epsilon_\tau}(\Delta_k)\mid\mathcal{F}_{k-1}]
-\Psi{\epsilon_\tau}(\Delta_{k-1})
&\geq
-\frac{1}{\Delta_{k-1}}
\left(
\mathbb{E}[\Delta_k\mid\mathcal{F}_{k-1}]
-\Delta_{k-1}
\right)\\
&\geq
2\alpha\theta\mu
-\frac{\alpha\omega+\alpha^2\nu}{\Delta_{k-1}}\\
&\geq
2\alpha\theta\mu
-\frac{\alpha\omega+\alpha^2\nu}{\epsilon_\tau}.
\end{aligned}
\]
By the choice of $\epsilon_\tau$,
\[
\epsilon_\tau
\geq
\frac{\omega+\alpha\nu}{2\theta\mu\tau^2},
\]
and hence
\[
\frac{\alpha\omega+\alpha^2\nu}{\epsilon_\tau}
\leq
2\alpha\theta\mu\tau^2.
\]
Consequently,
\[
\mathbb{E}[\Psi_{\epsilon_\tau}(\Delta_k)\mid\mathcal{F}_{k-1}]
-\Psi{\epsilon_\tau}(\Delta_{k-1})
\geq
2\alpha\theta\mu(1-\tau^2).
\]

Now we define $Y_k=\Psi_{\epsilon_\tau}(\Delta_k)-\Psi_{\epsilon_\tau}(\Delta_{k-1})$ and ${\cal Y}=2\alpha\theta\mu(1-\tau^2)$. 
Thus in Theorem \ref{thm:walds_inequal} 
 for any $m\geq1$,
\[
\begin{aligned}
S_m
&=
\sum_{k=1}^{K_{\epsilon_\tau}\wedge m}Y_k
=
\Psi_{\epsilon_\tau}
\left(\Delta_{K_{\epsilon_\tau}\wedge m}\right)
-\Psi_{\epsilon_\tau}(\Delta_0)\\
&\leq
-\log\epsilon_\tau+1+\log\Delta_0
=
\log\left(
\frac{\phi(x_0)-\phi^*}{\epsilon_\tau}
\right)+1.
\end{aligned}
\]
Taking
\[
Z=
\log\left(
\frac{\phi(x_0)-\phi^*}{\epsilon_\tau}
\right)+1
\]
and $N=K_{\epsilon_\tau}$ we apply Theorem \ref{thm:walds_inequal} to conclude that  $K_{\epsilon_\tau}<\infty$ almost surely and
\[
\mathbb{E}[K_{\epsilon_\tau}]
\leq
\frac{1}
{2\alpha\theta\mu(1-\tau^2)}
\left[
\log\left(
\frac{\phi(x_0)-\phi^*}{\epsilon_\tau}
\right)+1
\right].
\]
\end{proof}

Here again, we observe that any $\epsilon_\tau> \frac{\omega}{2\theta \mu \tau^2}$ is achievable by setting $\alpha$ sufficiently small, specifically, choosing $\alpha \sim \frac{\epsilon_{\tau}^2}{\nu}$. 
The result of Theorem~\ref{thm:PL_cond_stoch_complex} gives
\[
	\E[K_{\epsilon_\tau}] \leq  {\cal O}\left( \max\left \{1, \frac{\nu}{\epsilon_\tau}\right \}\log\left (\frac{1}{\epsilon_\tau}\right )\right). 
\]

\subsection{Analysis of fixed-step algorithm under stochastic gradient oracle}\label{sec:stoch_fixed_step}
We now analyze Algorithm~\ref{alg:first_order_fixed_step_sgd}, a simple fixed-step-size algorithm. 
It uses either Oracle~\ref{def:oracle_stoch_1.1}, Oracle~\ref{def:oracle_stoch_1.2}, or Oracle~\ref{def:oracle_stoch_1.3} to obtain its gradient estimate $g(x_k, \xi_k)$, where $\xi_k$ represents all randomness in iteration $k$. 
For all three types of oracles, we will derive inequality \eqref{eq:key_stoch_reduction} and the expected complexity as well as expected rate for Algorithm~\ref{alg:first_order_fixed_step_sgd} with specific values of $\theta$, $\omega$ and $\nu$. 
One of our key observations is that relative error constants, $\kappa_g$ and $\sigma_v$ affect the complexity constants, while the absolute bias $\epsilon_g$ determines an accuracy floor. The additive variance parameter $\sigma_c$ determines the additional dependence of the complexity on the target accuracy through the step-size restriction.

As in the deterministic case, it follows from Assumption~\ref{assum:lip_cont} that
\begin{align*}
    \phi(X_{k+1}) &\leq \phi (X_k) - \alpha \nabla \phi (X_k)^\intercal g_k + \frac{\alpha^2 L_1}{2} \| g_k\|^2, 
\end{align*}
and taking conditional expectation with respect to $\xi_k$ given $\mathcal{F}_{k}$ on both sides gives 
\begin{equation} \label{eq:expected_Lipschitz_UB}
    \E[\phi(X_{k+1}) \mid \mathcal{F}_{k}] \leq \phi(X_k) - \alpha \nabla\phi (X_k)^\intercal \E[g_k \mid \mathcal{F}_{k}] + \frac{\alpha^2 L_1}{2} \E\left[\|g_k\|^2 \mid \mathcal{F}_{k}\right].
\end{equation}

\begin{algorithm}[H]
    \caption{~\textbf{Stochastic Gradient Method with Fixed Step Size}}
    \label{alg:first_order_fixed_step_sgd}
    \begin{algorithmic}
        \State \textbf{Inputs:} $x_0 \in \mathbb{R}^n$ (starting point), $\alpha > 0$ (step size)
        \State \textbf{Requirements:} Oracle \ref{def:oracle_stoch_1.1}, \ref{def:oracle_stoch_1.2}  or \ref{def:oracle_stoch_1.3}.  
        \State \textbf{1:} \textbf{for $k=0,1,2,\dots$ do}
        \State \text{2:}~~~~Compute $g_k = g(x_k,\xi_k)$ via Oracle~\ref{def:oracle_stoch_1.1}, \ref{def:oracle_stoch_1.2} or \ref{def:oracle_stoch_1.3}.
        \State \text{3:}~~~~Update iterate $x_{k+1}\gets x_k-\alpha g_k$.
        \State \textbf{4:} \textbf{end for}
    \end {algorithmic}
\end{algorithm}

At each iteration, the oracle randomness is generated independently of the past, conditional on the current iterate $X_k$. Thus for brevity we will omit conditioning on 
$  \mathcal{F}_{k}$ when taking expectations over $g_k$.  In other words, we use $\E[g_k] =\E[g_k \mid \mathcal{F}_{k}] $ and $\E\left[\|g_k\|^2 \right]=\E\left[\|g_k\|^2 \mid \mathcal{F}_{k}\right]$. 

\subsubsection{Descent inequality and expected complexity analysis}
\paragraph{Analysis of Algorithm~\ref{alg:first_order_fixed_step_sgd} under Oracle \ref{def:oracle_stoch_1.1}.} 
We start with the following lemma that gives specific values to $\theta$, $\omega$ and $\nu$ in \eqref{eq:key_stoch_reduction}. 
\begin{lemma} \label{lem:descent_stoch_Case1}
    Let Assumption~\ref{assum:lip_cont} hold. 
    In any iteration $k \in \mathbb{N}$ of Algorithm~\ref{alg:first_order_fixed_step_sgd}, if $g_k$ comes from Oracle~\ref{def:oracle_stoch_1.1} with $\kappa_g<1$ and $\alpha \le 1/L_1$, the following inequality holds:
    \begin{equation*} 
        \E[\phi(X_{k+1}) \mid \mathcal{F}_k] \leq \phi (X_k) - (1-\kappa_g^2)\frac{\alpha}{2} \|\nabla \phi(X_k)\|^2 + \alpha\frac{\epsilon^2_g}{2} 
    \end{equation*}
\end{lemma}

\begin{proof}
    Since $g_k$ comes from Oracle~\ref{def:oracle_stoch_1.1}, it follows that 
    \begin{equation*} 
        \E\left[\|g_k - \nabla\phi(X_k)\|^2\right] \leq \kappa_g^2 \|\nabla \phi(X_k)\|^2+ \epsilon_g^2, 
    \end{equation*} 
    which implies 
    \begin{align*}
        -2 \nabla\phi (X_k)^\intercal \E[g_k] \leq \epsilon^2_g - \left(1-\kappa_g^2\right) \|\nabla \phi(X_k)\|^2 - \E\left[\|g_k\|^2\right]. 
    \end{align*}
    With the above inequality and $\alpha \leq 1/L_1$, it follows from \eqref{eq:expected_Lipschitz_UB} that
    \begin{align*}
        \E[ \phi(X_{k+1}) \mid \mathcal{F}_{k}]
        &\leq \phi (X_k) + \frac{\alpha}{2} \left( \epsilon^2_g - \left(1-\kappa_g^2\right) \|\nabla \phi(X_k)\|^2 -\E[ \|g_k\|^2]\right) + \frac{\alpha^2 L_1}{2} \E\left[\|g_k\|^2\right]\\
        &= \phi (X_k) - \frac{\alpha}{2}\left(1-\kappa_g^2\right) \|\nabla \phi(X_k)\|^2   - \frac{\alpha}{2}\left(1 - \alpha L_1\right) \E\left[\|g_k\|^2\right] + \frac{\alpha}{2}  \epsilon^2_g\\
        &\leq \phi (X_k) - \frac{\alpha}{2}\left(1-\kappa_g^2\right) \|\nabla \phi(X_k)\|^2 + \frac{\alpha}{2}  \epsilon^2_g. 
    \end{align*} 
\end{proof}

The above lemma shows, with Oracle~\ref{def:oracle_stoch_1.1}, Theorems~\ref{thm:nonconvex_exp_bound},  \ref{thm:KL_cond_stoch_complex} and \ref{thm:PL_cond_stoch_complex} apply with $\theta=\frac{1-\kappa_g^2}{2}$, $\omega = \frac{\epsilon_g^2}{2}$, and $\nu=0$ for Lipschitz smooth functions. 
Substituting  these constants back in these theorems, and we have the following complexity result for Algorithm~\ref{alg:first_order_fixed_step_sgd}. 

\begin{theorem}\label{thm:Case1_exp_complex}
    Suppose Assumptions~\ref{assum:low_bound} and \ref{assum:lip_cont} hold for the objective function $\phi$. 
    In Algorithm~\ref{alg:first_order_fixed_step_sgd}, set $\alpha\leq 1/L_1$ and compute $g_k$ using Oracle~\ref{def:oracle_stoch_1.1} with $\kappa_g<1$. 
    Then the following statements hold true. 

    \begin{enumerate}
    \item 
        For any $\tau\in (0,1)$ and $\epsilon_\tau\geq \frac{\epsilon_g}{\tau \sqrt{1-\kappa_g^2}}$, the stopping time $K_{\epsilon_\tau}=\min\{k \in \mathbb{N}: \|\nabla\phi (X_k)\| \leq \epsilon_\tau\}$ is bounded as 
        \begin{align*}
            \E[K_{\epsilon_\tau}] \leq \frac{2(\phi(x_0) - \hat \phi)}{ \alpha(1-\kappa_g^2)(1 - \tau^2)\epsilon_\tau^2}. 
        \end{align*}
    \item 
        Suppose in addition  Assumption \ref{assum:KLcond} holds with a constant $c=\frac{1}{D}$, for some $D>0$.
        Then for any $\tau\in (0,1)$ and $\epsilon_\tau\geq \frac{\epsilon_g D}{\tau \sqrt{1-\kappa_g^2}}$, the stopping time 
        $K_{\epsilon_\tau} = \min\{k \in \mathbb{N}: \phi (X_k) - \phi^*  \leq \epsilon_\tau\}$ is bounded as 
        \begin{align*}
            \E\left[K_{\epsilon_\tau}\right] \leq \frac{4D^2}{\alpha (1-\kappa_g^2) (1-\tau^2)\epsilon_\tau}. 
        \end{align*}
    \item 
        If Assumption~\ref{assum:PLcond} holds instead of \ref{assum:KLcond}, then for any $\tau\in (0,1)$ and $\epsilon_\tau \geq \frac{\epsilon_g^2}{ 2(1-\kappa_g^2)\mu\tau^2}$, the stopping time 
        $K_{\epsilon_\tau}=\min\{k\in\mathbb{N}: \phi (X_k) - \phi^*  \leq \epsilon_\tau\}$ is bounded as 
        \begin{align*}
            \E[K_{\epsilon_\tau}] 
            \le \frac{1}{\alpha(1-\kappa_g^2) \mu(1-\tau^2)}\left (\log\left(\frac{\phi(x_0) - \phi^*}{\epsilon_\tau}\right) + 1\right ). 
        \end{align*}
    \end{enumerate}
\end{theorem}

\paragraph{Analysis of Algorithm~\ref{alg:first_order_fixed_step_sgd} under Oracle \ref{def:oracle_stoch_1.2}.}

The next lemma derives the descent inequality \eqref{eq:key_stoch_reduction} for Case 2 of the stochastic oracle.

\begin{lemma}\label{lem:descent_stoch_Case2}
    Let Assumption~\ref{assum:lip_cont} hold. 
    In any iteration $k \in \mathbb{N}$ of Algorithm~\ref{alg:first_order_fixed_step_sgd}, if $g_k$ comes from Oracle~\ref{def:oracle_stoch_1.2} with $\kappa_g<1$ and $\alpha \le 1/L_1$, the following inequality holds:
    \begin{equation*} 
        \E[\phi(X_{k+1}) \mid \mathcal{F}_k] \leq \phi (X_k) - \left(1-\kappa_g^2-\alpha L_1 \sigma_v^2\right)\frac{\alpha}{2} \|\nabla \phi(X_k)\|^2 + \alpha\frac{\epsilon^2_g}{2} +\alpha^2\frac{L_1\sigma_c^2}{2}
    \end{equation*}
\end{lemma}
\begin{proof} 
    By the definition of Oracle~\ref{def:oracle_stoch_1.2}, it follows that
    \begin{align*}
      - 2 \nabla\phi (X_k)^\intercal \E[g_k] \leq - (1-\kappa_g^2) \|\nabla \phi(X_k)\|^2 + \epsilon^2_g -  \|\E[g_k]\|^2 
    \end{align*}
    and 
    \begin{align*}
        \E[\|g_k\|^2] \leq \|\E[g_k]\|^2 + \sigma_c^2 + \sigma_v^2\|\nabla \phi(X_k)\|^2. 
    \end{align*}
    Thus, together with \eqref{eq:expected_Lipschitz_UB}, we have 
    \begin{align*}
       \E[&\phi(X_{k+1})\mid \mathcal{F}_k] \leq \phi(X_k) - \alpha \nabla\phi(X_k)^\intercal \E[g_k] + \frac{\alpha^2 L_1}{2} \E\left[\|g_k\|^2\right] \\
        &\leq \phi(X_k) + \frac{\alpha}{2}\Big(- (1-\kappa_g^2) \|\nabla \phi(X_k)\|^2 + \epsilon^2_g - \|\E[g_k]\|^2\Big) \\
        &\qquad\qquad + \frac{\alpha^2 L_1}{2} \Big(\|\E[g_k]\|^2 + \sigma_c^2 + \sigma_v^2\|\nabla \phi(X_k)\|^2\Big) \\
        &= \phi (X_k) - \left(1-\kappa_g^2-\alpha L_1\sigma_v^2\right)\frac{\alpha}{2} \|\nabla \phi(X_k)\|^2 -\frac{\alpha}{2} (1-\alpha L_1) \|\E[g_k]\|^2 + \alpha\frac{\epsilon^2_g}{2} + \alpha^2\frac{L_1\sigma_c^2}{2} \\
        &\leq \phi(X_k) - \left(1-\kappa_g^2-\alpha L_1\sigma_v^2\right)\frac{\alpha}{2} \|\nabla \phi(X_k)\|^2 + \alpha\frac{\epsilon^2_g}{2} + \alpha^2\frac{L_1\sigma_c^2}{2}.
    \end{align*}
\end{proof}

For this case, we conclude that if $\alpha\leq \frac{1}{L_1}\min\{1, \frac{1-\kappa_g^2}{2\sigma_v^2}\}$, then Theorems~\ref{thm:nonconvex_exp_bound}, \ref{thm:KL_cond_stoch_complex}, and \ref{thm:PL_cond_stoch_complex} apply with $\theta=\frac{1-\kappa_g^2}{4}$, $\omega =\frac{\epsilon_g^2}{2}$, and $\nu=\frac{L_1\sigma_c^2}{2}$.  
We note that the condition $\alpha \leq \frac{1-\kappa_g^2}{2L_1\sigma_v^2}$ can be generalized to $\alpha \leq \frac{1-\kappa_g^2}{qL_1\sigma_v^2}$ for any $q>1$, which changes $\theta$ to $\frac{q-1}{q} \frac{1-\kappa_g^2}{2}$. 
However, to avoid introducing notation for yet another constant, we will assume $q=2$ henceforth. 

We now combine consequences of Lemma \ref{lem:descent_stoch_Case2} and Theorems \ref{thm:nonconvex_exp_bound}, \ref{thm:KL_cond_stoch_complex} and \ref{thm:PL_cond_stoch_complex} to obtain expected complexity  bounds for  Algorithm~\ref{alg:first_order_fixed_step_sgd} with  Oracle \ref{def:oracle_stoch_1.2}. 
\begin{theorem}\label{thm:Case2_exp_complex}
    Suppose Assumptions~\ref{assum:low_bound} and \ref{assum:lip_cont} hold for the objective function $\phi$. 
    In Algorithm~\ref{alg:first_order_fixed_step_sgd}, set $\alpha \leq \frac{1}{L_1} \min\{1, \frac{1-\kappa_g^2}{2\sigma_v^2}\}$ and compute $g_k$ using Oracle~\ref{def:oracle_stoch_1.2} with $\kappa_g<1$. 
    Then the following statements hold true. 
    \begin{enumerate}
        \item 
        For any $\tau\in (0,1)$ and $\epsilon_\tau \geq \frac{\sqrt{2(\epsilon_g^2+\alpha L_1\sigma_c^2)}}{\tau \sqrt{1-\kappa_g^2}}$, the stopping time 
        $K_{\epsilon_\tau}=\min\{k\in\mathbb{N}: \|\nabla\phi(X_k)\| \leq \epsilon_\tau\}$ is bounded as 
        \begin{align*}
            \E[K_{\epsilon_\tau}] \leq \frac{4(\phi(x_0) - \hat{\phi})}{ \alpha(1-\kappa_g^2)(1 - \tau^2)\epsilon_\tau^2}. 
        \end{align*}
        \item 
        Suppose in addition  Assumption \ref{assum:KLcond} holds with a constant $c=\frac{1}{D}$, for some $D>0$.
        Then for any $\tau\in (0,1)$ and $\epsilon_\tau \geq \frac{\sqrt{2(\epsilon_g^2+\alpha L_1\sigma_c^2)}D}{\tau \sqrt{1-\kappa_g^2}}$, the stopping time 
        $K_{\epsilon_\tau} = \min\{k\in\mathbb{N}: \phi(X_k) - \phi^* \leq \epsilon_\tau\}$ is bounded as 
        \begin{align*}
            \E\left [K_{\epsilon_\tau}\right ] \leq \frac{1}{\epsilon_\tau} \frac{8D^2}{\alpha (1-\kappa_g^2) (1-\tau^2)}. 
        \end{align*}
        \item If Assumption~\ref{assum:PLcond} holds instead of \ref{assum:KLcond}, then for any $\tau\in (0,1)$ and $\epsilon_\tau \geq \frac{\epsilon_g^2+\alpha L_1\sigma_c^2}{ \tau^2\mu (1-\kappa_g^2)}$,  the stopping time 
        $K_{\epsilon_\tau}=\min\{k\in\mathbb{N}: \phi(X_k) - \phi^* \leq \epsilon_\tau\}$ is bounded as 
        \begin{align*}
            \E[K_{\epsilon_\tau}] \leq \frac{2}{\alpha(1-\kappa_g^2) \mu(1-\tau^2)}\left (\log \left(\frac{\phi(x_0) - \phi^*}{\epsilon_\tau}\right) + 1\right)
        \end{align*}
    \end{enumerate}
\end{theorem}

\paragraph{Analysis of Algorithm~\ref{alg:first_order_fixed_step_sgd} under Oracle \ref{def:oracle_stoch_1.3}.}
Finally, we present the lemma deriving the descent inequality for  Case 3 of the stochastic oracle. This lemma statement is essentially identical to Lemma \ref{lem:descent_stoch_Case2} for Case 2 and the derivation follows the same steps. 

\begin{lemma} \label{lem:descent_stoch_Case3}
    Let Assumption~\ref{assum:lip_cont} hold. 
    In any iteration $k \in \mathbb{N}$ of Algorithm~\ref{alg:first_order_fixed_step_sgd}, if $g_k$ comes from Oracle~\ref{def:oracle_stoch_1.3} with $\kappa_g<1$ and $\alpha \le 1/L_1$, the following inequality holds:
    \begin{equation*}
        \E[\phi(X_{k+1}) \mid \mathcal{F}_k] \leq \phi(X_k) - {(1-\kappa_g^2-\alpha L_1 \sigma_v^2)}\frac{\alpha}{2} \|\nabla\phi(X_k)\|^2 + \alpha\frac{\epsilon^2_g}{2} + \alpha^2\frac{L_1\sigma_c^2}{2}. 
    \end{equation*}
\end{lemma}

\begin{proof}
    In this case, $\xi_k$ consists of two random variables $\zeta_1$ and $\zeta_2$ (the iteration index $k$ is omitted). The inequality \eqref{eq:expected_Lipschitz_UB} can be written as 
    \begin{align*}
       \E_{\zeta_1,\zeta_2}[\phi(X_{k+1}) \mid \mathcal{F}_{k}] &\leq \phi (X_k) - \alpha \nabla\phi(X_k)^\intercal \E_{\zeta_1,\zeta_2}[g_k] + \frac{\alpha^2 L_1}{2} \E_{\zeta_1,\zeta_2}\left[\|g_k\|^2\right]. 
    \end{align*} 
    By the properties of the oracle, it follows that
    \begin{align*}
      - 2 \nabla\phi (X_k)^\intercal \E_{\zeta_1,\zeta_2}[g_k] 
      \leq - (1-\kappa_g^2) \|\nabla \phi(X_k)\|^2 + \epsilon^2_g -  \E_{\zeta_1}\|\E_{\zeta_2}[g_k\mid \zeta_1]\|^2 
    \end{align*}
    and 
    \begin{align*}
        \E_{\zeta_1,\zeta_2}\left[\|g_k\|^2\right] \leq \E_{\zeta_1}\left[\|\E_{\zeta_2}[g_k\mid \zeta_1]\|^2\right] +\sigma_c^2+\sigma_v^2\|\nabla \phi(X_k)\|^2. 
    \end{align*}
    Thus, 
    \begin{align*}
        &\E_{\zeta_1,\zeta_2}[\phi(X_{k+1}) \mid \mathcal{F}_k]  \leq \phi (X_k) - \alpha \nabla\phi(X_k)^\intercal \E_{\zeta_1,\zeta_2}[g_k] + \frac{\alpha^2 L_1}{2} \E_{\zeta_1,\zeta_2}\left[\|g_k\|^2\right] \\
        &\leq \phi(X_k) + \frac{\alpha}{2} \Big(- (1-\kappa_g^2) \|\nabla \phi(X_k)\|^2 + \epsilon^2_g -  \E_{\zeta_1}\left[\|\E_{\zeta_2}[g_k\mid \zeta_1]\|^2\right]\Big) \\
        &\qquad\qquad + \frac{\alpha^2 L_1}{2} \Big(\E_{\zeta_1}\left[\|\E_{\zeta_2}[g_k\mid \zeta_1]\|^2\right] +\sigma_c^2+\sigma_v^2\|\nabla \phi(X_k)\|^2\Big) \\
        &= \phi (X_k) - \left(1-\kappa_g^2-\alpha L_1\sigma_v^2\right) \frac{\alpha}{2} \|\nabla \phi(X_k)\|^2 -\frac{\alpha}{2} (1-\alpha L_1) \E_{\zeta_1}\left[\|\E_{\zeta_2} [g_k\mid \zeta_1]\|^2\right] + \alpha\frac{\epsilon^2_g}{2} +\alpha^2\frac{L_1\sigma_c^2}{2} \\
        &\leq \phi (X_k) - (1-\kappa_g^2-\alpha L_1\sigma_v^2)\frac{\alpha}{2} \|\nabla \phi(X_k)\|^2 + \alpha\frac{\epsilon^2_g}{2} +\alpha^2\frac{L_1\sigma_c^2}{2}.
    \end{align*}
\end{proof}

The expected descent inequality for Case 2 and Case 3 has the same form, thus Theorem 
\ref{thm:Case2_exp_complex} can be restated with identical complexity bounds for  Oracle \ref{def:oracle_stoch_1.3}. 

\subsubsection{Expected rate analysis of Algorithm~\ref{alg:first_order_fixed_step_sgd}} 
Let Assumptions~\ref{assum:low_bound} and \ref{assum:lip_cont} hold. 
Given that we have derived inequality \eqref{eq:key_stoch_reduction} for Algorithm~\ref{alg:first_order_fixed_step_sgd} for all three cases of the stochastic oracle in Lemmas~\ref{lem:descent_stoch_Case1}, \ref{lem:descent_stoch_Case2}, or \ref{lem:descent_stoch_Case3}, 
we can apply Theorems~\ref{thm:nonconvex_exp_rate} and \ref{thm:PL_exp_rate} to obtaine expected convergence rates for general smooth functions and for functions under the P{\L} condition, respectively.

When using Oracle~\ref{def:oracle_stoch_1.1} with $\alpha \leq \frac{1}{L_1}$, we have $\theta=\frac{1-\kappa_g^2}{2}$, $\omega =\frac{\epsilon_g^2}{2}$, and $\nu=0$ from Lemma~\ref{lem:descent_stoch_Case1}. Thus, by Theorem~\ref{thm:nonconvex_exp_rate}, it holds for any $K\in\mathbb{N}$ that 
\begin{equation}\label{eq:nonconvex_exp_rate_case1} 
    \frac{1}{K} \sum_{k=0}^{K-1} \E\left[\|\nabla\phi(X_k)\|^2\right]
    \le \frac{2}{1-\kappa_g^2} \left( \frac{\phi(x_0) - \hat{\phi}}{\alpha K} + \frac{\epsilon_g^2}{2}\right), 
\end{equation}
and, if $\phi$ satisfies the P{\L} condition, we have
\begin{equation}\label{eq:PL_exp_rate_Case1}
    \E[\phi(X_{K})] - \phi^* \leq \left( 1 - \frac{1-\kappa_g^2}{2} \alpha \mu \right)^K \left(\phi (x_0) - \phi^*\right) + \frac{\epsilon_g^2}{2(1-\kappa_g^2) \mu}. 
\end{equation}

When using Oracle~\ref{def:oracle_stoch_1.2} or \ref{def:oracle_stoch_1.3} with $\alpha \leq \frac{1}{L_1} \min\{1, \frac{1-\kappa_g^2}{2\sigma_v^2}\}$, we have $\theta=\frac{1-\kappa_g^2}{4}$, $\omega =\frac{\epsilon_g^2}{2}$, and $\nu=\frac{L_1\sigma_c^2}{2}$ from Lemma~\ref{lem:descent_stoch_Case2} or Lemma~\ref{lem:descent_stoch_Case3}. 
Then, it holds for any $K\in\mathbb{N}$ that 
\begin{equation}\label{eq:nonconvex_exp_rate_case2_3} 
    \frac{1}{K} \sum_{k=0}^{K-1} \E\left[\|\nabla\phi(X_k)\|^2\right]
    \le \frac{4}{1-\kappa_g^2} \left( \frac{\phi(x_0) - \hat{\phi}}{\alpha K} + \frac{\epsilon_g^2}{2} +\alpha \frac{L_1\sigma_c^2}{2}\right), 
\end{equation}
and, if $\phi$ satisfies the P{\L} condition, we have
\begin{equation}\label{eq:PL_exp_rate_Case2_3}
    \E[\phi(X_{K})] - \phi^* \leq \left( 1 - \frac{1-\kappa_g^2}{4} \alpha \mu \right)^K \left(\phi (x_0) - \phi^*\right) + \frac{\epsilon_g^2+\alpha L_1\sigma_c^2}{(1-\kappa_g^2) \mu}. 
\end{equation}


As mentioned earlier, expected rate for convex functions is usually derived by directly using convexity and applying   {\em telescopic sum}  technique. We present this analysis in the next subsection. 
	 
\subsubsection{Telescopic sum analysis for convex functions under Oracle~\ref{def:oracle_stoch_1.2}}\label{sec:telescopic_stoch}

Our general assumptions on the stochastic first-order oracle \ref{def:oracle_stoch_1.2} have the form (omitting random variables and other dependencies for brevity): 
\begin{equation*} 
    \|\E [g(x)] - \nabla \phi(x)\|^2 \leq \epsilon_g^2+\kappa_g^2\|\nabla \phi(x)\|^2, 
\end{equation*} 
\begin{equation*} 
    \E \left[\|g(x) - \E[g(x)]\|^2\right] \leq \sigma_c^2+\sigma_v^2\|\nabla \phi(x)\|^2. 
\end{equation*}

In Section \ref{sec:telescopic}, we saw that the analysis of convex functions via telescopic sum suffers when $\kappa_g>0$.  
The stochastic case is more complicated and is typically analyzed under the assumption $\epsilon_g=\kappa_g=0$, in other words, when the stochastic first-order oracle is unbiased. 
Under this assumption and convexity, we can get an expected rate result without relying only on the generalized K{\L} conditions, as presented in the following theorem. 

\begin{theorem}\label{thm:stoch_convex_telescopic}
    Suppose the objective function $\phi$ is convex,  satisfies Assumption~\ref{assum:lip_cont} and attains it minimum at some point $x^*$. 
    In Algorithm~\ref{alg:first_order_fixed_step_sgd}: let $g_k$ be computed via Oracle~\ref{def:oracle_stoch_1.2} with $\epsilon_g=\kappa_g=0$ and $\alpha \le \frac{1}{(\sigma_v^2+1) L_1}$. 
    Then, the following inequality holds for any $K \in \mathbb{N}$: 
    \begin{equation}
        \frac{1}{K-1} \sum_{k=0}^K \mathbb{E}\big[\phi(X_k) - \phi^* \big]
        \le \frac{\|x_0 - x^*\|^2}{2\alpha K} + \frac{\alpha}{2}\sigma_c^2
    \end{equation}
\end{theorem}

\begin{proof}
    Let $g_k=g(X_k)$. It holds for $k=0, \ldots, K-1$ that 
    \begin{equation*} 
        \|X_{k+1} - x^*\|^2
        = \|X_k - x^*\|^2 - 2\alpha g_k^T (X_k - x^*)+ \alpha^2 \|g_k\|^2.
    \end{equation*}

    Taking conditional expectation we have
    \begin{equation*}
        \mathbb{E}\!\left[ \|X_{k+1} - x^*\|^2 \mid \mathcal{F}_k\right]
        = \|X_k - x^*\|^2 - 2\alpha \nabla\phi(X_k)^T (X_k - x^*) + \alpha^2 \mathbb{E}\left[ \|g_k\|^2 \mid \mathcal{F}_k\right].
    \end{equation*}
    The inequality 
    \[
        \nabla \phi(X_k)^T(X_k - x^*) \ge \phi(X_k) - \phi^* + \frac{1}{2L_1} \|\nabla\phi(X_k)\|^2
    \] 
    holds for convex functions satisfying Assumption~\ref{assum:lip_cont}, and the properties of the first-order oracle implies 
    \[ 
        \E\left[\|g(x)\|^2\right] 
        = \E \left[\|g(x) - \E[g(x)]\|^2\right] + \|\nabla\phi(x)\|^2
        \le \sigma_c^2 + (\sigma_v^2+1) \|\nabla\phi(X_k)\|^2. 
    \]
    Plug these into the equation above to obtain 
    \begin{equation*} \begin{aligned}
        &\mathbb{E}\!\left[\|X_{k+1} - x^*\|^2 \mid \mathcal{F}_k\right] \\
        &\le \|X_k - x^*\|^2 - 2\alpha\Big(\phi(X_k) - \phi^* + \frac{1}{2L_1} \|\nabla\phi(X_k)\|^2\Big) + \alpha^2 \sigma_c^2 + \alpha^2 (\sigma_v^2 + 1) \|\nabla \phi(X_k)\|^2 \\
        &= \|X_k - x^*\|^2 - 2 \alpha (\phi(X_k) - \phi^* )+ \alpha^2 \sigma_c^2 + \Big(\alpha^2 (\sigma_v^2+1) - \frac{\alpha}{L_1}\Big) \|\nabla \phi(X_k)\|^2 \\
        &\leq \|X_k - x^*\|^2 - 2 \alpha (\phi(X_k) - \phi^* )+ \alpha^2 \sigma_c^2, 
    \end{aligned} \end{equation*}
    where the last inequality holds because $\alpha \le \frac{1}{(\sigma_v^2+1) L_1}$ implies $\alpha^2 (\sigma_v^2+1) - \frac{\alpha}{L_1} \le 0$. 

    Summing the above inequality over $k = 0,\dots,K-1$, applying total expectation, and canceling telescoping terms yields 
    \begin{equation*}
        2\alpha \sum_{k=0}^{K-1} \mathbb{E}\big[\phi(X_k) - \phi^* \big]
        \le \|x_0 - x^*\|^2 - \E\|X_K - x^*\|^2 + K \alpha^2 \sigma_c^2,
    \end{equation*}
    which gives the result of the theorem after dropping the negative term on the right and dividing by $2\alpha K$. 
\end{proof}

Theorem~\ref{thm:stoch_convex_telescopic} implies that, to guarantee
\begin{equation*}
    \frac{1}{K-1} \sum_{k=0}^K \mathbb{E}\big[\phi(X_k) - \phi^* \big]\leq \epsilon
\end{equation*}
for a given $\epsilon>0$, one can choose $\alpha= \min\{ \frac{1}{(\sigma_v^2+1) L_1}, \frac{\epsilon}{\sigma_c^2}\}$ and 
\[
    K \geq \frac{\|x_0 - x^*\|^2}{\alpha\epsilon}. 
\]
This clearly implies that if the additive gradient variance $\sigma_c$ is very close to zero, then the complexity is $K \sim \frac{(\sigma_v^2+1)L_1 \|x_0 - x^*\|^2}{\epsilon}$, which is consistent with the exact gradient descent up to the constant $\sigma_v^2+1$. 
When $\sigma_c$ is larger, the complexity becomes $k\sim {\cal O}(\frac{\sigma_c^2}{\epsilon^2})$, which is the best known complexity for  stochastic gradient methods applied to convex functions. 

\paragraph{Benign bias}
There is a specific kind of biased oracle that does not hurt the telescopic sum analysis. Specifically, instead of assuming $\mathbb{E}[g_k] = \nabla\phi(x_k)$, it is sufficient to assume that $\mathbb{E}[g_k] = c\nabla \phi(x_k)$ for some constant $c>0$. Then the above derivations change as follows. 

The assumptions on the oracle now implie
\begin{equation*} 
    \E [\|g(x)\|]^2 
    = \E \left[\|g(x) - \E[g(x)]\|^2\right] + c^2 \|\nabla\phi(x)\|^2
    \leq \sigma_c^2+(\sigma_v^2+c^2)\|\nabla \phi(x)\|^2. 
\end{equation*} 
Then, following the same logic as in the proof above, we have when $\alpha \le \frac{c}{(\sigma_v^2+c^2) L_1}$ that 
\begin{equation*} \begin{aligned}
    \mathbb{E}\!&\left[\|X_{k+1} - x^*\|^2 \mid \mathcal{F}_k\right]
    = \|X_k - x^*\|^2 - 2\alpha c \nabla\phi(X_k)^T (X_k - x^*)+ \alpha^2 \mathbb{E}\left[ \|g_k\|^2 \mid \mathcal{F}_k\right] \\
    &\le \|X_k - x^*\|^2 - 2\alpha c \Big(\phi(X_k) - \phi^* + \frac{1}{2L_1} \|\nabla\phi(X_k)\|^2\Big) + \alpha^2 \sigma_c^2 + \alpha^2 (\sigma_v^2 + c^2) \|\nabla \phi(X_k)\|^2 \\
    &\leq \|X_k - x^*\|^2 - 2\alpha c (\phi(X_k) - \phi^* )+ \alpha^2 \sigma_c^2. 
\end{aligned} \end{equation*}
After telescoping, we again have the result 
\begin{equation}
    \frac{1}{K} \sum_{k=0}^{K-1} \mathbb{E}\big[\phi(X_k) - \phi^* \big]
    \le \frac{\|x_0 - x^*\|^2}{2\alpha c K} + \frac{\alpha}{2c}\sigma_c^2. 
\end{equation}

This observation is important as it shows that $\E[g(x)]$ needs to coincide with the gradient direction in expectation, but does not have to preserve the length. This is particularly useful, for example, when analyzing randomized coordinate descent (or more generally  when $g(x)$ is a projection onto a random subspace). 
In those cases, as shown in Section \ref{sec:gradproj}, $c=\frac{s}{n}$ (where $s$ is the dimension of the subspace), $\sigma_c=0$, and $\sigma_v^2= (n-s)s/n^2$ and thus $\sigma_v^2 + c^2 = \frac{s}{n}$, which recovers rate
\begin{equation}
    \frac{1}{K} \sum_{k=0}^{K-1} \mathbb{E}\big[\phi(X_k) - \phi^* \big]
    \le \frac{nL_1\|x_0 - x^*\|^2}{2sK}
\end{equation}
assuming $\alpha = \frac{c}{(\sigma_v^2+c^2) L_1} = \frac{1}{L_1}$.

\subsection{Total complexity bounds for Algorithm  \ref{alg:first_order_fixed_step_sgd} under specific examples of stochastic first order oracles}

We now show how the results of Section \ref{sec:stoch_fixed_step} apply to the specific examples of the stochastic first-order oracles discussed in Section~\ref{sec:stoch_first_order_ex}. 
By specifying the parameters of the oracles, this analytical framework immediately yields convergence or complexity guarantees for optimization algorithms using such oracles.
This approach recovers many well-known results from the literature and also leads to several new ones, demonstrating the utility of the stochastic oracle definitions together with the general analysis developed above.

As before $K_{\epsilon_\tau}$ will denote the number of iterations of Algorithm~\ref{alg:first_order_fixed_step_sgd} that achieve $\epsilon_{\tau}$ accuracy, defined appropriately depending on the properties of $\phi(x)$. We will use $T_{\epsilon_\tau}$ to denote the total oracle cost to reach such accuracy. 
The oracle costs will be defined according to the oracles in use. 

\subsubsection{Minibatch gradient descent in ERM} 
Consider the very popular setting of stochastic gradient descent in empirical risk minimization as discussed in Section~\ref{sec:erm}. Recall that given positive integer $\lambda$ as input, the oracle generates $\xi_1,\ldots,\xi_\lambda$ as i.i.d.\ samples from $\Xi$ and computes
\begin{equation} \label{eq:erm_oracle2}
    g(x, \xi^{(1)}) = \frac{1}{\lambda}\sum_{j=1}^\lambda \nabla_x \ell(x,\xi_j), 
\end{equation}
with    $\xi^{(1)}:=(\xi_1,\ldots,\xi_\lambda)$
As shown, in this case we have Oracle~\ref{def:oracle_stoch_1.2} with $\epsilon_g = 0$, $\kappa_g = 0$, $\sigma_c^2 = M_c^2 / \lambda$, and $\sigma_v^2 = M_v^2 / \lambda$, where $M_c$ and $M_v$ are absolute and relative variance constants for a single sample. 
By applying Theorem~\ref{thm:Case2_exp_complex}, we obtain the following complexity results. 

\begin{theorem}\label{thm:Case2_exp_complex_erm}
    Suppose Assumptions~\ref{assum:low_bound} and \ref{assum:lip_cont} hold. In Algorithm~\ref{alg:first_order_fixed_step_sgd}: let 
    $g_k$ be computed as \eqref{eq:erm_oracle1} and $\alpha\leq \frac{1}{L_1} \min \{1, \frac{\lambda}{2M_v^2}\}$.  
    \begin{enumerate}
        \item 
            For any $\tau\in (0,1)$
            and $\epsilon_\tau\geq \sqrt{\frac{\alpha L_1M_c^2}{ \tau^2\lambda}}$, the stopping time 
            $K_{\epsilon_\tau}=\min\{k\in\mathbb{N}:  \| \nabla\phi (X_k) \| \leq \epsilon_\tau\}$ is bounded as 
            \begin{align*}
                \E[K_{\epsilon_\tau}] \leq  \frac{4(\phi(x_0) - \hat{\phi})}{ \alpha(1 - \tau^2) \epsilon_\tau^2}. 
            \end{align*}
        \item 
            Suppose in addition  Assumption \ref{assum:KLcond} holds with a constant $c=\frac{1}{D}$, for some $D>0$.
            For any $\tau\in (0,1)$ and $\epsilon_\tau\geq \sqrt{\frac{2\alpha L_1M_c^2D^2}{\tau^2 \lambda}}$, the stopping time 
            $K_{\epsilon_\tau}=\min\{k\in\mathbb{N}: \phi (X_k) - \phi^* \leq \epsilon_\tau\}$ is bounded as 
            \begin{align*}
                \E\left [K_{\epsilon_\tau}\right ] \leq \frac{8D^2}{\alpha (1-\tau^2) \epsilon_\tau} 
            \end{align*}
        \item 
            If Assumption~\ref{assum:PLcond} holds instead of \ref{assum:KLcond}, then for any $\tau\in (0,1)$
            and $\epsilon_\tau \geq \frac{\alpha L_1M_c^2}{ \tau^2\lambda\mu}$,  the stopping time 
            $K_{\epsilon_\tau}=\min\{k\in\mathbb{N}: \phi (X_k) - \phi^* \leq \epsilon_\tau\}$ is bounded as 
            \begin{align*}
                \E[K_{\epsilon_\tau}] \le \frac{2}{\alpha\mu(1-\tau^2)}\left (\log\left(\frac{\phi(x_0) - \hat{\phi}}{\epsilon_\tau}\right) + 1\right ). 
            \end{align*}
    \end{enumerate}
\end{theorem}

The first result in Theorem~\ref{thm:Case2_exp_complex_erm} shows that Algorithm~\ref{alg:first_order_fixed_step_sgd} can find an $x_k$ such that $\|\nabla\phi(x_k)\| \le \epsilon_\tau$ for any $\epsilon_\tau\geq \sqrt{\frac{\alpha L_1M_c^2}{ \tau^2\lambda}}$ with any $\tau \in (0,1)$.  
Thus, given any $\epsilon>0$, if we choose $\alpha\leq \frac{\lambda\epsilon^2}{2L_1M_c^2}$, then since it ensures $\epsilon$ is no less than $\sqrt{\frac{2\alpha L_1M_c^2}{\lambda}}$, the lower bound for $\epsilon_\tau$ in Theorem~\ref{thm:Case2_exp_complex_erm} when $\tau^2=\frac{1}{2}$, Algorithm~\ref{alg:first_order_fixed_step_sgd} can find an $x_k$ such that $\|\nabla\phi(x_k)\| \le \epsilon$. 
Then, by Theorem \ref{thm:Case2_exp_complex_erm}, with $\alpha = \frac{1}{L_1}\min \{1, \frac{\lambda}{2M_v^2}, \frac{\lambda\epsilon^2}{2M_c^2}\}$, the expected stopping time $K_{\epsilon}=\min\{k\in\mathbb{N}:  \| \nabla\phi (X_k) \| \leq \epsilon\}$ is bounded as 
\begin{align*}
	\E[K_{\epsilon}] \leq  \frac{8L_1 (\phi(x_0) - \hat{\phi})}{\epsilon^2}  \max \left\{1, \frac{2M_v^2}{\lambda}, \frac{2M_c^2}{\lambda\epsilon^2}\right\}. 
\end{align*}

When $M_v=0$ and $M_c > 0$, this essentially recovers the well known ${\cal O}(\epsilon^{-4})$ complexity of stochastic gradient descent in the case of simple bounded variance \cite{ghadimi2013stochastic}. 
It also recovers the results in \cite{bottou2018optimization} where $M_c$ and $M_v$ both can be positive. 
Most importantly, it clearly shows the effect on complexity of these two constants, specifically, if $M_c$ (and thus $\sigma_c$) is zero (or very small) --- the case referred to as ``strong growth condition'', discussed in Sections \ref{strong_growth_reg} and \ref{strong_growth_class} and analyzed for example in \cite{schmidt2013fast, vaswani2019fast} -  then step size $\alpha$ can be chosen to be ``large'', i.e. independent of $\epsilon$, and complexity is comparable to that of the deterministic algorithm, up to the constant $\sigma_v^2$. 
On the other hand, if $\sigma_c$ is significant then the standard stochastic gradient complexity ${\cal O}(\epsilon^{-4})$  dominates.

The second result of Theorem \ref{thm:Case2_exp_complex_erm} shows that Algorithm~\ref{alg:first_order_fixed_step_sgd} with $\alpha\leq \frac{1}{L_1} \min \{1, \frac{\lambda}{2M_v^2}\}$ can reduce the optimality gap $\phi(x_k) - {\phi}^*$ below any $\epsilon_\tau \geq \sqrt{\frac{2\alpha L_1M_c^2D^2}{\lambda \tau^2}}$ with any $\tau\in (0,1)$.   
For any $\epsilon>0$, choosing $\alpha\leq \frac{\lambda\epsilon^2}{4L_1M_c^2D^2}$ ensures $\epsilon \geq \sqrt{\frac{2\alpha L_1M_c^2D^2}{\lambda \tau^2}}$ when $\tau^2=\frac{1}{2}$. 
Then, by Theorem \ref{thm:Case2_exp_complex}, with $\alpha = \frac{1}{L_1} \min \{1, \frac{\lambda}{2M_v^2}, \frac{\lambda\epsilon^2}{4M_c^2D^2}\}$, the stopping time $K_{\epsilon}=\min\{k\in\mathbb{N}: \phi(X_k) - \phi^* \leq \epsilon\}$ is bounded as 
\begin{equation}\label{eq:erm_KL_complex}
    \E\left[ K_{\epsilon} \right] \leq \frac{16L_1 D^2}{\epsilon} \max \left\{1, \frac{2M_v^2}{\lambda}, \frac{4M_c^2D^2}{\lambda\epsilon^2}\right\}. 
\end{equation}

If, instead of Assumption~\ref{assum:KLcond}, we assume that $\phi(x)$ is convex, then we can apply Theorem~\ref{thm:stoch_convex_telescopic}, because in this case $\kappa_g=\epsilon_g=0$. 
The result shows, when $\alpha \le \frac{1}{(\sigma_v^2+1) L_1} = \frac{\lambda}{(M_v^2+\lambda) L_1}$, it holds for any $K\in\mathbb{N}$ that 
\begin{equation*}
    \frac{1}{K-1} \sum_{k=0}^K \mathbb{E}\big[\phi(X_k) - \phi^*\big]
    \leq \frac{\|x_0 - x^*\|^2}{2\alpha K} + \frac{\alpha M_c^2}{2\lambda}. 
\end{equation*} 
We can further conclude from this result that, given any $\epsilon > 0$, if $\alpha = \min\{ \frac{\lambda}{(M_v^2+\lambda) L_1}, \frac{\epsilon\lambda}{M_c^2}\}$ and 
\begin{equation}\label{eq:erm_cvx_complex}
    K \geq \frac{\|x_0 - x^*\|^2}{\alpha \epsilon}
    = \frac{\|x_0 - x^*\|^2}{\lambda} \max\left\{\frac{(M_v^2+\lambda)L_1}{\epsilon},  \frac{M_c^2}{\epsilon^2}\right\}, 
\end{equation}
then $\frac{1}{K} \sum_{k=0}^K \mathbb{E}\big[\phi(X_k) - \phi^*\big] \le \epsilon$. 

Here we note the difference between dependence on $\epsilon$ in \eqref{eq:erm_KL_complex} and \eqref{eq:erm_cvx_complex}. If $M_c$ is positive, the complexity obtained under convexity is ${\cal O}(\epsilon^{-2})$ and 
that obtained under Assumption~\ref{assum:KLcond} is ${\cal O}(\epsilon^{-2})$. This is consistent with other known results in the literature \cite{fatkhullin2022sharp,Bubeck2015}

Finally, let us consider the case when Assumption \ref{assum:PLcond} holds. 
Applying the same logic as in the cases above, for any $\epsilon>0$, choosing $\alpha\leq \frac{\mu \lambda\epsilon}{2L_1M^2_c}$ ensures $\epsilon\geq \frac{2\alpha L_1M_c^2}{\lambda\mu}$. 
Thus, by Theorem \ref{thm:Case2_exp_complex_erm}, with $\alpha= \frac{1}{L_1} \min\{1, \frac{\lambda}{2M_v^2}, \frac{\lambda\mu \epsilon}{2M_c^2}\}$, the stopping time $K_{\epsilon}=\min\{k\in\mathbb{N}: \phi (X_k) - \phi^* \leq \epsilon\}$ is bounded as 
\[
    \E\left[ K_{\epsilon} \right] \le \frac{4L_1}{\mu} \max \left\{1, \frac{2M_v^2}{\lambda}, \frac{2M_c^2}{\lambda \mu \epsilon}\right\} \left(\log\left(\frac{\phi(x_0) - {\phi}^*}{\epsilon}\right) + 1\right ). 
\]

In the case of this stochastic first-order oracle, the total oracle complexity (also known as a sample complexity)
 $T_{\epsilon}= \lambda K_{\epsilon}$.  It is then clear that the bound on $\E[T_{\epsilon}]$  does not depend on the size of the minibatch,  unless $\lambda$ is so large that the last two terms in the minimum defining $\alpha$ are larger than $1$,  since the iteration complexity bounds  on $|K_{\epsilon}]$ (in all considered cases of $\phi(x)$)  scale linearly with $\frac{1}{\lambda}$ (via $\alpha$) while the cost of the first order oracle, scales linearly with $\lambda$. 

\subsubsection{Randomized forward finite difference methods with orthonormal directions and  deterministically bounded zeroth-order oracle}\label{sec:RFFD_orth_complex_epsilonf} 

Let us now turn to first-order oracle based on finite differences, as discussed in Section \ref{sec:RFFD_orth}, but restricted to the case where the  zeroth order oracle error is deterministically bounded. 
Specifically, we consider oracle implemented using \eqref{eq:RFFD}, where random directions ${\cal U}$ satisfy Assumption \ref{ass:haar}
and where function estimates  are computed using a special case of zeroth-order single-point   Oracle  \ref{def:oracle_stoch_0.1}  
 $\epsilon_f\geq 0$ and   $\sigma_f=0$. 

In this case, gradient estimate computed via randomized finite differences constitutes a case of Oracle~\ref{def:oracle_stoch_1.1} with 
$\epsilon_g^2 = s \big(\frac{L_1 {h}}{2} + \frac{2\epsilon_f}{h}\big)^2$ and $\kappa_g^2 = \frac{n-s}{n}$.  Let us assume the "ideal" conditions, 
where we can choose $\alpha=\frac{1}{L_1}$ and  $h$ so that $\epsilon_g^2$ is minimized. In this case $h=2\sqrt{\frac{\epsilon_f}{L_1}}$ and $\epsilon_g^2 = 4s L_1\epsilon_f$. 

To analyze  Algorithm~\ref{alg:first_order_fixed_step_sgd} we apply Theorem \ref{thm:Case1_exp_complex} and thus we have the following. 

\begin{theorem}\label{thm:Case1_exp_complex_rffd_det}
   Suppose Assumptions~\ref{assum:low_bound} and \ref{assum:lip_cont} hold. 
In Algorithm~\ref{alg:first_order_fixed_step_sgd} let $\alpha=1/L_1$ and  $g_k$ be computed by  \eqref{eq:RFFD} satisfying Assumption \ref{ass:haar}
and with $f(\cdot)$ computed by 
Oracle \ref{def:oracle_stoch_0.1} with $\sigma_f=0$, then 

\begin{enumerate}
\item 
 For any $\tau\in (0,1)$
and $\epsilon_\tau\geq 2\sqrt{\frac{nL_1\epsilon_f}{ \tau^2}}$, the stopping time 
$K_{\epsilon_\tau}=\min\{k=0,1, \ldots, :  \| \nabla \phi (X_k) \| \leq \epsilon_\tau\}$ is bounded as 
\begin{align*}
    \E[K_{\epsilon_\tau}] \leq  \frac{2nL_1}{s} \frac{(\phi(x_0) - \hat{\phi})}{ (1 - \tau^2)\epsilon_\tau^2}.
\end{align*}
\item 
Suppose in addition  Assumption \ref{assum:KLcond} holds with a constant $c=\frac{1}{D}$, for some $D>0$.
For any $\tau\in (0,1)$
and $\epsilon_\tau\geq 2\sqrt{\frac{nL_1\epsilon_fD^2}{\tau^2}}$  the stopping time 
$K_{\epsilon_\tau}=\min\{k\in {\mathbb N}:   \phi (X_k) - \phi^* \leq \epsilon_\tau\}$ is bounded as 
\begin{align*}
	\E\left [K_{\epsilon_\tau}\right ] \leq \frac{4nL_1 D^2}{(1-\tau^2)s\epsilon_\tau}.
\end{align*}
\item If instead of Assumption~\ref{assum:KLcond},  ~\ref{assum:PLcond} holds then, for any $\tau\in (0,1)$
and $\epsilon_\tau \geq  \frac{2nL_1 \epsilon_f}{\mu \tau^2}$,   the stopping time 
$K_{\epsilon_\tau}=\min\{k=0, 1, \ldots, :   \phi (X_k) - \phi^* \leq \epsilon_\tau\}$ is bounded as follows:
\begin{align*}
    \E[K_{\epsilon_\tau}]\leq  \frac{2nL_1}{s\mu(1-\tau^2)}\left (\log \frac{(\phi (x_0) - \phi^*)}{\epsilon_\tau}+1\right ). 
\end{align*}
\end{enumerate}
\end{theorem}

In the case when $L_1$ and $\epsilon_f$ are not known exactly we assume that  $\alpha\leq \frac{1}{L_1}$ and  $h\geq 2\sqrt{\frac{\epsilon_f}{L_1}}$. In this case the only change is
$\epsilon_g^2 = s (L_1h)^2$ and the results above can be restated accordingly.

The most important observations that we derive from the above results is that the choice of $s$ which is the number of random orthogonal directions, and thus is the dimension of the random  subspace where the gradient estimate is constructed, has no effect on the convergence neighborhood. It also has virtually no effect on the total complexity $T_{\epsilon_\tau}$,  since iterations complexity $K_{\epsilon_\tau}$, in each case, scales as $\frac{n}{s}$ and the cost of the Oracle~\ref{def:oracle_stoch_1.1} is proportional $s+1$.

\subsubsection{Randomized forward finite difference methods with orthonormal directions and single-point stochastic  zeroth-order oracle}\label{sec:RFFD_orth_complex_sigmaf} 

Let us again consider \eqref{eq:RFFD} with ${\cal U}$ satisfying Assumption \ref{ass:haar} but now with access to the single-point Oracle  \ref{def:oracle_stoch_0.1} for which   $\sigma_f\geq 0$ and $\epsilon_f=0$. 
In other words, we have access to an unbiased stochastic function value estimator, 
 which is the setting in many applications of simulation optimization, stochastic approximation and machine learning. 
 
In this case the gradient estimate \eqref{eq:RFFD}, as analyzed in Section \ref{sec:RFFD_orth},   constitutes an example of Oracle~\ref{def:oracle_stoch_1.3} with 
$\epsilon_g^2 = s (\frac{L_1 {h}}{2} )^2$ and $\kappa_g^2 = \frac{n-s}{n}$, $\sigma_c^2=\frac{2s\sigma_f^2}{h^2}$ and $\sigma_v^2=0$.  The cost of this oracle is proportional to $s$.

To analyze  Algorithm~\ref{alg:first_order_fixed_step_sgd} we apply Theorem \ref{thm:Case2_exp_complex} to the specific parameters of this oracle. Note that 
since $\sigma_v=0$ then choosing $\alpha\leq \frac{1}{L_1}$ is sufficient. We have the following convergence guarantees. 

\begin{theorem}\label{thm:Case2_exp_complex_ffd_stoch}
    Suppose Assumptions~\ref{assum:low_bound} and \ref{assum:lip_cont} hold. In Algorithm~\ref{alg:first_order_fixed_step_sgd}: let $\alpha\leq \frac{1}{L_1}$ and  $g_k$ computed by \eqref{eq:RFFD}, under Assumption \ref{ass:haar} and  with access to Oracle  \ref{def:oracle_stoch_0.1} with $\epsilon_f=0$.
    
\begin{enumerate}
\item 
For any $\tau\in (0,1)$
and $\epsilon_\tau\geq \sqrt{\frac{2n}{\tau^2} \left((\frac{L_1 {h}}{2} )^2+\alpha L_1\frac{2\sigma_f^2}{h^2}\right )}$, the stopping time 
$K_{\epsilon_\tau}=\min\{k=0,1, \ldots, :  \| \nabla \phi (X_k) \| \leq \epsilon_\tau\}$ is bounded as 
\begin{align*}
	\E[K_{\epsilon_\tau}] \leq  \frac{1}{\epsilon_\tau^2} \frac{4n(\phi(x_0) - \hat{\phi})}{ \alpha s (1 - \tau^2)}.
\end{align*}
\item 
Suppose in addition  Assumption \ref{assum:KLcond} holds with a constant $c=\frac{1}{D}$, for some $D>0$.
For any $\tau\in (0,1)$
and $\epsilon_\tau\geq \sqrt{\frac{2nD^2}{\tau^2} \left((\frac{L_1 {h}}{2} )^2+\alpha L_1\frac{2\sigma_f^2}{h^2}\right )}$   the stopping time 
$K_{\epsilon_\tau}=\min\{k=0, 1, \ldots, :   \phi (X_k) - \phi^* \leq \epsilon_\tau\}$ is bounded as 
\begin{align*}
	\E\left [K_{\epsilon_\tau}\right ] \leq \frac{1}{\epsilon_\tau} \frac{8nD^2}{\alpha s (1-\tau^2)}.
\end{align*}
\item If Assumption~\ref{assum:PLcond} holds instead of \ref{assum:KLcond}, then, for any $\tau\in (0,1)$ and $\epsilon_\tau \geq \frac{n}{\mu \tau^2} \left((\frac{L_1 {h}}{2} )^2+\alpha L_1\frac{2\sigma_f^2}{h^2}\right)$, the stopping time $K_{\epsilon_\tau}=\min\{k=0, 1, \ldots, :   \phi (x_k) - \phi^* \leq \epsilon_\tau\}$ is bounded as 
\begin{align*}
    \E[K_{\epsilon_\tau}]\leq   \frac{2n}{\alpha s\mu(1-\tau^2)}\left (\log \frac{\phi(x_0)-\phi^*}{\epsilon_\tau} + 1\right ).
\end{align*}
\end{enumerate}
\end{theorem}

The important observation in this case is that for any $\epsilon>0$ it is possible to find $\epsilon_\tau\leq \epsilon$ by choosing $h$ and $\alpha$ appropriately. 
Specifically, given $\epsilon>0$ and a $\tau\in (0,1)$ let us set $\alpha$ and $h$ so that $(\frac{L_1 {h}}{2} )^2\leq \frac{\epsilon^2\tau^2}{4n}$
and $\alpha L_1\frac{2\sigma_f^2}{h^2}\leq \frac{\epsilon^2\tau^2}{4n}$, it is easy to see that this implies $\epsilon\geq  \sqrt{\frac{2n}{\tau^2} \left((\frac{L_1 {h}}{2} )^2+\alpha L_1\frac{2\sigma_f^2}{h^2}\right )}$. For concreteness let us pick $\tau^2=\frac{1}{2}$, then (assuming we know the relevant constants) we can pick  $h^2=\frac{\epsilon^2}{2nL_1^2}$ and 
$\alpha =  \min\{\frac{1}{L}, \frac{\epsilon^4}{32n^2L_1^3\sigma_f^2}\}$. 

 Then $K_{\epsilon}=\min\{k=0,1, \ldots, :  \| \nabla \phi (X_k) \| \leq \epsilon \}$ becomes bounded as (using ${\cal O}$ to absorb small constants) 
\begin{align}\label{eq:Keps_rffd_sigmaf}
	\E[K_{\epsilon}] \leq   {\cal O}\left ( \frac{n^3L_1^3}{s\epsilon^6}\right ).  
\end{align}
Recall that the cost of each iteration is equal to the cost of the first-order oracle. The finite difference based first-order oracle that evaluates $g_k$ using  $s+1$ calls to the zeroth-order oracle
has a cost of $s$ and thus the total oracle complexity is  
\begin{equation}\label{eq:Teps_rffd_sigmaf}
 \E[T_{\epsilon}]=(s+1) \E[K_{\epsilon}]  \leq {\cal O}\left ( \frac{n^3L_1^3}{\epsilon^6}\right ), 
\end{equation}
 which does not depend on $s$. 

Alternatively, it is often preferred to reduce the variance of the zeroth-order oracle, $\sigma_f^2$ 
by averaging over multiple samples (whenever possible). It is evident from our derivations that smaller $\sigma_f^2$ allows larger choice of  $\alpha$.
Specifically, if $\sigma_f^2\leq \frac{\epsilon^4}{8n^2L_1^2}$ then  $\alpha$ can be chosen as large as $\frac{1}{L_1}$. 
Assuming that variance of a function estimate based on one sample is not greater than $1$,  then averaging 
$\frac{32L_1^2n^2}{\epsilon^4}$ samples gives a zeroth order oracle with $\sigma_f^2 \leq \frac{\epsilon^4}{32L_1^2n^2}$. 
In this case,  $\E[K_{\epsilon}] \leq   {\cal O}\left ( \frac{nL_1}{s\epsilon^2}\right )$, but the total oracle complexity is still 
\[
\E[T_{\epsilon}]=(s+1)\frac{8n^2L_1^2}{\epsilon^4}\E[K_{\epsilon}] \leq {\cal O}\left ( \frac{n^3L_1^3}{\epsilon^6}\right ).
\]
Thus choosing small $\alpha$ or forcing small $\sigma_f$ result in the same total oracle complexity of the algorithm. 
This is consistent with results in simulation optimization literature \cite{HaShashaaniPasupathy2025}. 

Applying similar derivations to the case when Assumption \ref{assum:KLcond} holds we set
$h^2=\frac{\epsilon^2}{2nL_1^2D^2}$ and 
$\alpha = \min\{\frac{1}{L_1}, \frac{\epsilon^4}{32n^2D^4L_1^3\sigma_f^2}\}$.
 This gives the total oracle complexity 
of achieving $\phi (x_k) - \phi^* \leq \epsilon$ being bounded as (again using ${\cal O}$ to absorb small constants)
\begin{equation}\label{eq:Teps_rffd_sigmaf_cvx}
	\E[T_{\epsilon}]= (s+1) \E[K_{\epsilon}]\leq {\cal O}\left ({\frac{n^3D^6L_1^3}{\epsilon^5} }\right ).
\end{equation}
Finally, under Assumption \ref{assum:PLcond}, we choose $h^2=\frac{\epsilon\mu}{nL_1^2}$ and 
$\alpha =   \min\{\frac{1}{L_1},\frac{\mu^2\epsilon^2}{8n^2L_1^3\sigma_f^2}\}$.
which gives final complexity as  
\begin{equation}\label{eq:Teps_rffd_sigmaf_strcvx}
   \E[T_{\epsilon}]=(s+1) \E[K_{\epsilon}]\leq  {\cal O} \left ( \frac{n^3L_1^3 }{\mu^3\epsilon^2} \log \frac{1}{\epsilon}\right ). 
\end{equation}

\paragraph{Extension to the scaled version of randomized finite difference method with orthonormal directions.}
In Section \ref{sec:SRFFD} we considered first-order oracle generated by 
\eqref{eq:SRFFD} with access to the single-point zeroth-order Oracle \ref{def:oracle_stoch_0.1}.  If we again assume that $\epsilon_f=0$ in  Oracle \ref{def:oracle_stoch_0.1} and that ${\cal U}$ satisfies Assumption \ref{ass:haar}, then  the first-order oracle satisfies properties of Oracle~\ref{def:oracle_stoch_1.2} with $\epsilon_g^2=n \left( \frac{L_1 {h}}{2} \right)^2$, $\kappa_g=0$, $\sigma_c^2 = 
 \frac{2n^2 \sigma_f^2}{ s{h}^2}+\frac{n^2}{s} \left( \frac{L_1 {h}}{2}\right)^2$  and $\sigma_v^2 = \frac{(2n-s)(n-s)}{ns}$. 

For simplicity, let us assume that the $\left(\frac{L_1 {h}}{2}\right)^2$ term in the expression for $\sigma_c^2$ is negligibly small compared to $\frac{ \sigma_f^2}{{h}^2}$. 
In this case we can easily compare scaled and unscaled oracle properties where $\epsilon_g^2$ is the same as its counterpart for  the unscaled version multiplied by  $\frac{n}{s}$
and $\sigma_c^2$ gets scaled by  $\frac{n^2}{s^2}$, while $\kappa_g^2=0$ and instead $\sigma_v^2>0$ and roughly scales as $\frac{n}{s}$ for $s\ll n$. 

Thus Theorem \ref{thm:Case2_exp_complex_ffd_stoch} changes so that the step size has to satisfy $\alpha \leq \min\{\frac{1}{L_1}, \frac{ns}{2(2n-s)(n-s)L_1}\}$ and
for the case of general smooth $\phi(x)$ 
$\epsilon_\tau\geq  \sqrt{\frac{n}{\tau^2} \left((\frac{L_1 {h}}{2} )^2+\alpha L_1\frac{n}{s}\frac{2\sigma_f^2}{h^2}\right )}$. 
Compared to the unscaled case this floor bound has the same first (bias) term under the square root but the second term is scaled by $\frac{n}{s}$.

Thus is we want to ensure convergence to an $\epsilon$ stationary point for an given $\epsilon$ the optimal choice of $h$ and $\alpha$,  taking $\tau^2=\frac{1}{2}$,  should be the same for $h$ as in the unscaled case, i.e. 
  $h^2=\frac{\epsilon^2}{2nL_1^2}$ and  $\alpha$ scaled by $\frac{s}{n}$ compared to the unscaled case $\alpha =  \frac{s\epsilon^2h^2}{16L_1n^2\sigma_f^2}= \frac{s\epsilon^4}{32L_1^3n^3\sigma_f^2}$. Assuming this choice of $\alpha$ is smaller than $\min\{\frac{1}{L_1}, \frac{s}{2(2n-s)(n-s)L_1}\}$ (which is essentially w.l.o.g.)  this results in 
\begin{equation}\label{eq:Teps_rffd_sigmaf_scaled}
		\E[T_{\epsilon}]=(s+1)\E[K_{\epsilon_\tau}] \leq   {\cal O}\left ( \frac{n^3L_1^3\sigma_f^2}{\epsilon^6}\right )
		\end{equation}	
which is the same total complexity bound of the unscaled case. Since the choice of $\alpha$ in the unscaled case is larger that the alpha we get here by the factor of $\frac{n}{s}$ this essentially means that, if we  scale up  the gradient estimator we scale down the step size by the same amount and the resulting algorithm is basically  the same. 

Similar derivations can be carried out for the bounds under Assumptions \ref{assum:KLcond} and \ref{assum:PLcond} with the similar  complexity bounds as those obtained under unscaled oracle \eqref{eq:RFFD}.

In a later section we will show that for the cases of \eqref{eq:SRFFD} with ${\cal U}\sim {\cal N}$ and  ${\cal U}\sim {\cal S}$  that are analyzed via smoothing argument we obtain an improved dependence on $n$ for the complexity bounds. 

\subsubsection{Randomized forward finite difference methods with orthonormal directions and multi-point stochastic zeroth-order oracle.}\label{sec:rffd_crn_complex}
We now turn to the first-order oracle generated by  \eqref{eq:RFFD} with ${\cal U}$ satisfying Assumption \ref{ass:haar} 
and with the access to the multi-point Oracle  \ref{def:oracle_stoch_0.2} for which  we assume $\epsilon_f=0$. 
As analyzed in  Section \ref{sec:RFFD_orth}  this first-order oracle   constitutes another  example of Oracle~\ref{def:oracle_stoch_1.3} but this time
with  $\epsilon_g^2 = s (\frac{L_1 {h}}{2} )^2$, $\kappa_g^2 = \frac{n-s}{n}$, $\sigma_c^2=sL_c^2$ and $\sigma_v^2=sL_v^2$. The cost of this oracle is proportional to $s+1$.  

This setting includes  the case of an unbiased stochastic zeroth-order under  common random numbers. 
The use of the multi-point oracle and the resulting bounds on  the error of the gradient estimates, imply significantly improved complexity bounds. 
Below we show how these bounds naturally come out of our analysis.  We want to note that in most of the literature we are aware of $L_v$ is assumed to be $0$, which is a more restrictive assumption than what we have here, yet $L_v>0$ is easily incorporated into the results and its value has smaller effect on complexity than $L_c$. 

We again apply Theorem \ref{thm:Case2_exp_complex} and observe that  $\sigma_c$ is now a constant that does not depend on $h$. Substituting this into  Theorem \ref{thm:Case2_exp_complex} we have the following result.

\begin{theorem} 
    Suppose Assumptions~\ref{assum:low_bound} and \ref{assum:lip_cont} hold. In Algorithm~\ref{alg:first_order_fixed_step_sgd}: let $\alpha\leq \min\{\frac{1}{L_1}, \frac{1}{2nL_1L_v^2}\}$ and  $g_k$ computed by \eqref{eq:RFFD}, under Assumption \ref{ass:haar}  with access to the multi-point Oracle  \ref{def:oracle_stoch_0.2} with $\epsilon_f=0$.
    
\begin{enumerate}
\item 
For any $\tau\in (0,1)$
and $\epsilon_\tau\geq \sqrt{\frac{2n}{\tau^2} \left((\frac{L_1 {h}}{2} )^2+\alpha L_1L_c^2 \right )}$, the stopping time 
$K_{\epsilon_\tau}=\min\{k\in\mathbb{N}:  \| \nabla \phi (X_k) \| \leq \epsilon_\tau\}$ is bounded as 
\begin{align*}
	\E[K_{\epsilon_\tau}] \leq  \frac{1}{\epsilon_\tau^2} \frac{4n(\phi(x_0) - \hat{\phi})}{ \alpha s (1 - \tau^2)}
\end{align*}
\item 
Suppose in addition  Assumption \ref{assum:KLcond} holds with a constant $c=\frac{1}{D}$, for some $D>0$.
For any $\tau\in (0,1)$
and $\epsilon_\tau\geq \sqrt{\frac{2nD^2}{\tau^2} \left((\frac{L_1 {h}}{2} )^2+\alpha L_1L_c^2\right )}$   the stopping time 
$K_{\epsilon_\tau}=\min\{k\in\mathbb{N}:   \phi (X_k) - \phi^* \leq \epsilon_\tau\}$ is bounded as 
\begin{align*}
	\E\left [K_{\epsilon_\tau}\right ] \leq \frac{1}{\epsilon_\tau} \frac{8nD^2}{\alpha s (1-\tau^2)}.
\end{align*}
\item If  Assumption~\ref{assum:PLcond} holds instead of \ref{assum:KLcond}, then, for any $\tau\in (0,1)$ and $\epsilon_\tau \geq \frac{n}{\mu \tau^2} \left((\frac{L_1 {h}}{2} )^2+\alpha L_1L_c^2\right)$,  the stopping time $K_{\epsilon_\tau}=\min\{k\in\mathbb{N}:   \phi (X_k) - \phi^* \leq \epsilon_\tau\}$ is bounded as 
\begin{align*}
    \E[K_{\epsilon_\tau}]\leq   \frac{2n}{\alpha s\mu(1-\tau^2)}\left (\log \frac{\phi(x_0)-\phi^*}{\epsilon_\tau}+1\right )
\end{align*}
\end{enumerate}
\end{theorem}

Let us see  how the complexity bounds for reaching arbitrary accuracy $\epsilon>0$ change as compared to the use of single-point stochastic oracle.  Specifically, given $\epsilon>0$ and a $\tau^2=\frac{1}{2}$ let us again 
set $\alpha$ and $h$ so that $(\frac{L_1 {h}}{2} )^2=\frac{\epsilon^2}{8n}$
and $\alpha L_1L_c^2=\frac{\epsilon^2}{8n}$. 
In this case we again set $h^2= \frac{\epsilon^2}{2nL_1^2}$ but now with $\alpha= \min\{ \frac{\epsilon^2}{8nL_1L_c^2}, \frac{1}{2nL_1L_v^2}, \frac{1}{L_1}\}$. Let us assume that the smallest term defining $\alpha$ is $\frac{\epsilon^2}{8nL_1L_c^2}$ (as this would be the most common case in most applications). 
Thus,  the resulting total sample complexity for reaching $\epsilon$-stationary point in the  general case of smooth functions 
is bounded  
\begin{equation}\label{eq:Teps_rffd_crn}
	\E[T_{\epsilon}] \leq {\cal O}\left ({\frac{n^2L_1L_c^2}{\epsilon^4} }\right ). 
\end{equation}

Note that in this case we can choose $h$ to be arbitrarily small without a negative affect on complexity. 
However, no matter how small $h$ is $\alpha$ would still need to be chosen to be smaller than 
$ \frac{\epsilon^2}{2nL_1L_c^2}$, so the positive effect on complexity by allowing larger $\alpha$ is up to a factor of two at the most.  

For the case when Assumption \ref{assum:KLcond} holds, select 
$h^2= \frac{\epsilon^2}{2nD^2L_1^2}$   and  $\alpha=\frac{ \epsilon^2}{8nD^2L_1L_c^2}$
we get total complexity bound for $T_{\epsilon}=(s+1)K_{\epsilon}$ where 
$K_{\epsilon}=\min\{k\in\mathbb{N}: \phi (X_k) - \phi^* \leq \epsilon\}$  as 
\begin{equation}\label{eq:Teps_rffd_crn_cvx}
	\E\left [T_{\epsilon} \right ] \leq  {\cal O}\left ({\frac{n^2D^4L_1L_c^2}{\epsilon^3} }\right ).
\end{equation}
And finally, under Assumption \ref{assum:PLcond}, we can choose $h^2\sim \frac{\mu\epsilon}{nL_1^2}$   and  $\alpha=\frac{\mu \epsilon}{4nL_1L_c^2}$  which gives final complexity as
\begin{equation}\label{eq:Teps_rffd_crn_strcvx}
    \E[T_{\epsilon}]\leq  {\cal O} \left ( \frac{n^2L_1L_c^2}{\mu^2\epsilon} \log \frac{1}{\epsilon}\right ). 
\end{equation}

The complexity obtained under Assumption \ref{assum:KLcond} is worse than that for convex functions (since convexity is a stronger assumption). 
To recover complexity under convexity we can apply the analysis in Section \ref{sec:telescopic_stoch}.

\paragraph{Extension to the scaled version.}
Let us now consider the first-order oracle generated by  
\eqref{eq:SRFFD} with access to the multi-point zeroth-order Oracle \ref{def:oracle_stoch_0.2}.  If we again assume that $\epsilon_f=0$ in  Oracle \ref{def:oracle_stoch_0.2},  the first-order oracle satisfies properties of Oracle~\ref{def:oracle_stoch_1.2} with $\epsilon_g^2=n\left( \frac{L_1 {h}}{2} \right)^2$, $\kappa_g=0$, $\sigma_c^2 = 
\frac{n^2}{s}L_c^2+\frac{n^2}{s} \left( \frac{L_1 {h}}{2}\right)^2 $  and $\sigma_v^2 = \frac{(2n-s)(n-s)}{ns}+\frac{n^2}{s}L_v^2$. 

Stepsize condition becomes $\alpha\leq \min\{\frac{1}{L_1}, \frac{ns}{2L_1(2n-s)(n-s)+n^3L_v^2}\}$. 
For simplicity, let us assume that $\left(\frac{L_1 {h}}{2}\right)^2$ in the expression for $\sigma_c^2$ is negligibly small compared to $L_c^2$. We can perform derivations as in the previous sections. We  observe that
the noise floor for the general smooth case is 
$\epsilon_\tau\geq \sqrt{ \frac{n}{\tau^2}\left (\left (\frac{L_1 {h}}{2}\right )^{2}+\frac{n}{s}\alpha L_1L_c^2 \right )}$ which is the same as the noise floor in the unscaled case with the second (variance) term being multiplied by $\frac{n}{s}$. The noise floor under Assumptions \ref{assum:KLcond}  and \ref{assum:PLcond} also has the second term  
 scaled by $\frac{n}{s}$. As a consequence, $h$ will be selected to be the same as in the unscaled case but $\alpha$ will be scaled down by $\frac{s}{n}$ and the complexity bounds 
 \eqref{eq:Teps_rffd_crn},  \eqref{eq:Teps_rffd_crn_cvx} and  \eqref{eq:Teps_rffd_crn_strcvx} hold for their corresponding cases.

\subsubsection{Randomized forward finite difference oracle with independent directions  and access to stochastic zeroth-order oracles. }
We now apply the same approach as in the previous sections to obtain results for the  fixed step size stochastic gradient descent based on any randomized finite difference method with independent directions that satisfy Assumption~\ref{ass:normal}, in other words, those  considered in Sections \ref{sec:RFFD} and \ref{sec:SRFFD}.
 These methods are typically used in the regime $s\ll n$, with $s=1$ being a popular choice. Here we will show how our analysis of expected complexity and its dependence on the noise constant $\epsilon_g$, $\kappa_g$, $\sigma_c$ and $\sigma_v$ reveals the key properties of this class of methods. 

\paragraph{The unscaled case}

 Recall that for the unscaled oracle defined by  \eqref{eq:RFFD} under Assumption~\ref{ass:normal}
we have 
\begin{align*}
    \epsilon_g^2 &= \frac{(n+s-1)}{n} \cdot s\, \Big(\frac{(1+\mathcal{O}(n^{-1}))L_1 h}{2}+ \frac{2\epsilon_f}{h} \Big)^2, \\
    \kappa_g^2 &= 1 - \frac{s}{n} + \frac{2s^2}{n^2} + \frac{\mathcal{O}(n^{-1})s}{n},  
\end{align*}
while $\sigma_c^2 = 2s\sigma_f^2 / {h}^2$, $\sigma_v = 0$ when using Oracle~\ref{def:oracle_stoch_0.1}, or 
\begin{equation*}
    \sigma_c^2 = \Big(\frac{n+s-1}{n} + \mathcal{O}(n^{-1})\Big) \cdot sL_c^2, \quad 
    \sigma_v^2 = \Big(\frac{n+s-1}{n} + \mathcal{O}(n^{-1})\Big) \cdot sL_v^2 
\end{equation*}

We observe that when $s\ll n$ the term  $\frac{(n+s-1)}{n} + \mathcal{O}(n^{-1})\approx 1 $, and 
   $\kappa_g^2 \approx 1 - \frac{s}{n}$ and thus we essentially obtain the same oracle constants 
   and in the case Assumption~\ref{ass:haar} in Section \ref{sec:RFFD_orth}.  This is not surprising because 
random  vectors from distributions  that satisfy Assumption~\ref{ass:normal} (e.g. Gaussian) 
in high dimension are nearly orthogonal. 
Thus in the  $s\ll n$  regime 
the analyses of Sections \ref{sec:RFFD_orth_complex_epsilonf}, \ref{sec:RFFD_orth_complex_sigmaf} and \ref{sec:rffd_crn_complex}  apply (approximately). 

 For the  general case  with $s< n$  we  can carry out a similar analysis. We first observe that $\kappa_g^2 <1$ is guaranteed as long as  $s<n-1$ . Let us assume as before that $\epsilon_f=0$.  To simplify the  discussion further,
and without a significant loss of generality,  we can assume that   $\frac{(n+s-1)}{n} + \mathcal{O}(n^{-1})\leq 2 $ and thus $\epsilon_g^2 \leq 2sL_1^2h^2$. 
Similarly we can assume  $\kappa_g^2 \leq  1 - \frac{s}{2n}$. Thus  $\epsilon_g^2$ and $1-\kappa_g^2$, $\sigma_c^2$ 
and $\sigma_v^2$ in the case of independent directions  differ from their counterparts for the orthogonal 
case only by a small constant factor. Thus  the complexity bounds from Sections \ref{sec:RFFD_orth_complex_epsilonf}, \ref{sec:RFFD_orth_complex_sigmaf} and \ref{sec:rffd_crn_complex} apply. Specifically,  \eqref{eq:Teps_rffd_sigmaf} holds under single-point oracle for general $\phi(x)$, \eqref{eq:Teps_rffd_sigmaf_cvx} 
holds under single-point oracle and Assumption \eqref{assum:KLcond}, while \eqref{eq:Teps_rffd_sigmaf_strcvx} 
holds under single-point oracle and Assumption \eqref{assum:PLcond}. 
Under multi-point oracle  \eqref{eq:Teps_rffd_crn} holds  for general $\phi(x)$, \eqref{eq:Teps_rffd_crn_cvx} holds under Assumption \eqref{assum:KLcond}, and 
\eqref{eq:Teps_rffd_crn_strcvx}  holds under Assumption \eqref{assum:PLcond}.

In Section \ref{sec:RFFD_orth_complex_epsilonf}  we considered \eqref{eq:RFFD} with orthogonal directions and with access to Oracle  \ref{def:oracle_stoch_0.1}   with $\sigma_f=0$.
This  case  the first-order oracle was treated as the case of Oracle~\ref{def:oracle_stoch_1.1}. The same can be applied in the case of independent directions and a close analogue of 
Theorem \ref{thm:Case1_exp_complex_rffd_det} can be stated. Thus Algorithm~\ref{alg:first_order_fixed_step_sgd} that uses \eqref{eq:RFFD} with independent directions and with access to Oracle  \ref{def:oracle_stoch_0.1}   with $\sigma_f=0$ enjoys deterministic-like bounds on the expected complexity, i.e.,  ${\cal O}\left ({\frac{1}{\epsilon^2} }\right )$ for the general smooth case,  ${\cal O}\left ({\frac{1}{\epsilon} }\right )$ under Assumption  \ref{assum:KLcond}
and ${\cal O}\left (\log \left ({\frac{1}{\epsilon} }\right )\right )$ under Assumption  \ref{assum:PLcond}. Consistent results for Gaussian distributions specifically we first developed in \cite{nesterov2017random}.  

Note that because the true dependence of the $\epsilon_g$ and $\kappa_g$ constants on $s$ is somewhat more complicated than their simplified bounds, 
the tighter bounds on resulting convergence neighborhood and complexity will turn out to depend on $s$, but only in a minor way, thus we ignore this dependency. 

\paragraph{Scaled case} 
Let us now consider \eqref{eq:SRFFD} for independent directions. In the case of general directions satisfying 
Assumption~\ref{ass:normal} we have pretty complicated expressions for the oracle constants (see Table \ref{tab:scaledffd}). However, for $s<<n$ these constants approximately equal to their counterparts for orthogonal directions. And thus general complexity results derived in Sections \ref{sec:RFFD_orth_complex_sigmaf} and \ref{sec:rffd_crn_complex} for  \eqref{eq:SRFFD} approximately apply.

We now turn to utilizing the improved bounds on $\epsilon_g^2$ derived in Section \ref{sec:smoothing}
under specific distributions $\mathcal{N}$  and $\mathcal{S}$. Let us consider \eqref{eq:SRFFD} under access to
single-point  Oracle~\ref{def:oracle_stoch_0.1} and let us assume again that it is unbiased, i.e.,   $\epsilon_f=0$.

We remind the reader here that for  the scaled version which is the vector $g(x)$ computed by  \eqref{eq:SRFFD} under distribution $\mathcal{N}$  we have the bias constants (due to Section \ref{sec:smoothing} )
\begin{align}\label{eq:const_rnd_dirder_scaled}
    \epsilon_g^2 = & ( L_1 h )^2 \text{ and }   \kappa_g^2 = 0  
                 \end{align}
while the variance constants are (due to $\bar s=s$ for $\mathcal{N}$  and $\mathcal{S}$)
\[ \begin{aligned}
    \sigma_c^2 &= \frac{2n^2\sigma_f^2}{s h^2} + \frac{n+s-1}{n} \cdot \frac{n^2}{s} \Big(\frac{(1+\mathcal{O}(n^{-1})) L_1 h}{2}  \Big)^2, \quad 
    \sigma_v^2 = \frac{(2n-s)(n-s)}{ns}
\end{aligned} \]

As in our previous similar analyses we assume that $L_1 h$ is negligible compared to $\frac{\sigma_f^2}{h^2}$ 
and thus our oracle constants are similar to the ones considered in Section \ref{sec:RFFD_orth_complex_sigmaf}
for \eqref{eq:SRFFD} except that $\epsilon_g^2$ is scaled by $\frac{4}{n}$.

Repeating the steps on that section,  
the first bound on  $\alpha \leq \min\{\frac{1}{L_1}, \frac{ns}{2(2n-s)(n-s)L_1}\}$ from Theorem \ref{thm:Case2_exp_complex_ffd_stoch} is unchanged. While the lower bound on $\epsilon_\tau$ 
for the case of general smooth $\phi(x)$ becomes 
$\epsilon_\tau\geq  \sqrt{\frac{1}{\tau^2} \left(({L_1 {h}} )^2+\alpha L_1\frac{n^2}{s}\frac{2\sigma_f^2}{h^2}\right )}$. 

For a given $\epsilon$ the optimal choice of $h$ and $\alpha$,  with $\tau^2=\frac{1}{2}$, now become as follows
  $h^2=\frac{\epsilon^2}{8L_2^2}$ and  $\alpha =  \frac{s\epsilon^2h^2}{16L_1n^2\sigma_f^2}= \frac{s\epsilon^4}{128n^2L_1^3\sigma_f^2}$. Assuming, as before, that this choice of $\alpha$ is smaller than $\min\{\frac{1}{L_1}, \frac{s}{2(2n-s)(n-s)L_1}\}$  this results in 
\begin{equation}\label{eq:Teps_rffd_sigmaf_scaled}
		\E[T_{\epsilon}]=(s+1)\E[K_{\epsilon_\tau}] \leq   {\cal O}\left ( \frac{n^2L_1^3\sigma_f^2}{\epsilon^6}\right )
		\end{equation}	
which is a better complexity by ${\cal O}(n)$ than all other results we obtained so far for any randomized 
forward finite difference scheme with access to the single point Oracle~\ref{def:oracle_stoch_0.1}.

We can also apply the the derivations for the case of  multipoint Oracle~\ref{def:oracle_stoch_0.2}, however, 
in this case no improvement in complexity will be observed, since the choice of $h$  does not have any influence on the choice of $\alpha$ and the the complexity is driven by the variance term rather than the bias term.

\subsubsection{Randomized central finite difference methods with orthogonal directions}\label{sec:gradproj_cfd_complex_sigmaf} 

We now want to turn to oracle generated by  the central finite difference \eqref{eq:RCFD} and  to highlight the difference in resulting complexity bounds for Algorithm~\ref{alg:first_order_fixed_step_sgd} as compared to those derived using forward finite differences.  
 The key difference between the oracles under \eqref{eq:RFFD} and \eqref{eq:RCFD} is in the bias term 
$\epsilon_g$ and its dependence on $h$. As we observe from the analysis of Algorithm~\ref{alg:first_order_fixed_step_sgd} under \eqref{eq:RFFD} the complexity is affected only under the single-point Oracle  \ref{def:oracle_stoch_0.1} with  $\sigma_f> 0$ because in other cases the choice of $\alpha$ is not affected by the choice of $h$. 

Here we consider  \eqref{eq:RCFD} with  ${\cal U}$ satisfying Assumption \ref{ass:haar} with access to  single-point Oracle  \ref{def:oracle_stoch_0.1}  which is assumed to be unbiased, i.e, $\epsilon_f=0$.

In this case, and under Assumption \ref{assum:lip_cont2},  the gradient estimate \eqref{eq:RCFD}, as analyzed in Section \ref{sec:RCFD},    constitutes an example of Oracle~\ref{def:oracle_stoch_1.3} with 
the constants  $\epsilon_g^2 =  s \cdot \min \left\{\left( \frac{L_1 {h}}{2}  \right)^2, \left( \frac{L_2 {h}^2}{6}  \right)^2\right\}$,  $\kappa_g^2 = \frac{n-s}{n}$,  $\sigma_c^2 = s\sigma_f^2 / (2{h}^2)$, and $\sigma_v = 0$.  Let us assume that $\epsilon_g^2 =  s  \left( \frac{L_2 {h}^2}{6}  \right)^2$. 
Theorem \ref{thm:Case2_exp_complex_ffd_stoch} holds with replacing $\frac{L_1 h}{2}$ terms in the expressions for $\epsilon_\tau$ by 
$ \frac{L_2 {h}^2}{6} $ with all other terms remaining the same. 
Thus the accuracy floor bound is  
$\epsilon_\tau\geq \sqrt{\frac{2n}{\tau^2} \left(\frac{L_2^2 h^4}{36}+\alpha L_1\frac{\sigma_f^2}{2h^2}\right )}$.  

As in the case of \eqref{eq:RFFD}  it is possible to find $\epsilon_\tau\leq \epsilon$ for any  $\epsilon>0$ by choosing $h$ and $\alpha$ appropriately, this time the choice of $h$ is different. 
Specifically, given $\epsilon>0$ and $\tau^2=\frac{1}{2}$  let us set 
$h^2=\frac{3\epsilon}{\sqrt{2n}L_2}$ and 
$\alpha =  \frac{3\epsilon^3}{4\sqrt{2}n^{\frac{3}{2}}L_1L_2\sigma_f^2}$. 

 Then for general smooth $\phi(x)$, $K_{\epsilon}=\min\{k\in \mathbb{N}:  \| \nabla \phi (X_k) \| \leq \epsilon \}$ becomes bounded as (using ${\cal O}$ to absorb small constants) 
\begin{align}\label{eq:Keps_rcfd_sigmaf}
	\E[K_{\epsilon}] \leq   {\cal O}\left ( \frac{n^{\frac{5}{2}}L_1L_2\sigma_f^2}{s\epsilon^5}\right ).  
\end{align}
and thus
\begin{equation}\label{eq:Teps_rcfd_sigmaf}
 \E[T_{\epsilon}]=2s \E[K_{\epsilon}]  \leq {\cal O}\left ( \frac{n^{\frac{5}{2}}L_1L_2\sigma_f^2}{\epsilon^5}\right ). 
\end{equation}

The complexity bounds for $\phi(x)$ satisfying Assumptions  \ref{assum:KLcond}  and \ref{assum:PLcond} improve accordingly, specifically 
when $\phi(x)$ satisfies Assumptions  \ref{assum:KLcond} 
\begin{equation}\label{eq:Teps_rcfd_sigmaf}
 \E[T_{\epsilon}] \leq {\cal O}\left ( \frac{D^5n^{\frac{5}{2}}L_1L_2\sigma_f^2}{\epsilon^4}\right ). 
\end{equation}
and when $\phi(x)$ satisfies Assumptions  \ref{assum:PLcond} 
\begin{equation}\label{eq:Teps_rcfd_sigmaf}
 \E[T_{\epsilon}] \leq {\cal O}\left ( \frac{n^{\frac{5}{2}}L_1L_2\sigma_f^2}{\mu^{\frac{5}{2}}\epsilon^{\frac{3}{2}}}\log \left (\frac{1}{\epsilon}\right )\right ). 
\end{equation}

\section{Inexact Second Order Oracles}\label{sec:second_order}
Second-order information can substantially improve the performance of optimization algorithms by enabling curvature-aware steps, improving scaling, and supporting fast local convergence \cite{NW,Conn-Gould-Toint-00,dennis1974characterization}. In nonconvex optimization, it can also be used to exploit negative curvature and establish stronger convergence guarantees, while carefully incorporating approximate curvature information can be beneficial even in the presence of noise \cite{curtis2019exploiting,cartis2012complexity}. Despite these advantages, obtaining and using second-order information can be computationally expensive and challenging, particularly in large-scale optimization and in settings where derivatives are unavailable or expensive to evaluate or contaminated with noise. Forming and storing Hessian matrices may be prohibitive, and solving the resulting Newton systems can introduce substantial additional computational costs. These difficulties motivate the use of approximate Hessians, Hessian-vector products, subsampling, iterative linear solvers, and quasi-Newton approximations \cite{xu2020newton,yao2023inexact,bollapragada2018exact}. However, inexact second-order information introduces challenges that do not arise in the same form for first-order methods. In particular, errors in Hessian approximations can affect not only the quality of a search direction but also the curvature information used to construct it. An inaccurate Hessian approximation may fail to identify directions of negative curvature or may introduce spurious negative curvature, potentially affecting both algorithmic performance and convergence guarantees.

These challenges become more pronounced when second-order information is used together with inexact function values and gradients. The accuracy requirements on the different quantities are often interdependent and may need to be coordinated with the current step, trust-region radius, regularization parameter, or desired stationarity tolerance \cite{xu2020newton,bellavia2021adaptive,yao2023inexact}. Consequently, understanding the accuracy required of second-order information, and how that accuracy affects convergence and computational cost, is an important component of the analysis and design of inexact optimization algorithms.

In the preceding sections, we studied optimization methods employing inexact first-order oracles, focusing primarily on their convergence to approximate first-order stationary points. The availability of second-order information introduces additional possibilities. In particular, suitable curvature information can be used to establish stronger convergence guarantees involving approximate second-order stationarity, characterized by
\begin{align*}
    \|\nabla\phi(x_k)\| \leq \epsilon_{1},
    \quad \text{and} \quad
    \lambda_{\min}(\nabla^2\phi(x_k))
    \geq -\epsilon_{2},
\end{align*}
for $\epsilon_{1}, \epsilon_{2} >0$.  Increasingly accurate derivative approximations may also lead to faster local convergence. These objectives generally impose different requirements on the accuracy of the available information. For example, an approximate Hessian that is sufficient to construct a useful descent direction need not be sufficiently accurate to identify negative curvature or support rapid local convergence.

In this section we review inexact second-order oracles, i.e. oracles that,  given $x\in\mathbb{R}^n$, return either an approximation of the Hessian of the objective function or its action along specified directions. The latter is particularly important in large-scale settings, where forming or storing a full Hessian matrix may be prohibitively expensive. We consider both deterministic and stochastic second-order oracle models, present representative constructions and examples, and discuss how the accuracy of this information affects the convergence and complexity of second-order optimization methods. Because second-order information can be incorporated through several substantially different algorithmic frameworks, our treatment is necessarily selective and focuses on the most common oracle models and accuracy requirements that arise in practice. In contrast to the first-order setting, we do not attempt to provide a unified analysis of algorithms based on these oracle models, as doing so would require a substantial separate treatment. Instead, we highlight representative results from the literature that analyze optimization methods under particular forms of inexact second-order information.

\subsection{Deterministic Inexact Second-Order Oracles}

We first consider deterministic inexact second-order oracles. We introduce two representative models, corresponding to approximations of the full Hessian and to Hessian-vector products along specified directions. These two oracle models capture the main ways in which inexact second-order information is used in optimization algorithms and will serve as the basis for the examples and discussion that follow.

\begin{customoracle}{D-2a}[General Deterministic Inexact Second-Order Oracle]
\label{def:oracle_det_2a}
Let $\psi:\mathbb{S}^n\times\mathbb{S}^n\to\mathbb{R}_+$ be a prescribed error measure, where $\mathbb{S}^n$ denotes the set of $n\times n$ symmetric matrices. Given a point $x\in\mathbb{R}^n$ and a control parameter $\lambda$, the oracle generates an output $H(x)\in\mathbb{S}^n$ that satisfies
\begin{equation}
\label{eq:2a_oracle_det}
\psi\bigl(H(x),\nabla^2\phi(x)\bigr)
\leq \epsilon_{H\psi}(\lambda),
\end{equation}
where $\epsilon_{H\psi}(\lambda)$ is a nonnegative deterministic function of $\lambda$. The computational cost of this oracle is some function ${\cal C}_{H}(\lambda;\psi)$.
\end{customoracle}

Common choices of the error measure $\psi$ include a matrix norm of the difference between the approximate and true Hessians and the absolute difference between selected spectral quantities, such as their smallest eigenvalues \cite{eisenstat1998three}. In large-scale applications, forming or storing a full Hessian approximation may be impractical, and many algorithms only require  products of the Hessian with specified vectors \cite{pearlmutter1994fast}. This motivates our second oracle.


\begin{customoracle}{D-2b}[Deterministic Inexact Hessian-Vector Product Oracle]
\label{def:oracle_det_2b}
Given a point $x\in\mathbb{R}^n$, a vector $v\in\mathbb{R}^n$, and a control parameter $\lambda$, the oracle generates an output $h(x,v)\in\mathbb{R}^n$ that satisfies
\begin{equation} \label{eq:2b_oracle_det}
	\|h(x,v)-\nabla^2\phi(x)v\| \leq \kappa_{Hv}(\lambda)\|v\|^p,
\end{equation}
where $p\geq1$ and $\kappa_{Hv}(\lambda)$ is a nonnegative deterministic function of $\lambda$. The computational cost of this oracle is some function ${\cal C}_{Hv}(\lambda)$.
\end{customoracle}

Oracle~\ref{def:oracle_det_2b} can be weaker than Oracle~\ref{def:oracle_det_2a} in the sense that it controls the accuracy of second-order information only along the queried directions and does not require the existence of a single approximate Hessian matrix underlying all oracle outputs. This distinction is important in large-scale optimization, where Hessian-vector products may be available much more cheaply than explicit Hessian approximations. On the other hand, directional accuracy alone generally does not provide global spectral information about the Hessian approximation. In the special case where $h(x,v)=H(x)v$ for a common symmetric matrix $H(x)$ and the bound with $p=1$ holds uniformly for all $v$, Oracle~\ref{def:oracle_det_2b} reduces to an operator-norm bound on $H(x)-\nabla^2\phi(x)$ and therefore implies an instance of Oracle~\ref{def:oracle_det_2a}.

\subsection{Examples of Deterministic Inexact Second-Order Oracles}

We next present several representative examples of deterministic inexact second-order oracles that fit the frameworks of Oracles~\ref{def:oracle_det_2a} and \ref{def:oracle_det_2b}. These include finite-difference approximations based on zeroth- and first-order information, quadratic interpolation and simplex Hessian constructions, the Gauss--Newton approximation, and directional finite differences for Hessian-vector products. For each example, we characterize the resulting approximation error and discuss the associated computational cost. Since second-order oracle constructions are not the primary focus of this review, our treatment is intentionally brief and, where appropriate, detailed derivations are omitted.


\subsubsection{Forward Finite Differences (FFD) via Zeroth-Order Oracles}
  
When only zeroth-order information is available, an approximation of the Hessian can be constructed by applying forward finite differences directly to function values. Given access to Oracle~\ref{def:oracle_det_0}, for any $x\in\mathbb{R}^n$, choose $h>0$ and evaluate $f(y)$ for all $y\in{\cal Y}
=
\{x\}
\cup
\{x+he_i:i=1,\ldots,n\}
\cup
\{x+he_i+he_j:1\leq i\leq j\leq n\}$, where $e_i$ denotes the $i$th column of the identity matrix. A Hessian approximation is then given by
\begin{equation}
\label{eq:ffd2}
H(x)
=
\sum_{i=1}^n\sum_{j=1}^n
\frac{
f(x+he_i+he_j)-f(x+he_i)-f(x+he_j)+f(x)
}{h^2}
e_i e_j^T.
\end{equation}

The following proposition provides an error bound for the Hessian approximation in \eqref{eq:ffd2}.

\begin{proposition}
\label{prop:finite_diff2}
Let Assumption~\ref{assum:lip_cont2} hold, and suppose that the zeroth-order oracle evaluations used in \eqref{eq:ffd2} satisfy Oracle~\ref{def:oracle_det_0} with error bound $\epsilon_f$. Then, for any $x\in\mathbb{R}^n$ and $h>0$, the Hessian approximation $H(x)$ defined in \eqref{eq:ffd2} satisfies
\begin{align*}
\|H(x)-\nabla^2\phi(x)\|_2
\leq
\left(
1+
\frac{\left(n-1\right)\left(\sqrt{2}+1\right)}{3}
\right)L_2h
+
\frac{4n\epsilon_f}{h^2}.
\end{align*}
\end{proposition}
\begin{proof}
Let $\widehat H(x)$ denote the finite-difference Hessian approximation obtained from \eqref{eq:ffd2} by replacing each evaluation of $f$ with the corresponding exact function value $\phi$. For any $i,j=1,\ldots,n$,
\begin{align*}
	e_i^T\widehat H(x)e_j=\frac{\phi\left(x+he_i+he_j\right)-\phi\left(x+he_i\right)-\phi\left(x+he_j\right)+\phi(x)}{h^2}.
\end{align*}
Thus, to bound the $\left(i,j\right)$th entry of
$\widehat H(x)-\nabla^2\phi(x)$, we consider Taylor expansions of the three shifted function values appearing in the numerator above.

By Assumption~\ref{assum:lip_cont2}, for any $s\in\mathbb{R}^n$,
\begin{align*}
	\phi(x+s)=\phi(x)+\nabla\phi(x)^Ts+\frac{1}{2}s^T\nabla^2\phi(x)s+r(s),
\end{align*}
where $\left|r(s)\right| \leq \frac{L_2}{6}\|s\|^3$.

For $i\neq j$, we apply this expansion to the three displacement vectors $s=he_i+he_j$, $s=he_i$, and $s=he_j$, corresponding, to $\phi\left(x+he_i+he_j\right)$, $\phi\left(x+he_i\right)$, and $\phi\left(x+he_j\right)$, respectively. Substituting these three expansions into the expression for $e_i^T\widehat H(x)e_j$, the constant and linear terms cancel. Moreover, using the symmetry of $\nabla^2\phi(x)$, the quadratic terms satisfy
\begin{align*}
&\frac{1}{2}
\left(he_i+he_j\right)^T
\nabla^2\phi(x)
\left(he_i+he_j\right)
-
\frac{1}{2}h^2e_i^T\nabla^2\phi(x)e_i
-
\frac{1}{2}h^2e_j^T\nabla^2\phi(x)e_j
\\
&\qquad =
h^2e_i^T\nabla^2\phi(x)e_j.
\end{align*}
Therefore,
\begin{align*}
e_i^T
\left(\widehat H(x)-\nabla^2\phi(x)\right)
e_j
=
\frac{
r\left(he_i+he_j\right)
-r\left(he_i\right)
-r\left(he_j\right)
}{h^2}.
\end{align*}
Hence,
\begin{align*}
\left|
e_i^T
\left(\widehat H(x)-\nabla^2\phi(x)\right)
e_j
\right|
&\leq
\frac{L_2}{6h^2}
\left(
\|he_i+he_j\|^3
+
\|he_i\|^3
+
\|he_j\|^3
\right)
\\
&=
\frac{L_2}{6h^2}
\left(
2\sqrt{2}h^3+h^3+h^3
\right)
\\
&=
\frac{\sqrt{2}+1}{3}L_2h.
\end{align*}

For $i=j$, the corresponding diagonal entry of $\widehat H(x)$ is
\begin{align*}
e_i^T\widehat H(x)e_i
=
\frac{
\phi\left(x+2he_i\right)
-2\phi\left(x+he_i\right)
+\phi(x)
}{h^2}.
\end{align*}
By two applications of the fundamental theorem of calculus,
\begin{align*}
e_i^T\widehat H(x)e_i
=
\int_0^1\int_0^1
e_i^T
\nabla^2\phi\left(x+h\left(s+t\right)e_i\right)
e_i
\,ds\,dt.
\end{align*}
Therefore,
\begin{align*}
e_i^T
\left(\widehat H(x)-\nabla^2\phi(x)\right)e_i
=
\int_0^1\int_0^1
e_i^T
\left(
\nabla^2\phi\left(x+h\left(s+t\right)e_i\right)
-\nabla^2\phi(x)
\right)
e_i
\,ds\,dt.
\end{align*}
By Assumption~\ref{assum:lip_cont2}, it follows that
\begin{align*}
\left|
e_i^T
\left(\widehat H(x)-\nabla^2\phi(x)\right)e_i
\right|
&\leq
\int_0^1\int_0^1
\left\|
\nabla^2\phi\left(x+h\left(s+t\right)e_i\right)
-\nabla^2\phi(x)
\right\|_2
\,ds\,dt
\\
&\leq
L_2h
\int_0^1\int_0^1
\left(s+t\right)
\,ds\,dt
\\
&=
L_2h.
\end{align*}

From the preceding arguments, 
\begin{align*}
\left|
e_i^T
\left(\widehat H(x)-\nabla^2\phi(x)\right)e_j
\right|
\leq
\begin{cases}
L_2h, & i=j,\\[2mm]
\frac{\sqrt{2}+1}{3}L_2h, & i\neq j.
\end{cases}
\end{align*}
Since 
$\widehat H(x)-\nabla^2\phi(x)$ is symmetric,
\begin{align}\label{eq.bound1}
	\|\widehat H(x)-\nabla^2\phi(x)\|_2\leq\|\widehat H(x)-\nabla^2\phi(x)\|_\infty\leq\left(1+\frac{\left(n-1\right)\left(\sqrt{2}+1\right)}{3}\right)L_2h.
\end{align}

We next bound the error due to the inexact zeroth-order oracle. For any $i,j=1,\ldots,n$, the $\left(i,j\right)$th entry of $H(x)-\widehat H(x)$ satisfies
\begin{align*}
&
\left|
e_i^T
\left(H(x)-\widehat H(x)\right)
e_j
\right|
\\
&\leq
\frac{1}{h^2}
\left(
|f(x+he_i+he_j)-\phi(x+he_i+he_j)|
+
|f(x+he_i)-\phi(x+he_i)|
\right.
\\
&\hspace{3cm}
\left.
+
|f(x+he_j)-\phi(x+he_j)|
+
|f(x)-\phi(x)|
\right)
\\
&\leq
\frac{4\epsilon_f}{h^2}.
\end{align*}
Therefore,
\begin{align}\label{eq.bound2}
	\|H(x)-\widehat H(x)\|_2 \leq\|H(x)-\widehat H(x)\|_F\leq\frac{4n\epsilon_f}{h^2}.
\end{align}

Combining \eqref{eq.bound1} and \eqref{eq.bound2}
and applying the triangle inequality yields
\begin{align*}
\|H(x)-\nabla^2\phi(x)\|_2
&\leq
\|H(x)-\widehat H(x)\|_2
+
\|\widehat H(x)-\nabla^2\phi(x)\|_2
\\
&\leq
\frac{4n\epsilon_f}{h^2}
+
\left(
1+
\frac{\left(n-1\right)\left(\sqrt{2}+1\right)}{3}
\right)L_2h,
\end{align*}
which proves the result.
\end{proof}

We now discuss the result in Proposition~\ref{prop:finite_diff2} in terms of Oracle~\ref{def:oracle_det_2a}. In the simplest case where Oracle~\ref{def:oracle_det_0} is not adjustable, the finite-difference interval $h>0$ serves as the control parameter, so that $\lambda=h$. Taking $\psi\left(A,B\right)=\|A-B\|_2$, Proposition~\ref{prop:finite_diff2} implies that the approximation satisfies Oracle~\ref{def:oracle_det_2a} with
\begin{align*}
\epsilon_{H\psi}(h)
=
\left(
1+
\frac{\left(n-1\right)\left(\sqrt{2}+1\right)}{3}
\right)L_2h
+
\frac{4n\epsilon_f}{h^2}.
\end{align*}
The construction requires $\frac{\left(n+1\right)\left(n+2\right)}{2}$ distinct evaluations of Oracle~\ref{def:oracle_det_0}. Thus, if each zeroth-order oracle evaluation has cost ${\cal C}_f$, the total computational cost is ${\cal C}_H\left(h;\psi\right)=\frac{\left(n+1\right)\left(n+2\right)}{2}{\cal C}_f+ \mathcal{O}\left(n^2\right)$.

Furthermore, if $h$ is chosen to satisfy $h\geq\left(\frac{4n\epsilon_f}{\left(1+\frac{\left(n-1\right)\left(\sqrt{2}+1\right)}{3}\right)L_2}\right)^{1/3}$,
then $\epsilon_{H\psi}
\leq
2
\left(
1+
\frac{\left(n-1\right)\left(\sqrt{2}+1\right)}{3}
\right)L_2h$. This upper bound is minimized by choosing $h$ as the smallest value permitted by the above condition.

 \subsubsection{Forward Finite Differences (FFD) via First-Order Oracles}
 
 Given access to Oracle~\ref{def:oracle_det_1}, an approximation of the Hessian can be constructed by applying forward finite differences to the first-order oracle outputs. For any $x\in\mathbb{R}^n$, choose $h>0$ and evaluate $g(y)$ for all $y\in{\cal Y}=\{x\}\cup\{x+he_i=1,\ldots,n\}$, where $e_i$ denotes the $i$th column of the identity matrix. A symmetric Hessian approximation is then given by
 \begin{equation}
\label{eq:ffd2_grad}
\begin{aligned}
H(x)
&=
\sum_{i=1}^n\sum_{j=1}^n
\frac{
(g(x+he_i)-g(x))^Te_j
+
(g(x+he_j)-g(x))^Te_i
}{2h}
e_ie_j^T\\
&=
\sum_{i=1}^n
\left[
\left(\frac{g(x+he_i)-g(x)}{2h}\right)e_i^T
+
e_i\left(\frac{g(x+he_i)-g(x)}{2h}\right)^T
\right].
\end{aligned}
\end{equation}
The symmetrization in \eqref{eq:ffd2_grad} ensures that $H(x)\in\mathbb{S}^n$, even though the finite-difference approximation obtained directly from the inexact first-order oracle need not be symmetric.

The following proposition provides an error bound for the Hessian approximation in \eqref{eq:ffd2_grad}.

\begin{proposition}
\label{prop:hessian_ffd1}
Let Assumptions~\ref{assum:lip_cont} and~\ref{assum:lip_cont2} hold, and suppose that the first-order oracle evaluations used in \eqref{eq:ffd2_grad} satisfy Oracle~\ref{def:oracle_det_1} with parameters $\epsilon_g$ and $\kappa_g$. Then, for any $x\in\mathbb{R}^n$ and $h>0$, the Hessian approximation $H(x)$ defined in \eqref{eq:ffd2_grad} satisfies
\begin{align*}
\|H(x)-\nabla^2\phi(x)\|_2
\leq
\frac{2\sqrt{n}\left(\epsilon_g+\kappa_g\|\nabla\phi(x)\|\right)}{h}
+\sqrt{n}\kappa_gL_1
+\frac{\sqrt{n}L_2h}{2}.
\end{align*}
\end{proposition}
\begin{proof}
For brevity, we suppress the dependence of $\epsilon_g$ and $\kappa_g$ on the control parameter $\lambda_g$ of Oracle~\ref{def:oracle_det_1}. Since $\nabla^2\phi(x)$ is symmetric, by the properties of the Frobenius norm,
\begin{align*}  
		\left\|H(x) - \nabla^2\phi(x)\right\| \le{}& \left\|H(x) - \nabla^2\phi(x)\right\|_F \\
		\le{}& \frac{1}{2} \left\| \sum_{i=1}^d \frac{g(x+he_i) - g(x)}{h} e_i^T- \nabla^2\phi(x) \right\|_F \\
		&+ \frac{1}{2} \bigg\| \sum_{i=1}^d e_i \Big(\frac{g(x+he_i) - g(x)}{h}\Big)^T- \nabla^2\phi(x) \bigg\|_F \\
		={}& \left\| \sum_{i=1}^d \frac{g(x+he_i) - g(x)}{h} e_i^T- \nabla^2\phi(x) \right\|_F \\
		={}& \sqrt{ \sum_{i=1}^d \left\|\frac{g(x+he_i) - g(x)}{h} - \nabla^2\phi(x) e_i\right\|^2 }. 
	\end{align*}
	For every $i = 1,\dots,n$, 
	\begin{align*}
		&\left\|\frac{g(x+he_i) - g(x)}{h} - \nabla^2\phi(x) e_i\right\|^2 \\
		&= \frac{1}{h^2} \Big\| \left[g(x+he_i) - \nabla\phi(x+he_i)\right] - \left[g(x) - \nabla\phi(x)\right]  \\
		& \qquad +  \left[\nabla\phi(x+he_i) - \nabla\phi(x) - \nabla^2\phi(x) he_i\right]  \Big\|^2 \\
		&\leq \frac{1}{h^2} \Big( \left\| g(x+he_i) - \nabla\phi(x+he_i)\right\| +  \left\| g(x) - \nabla\phi(x) \right\| \\
		& \qquad +  \left\| \nabla\phi(x+he_i) - \nabla\phi(x) - \nabla^2\phi(x) he_i \right\| \Big)^2 \\
		&\leq \frac{1}{h^2} \Big( \epsilon_g+\kappa_g\|\nabla\phi(x+he_i)\| + \epsilon_g+\kappa_g\|\nabla\phi(x)\| + \frac{L_2 h^2} {2} \Big)^2 \\
		&\leq \frac{1}{h^2} \Big( 2\epsilon_g+2\kappa_g\|\nabla\phi(x)\| + \kappa_gL_1 h + \frac{L_2 h^2} {2} \Big)^2. 
	\end{align*} 
	The second inequality above holds by the property of Oracle~\ref{def:oracle_det_1} and Assumption~\ref{assum:lip_cont2}, and the last inequality follows by Assumption~\ref{assum:lip_cont}. 
	
	Combining the above two inequalities yields the result. 
\end{proof}

We now discuss the result in Proposition~\ref{prop:hessian_ffd1} in terms of Oracle~\ref{def:oracle_det_2a}. The construction uses $n+1$ outputs of Oracle~\ref{def:oracle_det_1}, evaluated at $x$ and $x+he_i$, $i=1,\ldots,n$, with $h>0$ serving as the finite-difference control parameter. In this case, the error measure is $\psi(A,B)=\|A-B\|_2$. When $\kappa_g>0$, the bound in Proposition~\ref{prop:hessian_ffd1} depends on $\|\nabla\phi(x)\|$, which may vary substantially depending on the proximity to stationarity. As seen in our analysis of first-order methods, separating the error term proportional to $\|\nabla\phi(x)\|$ from the constant error term can lead to sharper convergence neighborhoods. A similar distinction may also be beneficial in the analysis of second-order methods, although such gradient-dependent second-order error models have received limited attention. Moreover, it is unclear how frequently finite-difference Hessian approximations such as \eqref{eq:ffd2_grad} are constructed from first-order oracles with $\kappa_g>0$ in practice. Our definition of Oracle~\ref{def:oracle_det_2a} here does not treat this situation as special. Thus, when $\kappa_g>0$ we can simply consider a bound  $\|\nabla\phi(x)\|\leq G$ over the region of interest, for some $G>0$. Taking $\lambda=h$, the approximation satisfies Oracle~\ref{def:oracle_det_2a} with
\begin{align*}
\epsilon_{H\psi}\left(h\right) = \frac{2\sqrt{n}\left(\epsilon_g+\kappa_g G\right)}{h}
+\sqrt{n}\kappa_gL_1
+\frac{\sqrt{n}L_2h}{2},
\end{align*}
and the associated cost is ${\cal C}_{H}(h;\psi) = (n+1) {\cal C}_{g}(h) + \mathcal{O}(n^2)$. 

Furthermore, if $h$ is chosen to satisfy $h \geq
    2\sqrt{\frac{\epsilon_g+\kappa_g G}{L_2}}$, then $\epsilon_{H\psi}(h)
    \leq
    \sqrt{n}\kappa_gL_1+\sqrt{n}L_2h$. This upper bound is minimized by choosing $h$ as the smallest value permitted by the above condition. Moreover, if $h$ is chosen to satisfy $h \geq
    \max\left\{
        2\sqrt{\frac{\epsilon_g+\kappa_g G}{L_2}},
        \frac{\kappa_gL_1}{L_2}
    \right\}$, then $\epsilon_{H\psi}\leq 2\sqrt{n}L_2h$.

\subsubsection{Central Finite Differences (CFD) via Zeroth-Order and First-Order Oracles}

As in the first-order setting, central finite differences can be used to construct Hessian approximations from either zeroth-order or first-order oracle information. Compared with forward differences, central differencing generally yields higher-order accuracy in the finite-difference interval, at the expense of additional oracle evaluations. We briefly present both constructions, state representative error bounds, and discuss how they fit within the framework of Oracle~\ref{def:oracle_det_2a}.

Given access to Oracle~\ref{def:oracle_det_0}, for any $x\in\mathbb{R}^n$, choose $h>0$. The diagonal entries of the Hessian can be approximated using the central second-difference formula, while the off-diagonal entries can be approximated using central mixed differences. This yields the symmetric Hessian approximation
\begin{align}
H(x)
={}&
\sum_{i=1}^n
\frac{
f\left(x+he_i\right)
-2f(x)
+f\left(x-he_i\right)
}{h^2}
e_i e_i^T\nonumber\\
&+
\frac{1}{4h^2}\sum_{i=1}^n\sum_{j=i+1}^n
\left(f\left(x+he_i+he_j\right)
-f\left(x+he_i-he_j\right)\right.
\label{eq:cfd_zeroth}\\
& \left.\qquad \qquad \qquad \qquad -f\left(x-he_i+he_j\right)
+f\left(x-he_i-he_j\right)\right)
\left(
e_i e_j^T+e_j e_i^T
\right).\nonumber
\end{align}

The following proposition provides an error bound for the central finite-difference Hessian approximation via zeroth-order information. To obtain an improved dependence on the finite-difference interval $h$, we impose the stronger assumption that $\nabla^3\phi(x)$ is Lipschitz continuous.

\begin{proposition}
\label{prop:cfd_hessian}
Suppose that $\nabla^3\phi(x)$ is $L_3$-Lipschitz continuous, and suppose that the zeroth-order oracle evaluations used in \eqref{eq:cfd_zeroth} satisfy Oracle~\ref{def:oracle_det_0} with error bound $\epsilon_f$. Then, for any $x\in\mathbb{R}^n$ and $h>0$, the Hessian approximation $H(x)$ defined in \eqref{eq:cfd_zeroth} satisfies
\begin{align*}
\|H(x)-\nabla^2\phi(x)\|_2
\leq
\frac{\left(2n-1\right)L_3h^2}{12}
+
\frac{\left(n+3\right)\epsilon_f}{h^2}.
\end{align*}
\end{proposition}

The bound in Proposition~\ref{prop:cfd_hessian} improves upon the corresponding forward finite-difference bound, with the truncation error scaling quadratically rather than linearly with the finite-difference interval $h$. This improvement requires the stronger assumption that $\nabla^3\phi$ is Lipschitz continuous. If only Assumption~\ref{assum:lip_cont2} holds, i.e., if the Hessian is Lipschitz continuous, then the truncation error is generally of order $h$ rather than $h^2$.

We now discuss the result in Proposition~\ref{prop:cfd_hessian} in terms of Oracle~\ref{def:oracle_det_2a}. In the simplest case where Oracle~\ref{def:oracle_det_0} is not adjustable, the finite-difference interval $h>0$ serves as the control parameter, so that $\lambda=h$. Taking $\psi\left(A,B\right)=\|A-B\|_2$, Proposition~\ref{prop:cfd_hessian} implies that the approximation satisfies Oracle~\ref{def:oracle_det_2a} with
\begin{align*}
\epsilon_{H\psi}\left(h\right)
=
\frac{\left(2n-1\right)L_3h^2}{12}
+
\frac{\left(n+3\right)\epsilon_f}{h^2}.
\end{align*}
If $h \geq \left(\frac{12\left(n+3\right)\epsilon_f}
{\left(2n-1\right)L_3}\right)^{1/4}$, then $\epsilon_{H\psi} \leq \frac{\left(2n-1\right)L_3h^2}{6}$. 
The construction requires $2n^2+1$ distinct evaluations of Oracle~\ref{def:oracle_det_0}. Thus, if each zeroth-order oracle evaluation has cost ${\cal C}_f$, the total computational cost is ${\cal C}_H\left(h;\psi\right)=\left(2n^2+1\right){\cal C}_f+\mathcal{O}\left(n^2\right)$.

Next, given access to Oracle~\ref{def:oracle_det_1}, for any $x\in\mathbb{R}^n$, choose $h>0$ and evaluate $g\left(x+he_i\right)$ and $g\left(x-he_i\right)$, $i=1,\ldots,n$. A symmetric Hessian approximation is then given by
\begin{align}\label{eq:cfd_first}
	H(x)=\sum_{i=1}^n\left[\frac{g\left(x+he_i\right)-g\left(x-he_i\right)}{4h}e_i^T+e_i\left(\frac{g\left(x+he_i\right)-g\left(x-he_i\right)}{4h}\right)^T\right].
\end{align}
The following proposition provides an error bound for the central finite-difference Hessian approximation via first-order information.

\begin{proposition}
\label{prop:cfd2_first}
Suppose that $\nabla^3\phi$ is $L_3$-Lipschitz continuous, and suppose that the first-order oracle evaluations used in \eqref{eq:cfd_first} satisfy Oracle~\ref{def:oracle_det_1} with error parameters $\epsilon_g$ and $\kappa_g$. Then, for any $x\in\mathbb{R}^n$ and $h>0$, the Hessian approximation $H(x)$ defined in \eqref{eq:cfd_first} satisfies
\begin{align*}
\|H(x)-\nabla^2\phi(x)\|_2
\leq
\frac{\sqrt{n}\left(\epsilon_g+\kappa_g\|\nabla\phi(x)\|\right)}{h}
+
\sqrt{n}\kappa_gL_1 + \frac{\sqrt{n}L_3h^2}{6}.
\end{align*}
\end{proposition}

The result in Proposition~\ref{prop:cfd2_first} improves upon the corresponding forward finite-difference bound, with the truncation error depending quadratically rather than linearly on the finite-difference interval $h$. As in the zeroth-order setting, this improvement relies on the stronger assumption that $\nabla^3\phi(x)$ is Lipschitz continuous; if only Assumption~\ref{assum:lip_cont2} holds, then the truncation term is generally of order $h$ rather than $h^2$.

To position Proposition~\ref{prop:cfd2_first} in terms of Oracle~\ref{def:oracle_det_2a} in case when $\kappa_g>0$, as in the previous section, we consider a bound  $\|\nabla\phi(x)\|\leq G$ over the region of interest, for some $G>0$. Taking $\lambda=h$ and $\psi\left(A,B\right)=\|A-B\|_2$, the approximation satisfies Oracle~\ref{def:oracle_det_2a} with
\begin{align*}
\epsilon_{H\psi}\left(h\right)
=
\frac{\sqrt{n}L_3h^2}{6}
+
\frac{\sqrt{n}\left(\epsilon_g+\kappa_g G\right)}{h}
+
\sqrt{n}\kappa_gL_1.
\end{align*}
Furthermore, if
$h\geq\left(\frac{6\left(\epsilon_g+\kappa_gG\right)}{L_3}\right)^{1/3}$,
then $\epsilon_{H\psi}\left(h\right)\leq\frac{\sqrt{n}L_3h^2}{3}+\sqrt{n}\kappa_gL_1$. The construction requires $2n$ evaluations of Oracle~\ref{def:oracle_det_1}; thus, if each first-order oracle evaluation has cost ${\cal C}_g$, then ${\cal C}_H\left(h;\psi\right)=2n{\cal C}_g+\mathcal{O}(n^2)$.

\subsubsection{Quadratic Interpolation via Zeroth-Order Oracles}

Quadratic interpolation extends linear interpolation by constructing a local quadratic model from zeroth-order information. Given access to Oracle~\ref{def:oracle_det_0}, let ${\cal U}=\{u_1,\ldots,u_p\}\subset\mathbb{R}^n$, where $p=\frac{(n+1)(n+2)}{2}-1$, and consider the interpolation set
$\{x\}\cup\{x+u_i:i=1,\ldots,p\}$.
Quadratic interpolation constructs a model
\begin{equation}
\label{eq:quadratic_model}
	m(y)=f(x)+g(x)^T(y-x)+\frac{1}{2}(y-x)^TH(y-x),
\end{equation}
where $g(x)\in\mathbb{R}^n$ and $H(x)\in\mathbb{S}^n$, by imposing the interpolation conditions
\begin{align}
\label{eq:quadratic_interp_cond}
	m(x+u_i)=f(x)+g(x)^Tu_i+\frac{1}{2}u_i^TH(x) u_i=f(x+u_i),
\qquad
i=1,\ldots,p.
\end{align}

To express these conditions as a linear system, define
\begin{align*}
\operatorname{vec}(H(x))
=
\left[
\frac{1}{\sqrt{2}}H_{11},
H_{12},\ldots,H_{1n},
\frac{1}{\sqrt{2}}H_{22},
H_{23},\ldots,
\frac{1}{\sqrt{2}}H_{nn}
\right]^T
\in
\mathbb{R}^{\frac{n(n+1)}{2}},
\end{align*}
and, similarly,
\begin{align*}
\operatorname{vec}(uu^T)
=
\left[
\frac{1}{\sqrt{2}}u_1^2,
u_1u_2,\ldots,u_1u_n,
\frac{1}{\sqrt{2}}u_2^2,
u_2u_3,\ldots,
\frac{1}{\sqrt{2}}u_n^2
\right]^T.
\end{align*}
With $v_i=\operatorname{vec}(u_i u_i^T)$, the interpolation conditions can be written as
\begin{equation}
\label{eq:quadratic_linear_system}
M
\begin{bmatrix}
g\\
\operatorname{vec}(H(x))
\end{bmatrix}
=
f_\Delta,
\end{equation}
where
\begin{align*}
M
=
\begin{bmatrix}
u_1^T& v_1^T\\
\vdots & \vdots\\
u_p^T& v_p^T
\end{bmatrix}
\in\mathbb{R}^{p\times p}, \qquad \text{and} \qquad
f_\Delta
=
\begin{bmatrix}
f(x+u_1)-f(x)\\
\vdots\\
f(x+u_p)-f(x)
\end{bmatrix}.
\end{align*}
We assume that the interpolation set is poised for quadratic interpolation, so that $M$ is nonsingular. Finally, suppose that $\max_{i=1,\ldots,p}\|u_i\|\leq h$. To separate the geometry of the interpolation set from its scale, define
\begin{align*}
\widehat M
=
M
\begin{bmatrix}
\frac{1}{h}I_{n\times n} & 0\\
0 &
\frac{1}{h^2}I_{\frac{n(n+1)}{2}\times\frac{n(n+1)}{2}}
\end{bmatrix}.
\end{align*}

The following proposition extends standard quadratic interpolation error bounds from the noise-free setting, $\epsilon_f=0$, studied in, e.g., \cite{ConnScheVice08c}, to the case of inexact zeroth-order information.
\begin{proposition}
\label{prop:quadratic_interp}
Let Assumption~\ref{assum:lip_cont2} hold, and suppose that the zeroth-order oracle evaluations used to construct the quadratic interpolation model \eqref{eq:quadratic_model} satisfy Oracle~\ref{def:oracle_det_0} with error bound $\epsilon_f$. Moreover, suppose that $\max_{i=1,\ldots,p}\|u_i\|\leq h$ where $p=\frac{(n+1)(n+2)}{2}-1$, and that the interpolation set is poised for quadratic interpolation. Then, the coefficients $g$ and $H$ of the quadratic interpolation model satisfy
\begin{align*}
	\|g(x)-\nabla\phi(x)\|&\leq\sqrt{p}\|\widehat M^{-1}\|_2\left(\frac{L_2h^2}{6}+\frac{2\epsilon_f}{h}\right),\\
	\|H(x)-\nabla^2\phi(x)\|_2&\leq\sqrt{2p}\|\widehat M^{-1}\|_2\left(\frac{L_2h}{6}+\frac{2\epsilon_f}{h^2}\right).
\end{align*}
\end{proposition}

When $\epsilon_f=0$, the bounds in Proposition~\ref{prop:quadratic_interp} recover the standard noise-free quadratic interpolation error bounds; see e.g., \cite{ConnScheVice08c}.

The quadratic interpolation model can be viewed as generating outputs of both Oracles~\ref{def:oracle_det_1} and~\ref{def:oracle_det_2a}, with the interpolation set serving as the input $\lambda$. Specifically, the gradient approximation satisfies Oracle~\ref{def:oracle_det_1} with
\begin{align*}
\epsilon_g
&=
\sqrt{p}\|\widehat M^{-1}\|_2
\left(
\frac{L_2h^2}{6}
+
\frac{2\epsilon_f}{h}
\right)
\qquad \text{and} \qquad 
\kappa_g=0,
\end{align*}
while, taking $\psi(A,B)=\|A-B\|_2$, the Hessian approximation satisfies Oracle~\ref{def:oracle_det_2a} with
\begin{align*}
\epsilon_{H\psi}(M)
=
\sqrt{2p}\|\widehat M^{-1}\|_2
\left(
\frac{L_2h}{6}
+
\frac{2\epsilon_f}{h^2}
\right).
\end{align*}

Since $p=\frac{(n+1)(n+2)}{2}-1$, constructing the model requires $p+1=\frac{(n+1)(n+2)}{2}$ evaluations of Oracle~\ref{def:oracle_det_0}. Thus, for a non-adjustable zeroth-order oracle, a naive implementation has total cost ${\cal C}_H = \left(\frac{(n+1)(n+2)}{2}\right){\cal C}_f
+
\mathcal{O}\left(p^3\right)$.


\subsubsection{Quadratic Interpolation via First-Order Oracles}

Quadratic interpolation can also be constructed from first-order information. Given access to Oracle~\ref{def:oracle_det_1}, with output $g(x) \in \mathbb{R}^n$, consider the quadratic model
\begin{align}
	m(y)=c+ g(x)^T(y-x) + \frac{1}{2}(y-x)^TH(x)(y-x),
\label{eq:quadratic_model_first}
\end{align}
where $c\in\mathbb{R}$ is arbitrary and $H\in\mathbb{S}^n$ serves as an approximation of $\nabla^2\phi(x)$. Since the constant term does not affect the Hessian approximation, $H$ can be constructed directly from first-order information.

Specifically, consider a set of $n$ linearly independent directions ${\cal U}=\{u_1,\ldots,u_n\}\subset\mathbb{R}^n$, 
and define $U=
\begin{bmatrix}
u_1 & \cdots & u_n
\end{bmatrix}$.
We first construct an auxiliary linear interpolation model of the gradient,
\begin{align}
m_g(y)
=
g(x)+B(y-x),
\label{eq:quadratic_gradient_model}
\end{align}
where $B\in\mathbb{R}^{n\times n}$. The vector-valued model $m_g$ approximates $\nabla\phi(y)$ locally and is determined by imposing the gradient interpolation conditions
\begin{align}
m_g(x+u_i)
=
g(x+u_i),
\qquad
i=1,\ldots,n.
\label{eq:quadratic_interp_first_cond}
\end{align}
Equivalently,
\begin{align*}
Bu_i
=
g(x+u_i)-g(x),
\qquad
i=1,\ldots,n.
\end{align*}
Thus, defining
\begin{align*}
G_\Delta
=
\begin{bmatrix}
g(x+u_1)-g(x) & \cdots & g(x+u_n)-g(x)
\end{bmatrix},
\end{align*}
we obtain
\begin{align}
B
=
G_\Delta U^{-1}.
\label{eq:quadratic_first_system}
\end{align}
The matrix $B$ need not be symmetric and therefore need not itself correspond to the Hessian of a scalar quadratic model. We consequently define
\begin{align}
H
=
\frac{1}{2}\left(B+B^T\right),
\label{eq:quadratic_first_hessian}
\end{align}
which yields the symmetric Hessian approximation used in \eqref{eq:quadratic_model_first}.

Suppose that $\max_{i=1,\ldots,n}\|u_i\|\leq h$, 
and define $\widehat U=\frac{1}{h}U$. The following proposition provides an error bound for the resulting Hessian approximation.

\begin{proposition}
\label{prop:quadratic_interp_first}
Let Assumptions~\ref{assum:lip_cont} and~\ref{assum:lip_cont2} hold, and suppose that the first-order oracle evaluations used to construct \eqref{eq:quadratic_first_hessian} satisfy Oracle~\ref{def:oracle_det_1} with error parameters $\epsilon_g$ and $\kappa_g$. Then,
\begin{align*}
\|H(x)-\nabla^2\phi(x)\|_2
\leq
\sqrt{n}\|\widehat U^{-1}\|_2
\left(
\frac{L_2h}{2}
+
\frac{2\left(\epsilon_g+\kappa_g\|\nabla\phi(x)\|\right)}{h}
+
\kappa_gL_1
\right).
\end{align*}
\end{proposition}

When $u_i=he_i$, $i=1,\ldots,n$, we have $\widehat U=I$, and the construction and error bound reduce to those of the forward finite-difference Hessian approximation via first-order information.

We now discuss Proposition~\ref{prop:quadratic_interp_first} in terms of Oracle~\ref{def:oracle_det_2a}. Suppose additionally that $\|\nabla\phi(x)\|\leq G$ over the region of interest, for some $G>0$. Taking $\psi(A,B)=\|A-B\|_2$,  the Hessian approximation satisfies Oracle~\ref{def:oracle_det_2a} with
\begin{align*}
\epsilon_{H\psi}(U)
=
\sqrt{n}\|\widehat U^{-1}\|_2
\left(
\frac{L_2h}{2}
+
\frac{2\left(\epsilon_g+\kappa_g G\right)}{h}
+
\kappa_gL_1
\right),
\end{align*}
where the interpolation set, and hence $h$ and $\widehat U$, serves as the input $\lambda$. A naive implementation requires $n+1$ evaluations of Oracle~\ref{def:oracle_det_1}, yielding the total cost ${\cal C}_H=(n+1){\cal C}_g+\mathcal{O}\left(n^3\right)$.

\subsubsection{Simplex Hessian}

The simplex-gradient construction can be extended to approximate second-order information from first-order oracle outputs. Given a set of possibly linearly dependent directions
\begin{align*}
{\cal U}=\{u_1,\ldots,u_p\}\subset\mathbb{R}^n,
\qquad
\max_{i=1,\ldots,p}\|u_i\|\leq h,
\end{align*}
define
\begin{align*}
U
=
\begin{bmatrix}
u_1^T\\
\vdots\\
u_p^T
\end{bmatrix},
\qquad
G_\Delta
=
\begin{bmatrix}
\left(g(x+u_1)-g(x)\right)^T\\
\vdots\\
\left(g(x+u_p)-g(x)\right)^T
\end{bmatrix}.
\end{align*}
The simplex Hessian approximation is given by
\begin{align}
B(x)=U^\dagger G_\Delta.
\label{eq:simplex_hess}
\end{align}
This construction generalizes the interpolation approach above to arbitrary collections of directions; see, e.g., \cite{hare2024matrix}. Since $B(x)$ need not be symmetric, we will subsequently symmetrize it to obtain a Hessian approximation $H(x)\in\mathbb{S}^n$. Let
$\widehat U=U/h$ and $\proj_U=U^\dagger U$.

The following proposition provides an error bound for the simplex Hessian approximation.

\begin{proposition}
\label{prop:simplex_hess}
Let Assumptions~\ref{assum:lip_cont} and~\ref{assum:lip_cont2} hold, and suppose that the first-order oracle evaluations used in \eqref{eq:simplex_hess} satisfy Oracle~\ref{def:oracle_det_1} with error parameters $\epsilon_g$ and $\kappa_g$. Then,
\begin{align*}
\left\|
\proj_U\left(B(x)-\nabla^2\phi(x)\right)
\right\|_2
&=
\left\|
B(x)-\proj_U\left(\nabla^2\phi(x)\right)
\right\|_2
\\
&\leq
\sqrt{p}\|\widehat U^\dagger\|_2
\left(
\frac{L_2h}{2}
+
\frac{2\left(\epsilon_g+\kappa_g\|\nabla\phi(x)\|\right)}{h}
+
\kappa_gL_1
\right).
\end{align*}
\end{proposition}

\begin{proof}
Using \eqref{eq:simplex_hess} and $\proj_U=U^\dagger U$,
\begin{align*}
	B(x)-\proj_U\left(\nabla^2\phi(x)\right)=U^\dagger\left(G_\Delta-U\nabla^2\phi(x)\right).
\end{align*}
For each $i=1,\ldots,p$, Assumptions~\ref{assum:lip_cont} and~\ref{assum:lip_cont2}, together with Oracle~\ref{def:oracle_det_1}, imply
\begin{align*}
\left\|
g(x+u_i)-g(x)-\nabla^2\phi(x)u_i
\right\|
\leq
\frac{L_2h^2}{2}
+
2\left(\epsilon_g+\kappa_g\|\nabla\phi(x)\|\right)
+
\kappa_gL_1h.
\end{align*}
Hence,
\begin{align*}
\left\|
G_\Delta-U\nabla^2\phi(x)
\right\|_2
\leq
\sqrt{p}
\left(
\frac{L_2h^2}{2}
+
2\left(\epsilon_g+\kappa_g\|\nabla\phi(x)\|\right)
+
\kappa_gL_1h
\right).
\end{align*}
The result follows from $\|U^\dagger\|_2=\|\widehat U^\dagger\|_2/h$ and
$\proj_U \left(B(x)\right)=B(x)$.
\end{proof}

When $U$ has full column rank, $\proj_U(\cdot)=I$, and Proposition~\ref{prop:simplex_hess} provides a bound with respect to the full Hessian. In this case, since $H(x)$ and $\nabla^2\phi(x)$ are symmetric,
\begin{align*}
\|H(x)-\nabla^2\phi(x)\|_2
\leq
\|B(x)-\nabla^2\phi(x)\|_2.
\end{align*}

In the simplest case where Oracle~\ref{def:oracle_det_1} is not adjustable, with  $\kappa_g=0$, the input $\lambda$ consists of the sample directions ${\cal U}$, which determine $h$ and $\widehat U$.
For $\psi(A,B)=\|A-B\|_2$,  $H(x)$ satisfies Oracle~\ref{def:oracle_det_2a} with
\begin{align*}
\epsilon_{H\psi}(U)
=
\sqrt{p}\|\widehat U^\dagger\|_2
\left(
\frac{L_2h}{2}
+
\frac{2\left(\epsilon_g\right)}{h}
+
\kappa_gL_1
\right).
\end{align*}
 The construction requires $p+1$ evaluations of Oracle~\ref{def:oracle_det_1}, with total cost ${\cal C}_H=(p+1){\cal C}_g+\mathcal{O}\left(pn^2\right)$.

If $\kappa_g>0$ and we need to assume a bound $\|\nabla\phi(x)\|\leq G$ over the region of interest, then this is 
equivalent to redefining Oracle~\ref{def:oracle_det_1} with constants $\epsilon_g^\prime = \epsilon_g+\kappa_gG$ and 
$\kappa_g^\prime=0$ and apply above.

\subsubsection{Gauss--Newton Hessian Approximation}

For nonlinear least-squares problems of the form
\begin{align*}
\phi(x)
=
\frac{1}{2}\|r(x)\|^2,
\qquad \text{where} \qquad
r(x)
=
\begin{bmatrix}
r_1(x)\\
\vdots\\
r_m(x)
\end{bmatrix}.
\end{align*}
Let $J(x)$ denote the Jacobian of $r(x)$. The Hessian of $\phi$ satisfies
\begin{align*}
\nabla^2\phi(x)
=
J(x)^TJ(x)
+
\sum_{i=1}^m
r_i(x)\nabla^2 r_i(x).
\end{align*}
The Gauss-Newton Hessian approximation neglects the second term and is given by
\begin{align*}
H(x)
=
J(x)^TJ(x).
\end{align*}
The following proposition bounds the error in the Gauss-Newton Hessian approximation.

\begin{proposition}
\label{prop:gauss_newton}
Suppose that each residual function $r_i$, $i=1,\ldots,m$, is twice differentiable and that, for some constant $L_r\geq0$, $\|\nabla^2 r_i(x)\|_2\leq L_r$ for $i=1,\ldots,m$. Then the Gauss-Newton approximation $H(x)$ satisfies
\begin{align*}
\|H(x)-\nabla^2\phi(x)\|_2
&\leq
L_r\|r(x)\|_1
\\
&\leq
\sqrt{m}L_r\|r(x)\|.
\end{align*}
\end{proposition}


Thus, if $\|r(x)\|\leq\epsilon_{\rm GN}$ over the region of interest and $\psi(A,B)=\|A-B\|_2$, then $H(x)$ satisfies Oracle~\ref{def:oracle_det_2a} with
\begin{align*}
\epsilon_{H\psi}
=
\sqrt{m}L_r\epsilon_{\rm GN}.
\end{align*}
Unlike the preceding constructions, the classical Gauss--Newton approximation has no natural tunable control parameter. Its accuracy instead depends directly on the current residual and improves as $\|r(x)\|$ decreases, while $H(x)\succeq0$ by construction. Exploiting this state-dependent accuracy more explicitly would require an extension of Oracle~\ref{def:oracle_det_2a} in which the error bound may depend on quantities evaluated at the current point, and the corresponding algorithmic analysis is beyond the scope of this review.

\subsubsection{Directional Finite Differences for Hessian-Vector Products}

We next present a representative example of Oracle~\ref{def:oracle_det_2b}, whose input includes a direction $v\in\mathbb{R}^n$ and whose output approximates the Hessian--vector product $\nabla^2\phi(x)v$. Given access to Oracle~\ref{def:oracle_det_1}, for $v\neq 0$ define
\begin{align}
H(x;v)
=
\frac{\|v\|}{h}
\left[
g\left(x+h\frac{v}{\|v\|}\right)-g(x)
\right],
\label{eq:directional_fd_hv}
\end{align}
and set $H(x;0)=0$. The normalization of $v$ ensures that the finite-difference displacement has length $h$, independently of $\|v\|$, and yields an approximation error that scales linearly with $\|v\|$, consistent with Oracle~\ref{def:oracle_det_2b}.

\begin{proposition}
\label{prop:directional_fd_hv}
Let Assumptions~\ref{assum:lip_cont} and~\ref{assum:lip_cont2} hold, and suppose that the first-order oracle evaluations used in \eqref{eq:directional_fd_hv} satisfy Oracle~\ref{def:oracle_det_1} with error parameters $\epsilon_g$ and $\kappa_g$. Then, for any $x,v\in\mathbb{R}^n$ and $h>0$,
\begin{align*}
\|H(x;v)-\nabla^2\phi(x)v\|
\leq
\left[
\frac{L_2h}{2}
+
\frac{2\left(\epsilon_g+\kappa_g\|\nabla\phi(x)\|\right)}{h}
+
\kappa_gL_1
\right]\|v\|.
\end{align*}
\end{proposition}

\begin{proof}
The result is immediate for $v=0$. For $v\neq0$, let $\widehat v=v/\|v\|$. Then,
\begin{align*}
\|H(x;v)-\nabla^2\phi(x)v\|
&=
\frac{\|v\|}{h}
\left\|
g(x+h\widehat v)-g(x)
-
h\nabla^2\phi(x)\widehat v
\right\|.
\end{align*}
By adding and subtracting the corresponding exact gradients, Assumption~\ref{assum:lip_cont2} gives a truncation error bounded by $L_2h^2/2$, while Oracle~\ref{def:oracle_det_1} and Assumption~\ref{assum:lip_cont} yield
\begin{align*}
\|g(x+h\widehat v)-\nabla\phi(x+h\widehat v)\|
+
\|g(x)-\nabla\phi(x)\|
\leq
2\left(\epsilon_g+\kappa_g\|\nabla\phi(x)\|\right)
+
\kappa_gL_1h.
\end{align*}
Combining these bounds and multiplying by $\|v\|/h$ yields the result.
\end{proof}

 In the simplest case where Oracle~\ref{def:oracle_det_1} is not adjustable, with $\kappa_g=0$, taking $\lambda=h$, the approximation satisfies Oracle~\ref{def:oracle_det_2b} with
\begin{align*}
\kappa_{Hv}(h)
=
\frac{L_2h}{2}
+
\frac{2\epsilon_g}{h}.
\end{align*}
The construction requires two evaluations of Oracle~\ref{def:oracle_det_1}, and hence ${\cal C}_{Hv}(h)=2{\cal C}_g+\mathcal{O}(n)$. If $g(x)$ is already available, only one additional first-order oracle evaluation is required.
If $\kappa_g>0$  we again  assume a bound $\|\nabla\phi(x)\|\leq G$ over the region of interest, then this is 
equivalent to redefining Oracle~\ref{def:oracle_det_1} with constants $\epsilon_g^\prime = \epsilon_g+\kappa_gG$ and 
$\kappa_g^\prime=0$ and apply above. 

A Hessian-vector product can also be approximated using only zeroth-order information. For $v\neq0$, letting $\widehat v=v/\|v\|$, define
\begin{align*}
H_0(x;v)
=
\frac{\|v\|}{h^2}
\sum_{j=1}^n
\Big[
&f(x+h\widehat v+he_j)
-f(x+h\widehat v) -f(x+he_j)
+f(x)
\Big]e_j.
\end{align*}
Under Assumption~\ref{assum:lip_cont2} and Oracle~\ref{def:oracle_det_0},
\begin{align*}
\|H_0(x;v)-\nabla^2\phi(x)v\|
\leq
\sqrt{n}
\left(
L_2h+\frac{4\epsilon_f}{h^2}
\right)\|v\|.
\end{align*}
Thus, this construction also satisfies Oracle~\ref{def:oracle_det_2b} with
\begin{align*}
\kappa_{Hv}(h)
=
\sqrt{n}
\left(
L_2h+\frac{4\epsilon_f}{h^2}
\right),
\end{align*}
but requires at most $2n+2$ evaluations of Oracle~\ref{def:oracle_det_0}.

The resulting bounds depend naturally on the accuracy of the underlying first- and zeroth-order oracles; notably, the zeroth-order construction exhibits the less favorable noise amplification $\epsilon_f/h^2$, compared with $\epsilon_g/h$ in the first-order construction. Here we have considered only forward finite-difference variants; central finite-difference constructions can be developed analogously, typically yielding improved dependence on $h$ under stronger smoothness assumptions and at the cost of additional oracle evaluations.

\subsection{Deterministic Inexact Second-Order Oracles in Optimization}\label{sec:det_second_order_optimization}

In this section, we review representative optimization methods that employ the deterministic inexact second-order oracles introduced in Definitions~\ref{def:oracle_det_2a} and~\ref{def:oracle_det_2b}. Our primary focus is on smooth nonconvex optimization, although we also discuss selected results for strongly convex problems when they provide useful insight into local convergence and the effects of noise. We consider several algorithmic frameworks, including Newton-type methods, trust-region methods, cubic regularization, Newton-CG and negative-curvature methods, adaptive sampling, and quasi-Newton methods. Particular attention is given to the different accuracy conditions imposed on approximate Hessians and Hessian-vector products, their relationship to the oracle definitions, and the convergence and complexity guarantees they support. Although our focus is on deterministic accuracy conditions, we also discuss selected methods that use randomized algorithms or sampling to construct approximations satisfying such conditions. Rather than providing an exhaustive treatment or reproducing complete algorithms and proofs, our objective is to illustrate the principal relationships between oracle accuracy, algorithmic design, and convergence guarantees.

\paragraph{General descent and Newton-type directions.}
Recall from Section~\ref{sec.other_oracles} that convergence of general descent methods can be established when the search directions satisfy suitable descent and angle conditions \eqref{eq:angle_condition}. Newton-type directions fit naturally within this framework. In particular, consider $d_k=-H_k^{-1}g_k$, where $g_k$ is an approximation of $\nabla\phi(x_k)$ and $H_k$ is symmetric positive definite. As discussed in Section~\ref{sec:quasi-newton}, such directions satisfy the required descent conditions when $H_k$ is positive definite and sufficiently well conditioned. Importantly, these conditions do not require $H_k$ to be a uniformly accurate approximation of $\nabla^2\phi(x_k)$.

Quasi-Newton methods provide an important example of this setting, since they construct positive-definite curvature approximations from gradient differences rather than direct Hessian information. When exact gradients are available, standard updates such as BFGS preserve positive definiteness under suitable curvature conditions \cite{NoceWrig06,dennis1974characterization}. With inexact gradients, however, the curvature difference $y_k=g_{k+1}-g_k$ may be contaminated by the difference between successive gradient errors, and care is required to prevent unreliable updates. For strongly convex objectives with bounded function and gradient errors, safeguarded BFGS variants can preserve useful curvature properties and retain convergence to a noise-dependent neighborhood \cite{xie2020analysis}. Thus, even when second-order information is constructed indirectly from an inexact first-order oracle, suitable safeguards can produce matrices satisfying the conditions required by general descent methods.

The positive-definiteness assumption, however, restricts the use of second-order information in nonconvex optimization, where negative curvature can provide information beyond that contained in the gradient. General frameworks can combine descent directions satisfying \eqref{eq:angle_condition} with directions of negative curvature, allowing algorithms to exploit both first- and second-order information \cite{curtis2019exploiting}. This motivates the more general second-order methods discussed below, where approximate Hessians may need to support not only descent but also reliable detection and exploitation of negative curvature.

\paragraph{Trust-region methods.}
As discussed in the preceding sections, trust-region methods can accommodate inexact function values and gradients. Here, we consider the additional use of approximate second-order information, with the quadratic model
\begin{align*}
    m_k(s)
    = f_k + g_k^T s + \frac{1}{2}s^T H_k s,
    \qquad \|s\| \leq \Delta_k,
\end{align*}
where $f_k$, $g_k$, and $H_k$ may all be inexact. With exact function values and gradients, uniformly bounded approximate Hessians, and sufficient model decrease, standard trust-region methods can establish convergence to first-order stationary points without requiring $H_k$ to accurately approximate the true Hessian \cite{Conn-Gould-Toint-00}.

This observation extends, in a weaker form, to settings with persistent noise. In particular, trust-region methods with bounded errors in function and gradient evaluations and uniformly bounded Hessian approximations can still provide meaningful convergence guarantees. 
A  version of such a result under Oracle~\ref{def:oracle_det_1} is given in Section~\ref{sec:tr}. 

To establish stronger convergence guarantees, however, additional accuracy conditions on the Hessian approximation are generally required. A common assumption is the operator-norm bound
\begin{align*}
    \|H_k-\nabla^2\phi(x_k)\|
    \leq \epsilon_{H},
\end{align*}
which is a direct instance of Oracle~\ref{def:oracle_det_2a}. When $\epsilon_{H}$ is sufficiently small relative to $\epsilon_{2}$, the approximate Hessian can reliably capture negative curvature at the desired tolerance, provided the trust-region subproblem is solved to sufficient accuracy, and thus second-order convergence and complexity can be established; see e.g., \cite{xu2020newton,cartis2012complexity}.

Uniform operator-norm accuracy is not always necessary. For example, the weaker directional condition
\begin{align*}
    \|(H_k-\nabla^2\phi(x_k))s_k\|
    \leq \epsilon_{H}(\lambda_k)\|s_k\| 
    \quad \text{and} \quad
    \|H_k\|\leq\kappa_H,
\end{align*}
where $s_k$ is the computed trust-region step, has been considered in \cite{xu2020newton,cartis2022evaluation}. In this setting, the Hessian accuracy is chosen in relation to the trust-region radius and the desired second-order stationarity tolerance, with $ \epsilon_{H}(\lambda_k) = \max\{c \epsilon_{2},\Delta_k\}$, with $\lambda_k=(\epsilon_2,\Delta_k)$  serving as the input to the oracle. The directional condition can be expressed through Oracle~\ref{def:oracle_det_2b}. Under suitable Hessian accuracy and subproblem-solution conditions, an iteration complexity bound of $\mathcal{O}\left(
        \max\left\{
            \epsilon_{1}^{-2}
            \epsilon_{2}^{-1},
            \epsilon_{2}^{-3}
        \right\}
    \right)$
can be established for approximate second-order stationarity; see e.g., \cite{xu2020newton}. Related results allow both the gradient and Hessian to be approximated while retaining worst-case iteration complexity guarantees comparable to those of their exact-derivative counterparts \cite{yao2023inexact}. More general complexity results for trust-region methods with
inexact function and derivative information are discussed in
\cite{cartis2022evaluation}.

\paragraph{Cubic and higher-order regularization.}
Cubic regularization methods provide another framework for incorporating approximate second-order information. Unlike trust-region methods, which explicitly constrain the step, these methods control its length through a regularization term. A typical cubic model is
\begin{align*}
    m_k(s)
    = f_k + g_k^T s
    + \frac{1}{2}s^T H_k s
    + \frac{\sigma_k}{3}\|s\|^3,
\end{align*}
where $\sigma_k>0$ is the regularization parameter and $f_k$, $g_k$, and $H_k$ may be inexact. With sufficiently accurate derivatives and appropriate conditions on the computed steps, cubic regularization methods achieve worst-case iteration complexity of $\mathcal{O}\left(
        \max\left\{
            \epsilon_{1}^{-3/2},
            \epsilon_{2}^{-3}
        \right\}
    \right)$  for approximate second-order stationarity; see e.g., \cite{cagotoPI,cagotoPII,cartis2012complexity}. 
    
When the gradient is exact but the Hessian is approximated, the accuracy requirements can be expressed in terms of the computed step. In particular, the directional accuracy condition
\begin{align*}
    \|(H_k-\nabla^2\phi(x_k))s_k\|
    \leq \kappa_H\|s_k\|^2
\end{align*}
allows the Hessian approximation error to decrease with the norm of the step, and can be sufficient to preserve first-order complexity guarantees under appropriate algorithmic conditions \cite{cagotoPI,cagotoPII}. The directional condition can be achieved by imposing stronger, sufficient condition as 
\begin{align*}
    \|H_k-\nabla^2\phi(x_k)\|
    \leq \epsilon_{H,k},
    \qquad
    \epsilon_{H,k}\leq\kappa_H\|s_k\|,
\end{align*}
which can be expressed as Oracle~\ref{def:oracle_det_2a}; see e.g., \cite{bellavia2021adaptive}.

When both the gradient and Hessian are inexact, their errors must be controlled jointly. A representative set of sufficient conditions is
\begin{align*}
    \|g_k-\nabla\phi(x_k)\|
    \leq \kappa_g\|s_k\|^2,
    \quad \text{and} \quad
    \|H_k-\nabla^2\phi(x_k)\|
    \leq \kappa_H\|s_k\|.
\end{align*}
These conditions ensure that the errors in the gradient of the cubic model are of order $\|s_k\|^2$, consistent with the regularization term. Since the step $s_k$ is not known in advance, practical adaptive schemes may use previous step lengths or dynamically refine the oracle accuracy during an iteration. Such approaches can preserve the  worst-case iteration complexity order of exact-derivative cubic regularization for attaining first-order stationarity; see e.g.,  \cite{wang2019note,bellavia2019adaptive}.

Stronger conditions are generally required to establish second-order stationarity. In particular, directional accuracy along $s_k$ alone, i.e., Oracle~\ref{def:oracle_det_2b}, does not guarantee that $H_k$ identifies negative curvature in other directions. An operator-norm condition such as
\begin{align*}
    \|H_k-\nabla^2\phi(x_k)\| \leq c_H\epsilon_{2},
\end{align*}
that can be expressed as Oracle~\ref{def:oracle_det_2a}, provides the additional spectral accuracy needed to relate the smallest eigenvalue of the approximate Hessian to that of the true Hessian. Combined with suitable gradient accuracy and approximate minimization conditions for the cubic model, this permits second-order complexity guarantees comparable to those obtained with exact derivatives; see e.g., \cite{cartis2012complexity,bellavia2019adaptive}.

These ideas extend to regularization methods based on higher-order Taylor models. Suppose that derivatives up to order $p\geq 2$ are available, possibly inexactly, and consider
\begin{align*}
    m_k(s)
    =
    f_k+
    \sum_{j=1}^{p}
    \frac{1}{j!}
    \overline{\nabla^j\phi}(x_k)[s]^j
    +
    \frac{\sigma_k}{p+\beta}\|s\|^{p+\beta},
\end{align*}
where $\beta\in(0,1]$ characterizes the H\"older continuity of the $p$th derivative and $\overline{\nabla^j\phi}(x_k)$ denotes an approximation of the $j$th derivative tensor. A general deterministic accuracy requirement takes the form
\begin{align*}
    \left\|
        \overline{\nabla^j\phi}(x_k)
        -
        \nabla^j\phi(x_k)
    \right\|
    \leq \epsilon_{j,k},
    \qquad j=1,\ldots,p,
\end{align*}
where the tolerances $ \epsilon_{j,k}$ are selected adaptively. 

Under appropriate smoothness and dynamic accuracy conditions, higher-order adaptive regularization methods with inexact evaluations retain general evaluation-complexity guarantees. For a $q$th-order approximate minimizer, with $1\leq q\leq p$, and the corresponding accuracy conditions measure, a representative bound is $
    \mathcal{O}\left(
        \epsilon^{-\frac{p+\beta}{p-q+\beta}}
    \right)$,
subject to the assumptions and accuracy requirements of the particular method; see e.g.,  \cite{bellavia2019adaptive,cartis2022evaluation}. These results demonstrate that dynamically controlled inexact derivatives can preserve the worst-case iteration complexity orders of higher-order regularization methods. As always, this complexity must be distinguished from the total computational cost of obtaining the required oracle accuracy.

\paragraph{Newton-CG and negative-curvature methods.}
Newton-CG methods provide a matrix-free approach to second-order optimization by approximately solving Newton systems using conjugate gradients (CG). Rather than forming or storing the Hessian, these methods require only Hessian-vector products, making them particularly suitable for large-scale problems \cite{steihaug1983conjugate}. When second-order information is inexact, the accuracy of these products and the stopping conditions of the iterative solver must be coordinated with the desired convergence guarantees; see e.g.,  \cite{yao2023inexact}.

A representative analysis considers inexact gradients and Hessians satisfying
\begin{align*}
    \|g_k-\nabla\phi(x_k)\|
    \leq \epsilon_g,
    \qquad \text{and} \qquad
    \|H_k-\nabla^2\phi(x_k)\|
    \leq \epsilon_H,
\end{align*}
which correspond to Oracle~\ref{def:oracle_det_1} with $\kappa_g=0$ and Oracle~\ref{def:oracle_det_2a} with the operator norm as the error measure, respectively. If $\epsilon_1$ and $\epsilon_2$ denote the desired first- and second-order stationarity tolerances, then a simple iteration-independent choice requires $\epsilon_g=\mathcal{O}(\epsilon_1)$ and $\epsilon_H=\mathcal{O}(\epsilon_2)$. Under these accuracy requirements, the inexact Newton--CG method
retains the corresponding approximate second-order stationarity
and worst-case iteration-complexity guarantees \cite{yao2023inexact}.  With $\epsilon_2=\sqrt{L_2\epsilon_1}$, the resulting inexact Newton-CG method achieves an iteration complexity of $\mathcal{O}\left(\epsilon_1^{-3/2}\right)$ for approximate second-order stationarity, matching the optimal complexity of the corresponding exact-derivative Newton--CG method \cite{royer2020newton}.

In nonconvex optimization, Newton--CG methods must also account for negative curvature. When CG encounters a direction $v_k$ satisfying $v_k^TH_k v_k<0$, this direction can be used to obtain additional decrease rather than continuing to solve the Newton system. More generally, methods that combine descent and negative-curvature directions can establish convergence to approximate second-order stationary points \cite{curtis2019exploiting,yao2023inexact}. With inexact Hessian information, however, negative curvature of $H_k$ may not accurately reflect the curvature of the true Hessian. 

A simple way to guarantee reliable negative-curvature detection is through spectral accuracy. Suppose that Oracle~\ref{def:oracle_det_2a} provides a Hessian approximation satisfying
\begin{align*}
    \|H_k-\nabla^2\phi(x_k)\|
    \leq \epsilon_2.
\end{align*}
Then, by Weyl's inequality,
\begin{align*}
    \left|
    \lambda_{\min}(H_k)
    -
    \lambda_{\min}(\nabla^2\phi(x_k))
    \right|
    \leq
    \epsilon_2.
\end{align*}
Hence, if $\lambda_{\min}(\nabla^2\phi(x_k))<-\epsilon_2$, then $\lambda_{\min}(H_k)<0$, so negative curvature of the true Hessian that exceeds the desired second-order stationarity tolerance is also present in the approximation. Thus, Hessian accuracy on the order of $\epsilon_2$ allows negative curvature to be detected at the accuracy required for establishing approximate second-order stationarity. Inexact Newton--CG methods combine such Hessian accuracy with suitable gradient accuracy, CG termination conditions, and negative-curvature detection while retaining worst-case complexity guarantees comparable to their exact-derivative counterparts \cite{yao2021inexact,yao2023inexact}.

Directional accuracy alone, as in Oracle~\ref{def:oracle_det_2b}, does not provide such a guarantee, since accurate Hessian-vector products along the computed directions do not rule out undetected negative curvature in other directions. Thus, directional conditions by themselves are generally insufficient for establishing second-order stationarity. They can, however, be combined with a separate spectral accuracy condition ensuring that sufficiently negative curvature is not missed; this permits weaker directional accuracy along the computed descent and negative-curvature directions while still yielding complexity guarantees for approximate second-order stationarity \cite{berahas2026exploiting}. Related methods combine inexact gradients and Hessians with randomized procedures for identifying negative curvature while retaining second-order complexity guarantees \cite{li2025randomized}.

\paragraph{Fast local convergence.}
The connection between Hessian accuracy and local convergence can be understood through classical inexact Newton and quasi-Newton theory. Suppose that $x_k$ converges to a solution $x^*$ at which the Hessian is positive definite. For a Newton-type step $s_k$, a directional accuracy condition of the form
\begin{align*}
    \frac{
        \|(H_k-\nabla^2\phi(x_k))s_k\|
    }{\|s_k\|}
    \longrightarrow 0
\end{align*}
requires the Hessian approximation to become asymptotically accurate along the computed steps. This condition can be interpreted directly through Oracle~\ref{def:oracle_det_2b}. In particular, if $h(x_k,s_k)=H_ks_k$ and $p=1$, then it is sufficient that $\kappa_{Hv}(\lambda_k)\to0$. More generally, when $p>1$, the oracle bound gives
\begin{align*}
    \frac{
        \|h(x_k,s_k)-\nabla^2\phi(x_k)s_k\|
    }{\|s_k\|}
    \leq
    \kappa_{Hv}(\lambda_k)\|s_k\|^{p-1},
\end{align*}
so the directional accuracy condition is also satisfied if $\kappa_{Hv}(\lambda_k)$ remains bounded and $\|s_k\|\to0$. With exact gradients and suitable regularity assumptions, such directional accuracy is closely related to the Dennis--Mor\'e condition for superlinear convergence \cite{dennis1974characterization}. When the gradient is also inexact, its accuracy must improve sufficiently rapidly as the iterates approach the solution.

Subsampled Newton methods provide a useful illustration of these requirements. For strongly convex stochastic optimization problems, a Newton method using a Hessian sample of fixed size and a gradient sample whose size increases geometrically can achieve a global linear rate of convergence in expectation \cite{bollapragada2018exact}. Faster local convergence requires increasingly accurate derivative information. Under additional regularity and moment assumptions, a unit-step Newton method achieves local R-superlinear convergence in expectation when the gradient sample size grows faster than geometrically and the Hessian sample size increases monotonically and tends to infinity \cite{bollapragada2018exact}. Related analyses of subsampled Newton methods establish local linear and superlinear rates under alternative Hessian and gradient sampling conditions \cite{roosta2019sub}.

\subsection{Stochastic second-order oracles} 

We next consider stochastic second-order oracles, in which the Hessian approximation or Hessian-vector product is random and its accuracy is characterized in expectation. As in the deterministic setting, we distinguish between oracles that return an approximation of the full Hessian and those that provide only its action along specified directions.

\begin{customoracle}{S-2a}[General Stochastic Second-Order Oracle]
\label{def:oracle_stoch_2a}
Let $\psi:\mathbb{S}^n\times\mathbb{S}^n\to\mathbb{R}_+$ be a prescribed error measure, where $\mathbb{S}^n$ denotes the set of $n\times n$ symmetric matrices. Given a point $x\in\mathbb{R}^n$ and a control parameter $\lambda$, the oracle generates an output $H(x,\xi^{(2)})\in\mathbb{S}^n$ that satisfies
\begin{equation}
\label{eq:2a_oracle_stoch_bias}
\psi\left(
\E_{\xi^{(2)}}[H(x,\xi^{(2)})],
\nabla^2\phi(x)
\right)
\leq
\epsilon^2_{H\psi}(\lambda),
\end{equation}
and
\begin{equation}
\label{eq:2a_oracle_stoch_var}
\E_{\xi^{(2)}}\left[
\psi\left(
H(x,\xi^{(2)}),
\E_{\xi^{(2)}}[H(x,\xi^{(2)})]
\right)
\right]
\leq
\sigma^2_{H\psi}(\lambda),
\end{equation}
where $\epsilon_{H\psi}(\lambda)$ and $\sigma_{H\psi}(\lambda)$ are nonnegative deterministic functions of $\lambda$. The computational cost of this oracle is some function ${\cal C}_{H}(\lambda;\psi)$.
\end{customoracle}

When $\psi(A,B)=\|A-B\|_F^2$, the two conditions above control the squared bias and variance of the Hessian approximation, respectively. In this case,
\begin{align*}
\E_{\xi^{(2)}}\left[
\|H(x,\xi^{(2)})-\nabla^2\phi(x)\|_F^2
\right]
&=
\left\|
\E_{\xi^{(2)}}[H(x,\xi^{(2)})]
-
\nabla^2\phi(x)
\right\|_F^2
\\
&\quad+
\E_{\xi^{(2)}}\left[
\left\|
H(x,\xi^{(2)})
-
\E_{\xi^{(2)}}[H(x,\xi^{(2)})]
\right\|_F^2
\right]
\\
&\leq
\epsilon^2_{H\psi}(\lambda)
+
\sigma^2_{H\psi}(\lambda).
\end{align*}
Thus, Oracle~\ref{def:oracle_stoch_2a} allows the systematic error and stochastic variability of the second-order approximation to be controlled separately.

\begin{customoracle}{S-2b}[Stochastic Inexact Hessian-Vector Product Oracle]
\label{def:oracle_stoch_2b}
Given a point $x\in\mathbb{R}^n$, a vector $v\in\mathbb{R}^n$, and a control parameter $\lambda$, the oracle generates an output $h(x,v,\xi^{(2)})\in\mathbb{R}^n$ that satisfies
\begin{equation}
\label{eq:2b_oracle_stoch}
\E_{\xi^{(2)}}\left[
\|h(x,v,\xi^{(2)})-\nabla^2\phi(x)v\|
\right]
\leq
\kappa_{Hv}(\lambda)\|v\|^p,
\end{equation}
where $p\geq1$ and $\kappa_{Hv}(\lambda)$ is a nonnegative deterministic function of $\lambda$. The computational cost of this oracle is some function ${\cal C}_{Hv}(\lambda)$.
\end{customoracle}

As in the deterministic setting, Oracle~\ref{def:oracle_stoch_2b} controls second-order accuracy only along the queried directions and does not require the oracle outputs to correspond to a single random Hessian matrix. Consequently, it can model matrix-free stochastic second-order computations for which constructing or storing a full Hessian approximation is unnecessary.

\subsection{Examples of Stochastic Second-Order Oracles}

We next present several representative examples of stochastic second-order oracles, including minibatch Hessian and Gauss--Newton approximations, stochastic Hessian-vector products, finite-difference constructions, and randomized or sketched Hessian approximations. For each example, we characterize the relevant approximation error, identify the corresponding stochastic second-order oracle framework, and discuss the associated computational cost. Detailed derivations are omitted where appropriate.

\subsubsection{Minibatch Stochastic Hessians and Hessian-Vector Products in EM/ERM}

For expected risk minimization problems \eqref{eq:ERM}, stochastic Hessian approximations can be constructed by averaging second-order information over a set of i.i.d.\ samples. Given a positive integer $\lambda$ as input, the oracle generates $\xi_1,\ldots,\xi_\lambda$ as i.i.d.\ samples from $\Xi$. Let $\xi^{(2)}=(\xi_1,\ldots,\xi_\lambda)$ and define
\begin{align*}
    H(x,\xi^{(2)})
    =
    \frac{1}{\lambda}
    \sum_{j=1}^{\lambda}
    \nabla_x^2\ell(x,\xi_j).
\end{align*}

In most cases, suitable properties of  $\nabla_x^2\ell(x,\xi)$ ensure that 
\begin{align*}
  \nabla^2\phi(x)=
    \E_{\xi\sim\Xi}
    \left[
        \nabla_x^2\ell(x,\xi)
    \right].
\end{align*}
 and  that there exists $M_H>0$,
\begin{align*}
    \E_{\xi\sim\Xi}
    \left[
        \|\nabla_x^2\ell(x,\xi)-\nabla^2\phi(x)\|_F^2
    \right]
    &\leq M_H^2.
\end{align*}

Thus this oracle is unbiased and 
 by independence of the samples,
\begin{align*}
    \E_{\xi^{(2)}}
    \left[
        \left\|
            H(x,\xi^{(2)})
            -
            \E_{\xi^{(2)}}[H(x,\xi^{(2)})]
        \right\|_F^2
    \right]
    \leq
    \frac{M_H^2}{\lambda}.
\end{align*}
Thus, taking $\psi(A,B)=\|A-B\|_F^2$, the minibatch Hessian satisfies Oracle~\ref{def:oracle_stoch_2a} with
\begin{align*}
    \epsilon_{H\psi}(\lambda)
    =0
    \qquad
    \sigma_{H\psi}(\lambda)^2
    =
    \frac{M_H^2}{\lambda}.
\end{align*}
Hence, increasing the minibatch size reduces the stochastic variability of the Hessian approximation. 

The same samples can also be used to construct a stochastic Hessian-vector product without explicitly forming the Hessian matrix. For any $v\in\mathbb{R}^n$, define
\begin{align*}
    h(x,v,\xi^{(2)})
    =
    \frac{1}{\lambda}
    \sum_{j=1}^{\lambda}
    \nabla_x^2\ell(x,\xi_j)v.
\end{align*}
Suppose now that the stochastic Hessian-vector product is unbiased, i.e.,
$\overline H(x)=\nabla^2\phi(x)$. Then
\begin{align*}
    \E_{\xi^{(2)}}[h(x,v,\xi^{(2)})]
    =
    \nabla^2\phi(x)v,
\end{align*}
and
\begin{align*}
    \E_{\xi^{(2)}}
    \left[
        \|h(x,v,\xi^{(2)})-\nabla^2\phi(x)v\|^2
    \right]
    \leq
    \frac{M_H^2}{\lambda}\|v\|^2.
\end{align*}
Consequently, by Jensen's inequality,
\begin{align*}
    \E_{\xi^{(2)}}
    \left[
        \|h(x,v,\xi^{(2)})-\nabla^2\phi(x)v\|
    \right]
    \leq
    \frac{M_H}{\sqrt{\lambda}}\|v\|.
\end{align*}
Hence, the minibatch Hessian-vector product satisfies Oracle~\ref{def:oracle_stoch_2b} with $p=1$ and
$\kappa_{Hv}(\lambda)=M_H/\sqrt{\lambda}$.

The computational cost of both constructions grows linearly with the minibatch size $\lambda$, while the Hessian-vector product formulation avoids explicitly forming or storing the Hessian matrix.

\subsubsection{Stochastic Gauss--Newton Hessian and Hessian-Vector Product Approximations}

Consider an expected least-squares problem of the form
\begin{align*}
    \phi(x)
    =
    \E_{\xi\sim\Xi}
    \left[
        \frac{1}{2}\|r(x,\xi)\|^2
    \right],
\end{align*}
where $J(x,\xi)$ denotes the Jacobian of $r(x,\xi)$. Given a positive integer $\lambda$ as input, the oracle generates $\xi_1,\ldots,\xi_\lambda$ as i.i.d.\ samples from $\Xi$. Let $\xi^{(2)}=(\xi_1,\ldots,\xi_\lambda)$ and define the stochastic Gauss--Newton approximation
\begin{align*}
    H(x,\xi^{(2)})
    =
    \frac{1}{\lambda}
    \sum_{j=1}^{\lambda}
    J(x,\xi_j)^T J(x,\xi_j).
\end{align*}
Defining
\begin{align*}
    G(x)
    =
    \E_{\xi\sim\Xi}
    \left[
        J(x,\xi)^T J(x,\xi)
    \right],
\end{align*}
we have $\E_{\xi^{(2)}}[H(x,\xi^{(2)})]=G(x)$. Suppose that, for some constants $M_{GN}>0$ and $\epsilon_{GN}\geq0$,
\begin{align*}
    \E_{\xi\sim\Xi}
    \left[
        \|J(x,\xi)^T J(x,\xi)-G(x)\|_F^2
    \right]
    \leq M_{GN}^2,
    \qquad
    \|G(x)-\nabla^2\phi(x)\|_F
    \leq \epsilon_{GN}.
\end{align*}
Then
\begin{align*}
    \left\|
        \E_{\xi^{(2)}}[H(x,\xi^{(2)})]
        -
        \nabla^2\phi(x)
    \right\|_F^2
    \leq
    \epsilon_{GN}^2,
\end{align*}
and, by independence of the samples,
\begin{align*}
    \E_{\xi^{(2)}}
    \left[
        \left\|
            H(x,\xi^{(2)})
            -
            \E_{\xi^{(2)}}[H(x,\xi^{(2)})]
        \right\|_F^2
    \right]
    \leq
    \frac{M_{GN}^2}{\lambda}.
\end{align*}
Thus, taking $\psi(A,B)=\|A-B\|_F^2$, the approximation satisfies Oracle~\ref{def:oracle_stoch_2a} with
\begin{align*}
    \epsilon^2_{H\psi}(\lambda)
    =
    \epsilon_{GN}^2,
    \qquad
    \sigma^2_{H\psi}(\lambda)
    =
    \frac{M_{GN}^2}{\lambda}.
\end{align*}
Equivalently,
\begin{align*}
    \E_{\xi^{(2)}}
    \left[
        \|H(x,\xi^{(2)})-\nabla^2\phi(x)\|_F^2
    \right]
    \leq
    \epsilon_{GN}^2+\frac{M_{GN}^2}{\lambda}.
\end{align*}

The same samples can be used to construct the Hessian-vector product approximation
\begin{align*}
    h(x,v,\xi^{(2)})
    =
    \frac{1}{\lambda}
    \sum_{j=1}^{\lambda}
    J(x,\xi_j)^T
    \left(J(x,\xi_j)v\right),
\end{align*}
without explicitly forming the Gauss--Newton matrix. Moreover,
\begin{align*}
    \E_{\xi^{(2)}}
    \left[
        \|h(x,v,\xi^{(2)})-\nabla^2\phi(x)v\|^2
    \right]
    \leq
    \left(
        \epsilon_{GN}^2+\frac{M_{GN}^2}{\lambda}
    \right)\|v\|^2.
\end{align*}
Hence, by Jensen's inequality,
\begin{align*}
    \E_{\xi^{(2)}}
    \left[
        \|h(x,v,\xi^{(2)})-\nabla^2\phi(x)v\|
    \right]
    \leq
    \sqrt{
        \epsilon_{GN}^2+\frac{M_{GN}^2}{\lambda}
    }\|v\|.
\end{align*}
Thus, the construction satisfies Oracle~\ref{def:oracle_stoch_2b} with $p=1$ and
$\kappa_{Hv}(\lambda)
=
\sqrt{\epsilon_{GN}^2+M_{GN}^2/\lambda}$.
Increasing $\lambda$ reduces the stochastic variability, while the intrinsic Gauss--Newton approximation error $\epsilon_{GN}$ remains.

\subsubsection{Finite Differences via Stochastic Zeroth-Order Oracles}

Given access to Oracle~\ref{def:oracle_stoch_0.1} (with a control parameter $\lambda_f$), a Hessian approximation can be constructed by applying forward finite differences to stochastic function evaluations. For any $x\in\mathbb{R}^n$, choose $h>0$ and generate independent oracle outputs at the points
\begin{align*}
{\cal Y}
=
\{x\}
\cup
\{x+he_i:i=1,\ldots,n\}
\cup
\{x+he_i+he_j:1\leq i\leq j\leq n\},
\end{align*}
where $e_i$ denotes the $i$th column of the identity matrix. A symmetric Hessian approximation is given by
\begin{align*}
H(x,\xi^{(2)})
=
\sum_{i=1}^n\sum_{j=1}^n
\frac{
f(x+he_i+he_j,\xi_{ij})
-f(x+he_i,\xi_i)
-f(x+he_j,\xi_j)
+f(x,\xi_0)
}{h^2}
e_ie_j^T,
\end{align*}
where $\xi^{(2)}$ collects the random variables associated with the stochastic function evaluations (with the convention that $\xi_{ij}=\xi_{ji}$).

Let
\begin{align*}
c_n
=
\left[
n+
n(n-1)
\left(
\frac{\sqrt{2}+1}{3}
\right)^2
\right]^{1/2}.
\end{align*}
Under Assumption~\ref{assum:lip_cont2} and Oracle~\ref{def:oracle_stoch_0.1}, the bias of the resulting Hessian approximation satisfies
\begin{align*}
\left\|
\E_{\xi^{(2)}}[H(x,\xi^{(2)})]
-\nabla^2\phi(x)
\right\|_F
\leq
c_nL_2h
+
\frac{4n\epsilon_f}{h^2}.
\end{align*}
Moreover, using independence of the stochastic function evaluations,
\begin{align*}
\E_{\xi^{(2)}}
\left[
\left\|
H(x,\xi^{(2)})
-
\E_{\xi^{(2)}}[H(x,\xi^{(2)})]
\right\|_F^2
\right]
\leq
\frac{2n(2n+1)\sigma_f^2}{h^4}.
\end{align*}

Thus, taking $\psi(A,B)=\|A-B\|_F^2$, the construction satisfies Oracle~\ref{def:oracle_stoch_2a} with
\begin{align*}
\epsilon^2_{H\psi}(h,\lambda_f)
&=
\left(
c_nL_2h
+
\frac{4n\epsilon_f}{h^2}
\right)^2,
\\
\sigma^2_{H\psi}(h,\lambda_f)
&=
\frac{2n(2n+1)\sigma_f^2}{h^4}.
\end{align*}
Thus, the finite-difference interval $h$ and the control parameter $\lambda_f$ of Oracle~\ref{def:oracle_stoch_0.1} jointly determine the bias and stochastic variability of the resulting second-order oracle.

The construction requires $\frac{(n+1)(n+2)}{2}$ stochastic zeroth-order oracle evaluations. Thus, if each evaluation of Oracle~\ref{def:oracle_stoch_0.1} with input $\lambda$ has cost ${\cal C}_f(\lambda)$, the total computational cost is
\begin{align*}
{\cal C}_{H}(h,\lambda_f;\psi)
=
\frac{(n+1)(n+2)}{2}{\cal C}_f(\lambda_f)
+
\mathcal{O}(n^2).
\end{align*}

Furthermore, if $L_2>0$ and $h$ is chosen to satisfy $h
\geq
\left(
\frac{4n\epsilon_f}
{c_nL_2}
\right)^{1/3}$, then
\begin{align*}
\epsilon^2_{H\psi}(h,\lambda_f)
\leq
4c_n^2L_2^2h^2.
\end{align*}
Unlike the deterministic setting, however, controlling the bias alone does not control the stochastic variability, since $\sigma^2_{H\psi}(h,\lambda_f)$ grows as $h^{-4}$. Hence, $h$ and $\lambda$ should be selected jointly to balance the finite-difference truncation error, the bias of the zeroth-order oracle, and its variance.

\subsubsection{Finite Differences via Stochastic First-Order Oracles}

Given access to Oracle~\ref{def:oracle_stoch_1.2}, a Hessian approximation can be constructed by applying forward finite differences to stochastic gradient evaluations. For any $x\in\mathbb{R}^n$, choose $h>0$ and generate independent oracle outputs at $x$ and $x+he_i$, $i=1,\ldots,n$, where $e_i$ denotes the $i$th column of the identity matrix. A symmetric Hessian approximation is given by
\begin{align*}
H(x,\xi^{(2)})
&=
\sum_{i=1}^n
\left[
\left(
\frac{g(x+he_i,\xi_i)-g(x,\xi_0)}{2h}
\right)e_i^T
+
e_i
\left(
\frac{g(x+he_i,\xi_i)-g(x,\xi_0)}{2h}
\right)^T
\right],
\end{align*}
where $\xi^{(2)}=(\xi_0,\xi_1,\ldots,\xi_n)$. The symmetrization ensures that $H(x,\xi^{(2)})\in\mathbb{S}^n$, even though the finite-difference approximation obtained directly from the stochastic first-order oracle need not be symmetric.

Suppose that Oracle~\ref{def:oracle_stoch_1.2} holds with parameters $\epsilon_g$, $\kappa_g$, $\sigma_c$, and $\sigma_v$, where we suppress their dependence on the control parameter $\lambda$ for brevity. Under Assumptions~\ref{assum:lip_cont} and~\ref{assum:lip_cont2}, the bias term of the resulting Hessian approximation satisfies
\begin{align*}
\left\|
\E_{\xi^{(2)}}[H(x,\xi^{(2)})]
-
\nabla^2\phi(x)
\right\|^2_F
\leq
n
\left(
\frac{2\left(\epsilon_g+
\kappa_g\|\nabla\phi(x)\| \right) }{h}
+\kappa_gL_1+\frac{L_2h}{2}
\right)^2.
\end{align*}
As in the deterministic setting, assume that $\|\nabla\phi(x)\|\leq G$ over the region of interest. Then
\begin{align*}
\left\|
\E_{\xi^{(2)}}[H(x,\xi^{(2)})]
-
\nabla^2\phi(x)
\right\|_F^2
\leq
n
\left(
\frac{
2\left(
\epsilon_g
+
\kappa_g G
\right)
}{h}
+
\kappa_gL_1
+
\frac{L_2h}{2}
\right)^2.
\end{align*}

The stochastic variability can be bounded separately. Using independence of the stochastic gradient evaluations and Assumption~\ref{assum:lip_cont},
\begin{align*}
&\E_{\xi^{(2)}}
\left[
\left\|
H(x,\xi^{(2)})
-
\E_{\xi^{(2)}}[H(x,\xi^{(2)})]
\right\|_F^2
\right]\leq
\frac{n}{h^2}
\left[
2\sigma_c^2
+
\sigma_v^2
\left(
G^2+(G+L_1h)^2
\right)
\right].
\end{align*}

Thus, taking $\psi(A,B)=\|A-B\|_F^2$, the construction satisfies Oracle~\ref{def:oracle_stoch_2a} with
\begin{align*}
\epsilon^2_{H\psi}(h)
&=
n
\left(
\frac{
2\left(
\epsilon_g
+
\kappa_g G
\right)
}{h}
+
\kappa_gL_1
+
\frac{L_2h}{2}
\right)^2,
\\
\sigma^2_{H\psi}(h)
&=
\frac{n}{h^2}
\left[
2\sigma_c^2
+
\sigma_v^2
\left(
G^2+(G+L_1h)^2
\right)
\right].
\end{align*}

The construction requires $n+1$ stochastic first-order oracle evaluations. Thus, if each evaluation of Oracle~\ref{def:oracle_stoch_1.2} has cost ${\cal C}_g$, the total computational cost is  ${\cal C}_{H}(h;\psi)
=
(n+1){\cal C}_g
+
\mathcal{O}(n^2)$.

The finite-difference interval $h$ and the control parameter $\lambda$ jointly determine the quality and cost of the resulting stochastic second-order oracle. Decreasing $h$ reduces the finite-difference truncation error, but amplifies both the bias and stochastic variability inherited from the stochastic first-order oracle. Consequently, $h$ and $\lambda$ should be selected jointly to balance these competing effects.

\subsubsection{Directional Finite Differences for Stochastic Hessian-Vector Products}

Given access to Oracle~\ref{def:oracle_stoch_1.2}, a Hessian-vector product can be approximated directly without constructing a Hessian matrix. For any $x\in\mathbb{R}^n$, $v\neq0$, and $h>0$, define
\begin{align*}
    h(x,v,\xi^{(2)})
    =
    \frac{\|v\|}{h}
    \left[
        g\left(x+h\frac{v}{\|v\|},\xi_1\right)
        -
        g(x,\xi_0)
    \right],
\end{align*}
where $\xi_0$ and $\xi_1$ denote independent oracle samples. For $v=0$, we set $h(x,0,\xi^{(2)})=0$.

Suppose that Oracle~\ref{def:oracle_stoch_1.2} holds with parameters $\epsilon_g$, $\kappa_g$, $\sigma_c$, and $\sigma_v$, where we suppress their dependence on the control parameter $\lambda$ for brevity. Under Assumptions~\ref{assum:lip_cont} and~\ref{assum:lip_cont2}, and assuming $\|\nabla\phi(x)\|\leq G$ over the region of interest, the bias of the directional finite-difference approximation satisfies
\begin{align*}
\left\|
    \E_{\xi^{(2)}}[h(x,v,\xi^{(2)})]
    -
    \nabla^2\phi(x)v
\right\|
\leq
\left(
    \frac{2(\epsilon_g+\kappa_gG)}{h}
    +
    \kappa_gL_1
    +
    \frac{L_2h}{2}
\right)\|v\|.
\end{align*}

The stochastic variability can also be bounded separately. Using independence of the stochastic gradient evaluations,
\begin{align*}
&\E_{\xi^{(2)}}
\left[
    \left\|
        h(x,v,\xi^{(2)})
        -
        \E_{\xi^{(2)}}[h(x,v,\xi^{(2)})]
    \right\|^2
\right]\leq
\frac{\|v\|^2}{h^2}
\left[
    2\sigma_c^2
    +
    \sigma_v^2
    \left(
        G^2+(G+L_1h)^2
    \right)
\right].
\end{align*}

Consequently, by the triangle inequality and Jensen's inequality,
\begin{align*}
&\E_{\xi^{(2)}}
\left[
    \|h(x,v,\xi^{(2)})-\nabla^2\phi(x)v\|
\right]
\\
&\qquad\leq
\Bigg[
    \frac{2(\epsilon_g+\kappa_gG)}{h}
    +
    \kappa_gL_1
    +
    \frac{L_2h}{2}+
    \frac{1}{h}
    \sqrt{
        2\sigma_c^2
        +
        \sigma_v^2
        \left(
            G^2+(G+L_1h)^2
        \right)
    }
\Bigg]\|v\|.
\end{align*}
Thus, the construction satisfies Oracle~\ref{def:oracle_stoch_2b} with $p=1$ and
\begin{align*}
\kappa_{Hv}(h)
=
    \frac{2(\epsilon_g+\kappa_gG)}{h}
    +
    \kappa_gL_1
    +
    \frac{L_2h}{2}
    +
    \frac{1}{h}
    \sqrt{
        2\sigma_c^2
        +
        \sigma_v^2
        \left(
            G^2+(G+L_1h)^2
        \right)
    }.
\end{align*}

Each Hessian-vector product requires two evaluations of Oracle~\ref{def:oracle_stoch_1.2}. The finite-difference interval $h$ and the control parameter $\lambda$ jointly determine the accuracy and cost of the resulting stochastic Hessian-vector product. Decreasing $h$ reduces the finite-difference truncation error but amplifies both the bias and stochastic variability inherited from the stochastic first-order oracle.

\subsubsection{Randomized Subspace and Sketched Hessian Approximations}

Randomized sketching provides another mechanism for constructing inexpensive second-order approximations. Suppose that the Hessian admits a factorization
\begin{align*}
    \nabla^2\phi(x)=B(x)^T B(x),
\end{align*}
as occurs, for example, in many least-squares and generalized linear models. Given a random sketching matrix $S\in\mathbb{R}^{s\times m}$, where $s\ll m$, define
\begin{align*}
    H(x,S)
    =
    B(x)^T S^T S B(x).
\end{align*}
The sketch dimension $s$ controls both the accuracy and computational cost of the approximation
\cite{martinsson2020randomized,pilanci2017newton,berahas2020investigation}.

As a simple example, suppose that the entries of $S$ are independent Gaussian random variables with mean zero and variance $1/s$. Then
\begin{align*}
    \E[S^T S]=I,
\end{align*}
and hence
$\E[H(x,S)]=\nabla^2\phi(x)$, so the Gaussian sketch is unbiased. Standard moment properties of Gaussian and Wishart matrices imply
\begin{align*}
\E\left[
\|H(x,S)-\nabla^2\phi(x)\|_F^2
\right]
\leq
\frac{n+1}{s}
\|\nabla^2\phi(x)\|_F^2
\end{align*}
\cite{muirhead1982aspects}. Consequently, if $\|\nabla^2\phi(x)\|_F\leq M_H$ over the region of interest, then
\begin{align*}
\E\left[
\|H(x,S)-\nabla^2\phi(x)\|_F^2
\right]
\leq
\frac{(n+1)M_H^2}{s}.
\end{align*}
Thus, taking $\psi(A,B)=\|A-B\|_F^2$, the Gaussian sketch satisfies Oracle~\ref{def:oracle_stoch_2a} with
\begin{align*}
    \epsilon_{H\psi}(s)=0,
    \qquad
    \sigma^2_{H\psi}(s)
    =
    \frac{(n+1)M_H^2}{s}.
\end{align*}
Hence, increasing the sketch dimension reduces the stochastic variability while the approximation remains unbiased.

Other common choices include subsampled randomized Hadamard transforms and row-sampling sketches
\cite{pilanci2017newton,berahas2020investigation}. In particular, Hessian subsampling can be viewed as a special case of sketching. Suppose that
\begin{align*}
    \phi(x)
    =
    \frac{1}{N}\sum_{i=1}^N\phi_i(x),
    \qquad
    \nabla^2\phi_i(x)=B_i(x)^T B_i(x),
\end{align*}
and form $B(x)$ by stacking the appropriately scaled matrices $B_i(x)$. Choosing $S$ to select and rescale blocks of rows of $B(x)$ produces
\begin{align*}
    H(x,S)
    =
    \frac{1}{s}
    \sum_{j=1}^{s}
    \nabla^2\phi_{i_j}(x),
\end{align*}
which is precisely a subsampled Hessian. Suppose that, for some $M_H>0$,
\begin{align*}
    \frac{1}{N}\sum_{i=1}^N
    \left\|
        \nabla^2\phi_i(x)-\nabla^2\phi(x)
    \right\|_F^2
    \leq M_H^2.
\end{align*}
If the indices are sampled independently and uniformly with replacement, then
\begin{align*}
    \E[H(x,S)]
    =
    \nabla^2\phi(x),
\end{align*}
and
\begin{align*}
    \E\left[
        \|H(x,S)-\nabla^2\phi(x)\|_F^2
    \right]
    \leq
    \frac{M_H^2}{s}.
\end{align*}
Thus, taking $\psi(A,B)=\|A-B\|_F^2$, uniform Hessian subsampling satisfies Oracle~\ref{def:oracle_stoch_2a} with
\begin{align*}
    \epsilon_{H\psi}(s)=0,
    \qquad
    \sigma^2_{H\psi}(s)
    \leq
    \frac{M_H^2}{s}.
\end{align*}

If instead $s$ indices are sampled uniformly without replacement, the estimator remains unbiased and the finite-population correction yields
\begin{align*}
    \sigma^2_{H\psi}(s)
    \leq
    \frac{N-s}{N-1}\frac{M_H^2}{s}.
\end{align*}
Uniform and nonuniform sampling strategies, including appropriately rescaled schemes based on matrix norms or leverage-score information, have been studied extensively in this setting
\cite{roosta2019sub,xu2016subsampled,bollapragada2018adaptive,erdogdu2015convergence}.

Thus, Hessian subsampling and more general random projections can be viewed within a common sketching framework. Row or block sampling is inexpensive to apply, whereas mixing-based sketches such as randomized Hadamard or Gaussian projections can provide stronger spectral approximation properties at additional computational cost
\cite{pilanci2017newton,berahas2020investigation}. In each case, the sketch dimension and sampling distribution act as control parameters that determine the stochastic variability and computational cost of the resulting second-order oracle.

\subsection{Stochastic Second-Order Oracles in Optimization}\label{sec:stoch_second_order_optimization}

The use of stochastic second-order information in optimization introduces an additional trade-off between the accuracy of the derivative approximations and the computational cost required to obtain them. Many of the algorithmic principles discussed in Section~\ref{sec:det_second_order_optimization} continue to apply, but the accuracy of the gradient, Hessian, or Hessian-vector products must now be controlled through sampling or other variance-reduction mechanisms. Rather than repeating the deterministic discussion, we briefly highlight several representative ways in which stochastic second-order oracles arise in optimization algorithms, with particular emphasis on subsampled Newton methods, adaptive sampling and Hessian averaging, and stochastic variants of second-order globalization methods.

\paragraph{Subsampled and sketched Newton methods.}
Subsampled Newton methods reduce the cost of second-order information by constructing Newton directions from approximations of the Hessian obtained using only a subset of the available curvature information. Early approaches demonstrated that relatively small curvature samples can be effective within Newton--CG and quasi-Newton frameworks \cite{byrd2011use}. Subsequent analyses established global and local convergence guarantees and related the required sample sizes to the accuracy with which the subsampled Hessian approximates the true Hessian \cite{roosta2019sub,erdogdu2015convergence}. Nonuniform sampling schemes based on matrix norms or leverage-score information can further reduce the amount of curvature information required \cite{xu2016subsampled}, while appropriately increasing the accuracy of subsampled gradients and Hessians can yield superlinear convergence in expectation, including when the resulting Newton systems are solved  approximately \cite{bollapragada2018exact}.

More generally, Hessian subsampling can be viewed as a special case of randomized sketching, where the sketching matrix selects and rescales rows or blocks of rows of an underlying Hessian factorization. Other sketching constructions compress the curvature information using randomized projections rather than direct sampling. The Newton sketch, for example, replaces the Hessian by such a randomized low-dimensional approximation and can retain fast local convergence under suitable assumptions \cite{pilanci2017newton}. In both subsampling and sketching, the sample or sketch dimension acts as an accuracy-control parameter that governs the trade-off between the cost of obtaining second-order information and the quality of the resulting Newton direction.

\paragraph{Stochastic globalization and nonconvex optimization.}
Classical  globalization mechanisms such as trust-region and cubic regularization methods can also be extended to settings in which the gradient and second-order information are stochastic. Stochastic trust-region methods construct quadratic models from stochastic gradient and Hessian estimates and can retain convergence and complexity guarantees under suitable conditions on the accuracy and variability of these estimates \cite{curtis2022fully}.

Cubic-regularization methods have similarly been developed using subsampled or stochastic first- and second-order information. Subsampled gradients and Hessians can be used while retaining convergence guarantees when their accuracies are appropriately controlled \cite{kohler2017subsampled}, and fully stochastic variants based on stochastic gradients and Hessian-vector products can achieve complexity guarantees for approximate second-order stationarity \cite{tripuraneni2018stochastic}. More general adaptive-regularization frameworks allow the accuracy of inexact derivative information to vary dynamically and extend deterministic complexity results to stochastic settings through adaptive sampling \cite{bellavia2019adaptive}. Recent stochastic adaptive cubic-regularization methods further establish iteration and sample-complexity guarantees for first- and second-order stationarity using stochastic zeroth-, first-, and second-order information \cite{scheinberg2026stochastic}. Overall, these results illustrate that deterministic second-order globalization mechanisms can be retained in stochastic settings when the accuracy of the stochastic oracles is appropriately coordinated with the progress of the algorithm.

\newpage

\paragraph{Acknowledgements:} The work on this manuscript was support in part by ONR award N00014-25-1-2088 and N00014-24-1-2638.



\clearpage
\phantomsection
\addcontentsline{toc}{section}{References}
\bibliographystyle{plain}
\bibliography{Katya}

\newpage
\appendix
\counterwithin{table}{section}

 
\section{Notation}
\label{app:notation}

This appendix summarizes the notation used throughout the paper. We group the notation according to its role in the problem formulation, oracle models, and optimization algorithms. Symbols that are used only locally are defined in the corresponding sections and are generally omitted here. Specifically, 
symbols such as $C_1$, $C_2$, and $c$ are used to represent constants within a particular section to facilitate proofs but do not retain that meaning in other sections. Unless otherwise stated, $\|\cdot\|$ denotes the Euclidean norm for vectors and the induced spectral norm for matrices. When a symbol is used in more than one oracle model, its interpretation is determined by the corresponding oracle definition.

\subsection{Problem and Basic Notation}

\begin{longtable}{@{}p{0.25\textwidth}p{0.69\textwidth}@{}}
\toprule
\textbf{Symbol} & \textbf{Description} \\
\midrule
\endhead

$n$
& Dimension of the decision variable. \\
\hdashline[1pt/2pt]

$x\in\mathbb{R}^n$
& Decision variable. \\
\hdashline[1pt/2pt]

$x_k\in\mathbb{R}^n$
& Iterate at iteration $k$. \\
\hdashline[1pt/2pt]

$k\in\{0,1,\ldots\}$
& Iteration counter. \\
\hdashline[1pt/2pt]

$\phi:\mathbb{R}^n\rightarrow\mathbb{R}$
& True objective function. \\
\hdashline[1pt/2pt]

$\nabla\phi:\mathbb{R}^n\rightarrow\mathbb{R}^n$
& Gradient of the true objective function. \\
\hdashline[1pt/2pt]

$\nabla^2\phi:\mathbb{R}^n\rightarrow\mathbb{R}^{n\times n}$
& Hessian of the true objective function. \\
\hdashline[1pt/2pt]

$\phi^*=\inf_{x\in\mathbb{R}^n}\phi(x)$
& Optimal objective function value (infimum). \\
\hdashline[1pt/2pt]

$x^*$
& Optimal solution, when the infimum is attained. \\
\hdashline[1pt/2pt]

$X^*$
& Set of optimal solutions, when nonempty. \\
\hdashline[1pt/2pt]

$\widehat{\phi}$
& Lower bound on the objective function. \\
\hdashline[1pt/2pt]

$e_i$
& $i$th column of the identity matrix. \\
\hdashline[1pt/2pt]

$B(x,\delta)$
& Euclidean ball centered at $x$ with radius
$\delta\geq 0$. \\

\bottomrule
\end{longtable}

\subsection{Problem Constants and Regularity Parameters}

\begin{longtable}{@{}p{0.25\textwidth}p{0.69\textwidth}@{}}
\toprule
\textbf{Symbol} & \textbf{Description} \\
\midrule
\endhead

$L_0$
& Lipschitz constant of the function of $\phi$. \\
\hdashline[1pt/2pt]

$L_1$
& Lipschitz constant of the gradient of $\phi$. \\
\hdashline[1pt/2pt]

$L_2$
& Lipschitz constant of the Hessian of $\phi$. \\
\hdashline[1pt/2pt]

$\mu$
& Strong convexity or Polyak--{\L}ojasiewicz (PL)
parameter, according to context. \\
\hdashline[1pt/2pt]


$D$
& Diameter of the relevant level set. \\
\hdashline[1pt/2pt]

$c$
& Constant appearing in the Kurdyka--{\L}ojasiewicz (KL)
condition. \\
\hdashline[1pt/2pt]

$G$
& Upper bound on $\|\nabla\phi(x)\|$ over the region
under consideration. \\

\bottomrule
\end{longtable}

\subsection{Oracle Notation}

\begin{longtable}{@{}p{0.25\textwidth}p{0.69\textwidth}@{}}
\toprule
\textbf{Symbol} & \textbf{Description} \\
\midrule
\endhead

$f$
& Inexact zeroth-order approximation of the objective
function $\phi$. \\
\hdashline[1pt/2pt]

$g$
& Inexact first-order approximation of the gradient
$\nabla\phi$. \\
\hdashline[1pt/2pt]

$H$
& Inexact matrix-valued second-order approximation
of the Hessian $\nabla^2\phi$. \\
\hdashline[1pt/2pt]

$h(x,v)$
& Inexact approximation of the Hessian-vector product
$\nabla^2\phi(x)v$. \\
\hdashline[1pt/2pt]

$\lambda$
& Oracle control parameter, or collection of control
parameters, governing oracle accuracy and
computational cost. \\
\hdashline[1pt/2pt]

$\psi(A,B)$
& Error measure used to quantify the discrepancy
between two matrix-valued second-order quantities. \\
\hdashline[1pt/2pt]

$p\geq 1$
& Exponent controlling the dependence of
Hessian-vector product approximation errors on
the direction norm $\|v\|$. \\
\hdashline[1pt/2pt]

$\mathcal{C}_i$
& Computational cost associated with oracle $i$. \\
\hdashline[1pt/2pt]

$\mathcal{C}_f(\lambda)$
& Computational cost of a zeroth-order oracle call. \\
\hdashline[1pt/2pt]

$\mathcal{C}_g(\lambda)$
& Computational cost of a first-order oracle call. \\
\hdashline[1pt/2pt]

$\mathcal{C}_H(\lambda)$
& Computational cost of a matrix-valued second-order
oracle call. \\
\hdashline[1pt/2pt]

$\mathcal{C}_{Hv}(\lambda)$
& Computational cost of a Hessian-vector product
oracle call. \\
\hdashline[1pt/2pt]

$\xi^{(0)}$
& Randomness associated with a stochastic zeroth-order
oracle. \\
\hdashline[1pt/2pt]

$\xi^{(1)}$
& Randomness associated with a stochastic first-order
oracle. \\
\hdashline[1pt/2pt]

$\xi^{(2)}$
& Randomness associated with a stochastic second-order
oracle. \\

\bottomrule
\end{longtable}

\subsection{Deterministic Oracle Accuracy Parameters}

\begin{longtable}{@{}p{0.27\textwidth}p{0.67\textwidth}@{}}
\toprule
\textbf{Symbol} & \textbf{Description} \\
\midrule
\endhead

$\epsilon_f(\lambda)$
& Additive error bound associated with the deterministic
zeroth-order oracle. \\
\hdashline[1pt/2pt]

$\epsilon_g(\lambda)$
& Additive error bound associated with the deterministic
first-order oracle. \\
\hdashline[1pt/2pt]

$\kappa_g(\lambda)$
& Multiplicative error parameter associated with the
deterministic first-order oracle. \\
\hdashline[1pt/2pt]

$\kappa_{g\tau}(\lambda)$
& Modified multiplicative first-order error parameter,
defined by
$\kappa_{g\tau}^2(\lambda)
=\tau^2+(1-\tau^2)\kappa_g^2(\lambda)$. \\
\hdashline[1pt/2pt]

$\epsilon_{H\psi}(\lambda)$
& Accuracy parameter for the deterministic
matrix-valued second-order oracle under the
error measure $\psi$. \\
\hdashline[1pt/2pt]

$\kappa_{Hv}(\lambda)$
& Accuracy parameter for the deterministic
Hessian-vector product oracle, controlling
the approximation error relative to
$\|v\|^p$. \\

\bottomrule
\end{longtable}

\subsection{Stochastic Oracle Accuracy Parameters}

\begin{longtable}{@{}p{0.27\textwidth}p{0.67\textwidth}@{}}
\toprule
\textbf{Symbol} & \textbf{Description} \\
\midrule
\endhead

$\epsilon_f^2(\lambda)$
& Squared bias bound associated with a stochastic
zeroth-order oracle. \\
\hdashline[1pt/2pt]

$\sigma_f^2(\lambda)$
& Variance bound associated with the single-point
stochastic zeroth-order oracle. \\
\hdashline[1pt/2pt]

$e(x,\xi^{(0)})$
& Centered stochastic zeroth-order error, defined by
$e(x,\xi^{(0)})
=f(x,\xi^{(0)})
-\E_{\xi^{(0)}}[f(x,\xi^{(0)})]$. \\
\hdashline[1pt/2pt]

$L_c(\lambda)$
& Additive parameter controlling the mean-square
variation of stochastic zeroth-order errors
between multiple query points. \\
\hdashline[1pt/2pt]

$L_v(\lambda)$
& Gradient-dependent parameter controlling the
mean-square variation of stochastic zeroth-order
errors between multiple query points. \\
\hdashline[1pt/2pt]

$\epsilon_g^2(\lambda)$
& Additive component of the mean-square error
or squared bias bound associated with a
stochastic first-order oracle. \\
\hdashline[1pt/2pt]

$\kappa_g^2(\lambda)$
& Gradient-dependent component of the mean-square
error or squared bias bound associated with
a stochastic first-order oracle. \\
\hdashline[1pt/2pt]

$\sigma_c^2(\lambda)$
& Additive component of the variance bound associated
with a stochastic first-order oracle. \\
\hdashline[1pt/2pt]

$\sigma_v^2(\lambda)$
& Gradient-dependent component of the variance
bound associated with a stochastic first-order
oracle. \\
\hdashline[1pt/2pt]

$\zeta_1,\zeta_2$
& Random variables used in the conditional stochastic
first-order oracle, with
$\xi^{(1)}=(\zeta_1,\zeta_2)$. \\
\hdashline[1pt/2pt]

$\epsilon_{H\psi}(\lambda)$
& Bias bound for the stochastic matrix-valued
second-order oracle, measured by $\psi$,
controlling the discrepancy between
$\E_{\xi^{(2)}}[H(x,\xi^{(2)})]$
and $\nabla^2\phi(x)$. \\
\hdashline[1pt/2pt]

$\sigma_{H\psi}(\lambda)$
& Variability bound for the stochastic matrix-valued
second-order oracle, measured by $\psi$,
controlling the expected discrepancy between
$H(x,\xi^{(2)})$ and its expectation. \\
\hdashline[1pt/2pt]

$\kappa_{Hv}(\lambda)$
& Accuracy parameter for the stochastic
Hessian-vector product oracle, controlling
the expected approximation error relative
to $\|v\|^p$. \\

\bottomrule
\end{longtable}

\subsection{Optimization and Algorithmic Parameters}

\begin{longtable}{@{}p{0.25\textwidth}p{0.69\textwidth}@{}}
\toprule
\textbf{Symbol} & \textbf{Description} \\
\midrule
\endhead

$\epsilon > 0$
&  Desired accuracy aimed by the algorithm. \\
\hdashline[1pt/2pt]

%
%

$\alpha_k$ (or $\alpha$)
& Step-size parameter. \\
\hdashline[1pt/2pt]

$\alpha_{\max}$
& Maximum allowable step size. \\
\hdashline[1pt/2pt]

$h>0$
& Finite-difference interval. \\
\hdashline[1pt/2pt]

$\Delta_k$
& Trust-region radius at iteration $k$. \\
\hdashline[1pt/2pt]

$c_1$
& Line-search sufficient-decrease parameter. \\
\hdashline[1pt/2pt]

$\rho\in(0,1)$
& Backtracking contraction factor. \\
\hdashline[1pt/2pt]

$r_k\geq 0$
& Line-search correction factor. \\
\hdashline[1pt/2pt]

$\eta_1\in(0,1)$, $\eta_2>0$
& Parameters used in trust-region acceptance
and radius-update rules. \\
\hdashline[1pt/2pt]

$\gamma\in(0,1)$
& Trust-region radius contraction factor. \\
\hdashline[1pt/2pt]

$0<\beta_1\leq\beta_2$
& Parameters appearing in the angle condition. \\
\hdashline[1pt/2pt]

$\theta\in(0,\frac{1}{2}]$
& Descent parameter appearing in the general descent
inequality. \\
\hdashline[1pt/2pt]

$\omega, \nu \geq 0$
& Error parameters appearing in the general descent
inequality. \\
\hdashline[1pt/2pt]

$\tau\in[0,1)$
& Parameter appearing in the neighborhood
and convergence analysis. \\
\hdashline[1pt/2pt]

$\epsilon_\tau$
& Neighborhood accuracy parameter associated
with $\tau$. \\

\bottomrule
\end{longtable}

\subsection{Other Notation}

\begin{longtable}{@{}p{0.25\textwidth}p{0.69\textwidth}@{}}
\toprule
\textbf{Symbol} & \textbf{Description} \\
\midrule
\endhead

$s$
& Subspace dimension, with $1\leq s\leq n$. \\
\hdashline[1pt/2pt]

$S\in\mathbb{R}^{n\times s}$
& Matrix with orthonormal columns defining
an $s$-dimensional subspace. \\
\hdashline[1pt/2pt]

$\mathcal{S}(x)$
& Set of cardinality $s$ corresponding to the
components of $\nabla\phi(x)$ with largest
absolute values. \\
\hdashline[1pt/2pt]

$\mathcal{S}$
& Subset of $\mathbb{R}^n$ on which a regularity
condition is imposed, when used in this context. \\
\hdashline[1pt/2pt]

$\mathcal{Y}$
& Finite set of points queried by a multi-point
oracle or finite-difference construction. \\
\hdashline[1pt/2pt]

$\kappa_B$
& Condition number of a matrix $B$. \\
\hdashline[1pt/2pt]

$\mathbb{S}^n$
& Space of $n\times n$ real symmetric matrices. \\
\hdashline[1pt/2pt]

$\lambda_{\min}(A)$
& Smallest eigenvalue of a symmetric matrix $A$. \\
\hdashline[1pt/2pt]

$\lambda_{\max}(A)$
& Largest eigenvalue of a symmetric matrix $A$. \\
\hdashline[1pt/2pt]

$\|\cdot\|_F$
& Frobenius norm of a matrix. \\
\hdashline[1pt/2pt]

$\langle u,v\rangle$
& Euclidean inner product of vectors $u$ and $v$. \\
\hdashline[1pt/2pt]

$\E[\cdot]$
& Expectation with respect to the indicated
randomness. \\
\hdashline[1pt/2pt]

$\mathbb{P}[\cdot]$
& Probability with respect to the indicated
randomness. \\
\hdashline[1pt/2pt]

$|\mathcal{S}|$
& Cardinality of a finite set $\mathcal{S}$. \\

\bottomrule
\end{longtable}

\newpage

\section{Summary Tables}\label{app.summary}

\setcounter{table}{0}

In this appendix, we collect summary tables of the oracle constructions and convergence and complexity results presented throughout the paper. The tables summarize the properties of different inexact oracles, including their error bounds and associated constants, as well as the convergence guarantees and iteration complexity bounds established for the corresponding optimization algorithms. The goal is to provide a convenient reference and facilitate comparisons across different oracle constructions, error models, and problem settings.

\subsection{Summary of Deterministic Inexact First-Order Oracle Examples}
\label{sec:summary_det_oracles}

Table~\ref{tbl.example_deterministic} summarizes the deterministic inexact first-order oracle constructions presented in Section~\ref{sec.examples_first}. For each example, we report the relevant input parameters and the corresponding absolute and relative error bounds. The table provides a convenient comparison of the different oracle constructions and illustrates how their accuracy depends on the underlying approximation mechanisms and associated parameters.

\begin{table}[htbp]
\centering
\renewcommand{\arraystretch}{1.2}
\setlength{\dashlinedash}{0.5pt}
\setlength{\dashlinegap}{1.5pt}

\caption{Summary of deterministic inexact first-order oracle constructions satisfying Oracle~\ref{def:oracle_det_1}. For each construction, we report the input parameters and the corresponding absolute and relative error bounds, $\epsilon_g(\lambda)$ and $\kappa_g(\lambda)$, respectively.}
\label{tbl.example_deterministic}

\resizebox{\linewidth}{!}{%
\begin{tabular}{|l|l|c|c|}
\hline

\textbf{First-Order Oracle Example}
& \textbf{Inputs ($\lambda$)}
& \textbf{Absolute Error ($\epsilon_g(\lambda)$)}
& \textbf{Relative Error ($\kappa_g(\lambda)$)}
\\
\hline

Forward Finite Differences (FFD)
& Finite difference interval $h>0$
& $\displaystyle \frac{\sqrt{n}L_1 h}{2}
  + \frac{2\sqrt{n}\epsilon_f}{h}$
& $0$
\\[8pt]
\hdashline

Central Finite Differences (CFD)
& Finite difference interval $h>0$
& $\displaystyle \frac{\sqrt{n}L_2 h^2}{6}
  + \frac{\sqrt{n}\epsilon_f}{h}$
& $0$
\\[8pt]
\hdashline

Block-Coordinate
& Subspace dimension $1\leq s\leq n$
& $0$
& $\displaystyle \sqrt{\frac{n-s}{n}}$
\\[5pt]
\hdashline

Orthogonal Projected Gradient
& Orthonormal matrix $S\in\mathbb{R}^{n\times s}$
& $0$
& $\sqrt{1-\theta^2}$, $0\leq\theta\leq 1$
\\
\hdashline

Block-Coordinate and FFD
& \shortstack[l]{
    Finite difference interval $h>0$\\
    Subspace dimension $1\leq s\leq n$
  }
& $\displaystyle \frac{\sqrt{s}L_1 h}{2}
  + \frac{2\sqrt{s}\epsilon_f}{h}$
& $\displaystyle \sqrt{\frac{n-s}{n}}$
\\[8pt]
\hdashline

\shortstack[l]{Orthogonal Projected Gradient\\and FFD}
& \shortstack[l]{
    Finite difference interval $h>0$\\
    Orthonormal matrix $S\in\mathbb{R}^{n\times s}$
  }
& $\displaystyle \frac{\sqrt{s}L_1 h}{2}
  + \frac{2\sqrt{s}\epsilon_f}{h}$
& $\sqrt{1-\theta^2}$, $0\leq\theta\leq 1$
\\[8pt]
\hdashline

Linear Interpolation
& Interpolation set (matrix) $M\in\mathbb{R}^{n\times n}$
& $\displaystyle \|\hat{M}^{-1}\|
  \left(\frac{\sqrt{n}L_1h}{2}
  +\frac{2\sqrt{n}\epsilon_f}{h}\right)$
& $0$
\\[8pt]
\hdashline

Simplex Gradients
& Sample matrix $M\in\mathbb{R}^{n\times p}$
& $\displaystyle \|\hat{M}^{\dagger}\|
  \left(\frac{\sqrt{p}L_1h}{2}
  +\frac{2\sqrt{p}\epsilon_f}{h}\right)$
& $\sqrt{1-\theta^2}$, $0\leq\theta\leq 1$
\\[8pt]
\hdashline

Gradient Quantization
& Quantization level
  $b\geq\lfloor\log_2(1+\sqrt{n})\rfloor+1$
& $0$
& $\displaystyle \sqrt{\frac{n}{(2^b-1)^2}}$
\\[8pt]

\hline
\end{tabular}%
}
\end{table}

\subsection{Summary of Convergence and Complexity Results for Fixed-Step Gradient Descent with Deterministic First-Order Oracles}
\label{sec:summary_det_fixed_step}

Tables~\ref{tab:nonconvex}-\ref{tab:PL} summarize the convergence neighborhoods and iteration complexity bounds established in Theorem~\ref{cor:all_complexity_fixed}, presented in Section~\ref{sec:fixed_step_det}, for the fixed-step gradient descent method (Algorithm~\ref{alg:first_order_fixed_step}). The tables report these guarantees for the different deterministic inexact first-order oracle constructions considered in the paper, under the general nonconvex, KL, and PL settings.

\begin{table}[htbp]
\centering
\renewcommand{\arraystretch}{1.2}
\setlength{\dashlinedash}{0.5pt}
\setlength{\dashlinegap}{1.5pt}

\caption{Summary of the convergence neighborhoods and iteration complexity bounds established in Theorem~\ref{cor:all_complexity_fixed} for fixed-step gradient descent (Algorithm~\ref{alg:first_order_fixed_step}) using deterministic inexact first-order oracles satisfying Oracle~\ref{def:oracle_det_1}, with constant stepsize $0<\alpha\leq 1/L_1$. The results are presented for the general nonconvex, KL, and PL settings. The oracle constructions and their input parameters are summarized in Table~\ref{tbl.example_deterministic}.}
\label{tbl.example_deterministic_constantstepsize}

\begin{subtable}[t]{\textwidth}
\centering
\caption{General nonconvex setting.}
\label{tab:nonconvex}
\resizebox{\linewidth}{!}{%
\begin{tabular}{|l|c|c|c|c|}
\hline
\textbf{First-Order Oracle Example}
& \textbf{Absolute Error ($\epsilon_g$)}
& \textbf{Relative Error ($\kappa_g$)}
& \shortstack[c]{\textbf{Convergence Neighborhood}\\ $\boldsymbol{\epsilon_\tau}$}
& \shortstack[c]{\textbf{Iteration Complexity}\\ $\boldsymbol{K_{\epsilon_\tau}}$}
\\
\hline

Forward Finite Differences (FFD)& $\displaystyle \frac{\sqrt{n}L_1h}{2}+\frac{2\sqrt{n}\epsilon_f}{h}$& $\displaystyle 0$& $\displaystyle \frac{1}{\tau}\left(\frac{\sqrt{n}L_1h}{2}+\frac{2\sqrt{n}\epsilon_f}{h}\right)$& $\displaystyle \left\lceil\frac{1}{\epsilon_\tau^2}\frac{2(\phi(x_0)-\hat{\phi})}{\alpha(1-\tau^2)}\right\rceil$
\\[8pt]
\hdashline

Central Finite Differences (CFD)& $\displaystyle \frac{\sqrt{n}L_2h^2}{6}+\frac{\sqrt{n}\epsilon_f}{h}$& $\displaystyle 0$& $\displaystyle \frac{1}{\tau}\left(\frac{\sqrt{n}L_2h^2}{6}+\frac{\sqrt{n}\epsilon_f}{h}\right)$& $\displaystyle \left\lceil\frac{1}{\epsilon_\tau^2}\frac{2(\phi(x_0)-\hat{\phi})}{\alpha(1-\tau^2)}\right\rceil$
\\[8pt]
\hdashline

Block-Coordinate& $\displaystyle 0$& $\displaystyle \sqrt{\frac{n-s}{n}}$& $\displaystyle 0$& $\displaystyle \left\lceil\frac{1}{\epsilon_\tau^2}\frac{2(\phi(x_0)-\hat{\phi})}{\alpha\left(\frac{s}{n}\right)(1-\tau^2)}\right\rceil$
\\[8pt]
\hdashline

\shortstack[l]{Orthogonal Projected\\Gradient}& $\displaystyle 0$& $\displaystyle \sqrt{1-\theta^2},\quad 0\leq\theta\leq 1$& $\displaystyle 0$& $\displaystyle \left\lceil\frac{1}{\epsilon_\tau^2}\frac{2(\phi(x_0)-\hat{\phi})}{\alpha\theta^2(1-\tau^2)}\right\rceil$
\\[8pt]
\hdashline

Block-Coordinate and FFD& $\displaystyle \frac{\sqrt{s}L_1h}{2}+\frac{2\sqrt{s}\epsilon_f}{h}$& $\displaystyle \sqrt{\frac{n-s}{n}}$& $\displaystyle \frac{1}{\tau}\sqrt{\frac{n}{s}}\left(\frac{\sqrt{s}L_1h}{2}+\frac{2\sqrt{s}\epsilon_f}{h}\right)$& $\displaystyle \left\lceil\frac{1}{\epsilon_\tau^2}\frac{2(\phi(x_0)-\hat{\phi})}{\alpha\left(\frac{s}{n}\right)(1-\tau^2)}\right\rceil$
\\[8pt]
\hdashline

\shortstack[l]{Orthogonal Projected Gradient\\and FFD}& $\displaystyle \frac{\sqrt{s}L_1h}{2}+\frac{2\sqrt{s}\epsilon_f}{h}$& $\displaystyle \sqrt{1-\theta^2},\quad 0\leq\theta\leq 1$& $\displaystyle \frac{1}{\tau\theta}\left(\frac{\sqrt{s}L_1h}{2}+\frac{2\sqrt{s}\epsilon_f}{h}\right)$& $\displaystyle \left\lceil\frac{1}{\epsilon_\tau^2}\frac{2(\phi(x_0)-\hat{\phi})}{\alpha\theta^2(1-\tau^2)}\right\rceil$
\\[8pt]
\hdashline

Linear Interpolation& $\displaystyle \|\hat{M}^{-1}\|\left(\frac{\sqrt{n}L_1h}{2}+\frac{2\sqrt{n}\epsilon_f}{h}\right)$& $\displaystyle 0$& $\displaystyle \frac{1}{\tau}\|\hat{M}^{-1}\|\left(\frac{\sqrt{n}L_1h}{2}+\frac{2\sqrt{n}\epsilon_f}{h}\right)$& $\displaystyle \left\lceil\frac{1}{\epsilon_\tau^2}\frac{2(\phi(x_0)-\hat{\phi})}{\alpha(1-\tau^2)}\right\rceil$
\\[8pt]
\hdashline

Simplex Gradients& $\displaystyle \|\hat{M}^{\dagger}\|\left(\frac{\sqrt{p}L_1h}{2}+\frac{2\sqrt{p}\epsilon_f}{h}\right)$& $\displaystyle \sqrt{1-\theta^2},\quad 0\leq\theta\leq 1$& $\displaystyle \frac{1}{\tau\theta}\|\hat{M}^{\dagger}\|\left(\frac{\sqrt{p}L_1h}{2}+\frac{2\sqrt{p}\epsilon_f}{h}\right)$& $\displaystyle \left\lceil\frac{1}{\epsilon_\tau^2}\frac{2(\phi(x_0)-\hat{\phi})}{\alpha\theta^2(1-\tau^2)}\right\rceil$
\\[8pt]
\hdashline

Gradient Quantization& $\displaystyle 0$& $\displaystyle \sqrt{\frac{n}{(2^b-1)^2}}$& $\displaystyle 0$& $\displaystyle \left\lceil\frac{1}{\epsilon_\tau^2}\frac{2(\phi(x_0)-\hat{\phi})}{\alpha\left(1-\frac{n}{(2^b-1)^2}\right)(1-\tau^2)}\right\rceil$
\\[8pt]
\hline

\end{tabular}%
}
\end{subtable}
\end{table}

\begin{table}[htbp]
\ContinuedFloat
\centering
\renewcommand{\arraystretch}{1.2}
\setlength{\dashlinedash}{0.5pt}
\setlength{\dashlinegap}{1.5pt}

\caption[]{Summary of convergence neighborhoods and iteration complexity bounds (continued).}

\begin{subtable}[t]{\textwidth}
\centering
\caption{KL setting.}
\label{tab:KL}
\resizebox{\linewidth}{!}{%
\begin{tabular}{|l|c|c|c|c|}
\hline
\textbf{First-Order Oracle Example}
& \textbf{Absolute Error ($\epsilon_g$)}
& \textbf{Relative Error ($\kappa_g$)}
& \shortstack[c]{\textbf{Convergence Neighborhood}\\ $\boldsymbol{\epsilon_\tau}$}
& \shortstack[c]{\textbf{Iteration Complexity}\\ $\boldsymbol{K_{\epsilon_\tau}}$}
\\
\hline

Forward Finite Differences (FFD)& $\displaystyle \frac{\sqrt{n}L_1h}{2}+\frac{2\sqrt{n}\epsilon_f}{h}$& $\displaystyle 0$& $\displaystyle \frac{D}{\tau}\left(\frac{\sqrt{n}L_1h}{2}+\frac{2\sqrt{n}\epsilon_f}{h}\right)$& $\displaystyle \left\lceil\frac{1}{\epsilon_\tau}\frac{4D^2}{\alpha(1-\tau^2)}\right\rceil$
\\[8pt]
\hdashline

Central Finite Differences (CFD)& $\displaystyle \frac{\sqrt{n}L_2h^2}{6}+\frac{\sqrt{n}\epsilon_f}{h}$& $\displaystyle 0$& $\displaystyle \frac{D}{\tau}\left(\frac{\sqrt{n}L_2h^2}{6}+\frac{\sqrt{n}\epsilon_f}{h}\right)$& $\displaystyle \left\lceil\frac{1}{\epsilon_\tau}\frac{4D^2}{\alpha(1-\tau^2)}\right\rceil$
\\[8pt]
\hdashline

Block-Coordinate& $\displaystyle 0$& $\displaystyle \sqrt{\frac{n-s}{n}}$& $\displaystyle 0$& $\displaystyle \left\lceil\frac{1}{\epsilon_\tau}\frac{4D^2}{\alpha\left(\frac{s}{n}\right)(1-\tau^2)}\right\rceil$
\\[8pt]
\hdashline

\shortstack[l]{Orthogonal Projected\\Gradient}& $\displaystyle 0$& $\displaystyle \sqrt{1-\theta^2},\quad 0\leq\theta\leq 1$& $\displaystyle 0$& $\displaystyle \left\lceil\frac{1}{\epsilon_\tau}\frac{4D^2}{\alpha\theta^2(1-\tau^2)}\right\rceil$
\\[8pt]
\hdashline

Block-Coordinate and FFD& $\displaystyle \frac{\sqrt{s}L_1h}{2}+\frac{2\sqrt{s}\epsilon_f}{h}$& $\displaystyle \sqrt{\frac{n-s}{n}}$& $\displaystyle \frac{D}{\tau}\sqrt{\frac{n}{s}}\left(\frac{\sqrt{s}L_1h}{2}+\frac{2\sqrt{s}\epsilon_f}{h}\right)$& $\displaystyle \left\lceil\frac{1}{\epsilon_\tau}\frac{4D^2}{\alpha\left(\frac{s}{n}\right)(1-\tau^2)}\right\rceil$
\\[8pt]
\hdashline

\shortstack[l]{Orthogonal Projected Gradient\\and FFD}& $\displaystyle \frac{\sqrt{s}L_1h}{2}+\frac{2\sqrt{s}\epsilon_f}{h}$& $\displaystyle \sqrt{1-\theta^2},\quad 0\leq\theta\leq 1$& $\displaystyle \frac{D}{\tau\theta}\left(\frac{\sqrt{s}L_1h}{2}+\frac{2\sqrt{s}\epsilon_f}{h}\right)$& $\displaystyle \left\lceil\frac{1}{\epsilon_\tau}\frac{4D^2}{\alpha\theta^2(1-\tau^2)}\right\rceil$
\\[8pt]
\hdashline

Linear Interpolation& $\displaystyle \|\hat{M}^{-1}\|\left(\frac{\sqrt{n}L_1h}{2}+\frac{2\sqrt{n}\epsilon_f}{h}\right)$& $\displaystyle 0$& $\displaystyle \frac{D}{\tau}\|\hat{M}^{-1}\|\left(\frac{\sqrt{n}L_1h}{2}+\frac{2\sqrt{n}\epsilon_f}{h}\right)$& $\displaystyle \left\lceil\frac{1}{\epsilon_\tau}\frac{4D^2}{\alpha(1-\tau^2)}\right\rceil$
\\[8pt]
\hdashline

Simplex Gradients& $\displaystyle \|\hat{M}^{\dagger}\|\left(\frac{\sqrt{p}L_1h}{2}+\frac{2\sqrt{p}\epsilon_f}{h}\right)$& $\displaystyle \sqrt{1-\theta^2},\quad 0\leq\theta\leq 1$& $\displaystyle \frac{D}{\tau\theta}\|\hat{M}^{\dagger}\|\left(\frac{\sqrt{p}L_1h}{2}+\frac{2\sqrt{p}\epsilon_f}{h}\right)$& $\displaystyle \left\lceil\frac{1}{\epsilon_\tau}\frac{4D^2}{\alpha\theta^2(1-\tau^2)}\right\rceil$
\\[8pt]
\hdashline

Gradient Quantization& $\displaystyle 0$& $\displaystyle \sqrt{\frac{n}{(2^b-1)^2}}$& $\displaystyle 0$& $\displaystyle \left\lceil\frac{1}{\epsilon_\tau}\frac{4D^2}{\alpha\left(1-\frac{n}{(2^b-1)^2}\right)(1-\tau^2)}\right\rceil$
\\[8pt]
\hline

\end{tabular}%
}
\end{subtable}
\end{table}

\begin{table}[htbp]
\ContinuedFloat
\centering
\renewcommand{\arraystretch}{1.2}
\setlength{\dashlinedash}{0.5pt}
\setlength{\dashlinegap}{1.5pt}

\caption[]{Summary of convergence neighborhoods and iteration complexity bounds (continued).}

\begin{subtable}[t]{\textwidth}
\centering
\caption{PL setting.}
\label{tab:PL}
\resizebox{\linewidth}{!}{%
\begin{tabular}{|l|c|c|c|c|}
\hline
\textbf{First-Order Oracle Example}
& \textbf{Absolute Error ($\epsilon_g$)}
& \textbf{Relative Error ($\kappa_g$)}
& \shortstack[c]{\textbf{Convergence Neighborhood}\\ $\boldsymbol{\epsilon_\tau}$}
& \shortstack[c]{\textbf{Iteration Complexity}\\ $\boldsymbol{K_{\epsilon_\tau}}$}
\\
\hline

Forward Finite Differences (FFD)& $\displaystyle \frac{\sqrt{n}L_1h}{2}+\frac{2\sqrt{n}\epsilon_f}{h}$& $\displaystyle 0$& $\displaystyle \frac{1}{2\mu\tau^2}\left(\frac{\sqrt{n}L_1h}{2}+\frac{2\sqrt{n}\epsilon_f}{h}\right)^2$& $\displaystyle \left\lceil\frac{1}{\alpha\mu(1-\tau^2)}\log\left(\frac{\phi(x_0)-\phi^*}{\epsilon_\tau}\right)\right\rceil$
\\[8pt]
\hdashline

Central Finite Differences (CFD)& $\displaystyle \frac{\sqrt{n}L_2h^2}{6}+\frac{\sqrt{n}\epsilon_f}{h}$& $\displaystyle 0$& $\displaystyle \frac{1}{2\mu\tau^2}\left(\frac{\sqrt{n}L_2h^2}{6}+\frac{\sqrt{n}\epsilon_f}{h}\right)^2$& $\displaystyle \left\lceil\frac{1}{\alpha\mu(1-\tau^2)}\log\left(\frac{\phi(x_0)-\phi^*}{\epsilon_\tau}\right)\right\rceil$
\\[8pt]
\hdashline

Block-Coordinate& $\displaystyle 0$& $\displaystyle \sqrt{\frac{n-s}{n}}$& $\displaystyle 0$& $\displaystyle \left\lceil\frac{1}{\alpha\mu\left(\frac{s}{n}\right)(1-\tau^2)}\log\left(\frac{\phi(x_0)-\phi^*}{\epsilon_\tau}\right)\right\rceil$
\\[8pt]
\hdashline

\shortstack[l]{Orthogonal Projected\\Gradient}& $\displaystyle 0$& $\displaystyle \sqrt{1-\theta^2},\quad 0\leq\theta\leq 1$& $\displaystyle 0$& $\displaystyle \left\lceil\frac{1}{\alpha\mu\theta^2(1-\tau^2)}\log\left(\frac{\phi(x_0)-\phi^*}{\epsilon_\tau}\right)\right\rceil$
\\[8pt]
\hdashline

Block-Coordinate and FFD& $\displaystyle \frac{\sqrt{s}L_1h}{2}+\frac{2\sqrt{s}\epsilon_f}{h}$& $\displaystyle \sqrt{\frac{n-s}{n}}$& $\displaystyle \frac{1}{2\mu\tau^2\left(\frac{s}{n}\right)}\left(\frac{\sqrt{s}L_1h}{2}+\frac{2\sqrt{s}\epsilon_f}{h}\right)^2$& $\displaystyle \left\lceil\frac{1}{\alpha\mu\left(\frac{s}{n}\right)(1-\tau^2)}\log\left(\frac{\phi(x_0)-\phi^*}{\epsilon_\tau}\right)\right\rceil$
\\[8pt]
\hdashline

\shortstack[l]{Orthogonal Projected Gradient\\and FFD}& $\displaystyle \frac{\sqrt{s}L_1h}{2}+\frac{2\sqrt{s}\epsilon_f}{h}$& $\displaystyle \sqrt{1-\theta^2},\quad 0\leq\theta\leq 1$& $\displaystyle \frac{1}{2\mu\tau^2\theta^2}\left(\frac{\sqrt{s}L_1h}{2}+\frac{2\sqrt{s}\epsilon_f}{h}\right)^2$& $\displaystyle \left\lceil\frac{1}{\alpha\mu\theta^2(1-\tau^2)}\log\left(\frac{\phi(x_0)-\phi^*}{\epsilon_\tau}\right)\right\rceil$
\\[8pt]
\hdashline

Linear Interpolation& $\displaystyle \|\hat{M}^{-1}\|\left(\frac{\sqrt{n}L_1h}{2}+\frac{2\sqrt{n}\epsilon_f}{h}\right)$& $\displaystyle 0$& $\displaystyle \frac{\|\hat{M}^{-1}\|^2}{2\mu\tau^2}\left(\frac{\sqrt{n}L_1h}{2}+\frac{2\sqrt{n}\epsilon_f}{h}\right)^2$& $\displaystyle \left\lceil\frac{1}{\alpha\mu(1-\tau^2)}\log\left(\frac{\phi(x_0)-\phi^*}{\epsilon_\tau}\right)\right\rceil$
\\[8pt]
\hdashline

Simplex Gradients& $\displaystyle \|\hat{M}^{\dagger}\|\left(\frac{\sqrt{p}L_1h}{2}+\frac{2\sqrt{p}\epsilon_f}{h}\right)$& $\displaystyle \sqrt{1-\theta^2},\quad 0\leq\theta\leq 1$& $\displaystyle \frac{\|\hat{M}^{\dagger}\|^2}{2\mu\tau^2\theta^2}\left(\frac{\sqrt{p}L_1h}{2}+\frac{2\sqrt{p}\epsilon_f}{h}\right)^2$& $\displaystyle \left\lceil\frac{1}{\alpha\mu\theta^2(1-\tau^2)}\log\left(\frac{\phi(x_0)-\phi^*}{\epsilon_\tau}\right)\right\rceil$
\\[8pt]
\hdashline

Gradient Quantization& $\displaystyle 0$& $\displaystyle \sqrt{\frac{n}{(2^b-1)^2}}$& $\displaystyle 0$& $\displaystyle \left\lceil\frac{1}{\alpha\mu\left(1-\frac{n}{(2^b-1)^2}\right)(1-\tau^2)}\log\left(\frac{\phi(x_0)-\phi^*}{\epsilon_\tau}\right)\right\rceil$
\\[8pt]
\hline

\end{tabular}%
}
\end{subtable}
\end{table}

\pagebreak

\subsection{Summary of Randomized Finite-Difference Oracle Parameters}
\label{sec:summary_rand_ffd}

Tables~\ref{tab:randomized_ffd_unscaled}-\ref{tab:randomized_cfd_scaled} summarize the accuracy and variance parameters of the stochastic first-order oracles associated with the randomized finite-difference constructions discussed in this section. We present the results for forward and central finite differences, both with and without scaling, highlighting their dependence on the distribution of the sampled directions and the assumptions on the underlying stochastic zeroth-order oracle. For forward finite differences, we consider both orthonormal and independent sampling directions, whereas the central finite-difference results are presented under the orthonormal sampling assumption.

\begin{table}[htbp]
\centering
\renewcommand{\arraystretch}{1.2}
\setlength{\dashlinedash}{0.5pt}
\setlength{\dashlinegap}{1.5pt}

\caption{The constants in Oracle~\ref{def:oracle_stoch_1.3}
for randomized forward finite difference without scaling (Sections~\ref{sec:RFFD_orth} and \ref{sec:RFFD}).}
\label{tab:randomized_ffd_unscaled}

\resizebox{\linewidth}{!}{%
\begin{tabular}{|l|l|c:c|}
\hline

\multicolumn{2}{|l|}{\textbf{Oracle/Parameter}}
& Assumption~\ref{ass:haar} (orthonormal)
& Assumption~\ref{ass:normal} (independent)
\\
\hline

\multicolumn{2}{|l|}{$\boldsymbol{\epsilon_g^2}$}
& $\displaystyle s \left(\frac{L_1 h}{2}+\frac{2\epsilon_f}{h}\right)^2$
& $\displaystyle \frac{2(n+s-1)}{n}\cdot s
\left(
\frac{(1+\mathcal{O}(n^{-1}))L_1h}{2}
+\frac{2\epsilon_f}{h}
\right)^2$
\\[8pt]
\hline

\multicolumn{2}{|l|}{$\boldsymbol{\kappa_g^2}$}
& $\displaystyle \frac{n-s}{n}$
& $\displaystyle 1-\frac{s}{n}
+\frac{2s^2}{n^2}
+\frac{\mathcal{O}(n^{-1})s}{n}$
\\[8pt]
\hline

\multirow{2}{*}{\textbf{Oracle~\ref{def:oracle_stoch_0.1}}}
& $\boldsymbol{\sigma_c^2}$
& $\displaystyle \frac{2s\sigma_f^2}{h^2}$
& $\displaystyle \frac{2s\sigma_f^2}{h^2}$
\\[8pt]
\cline{2-4}

& $\boldsymbol{\sigma_v^2}$
& $0$
& $0$
\\[5pt]
\hline

\multirow{2}{*}{\textbf{Oracle~\ref{def:oracle_stoch_0.2}}}
& $\boldsymbol{\sigma_c^2}$
& $\displaystyle sL_c^2$
& $\displaystyle
\left(\frac{n+s-1}{n}+\mathcal{O}(n^{-1})\right)\cdot sL_c^2$
\\[8pt]
\cline{2-4}

& $\boldsymbol{\sigma_v^2}$
& $\displaystyle sL_v^2$
& $\displaystyle
\left(\frac{n+s-1}{n}+\mathcal{O}(n^{-1})\right)\cdot sL_v^2$
\\[8pt]
\hline

\end{tabular}%
}
\end{table}

\begin{table}[htbp]
\centering
\renewcommand{\arraystretch}{1.2}
\setlength{\dashlinedash}{0.5pt}
\setlength{\dashlinegap}{1.5pt}

\caption{The constants in Oracle~\ref{def:oracle_stoch_1.2}
for randomized forward finite difference with scaling 
(Section~\ref{sec:SRFFD}).
The constant $\epsilon_g^2$ improves to
$\Big(L_1 h + \frac{\sqrt{n}\epsilon_f}{h}\Big)^2$
when $\mathcal{U} \sim \mathcal{N}$ or $\mathcal{S}$ (Section~\ref{sec:smoothing}).}
\label{tab:scaledffd}

\resizebox{\linewidth}{!}{%
\begin{tabular}{|l|l|c:c|}
\hline

\multicolumn{2}{|l|}{\textbf{Oracle/Parameter}}
& Assumption~\ref{ass:haar} (orthonormal)
& Assumption~\ref{ass:normal} (independent)
\\
\hline

\multicolumn{2}{|l|}{$\boldsymbol{\epsilon_g^2}$}
& $\displaystyle n \left(
\frac{L_1 h}{2}+\frac{2\epsilon_f}{h}
\right)^2$
& $\displaystyle n \left(
\frac{(1+\mathcal{O}(n^{-1}))L_1 h}{2}
+\frac{2\epsilon_f}{h}
\right)^2$
\\[8pt]
\hline

\multicolumn{2}{|l|}{$\boldsymbol{\kappa_g^2}$}
& $0$
& $0$
\\[8pt]
\hline

\multirow{2}{*}{\textbf{Oracle~\ref{def:oracle_stoch_0.1}}}
& $\boldsymbol{\sigma_c^2}$
& $\displaystyle
\frac{2n^2\sigma_f^2}{s h^2}
+\frac{n^2}{s}
\left(
\frac{L_1 h}{2}+\frac{2\epsilon_f}{h}
\right)^2$
& $\displaystyle
\frac{2n^2\sigma_f^2}{s h^2}
+\frac{n+s-1}{n}\cdot\frac{n^2}{s}
\left(
\frac{(1+\mathcal{O}(n^{-1}))L_1 h}{2}
+\frac{2\epsilon_f}{h}
\right)^2$
\\[8pt]
\cline{2-4}

& $\boldsymbol{\sigma_v^2}$
& $\displaystyle
\frac{(2n-s)(n-s)}{ns}$
& $\displaystyle
\frac{n-\bar{s}}{n}
+2\left(
\frac{(1+\mathcal{O}(n^{-1}))n}{s}
-2+\frac{\bar{s}}{n}
\right)$
\\[8pt]
\hline

\multirow{2}{*}{\textbf{Oracle~\ref{def:oracle_stoch_0.2}}}
& $\boldsymbol{\sigma_c^2}$
& $\displaystyle
\frac{n^2 L_c^2}{s}
+\frac{n^2}{s}
\left(
\frac{L_1 h}{2}+\frac{2\epsilon_f}{h}
\right)^2$
& $\displaystyle
\left(
\frac{n+s-1}{n}+\mathcal{O}(n^{-1})
\right)\cdot\frac{n^2 L_c^2}{s}
+\frac{n+s-1}{n}\cdot\frac{n^2}{s}
\left(
\frac{(1+\mathcal{O}(n^{-1}))L_1 h}{2}
+\frac{2\epsilon_f}{h}
\right)^2$
\\[8pt]
\cline{2-4}

& $\boldsymbol{\sigma_v^2}$
& $\displaystyle
\frac{n^2 L_v^2}{s}
+\frac{(2n-s)(n-s)}{ns}$
& $\displaystyle
\left(
\frac{n+s-1}{n}+\mathcal{O}(n^{-1})
\right)\cdot\frac{n^2 L_v^2}{s}
+\frac{n-\bar{s}}{n}
+2\left(
\frac{(1+\mathcal{O}(n^{-1}))n}{s}
-2+\frac{\bar{s}}{n}
\right)$
\\[8pt]
\hline

\end{tabular}%
}
\end{table}

\begin{table}[htbp]
\centering
\renewcommand{\arraystretch}{1.2}
\setlength{\dashlinedash}{0.5pt}
\setlength{\dashlinegap}{1.5pt}

\caption{The constants in Oracle~\ref{def:oracle_stoch_1.3}
for randomized central finite difference without scaling
under Assumption~\ref{ass:haar} (Section~\ref{sec:RCFD}).}
\label{tab:randomized_cfd_unscaled}

\scalebox{0.9}{%
\begin{tabular}{|l|l|c|}
\hline

\multicolumn{2}{|l|}{\textbf{Oracle/Parameter}}
& Assumption~\ref{ass:haar} (orthonormal)
\\
\hline

\multicolumn{2}{|l|}{$\boldsymbol{\epsilon_g^2}$}
& $\displaystyle
s \cdot \min \left\{
\left(\frac{L_1 h}{2}+\frac{2\epsilon_f}{h}\right)^2,
\left(\frac{L_2 h^2}{6}+\frac{\epsilon_f}{h}\right)^2
\right\}$
\\[8pt]
\hline

\multicolumn{2}{|l|}{$\boldsymbol{\kappa_g^2}$}
& $\displaystyle \frac{n-s}{n}$
\\[8pt]
\hline

\multirow{2}{*}{\textbf{Oracle~\ref{def:oracle_stoch_0.1}}}
& $\boldsymbol{\sigma_c^2}$
& $\displaystyle \frac{s\sigma_f^2}{2h^2}$
\\[8pt]
\cline{2-3}

& $\boldsymbol{\sigma_v^2}$
& $0$
\\[8pt]
\hline

\multirow{2}{*}{\textbf{Oracle~\ref{def:oracle_stoch_0.2}}}
& $\boldsymbol{\sigma_c^2}$
& $\displaystyle sL_c^2$
\\[8pt]
\cline{2-3}

& $\boldsymbol{\sigma_v^2}$
& $\displaystyle sL_v^2$
\\[8pt]
\hline

\end{tabular}%
}
\end{table}

\begin{table}[htbp]
\centering
\renewcommand{\arraystretch}{1.2}
\setlength{\dashlinedash}{0.5pt}
\setlength{\dashlinegap}{1.5pt}

\caption{The constants in Oracle~\ref{def:oracle_stoch_1.2}
for randomized central finite difference with scaling
under Assumption~\ref{ass:haar} (Section~\ref{sec:SRCFD}).
The constant $\epsilon_g^2$ improves to
$\min\left\{
\left(L_1 h+\frac{\sqrt{n}\epsilon_f}{h}\right)^2,
\left(\frac{L_2 h^2}{2}+\frac{\sqrt{n}\epsilon_f}{h}\right)^2
\right\}$
when $\mathcal{U}\sim\mathcal{N}$ or $\mathcal{S}$
(Section~\ref{sec:smoothing}).}
\label{tab:randomized_cfd_scaled}

\scalebox{0.9}{%
\begin{tabular}{|l|l|c|}
\hline

\multicolumn{2}{|l|}{\textbf{Oracle/Parameter}}
& Assumption~\ref{ass:haar} (orthonormal)
\\
\hline

\multicolumn{2}{|l|}{$\boldsymbol{\epsilon_g^2}$}
& $\displaystyle
n\cdot\min\left\{
\left(\frac{L_1 h}{2}+\frac{2\epsilon_f}{h}\right)^2,
\left(\frac{L_2 h^2}{6}+\frac{\epsilon_f}{h}\right)^2
\right\}$
\\[8pt]
\hline

\multicolumn{2}{|l|}{$\boldsymbol{\kappa_g^2}$}
& $0$
\\[8pt]
\hline

\multirow{2}{*}{\textbf{Oracle~\ref{def:oracle_stoch_0.1}}}
& $\boldsymbol{\sigma_c^2}$
& $\displaystyle
\frac{n^2\sigma_f^2}{2s h^2}
+\frac{n^2}{s}
\min\left\{
\left(\frac{L_1 h}{2}+\frac{2\epsilon_f}{h}\right)^2,
\left(\frac{L_2 h^2}{6}+\frac{\epsilon_f}{h}\right)^2
\right\}$
\\[8pt]
\cline{2-3}

& $\boldsymbol{\sigma_v^2}$
& $\displaystyle
\frac{(2n-s)(n-s)}{ns}$
\\[8pt]
\hline

\multirow{2}{*}{\textbf{Oracle~\ref{def:oracle_stoch_0.2}}}
& $\boldsymbol{\sigma_c^2}$
& $\displaystyle
\frac{n^2 L_c^2}{s}
+\frac{n^2}{s}
\min\left\{
\left(\frac{L_1 h}{2}+\frac{2\epsilon_f}{h}\right)^2,
\left(\frac{L_2 h^2}{6}+\frac{\epsilon_f}{h}\right)^2
\right\}$
\\[8pt]
\cline{2-3}

& $\boldsymbol{\sigma_v^2}$
& $\displaystyle
\frac{n^2 L_v^2}{s}
+\frac{(2n-s)(n-s)}{ns}$
\\[8pt]
\hline

\end{tabular}%
}
\end{table}

\newcommand{\paddedfrac}[2]{%
  \vphantom{\dfrac{A^{A^A}}{A_{A_A}}}%
  \dfrac{#1}{#2}%
  \vphantom{\dfrac{A^A}{A_A}}%
}

\newpage
\subsection{Summary of Convergence and Complexity Results for Fixed-Step Gradient Descent with Stochastic First-Order Oracles}
\label{sec:summary_stoch_fixed_step}

Tables~\ref{tab:stoch_fixed_step_oracle_11}--\ref{tab:stoch_fixed_step_oracle_13} summarize the convergence and iteration complexity guarantees for the fixed-step stochastic gradient descent method presented in Algorithm~\ref{alg:first_order_fixed_step_sgd}, using the stochastic first-order oracle models described in Oracles~\ref{def:oracle_stoch_1.1}--\ref{def:oracle_stoch_1.3}, respectively. Specifically, we specialize the general results to the oracle constructions and settings discussed in Sections~\ref{sec:erm}, \ref{strong_growth_reg}, \ref{sec:Byzantine}, \ref{sec:gradproj}, \ref{sec:rand_dirder}, \ref{sec:SRFFD}, and \ref{sec:SRCFD}, and report the corresponding convergence neighborhoods and expected iteration complexity bounds under the general nonconvex, KL, and PL settings.

\begin{table}[htbp]
\centering
\renewcommand{\arraystretch}{1.2}
\setlength{\dashlinedash}{0.5pt}
\setlength{\dashlinegap}{1.5pt}

\caption{Summary of complexity results for
Algorithm~\ref{alg:first_order_fixed_step_sgd}
using Oracle~\ref{def:oracle_stoch_1.1}.
Here, $L_1$ denotes Assumption~\ref{assum:lip_cont},
KL denotes Assumptions~\ref{assum:lip_cont}
and~\ref{assum:KLcond}, and PL denotes
Assumptions~\ref{assum:lip_cont}
and~\ref{assum:PLcond}.}
\label{tab:stoch_fixed_step_oracle_11}

\resizebox{\linewidth}{!}{%
\begin{tabular}{|l|l|c|c|c|}
\hline

\textbf{Setting}
& \multicolumn{2}{c|}{\textbf{Assumptions}}
& \makecell{\textbf{Convergence} \\ \textbf{Neighborhood} $\epsilon_\tau \geq \mathcal{O}(\cdot)$}
& \makecell{\textbf{iteration}\\  \textbf{Complexity} $\E[K_{\epsilon_\tau}] \leq \mathcal{O}(\cdot)$}
\\
\hline

\multirow{3}{*}{
    \shortstack[l]{
        Minibatch stochastic\\
        gradient descent\\
        for ERM\\
        (Sec~\ref{sec:erm})
    }
}
& \multirow{3}{*}{
    \shortstack[l]{
        $\alpha \leq \frac{1}{L_1}$,\\
        \eqref{BCN}, and\\
        minibatch size\\
        $|\mathcal{S}|>M_v^2$
    }
}
& $L_1$
& $\displaystyle
\frac{M_c}{\tau\sqrt{|\mathcal{S}|-M_v^2}}$
& $\displaystyle
\frac{\phi(x_0)-\hat{\phi}}
{\alpha(1-M_v^2/|\mathcal{S}|)(1-\tau^2)\epsilon_\tau^2}$
\\[8pt]
\cdashline{3-5}

& & KL
& $\displaystyle
\frac{M_cD}{\tau\sqrt{|\mathcal{S}|-M_v^2}}$
& $\displaystyle
\frac{D^2}
{\alpha(1-M_v^2/|\mathcal{S}|)(1-\tau^2)\epsilon_\tau}$
\\[8pt]
\cdashline{3-5}

& & PL
& $\displaystyle
\frac{M_c^2}{\mu\tau^2(|\mathcal{S}|-M_v^2)}$
& $\displaystyle
\frac{\log((\phi(x_0)-\hat{\phi})/\epsilon_\tau)}
{\alpha\mu(1-M_v^2/|\mathcal{S}|)(1-\tau^2)}$
\\[8pt]
\hdashline

\multirow{3}{*}{
    \shortstack[l]{
        Minibatch stochastic\\
        gradient descent\\
        for least squares\\
        (Sec~\ref{strong_growth_reg})
    }
}
& \multirow{3}{*}{
    \shortstack[l]{
        $\alpha \leq \frac{1}{L_1}$,\\
        \eqref{assum:M1}, and\\
        minibatch size\\
        $|\mathcal{S}|>2M_1$
    }
}
& $L_1$
& $\displaystyle
\frac{1}{\tau}
\sqrt{
\frac{\E_{\xi\sim\Xi}
[\|\nabla\ell(x^\star,\xi)\|^2]}
{|\mathcal{S}|-2M_1}
}$
& $\displaystyle
\frac{\phi(x_0)-\hat{\phi}}
{\alpha(1-2M_1/|\mathcal{S}|)(1-\tau^2)\epsilon_\tau^2}$
\\[8pt]
\cdashline{3-5}

& & KL
& $\displaystyle
\frac{D}{\tau}
\sqrt{
\frac{\E_{\xi\sim\Xi}
[\|\nabla\ell(x^\star,\xi)\|^2]}
{|\mathcal{S}|-2M_1}
}$
& $\displaystyle
\frac{D^2}
{\alpha(1-2M_1/|\mathcal{S}|)(1-\tau^2)\epsilon_\tau}$
\\[8pt]
\cdashline{3-5}

& & PL
& $\displaystyle
\frac{\E_{\xi\sim\Xi}
[\|\nabla\ell(x^\star,\xi)\|^2]}
{\mu\tau^2(|\mathcal{S}|-2M_1)}$
& $\displaystyle
\frac{\log((\phi(x_0)-\hat{\phi})/\epsilon_\tau)}
{\alpha\mu(1-2M_1/|\mathcal{S}|)(1-\tau^2)}$
\\[8pt]
\hdashline

\multirow{3}{*}{
    \shortstack[l]{
        Robust gradient descent\\
        in Byzantine machine\\
        setting with $m$ machines\\
        (Sec~\ref{sec:Byzantine})
    }
}
& \multirow{3}{*}{
    \shortstack[l]{
        $\alpha \leq \frac{1}{L_1}$,\\
        \eqref{BCN} holds for\\
        at least $(1-\theta)m$\\
        machines, and\\
        minibatch size\\
        $|\mathcal{S}|>
        \frac{(1-\theta)^3mM_v^2}{1-2\theta}$
    }
}
& $L_1$
& $\displaystyle
\frac{M_c}
{\tau\sqrt{
\frac{(1-2\theta)|\mathcal{S}|}{(1-\theta)^3m}
-M_v^2
}}$
& $\displaystyle
\frac{\phi(x_0)-\hat{\phi}}
{\alpha
\left(
1-\frac{(1-\theta)^3mM_v^2}
{(1-2\theta)|\mathcal{S}|}
\right)
(1-\tau^2)\epsilon_\tau^2}$
\\[8pt]
\cdashline{3-5}

& & KL
& $\displaystyle
\frac{M_cD}
{\tau\sqrt{
\frac{(1-2\theta)|\mathcal{S}|}{(1-\theta)^3m}
-M_v^2
}}$
& $\displaystyle
\frac{D^2}
{\alpha
\left(
1-\frac{(1-\theta)^3mM_v^2}
{(1-2\theta)|\mathcal{S}|}
\right)
(1-\tau^2)\epsilon_\tau}$
\\[8pt]
\cdashline{3-5}

& & PL
& $\displaystyle
\frac{M_c^2}
{\mu\tau^2
\left(
\frac{(1-2\theta)|\mathcal{S}|}{(1-\theta)^3m}
-M_v^2
\right)}$
& $\displaystyle
\frac{\log((\phi(x_0)-\hat{\phi})/\epsilon_\tau)}
{\alpha\mu
\left(
1-\frac{(1-\theta)^3mM_v^2}
{(1-2\theta)|\mathcal{S}|}
\right)
(1-\tau^2)}$
\\[8pt]
\hdashline

\multirow{3}{*}{
    \shortstack[l]{
        Gradient descent in\\
        random subspaces\\
        (Sec~\ref{sec:gradproj}) or\\
        random directional\\
        derivative descent\\
        (Sec~\ref{sec:rand_dirder})
    }
}
& \multirow{3}{*}{
    \shortstack[l]{
        $\alpha\leq\frac{1}{L_1}$,\\
        subspace dimension\\
        $1\leq s\leq n$ or\\
        number of directions\\
        $s\geq 1$
    }
}
& $L_1$
& $0$
& $\displaystyle
\frac{n}{s}\cdot
\frac{\phi(x_0)-\hat{\phi}}
{\alpha(1-\tau^2)\epsilon_\tau^2}$
\\[8pt]
\cdashline{3-5}

& & KL
& $0$
& $\displaystyle
\frac{n}{s}\cdot
\frac{D^2}
{\alpha(1-\tau^2)\epsilon_\tau}$
\\[8pt]
\cdashline{3-5}

& & PL
& $0$
& $\displaystyle
\frac{n}{s}\cdot
\frac{\log((\phi(x_0)-\hat{\phi})/\epsilon_\tau)}
{\alpha\mu(1-\tau^2)}$
\\[8pt]

\hline
\end{tabular}%
}
\end{table}

\begin{table}[htbp]
\centering
\renewcommand{\arraystretch}{1.2}
\setlength{\dashlinedash}{0.5pt}
\setlength{\dashlinegap}{1.5pt}

\caption{Summary of complexity results for
Algorithm~\ref{alg:first_order_fixed_step_sgd}
using Oracle~\ref{def:oracle_stoch_1.2}.
Here, $L_1$ denotes Assumption~\ref{assum:lip_cont},
KL denotes Assumptions~\ref{assum:lip_cont}
and~\ref{assum:KLcond}, and PL denotes
Assumptions~\ref{assum:lip_cont}
and~\ref{assum:PLcond}.}
\label{tab:stoch_fixed_step_oracle_12}

\resizebox{\linewidth}{!}{%
\begin{tabular}{|l|l|c|c|c|}
\hline

\textbf{Setting}
& \multicolumn{2}{c|}{\textbf{Assumptions}}
& \makecell{\textbf{Convergence} \\ \textbf{Neighborhood} $\epsilon_\tau \geq \mathcal{O}(\cdot)$}
& \makecell{\textbf{iteration}\\  \textbf{Complexity} $\E[K_{\epsilon_\tau}] \leq \mathcal{O}(\cdot)$}
\\
\hline


\multirow{3}{*}{
    \shortstack[l]{
        Minibatch gradient\\
        descent in ERM\\
        (Sec~\ref{sec:erm})
    }
}
& \multirow{3}{*}{
    \shortstack[l]{
        \eqref{BCN} and\\
        $\displaystyle \alpha\leq\frac{1}{L_1}\cdot\min\left\{1,
        \frac{|\mathcal{S}|}{2M_v^2}\right\}$
    }
}
& $L_1$
& $\displaystyle
\frac{M_c}{\tau}
\sqrt{\frac{\alpha L_1}{|\mathcal{S}|}}$
& $\displaystyle
\frac{\phi(x_0)-\hat{\phi}}
{\alpha(1-\tau^2)\epsilon_\tau^2}$
\\[8pt]
\cdashline{3-5}

& & KL
& $\displaystyle
\frac{M_cD}{\tau}
\sqrt{\frac{\alpha L_1}{|\mathcal{S}|}}$
& $\displaystyle
\frac{D^2}
{\alpha(1-\tau^2)\epsilon_\tau}$
\\[8pt]
\cdashline{3-5}

& & PL
& $\displaystyle
\frac{\alpha L_1M_c^2}
{\tau^2\mu|\mathcal{S}|}$
& $\displaystyle
\frac{\log((\phi(x_0)-\hat{\phi})/\epsilon_\tau)}
{\alpha\mu(1-\tau^2)}$
\\[8pt]
\hdashline


\multirow{3}{*}{
    \shortstack[l]{
        Minibatch stochastic\\
        gradient descent\\
        for least squares\\
        (Sec~\ref{strong_growth_reg})
    }
}
& \multirow{3}{*}{
    \shortstack[l]{
        \eqref{assum:M1} and\\
        $\displaystyle\alpha\leq\frac{1}{L_1}\cdot\min\left\{1,
        \frac{|\mathcal{S}|}{4M_1}\right\}$
    }
}
& $L_1$
& $\displaystyle
\frac{1}{\tau}
\sqrt{
\frac{\alpha L_1\E_{\xi\sim\Xi}
[\|\nabla\ell(x^\star,\xi)\|^2]}
{|\mathcal{S}|}
}$
& $\displaystyle
\frac{\phi(x_0)-\hat{\phi}}
{\alpha(1-\tau^2)\epsilon_\tau^2}$
\\[8pt]
\cdashline{3-5}

& & KL
& $\displaystyle
\frac{D}{\tau}
\sqrt{
\frac{\alpha L_1\E_{\xi\sim\Xi}
[\|\nabla\ell(x^\star,\xi)\|^2]}
{|\mathcal{S}|}
}$
& $\displaystyle
\frac{D^2}
{\alpha(1-\tau^2)\epsilon_\tau}$
\\[8pt]
\cdashline{3-5}

& & PL
& $\displaystyle
\frac{\alpha L_1\E_{\xi\sim\Xi}
[\|\nabla\ell(x^\star,\xi)\|^2]}
{\tau^2\mu|\mathcal{S}|}$
& $\displaystyle
\frac{\log((\phi(x_0)-\hat{\phi})/\epsilon_\tau)}
{\alpha\mu(1-\tau^2)}$
\\[8pt]
\hdashline


\multirow{3}{*}{
    \shortstack[l]{
        Gradient descent in\\
        random subspaces\\
        (Sec~\ref{sec:gradproj}) or\\
        random directional\\
        derivative descent\\
        (Sec~\ref{sec:rand_dirder})
    }
}
& \multirow{3}{*}{
    \shortstack[l]{
        (Subspace)\\
        $\displaystyle\alpha\leq\frac{1}{L_1}$\\
        (Derivative)\\
        $\displaystyle\alpha\leq
        \frac{1-s/(2n)}
        {(1+\mathcal{O}(n^{-1}))L_1}$
    }
}
& $L_1$
& $0$
& $\displaystyle
\frac{n}{s}\cdot
\frac{\phi(x_0)-\hat{\phi}}
{\alpha(1-\tau^2)\epsilon_\tau^2}$
\\[8pt]
\cdashline{3-5}

& & KL
& $0$
& $\displaystyle
\frac{n}{s}\cdot
\frac{D^2}
{\alpha(1-\tau^2)\epsilon_\tau}$
\\[8pt]
\cdashline{3-5}

& & PL
& $0$
& $\displaystyle
\frac{n}{s}\cdot
\frac{\log((\phi(x_0)-\hat{\phi})/\epsilon_\tau)}
{\alpha\mu(1-\tau^2)}$
\\[8pt]
\hdashline


\multirow{3}{*}{
    \shortstack[l]{
        Previous one with\\
        scaling $n/s$\\
        (Sec~\ref{sec:gradproj}\\
        and~\ref{sec:rand_dirder})
    }
}
& \multirow{3}{*}{
    \shortstack[l]{
        (Subspace)\\
        $\displaystyle\alpha\leq\frac{1}{L_1}\cdot\min\left\{1,
        \frac{s}{2n-2s}\right\}$\\
        (Derivative)\\
        $\displaystyle\alpha\leq
        \frac{s}
        {(1+\mathcal{O}(n^{-1}))nL_1}$
    }
}
& $L_1$
& $0$
& $\displaystyle
\frac{\phi(x_0)-\hat{\phi}}
{\alpha(1-\tau^2)\epsilon_\tau^2}$
\\[8pt]
\cdashline{3-5}

& & KL
& $0$
& $\displaystyle
\frac{D^2}
{\alpha(1-\tau^2)\epsilon_\tau}$
\\[8pt]
\cdashline{3-5}

& & PL
& $0$
& $\displaystyle
\frac{\log((\phi(x_0)-\hat{\phi})/\epsilon_\tau)}
{\alpha\mu(1-\tau^2)}$
\\[8pt]
\hdashline


\multirow{3}{*}{
    \shortstack[l]{
        Scaled randomized\\
        forward finite\\
        difference methods\\
        with Oracle~\ref{def:oracle_stoch_0.1}\\
        (Sec~\ref{sec:SRFFD})
    }
}
& \multirow{3}{*}{
    \shortstack[l]{
        $\displaystyle\alpha\leq\frac{1}{L_1}\cdot\min\left\{1,
        \mathcal{O}\left(\frac{s}{n}\right)\right\}$
    }
}
& $L_1$
& $\displaystyle
\frac{n}{\tau\sqrt{s}}
\sqrt{
\left(\frac{s}{n}+\alpha L_1\right)
\left(\frac{L_1h}{2}+\frac{2\epsilon_f}{h}\right)^2
+\frac{\alpha L_1\sigma_f^2}{h^2}
}$
& $\displaystyle
\frac{\phi(x_0)-\hat{\phi}}
{\alpha(1-\tau^2)\epsilon_\tau^2}$
\\[8pt]
\cdashline{3-5}

& & KL
& $\displaystyle
\frac{nD}{\tau\sqrt{s}}
\sqrt{
\left(\frac{s}{n}+\alpha L_1\right)
\left(\frac{L_1h}{2}+\frac{2\epsilon_f}{h}\right)^2
+\frac{\alpha L_1\sigma_f^2}{h^2}
}$
& $\displaystyle
\frac{D^2}
{\alpha(1-\tau^2)\epsilon_\tau}$
\\[8pt]
\cdashline{3-5}

& & PL
& $\displaystyle
\frac{n^2}{\mu\tau^2s}
\left(
\left(\frac{s}{n}+\alpha L_1\right)
\left(\frac{L_1h}{2}+\frac{2\epsilon_f}{h}\right)^2
+\frac{\alpha L_1\sigma_f^2}{h^2}
\right)$
& $\displaystyle
\frac{\log((\phi(x_0)-\hat{\phi})/\epsilon_\tau)}
{\alpha\mu(1-\tau^2)}$
\\[8pt]
\hdashline


\multirow{3}{*}{
    \shortstack[l]{
        Scaled randomized\\
        forward finite\\
        difference methods\\
        with Oracle~\ref{def:oracle_stoch_0.2}\\
        (Sec~\ref{sec:SRFFD})
    }
}
& \multirow{3}{*}{
    \shortstack[l]{
        $\displaystyle\alpha\leq\frac{1}{L_1}$\\
        $\displaystyle\cdot\min\left\{1,
        \mathcal{O}\left(
        \frac{s}{n+n^2L_v^2}
        \right)\right\}$
    }
}
& $L_1$
& $\displaystyle
\frac{n}{\tau\sqrt{s}}
\sqrt{
\left(\frac{s}{n}+\alpha L_1\right)
\left(\frac{L_1h}{2}+\frac{2\epsilon_f}{h}\right)^2
+\alpha L_1L_c^2
}$
& $\displaystyle
\frac{\phi(x_0)-\hat{\phi}}
{\alpha(1-\tau^2)\epsilon_\tau^2}$
\\[8pt]
\cdashline{3-5}

& & KL
& $\displaystyle
\frac{nD}{\tau\sqrt{s}}
\sqrt{
\left(\frac{s}{n}+\alpha L_1\right)
\left(\frac{L_1h}{2}+\frac{2\epsilon_f}{h}\right)^2
+\alpha L_1L_c^2
}$
& $\displaystyle
\frac{D^2}
{\alpha(1-\tau^2)\epsilon_\tau}$
\\[8pt]
\cdashline{3-5}

& & PL
& $\displaystyle
\frac{n^2}{\mu\tau^2s}
\left(
\left(\frac{s}{n}+\alpha L_1\right)
\left(\frac{L_1h}{2}+\frac{2\epsilon_f}{h}\right)^2
+\alpha L_1L_c^2
\right)$
& $\displaystyle
\frac{\log((\phi(x_0)-\hat{\phi})/\epsilon_\tau)}
{\alpha\mu(1-\tau^2)}$
\\[8pt]
\hdashline


\multirow{3}{*}{
    \shortstack[l]{
        Scaled randomized\\
        central finite\\
        difference methods\\
        with Oracle~\ref{def:oracle_stoch_0.1}\\
        (Sec~\ref{sec:SRCFD})
    }
}
& \multirow{3}{*}{
    \shortstack[l]{
        $\displaystyle\alpha\leq\frac{1}{L_1}\cdot\min\left\{1,
        \mathcal{O}\left(\frac{s}{n}\right)\right\}$
    }
}
& $L_1$
& $\displaystyle
\frac{n}{\tau\sqrt{s}}
\sqrt{
\left(\frac{s}{n}+\alpha L_1\right)
\left(\frac{L_2h^2}{6}+\frac{\epsilon_f}{h}\right)^2
+\frac{\alpha L_1\sigma_f^2}{h^2}
}$
& $\displaystyle
\frac{\phi(x_0)-\hat{\phi}}
{\alpha(1-\tau^2)\epsilon_\tau^2}$
\\[8pt]
\cdashline{3-5}

& & KL
& $\displaystyle
\frac{nD}{\tau\sqrt{s}}
\sqrt{
\left(\frac{s}{n}+\alpha L_1\right)
\left(\frac{L_2h^2}{6}+\frac{\epsilon_f}{h}\right)^2
+\frac{\alpha L_1\sigma_f^2}{h^2}
}$
& $\displaystyle
\frac{D^2}
{\alpha(1-\tau^2)\epsilon_\tau}$
\\[8pt]
\cdashline{3-5}

& & PL
& $\displaystyle
\frac{n^2}{\mu\tau^2s}
\left(
\left(\frac{s}{n}+\alpha L_1\right)
\left(\frac{L_2h^2}{6}+\frac{\epsilon_f}{h}\right)^2
+\frac{\alpha L_1\sigma_f^2}{h^2}
\right)$
& $\displaystyle
\frac{\log((\phi(x_0)-\hat{\phi})/\epsilon_\tau)}
{\alpha\mu(1-\tau^2)}$
\\[8pt]
\hdashline


\multirow{3}{*}{
    \shortstack[l]{
        Scaled randomized\\
        central finite\\
        difference methods\\
        with Oracle~\ref{def:oracle_stoch_0.2}\\
        (Sec~\ref{sec:SRCFD})
    }
}
& \multirow{3}{*}{
    \shortstack[l]{
        $\displaystyle\alpha\leq\frac{1}{L_1}$\\
        $\displaystyle\cdot\min\left\{1,
        \mathcal{O}\left(
        \frac{s}{n+n^2L_v^2}
        \right)\right\}$
    }
}
& $L_1$
& $\displaystyle
\frac{n}{\tau\sqrt{s}}
\sqrt{
\left(\frac{s}{n}+\alpha L_1\right)
\left(\frac{L_2h^2}{6}+\frac{\epsilon_f}{h}\right)^2
+\alpha L_1L_c^2
}$
& $\displaystyle
\frac{\phi(x_0)-\hat{\phi}}
{\alpha(1-\tau^2)\epsilon_\tau^2}$
\\[8pt]
\cdashline{3-5}

& & KL
& $\displaystyle
\frac{nD}{\tau\sqrt{s}}
\sqrt{
\left(\frac{s}{n}+\alpha L_1\right)
\left(\frac{L_2h^2}{6}+\frac{\epsilon_f}{h}\right)^2
+\alpha L_1L_c^2
}$
& $\displaystyle
\frac{D^2}
{\alpha(1-\tau^2)\epsilon_\tau}$
\\[8pt]
\cdashline{3-5}

& & PL
& $\displaystyle
\frac{n^2}{\mu\tau^2s}
\left(
\left(\frac{s}{n}+\alpha L_1\right)
\left(\frac{L_2h^2}{6}+\frac{\epsilon_f}{h}\right)^2
+\alpha L_1L_c^2
\right)$
& $\displaystyle
\frac{\log((\phi(x_0)-\hat{\phi})/\epsilon_\tau)}
{\alpha\mu(1-\tau^2)}$
\\[8pt]

\hline
\end{tabular}%
}
\end{table}

\begin{table}[htbp]
\centering
\renewcommand{\arraystretch}{1.2}
\setlength{\dashlinedash}{0.5pt}
\setlength{\dashlinegap}{1.5pt}

\caption{Summary of complexity results for
Algorithm~\ref{alg:first_order_fixed_step_sgd}
using Oracle~\ref{def:oracle_stoch_1.3}.
The random vectors $\mathcal{U}$ satisfy either
Assumption~\ref{ass:haar} or Assumption~\ref{ass:normal}.
The results under Assumption~\ref{ass:normal}
apply only when $s\ll n$.
Here, $L_1$ denotes Assumption~\ref{assum:lip_cont},
KL denotes Assumptions~\ref{assum:lip_cont}
and~\ref{assum:KLcond}, and PL denotes
Assumptions~\ref{assum:lip_cont}
and~\ref{assum:PLcond}.}
\label{tab:stoch_fixed_step_oracle_13}

\resizebox{\linewidth}{!}{%
\begin{tabular}{|l|l|c|c|c|}
\hline

\textbf{Setting}
& \multicolumn{2}{c|}{\textbf{Assumptions}}
& \makecell{\textbf{Convergence} \\ \textbf{Neighborhood} $\epsilon_\tau \geq \mathcal{O}(\cdot)$}
& \makecell{\textbf{iteration}\\  \textbf{Complexity} $\E[K_{\epsilon_\tau}] \leq \mathcal{O}(\cdot)$}
\\
\hline


\multirow{3}{*}{
    \shortstack[l]{
        Randomized forward\\
        finite difference\\
        methods with\\
        Oracle~\ref{def:oracle_stoch_0.1}\\
        (Sec~\ref{sec:SRFFD})
    }
}
& \multirow{3}{*}{
    \shortstack[l]{
        $\displaystyle\alpha\leq\frac{1}{L_1}$
    }
}
& $L_1$
& $\displaystyle
\frac{1}{\tau}
\sqrt{
n\left(
\frac{L_1h}{2}+\frac{2\epsilon_f}{h}
\right)^2
+\frac{\alpha L_1n\sigma_f^2}{sh^2}
}$
& $\displaystyle
\frac{n}{s}\cdot
\frac{\phi(x_0)-\hat{\phi}}
{\alpha(1-\tau^2)\epsilon_\tau^2}$
\\[8pt]
\cdashline{3-5}

& & KL
& $\displaystyle
\frac{D}{\tau}
\sqrt{
n\left(
\frac{L_1h}{2}+\frac{2\epsilon_f}{h}
\right)^2
+\frac{\alpha L_1n\sigma_f^2}{sh^2}
}$
& $\displaystyle
\frac{n}{s}\cdot
\frac{D^2}
{\alpha(1-\tau^2)\epsilon_\tau}$
\\[8pt]
\cdashline{3-5}

& & PL
& $\displaystyle
\frac{1}{\mu\tau^2}
\left(
n\left(
\frac{L_1h}{2}+\frac{2\epsilon_f}{h}
\right)^2
+\frac{\alpha L_1n\sigma_f^2}{sh^2}
\right)$
& $\displaystyle
\frac{n}{s}\cdot
\frac{\log((\phi(x_0)-\hat{\phi})/\epsilon_\tau)}
{\alpha\mu(1-\tau^2)}$
\\[8pt]
\hdashline


\multirow{3}{*}{
    \shortstack[l]{
        Randomized forward\\
        finite difference\\
        methods with\\
        Oracle~\ref{def:oracle_stoch_0.2}\\
        (Sec~\ref{sec:SRFFD})
    }
}
& \multirow{3}{*}{
    \shortstack[l]{
        $\displaystyle\alpha\leq\frac{1}{L_1}$\\
        $\displaystyle\cdot\min\left\{1,
        \mathcal{O}\left(
        \frac{1}{nL_v^2}
        \right)\right\}$
    }
}
& $L_1$
& $\displaystyle
\frac{1}{\tau}
\sqrt{
n\left(
\frac{L_1h}{2}+\frac{2\epsilon_f}{h}
\right)^2
+\frac{\alpha L_1nL_c^2}{s}
}$
& $\displaystyle
\frac{n}{s}\cdot
\frac{\phi(x_0)-\hat{\phi}}
{\alpha(1-\tau^2)\epsilon_\tau^2}$
\\[8pt]
\cdashline{3-5}

& & KL
& $\displaystyle
\frac{D}{\tau}
\sqrt{
n\left(
\frac{L_1h}{2}+\frac{2\epsilon_f}{h}
\right)^2
+\frac{\alpha L_1nL_c^2}{s}
}$
& $\displaystyle
\frac{n}{s}\cdot
\frac{D^2}
{\alpha(1-\tau^2)\epsilon_\tau}$
\\[8pt]
\cdashline{3-5}

& & PL
& $\displaystyle
\frac{1}{\mu\tau^2}
\left(
n\left(
\frac{L_1h}{2}+\frac{2\epsilon_f}{h}
\right)^2
+\frac{\alpha L_1nL_c^2}{s}
\right)$
& $\displaystyle
\frac{n}{s}\cdot
\frac{\log((\phi(x_0)-\hat{\phi})/\epsilon_\tau)}
{\alpha\mu(1-\tau^2)}$
\\[8pt]
\hdashline


\multirow{3}{*}{
    \shortstack[l]{
        Randomized central\\
        finite difference\\
        methods with\\
        Oracle~\ref{def:oracle_stoch_0.1}\\
        (Sec~\ref{sec:SRCFD})
    }
}
& \multirow{3}{*}{
    \shortstack[l]{
        $\displaystyle\alpha\leq\frac{1}{L_1}$
    }
}
& $L_1$
& $\displaystyle
\frac{1}{\tau}
\sqrt{
n\left(
\frac{L_2h^2}{6}+\frac{\epsilon_f}{h}
\right)^2
+\frac{\alpha L_1n\sigma_f^2}{sh^2}
}$
& $\displaystyle
\frac{n}{s}\cdot
\frac{\phi(x_0)-\hat{\phi}}
{\alpha(1-\tau^2)\epsilon_\tau^2}$
\\[8pt]
\cdashline{3-5}

& & KL
& $\displaystyle
\frac{D}{\tau}
\sqrt{
n\left(
\frac{L_2h^2}{6}+\frac{\epsilon_f}{h}
\right)^2
+\frac{\alpha L_1n\sigma_f^2}{sh^2}
}$
& $\displaystyle
\frac{n}{s}\cdot
\frac{D^2}
{\alpha(1-\tau^2)\epsilon_\tau}$
\\[8pt]
\cdashline{3-5}

& & PL
& $\displaystyle
\frac{1}{\mu\tau^2}
\left(
n\left(
\frac{L_2h^2}{6}+\frac{\epsilon_f}{h}
\right)^2
+\frac{\alpha L_1n\sigma_f^2}{sh^2}
\right)$
& $\displaystyle
\frac{n}{s}\cdot
\frac{\log((\phi(x_0)-\hat{\phi})/\epsilon_\tau)}
{\alpha\mu(1-\tau^2)}$
\\[8pt]
\hdashline


\multirow{3}{*}{
    \shortstack[l]{
        Randomized central\\
        finite difference\\
        methods with\\
        Oracle~\ref{def:oracle_stoch_0.2}\\
        (Sec~\ref{sec:SRCFD})
    }
}
& \multirow{3}{*}{
    \shortstack[l]{
        $\displaystyle\alpha\leq\frac{1}{L_1}$\\
        $\displaystyle\cdot\min\left\{1,
        \mathcal{O}\left(
        \frac{1}{nL_v^2}
        \right)\right\}$
    }
}
& $L_1$
& $\displaystyle
\frac{1}{\tau}
\sqrt{
n\left(
\frac{L_2h^2}{6}+\frac{\epsilon_f}{h}
\right)^2
+\frac{\alpha L_1nL_c^2}{s}
}$
& $\displaystyle
\frac{n}{s}\cdot
\frac{\phi(x_0)-\hat{\phi}}
{\alpha(1-\tau^2)\epsilon_\tau^2}$
\\[8pt]
\cdashline{3-5}

& & KL
& $\displaystyle
\frac{D}{\tau}
\sqrt{
n\left(
\frac{L_2h^2}{6}+\frac{\epsilon_f}{h}
\right)^2
+\frac{\alpha L_1nL_c^2}{s}
}$
& $\displaystyle
\frac{n}{s}\cdot
\frac{D^2}
{\alpha(1-\tau^2)\epsilon_\tau}$
\\[8pt]
\cdashline{3-5}

& & PL
& $\displaystyle
\frac{1}{\mu\tau^2}
\left(
n\left(
\frac{L_2h^2}{6}+\frac{\epsilon_f}{h}
\right)^2
+\frac{\alpha L_1nL_c^2}{s}
\right)$
& $\displaystyle
\frac{n}{s}\cdot
\frac{\log((\phi(x_0)-\hat{\phi})/\epsilon_\tau)}
{\alpha\mu(1-\tau^2)}$
\\[8pt]

\hline
\end{tabular}%
}
\end{table}

\end{document}